\documentclass[twoside,12pt]{book}
\usepackage{epsfig,mathptm,float}
\pdfoutput=1
\begin{document}

\newcommand{\ben}{\begin{enumerate}}
\newcommand{\een}{\end{enumerate}}
\newcommand{\bit}{\begin{itemize}}
\newcommand{\eit}{\end{itemize}}
\newcommand{\beq}{\begin{equation}}
\newcommand{\eeq}{\end{equation}}
\newcommand{\beqar}{\begin{eqnarray*}}
\newcommand{\eeqar}{\end{eqnarray*}}
\newcommand{\ul}{\underline}
\newcommand{\no}{{\it no answer provided for this exercise}}
\newtheorem{D}{Definition}
\newtheorem{theorem}{Theorem}
\newtheorem{property}{Property}
\newtheorem{C}{Conclusion}
\newcommand{\bec}{\begin{C}}
\newcommand{\eec}{\end{C}}
\newlength{\figurewidth}
\setlength{\figurewidth}{0.45\textwidth}
\pagestyle{headings}

\begin{titlepage}
\begin{center}
\vspace*{2cm}
{\Huge \bf QUALITATIVE ANALYSIS OF 

\vskip 1pc
DIFFERENTIAL EQUATIONS }

\vspace{.8cm}
{\Large Alexander Panfilov}
\end{center}
\begin{figure}[h]
\centerline{
\psfig{type=pdf,ext=.pdf,read=.pdf,figure=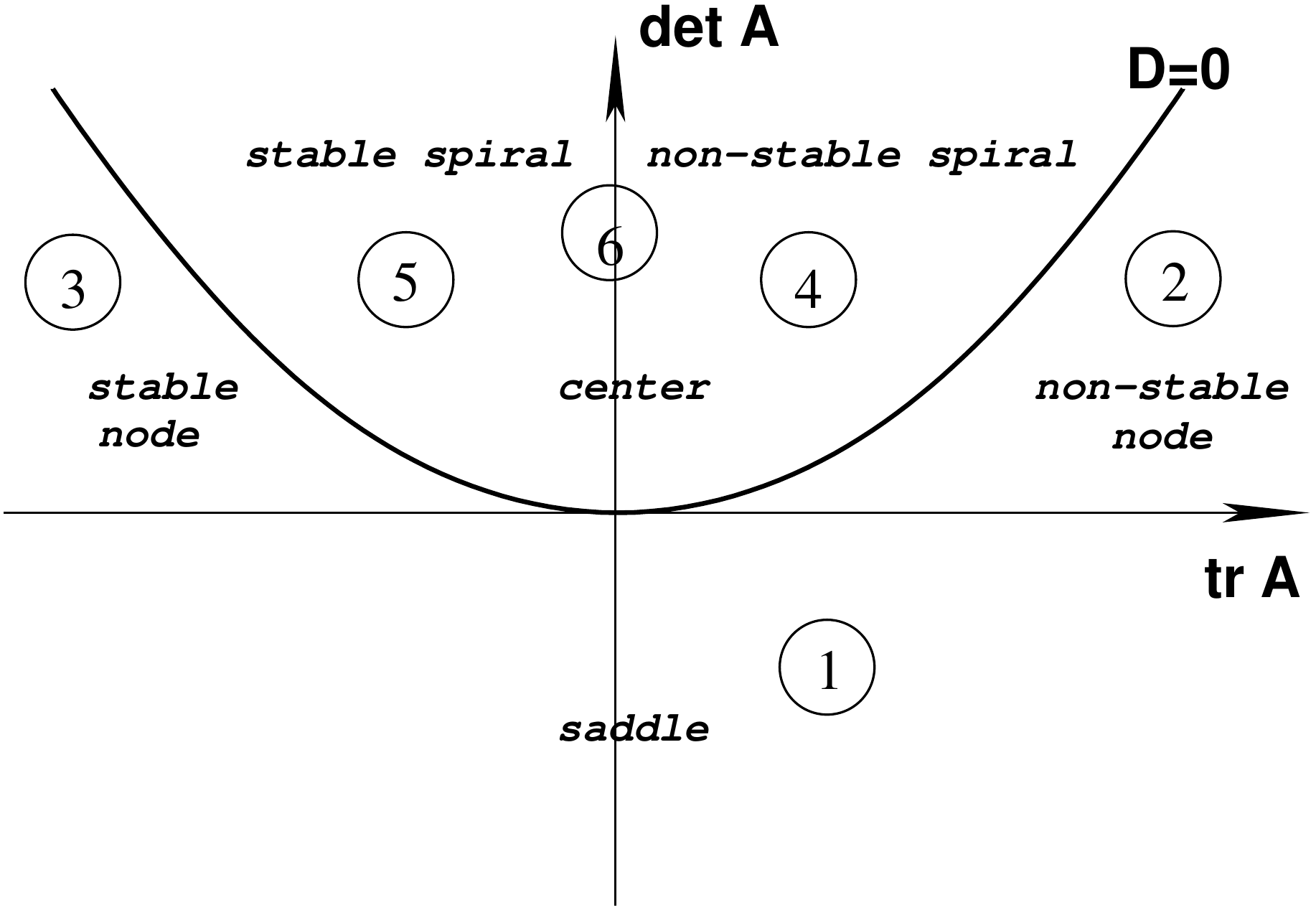,width=12cm}
}
\end{figure}
\end{titlepage}

\title{QUALITATIVE ANALYSIS OF DIFFERENTIAL EQUATIONS}
\author{Alexander Panfilov \\
Theoretical Biology, Utrecht University, Utrecht\\ {\bf \copyright
2010}\\
\date{2010}
}
\maketitle



\tableofcontents


\chapter{Preliminaries }
\setlength{\parskip}{0.65\baselineskip}
\section{Basic algebra  \label{secBasicAlgebra}}

\subsection{Algebraic expressions}
Algebraic expressions are  formed from numbers, letters and   arithmetic operations. The letters may represent unknown variables, which should be found from solutions of equations,  or parameters (unknown numbers) on which the  solutions depend. 

Below,  we review  examples several basic operations which help us to
work with algebraic express of ions.

One of the most basic algebraic operations is opening of  parentheses, or simplification of   expressions. For that we use the following  rule:
$$ (a+b)*(c+d)=ac+ad+bc+bd $$ 
note, that here $ac$ means $a*c$, etc.,  as in algebra the
multiplication is often    omitted.

{\it Example (open parenthesis):}
$(4x+2a)(2-3x)=8x-12x^2+4a-6ax$

Sometimes we use parentheses to factor expressions:

{\it Example (factor expression):}
$ 9x^3+3x^2-6a^4x^2=x^2(9x+3-6a^4)=3x^2(3x+1-2a^4)$

In many cases we also need to work  with fractions:

{\it Example (the same denominator):} ${a^2-ca \over 3}= {a^2 \over 3}- {ca \over 3} $ 

{\it Example (different  denominators):} ${a^2 \over b} +{b \over a}=
{a^2 \over b}{a \over a} +{b \over a}{b \over b}
={a^2*a \over ab} +{b*b \over ab} = {a^3+b^2 \over ab} $

{\it Example (fractions simplifications):} ${a^2-ca \over 3a}= {a-c \over 3} $

To divide a fraction  ${a \over b}$ by another fraction ${c \over d}$ we just need to multiply it by its  inverse ${d \over c}$: 
$ {a \over b} : {c \over d}={a \over b} * {d \over c}={ad \over bc} $, 

{\it Example (division):} ${4 \over 7} : {2a \over 5}={4 \over 7} {5 \over
2a}={20 \over 14a}={10 \over 7a}$

Also note that: ${ x+y \over z+d} \mathbf{\ne}{ x \over z}+ { y \over d}$ 

\subsection{Limits \label{sec_limits}} 

We call $A$ a  limit of the function $f(x)$ when  $x$ approaches
$a$, if the value of
$f(x)$ get closer and closer to $A$ when $x$ takes values closer and closer to $a$. We write it formally as:

\beq
\label{limit1}
\mathop {\lim }\limits_{x \to a }    f(x) =A
\eeq

In many cases,  finding the  limit is trivial: we just need  to substitute
the value of $x=a$ into our function:  

\beq
\label{limit2}{\lim_{x \rightarrow a }}
f(x) =f(a)
\eeq

{\it Example:} $\mathop {\lim }\limits_{x \to 2 }  \quad x^3 =2^3=8$ 

Functions which have such   property are called continuous and most of
the functions  used in biology are
continuous. However, there are several important exception.

The first case, which will be the most important for us, is finding
of  limit of the function when $x \rightarrow \infty$. Finding  such
limits is important as it gives an asymptotic behaviour of our system
when the size of a population becomes very large. Unfortunately, there
is no such number '$\infty$' which we can substitute into our function
to find a limit using formula (\ref{limit2}). For functions without parameters, we can guess the limit by substituting large values to (\ref{limit2}), e.g. $x=10000, 20000, etc$, but what to do for functions with parameters?

Let us discuss this problem  for a special class of
functions,  which are the most relevant to  our course, the  so called
rational functions $f(x)={p(x) \over g(x)}$, where $p(x)$ and $g(x)$
are polynomials.  In that case we can always find the limit using the
following property of the power function:

\beq
\label{limit3}
{\lim_{x \rightarrow \infty }}
{C \over x^\alpha}  =0
\eeq
where $C$ is an arbitrary constant  and $\alpha>0$.

To prove it note that if $x$ approaches $\infty$ (becomes
larger and larger),  the power function $x^\alpha$ with  $\alpha>0 $
also becomes larger and larger and therefore ${C \over x^\alpha}$ will
be closer and closer to zero, thus in accordance with the definition
$\mathop {\lim }\limits_{x \to \infty }  {C \over x^\alpha} =0$

To find the limit using this rule we need to do the following: (1) 
 find the highest power of $x$ in our expression ${p(x) \over
g(x)}$, (2)  divide each term in our function by $x$ in that power,  and
(3)  find the limit of each term using property (\ref{limit3}). Let us consider three typical examples:

{\it Example (find the limit):}  $\mathop {\lim }\limits_{N \to \infty } 
{aN^2-3N\over 3-2N^2}$. 

The highest power is $N^2$, division gives: ${{aN^2 \over N^2}-3{N
\over N^2} \over {3 \over N^2}-{2 N^2 \over N^2}}= {{a}-{3 \over N}
\over {3 \over N^2}-{2 }}$. The limits of the individual terms are:
${{a}- 0 \over 0 -2 }=-{a \over 2}$

{\it Example (find the limit):}  $\mathop {\lim }\limits_{P \to \infty } 
{aP-3bP^3\over cP-dP^2}$, $a,b,c,d \ne 0$.

Similar steps give us: ${aP-3bP^3\over cP-dP^2}={{aP \over P^3}
-{3bP^3 \over P^3} \over {cP \over P^3} -{dP^2 \over P^3}}={{a \over
P^2} -{3b} \over {c \over P^2} -{d \over P}}={0 -{3b} \over 0
-0}={-{3b} \over 0} $. This expression does not have sense as we cannot divide by zero  and we do not have a finite limit for this
function.

{\it Example (find the limit):}  $\mathop {\lim }\limits_{x \to \infty } 
{ax^3-bx^2+c \over a x^4 -b}$, \hskip 2pc $a,b,c \ne 0$ . 

${ax^3-bx^2+c \over a x^4 -b}={{ax^3 \over x^4 } -{bx^2\over x^4
}+{c\over x^4 } \over { a x^4 \over x^4 }-{ b\over x^4 }}={{a \over x
} -{b\over x^2 }+{c\over x^4 } \over { a }-{ b\over x^4 }}={0 -0+0
\over { a }-0}= {0 \over a}=0$.

Another non-trivial situation occurs when the denominator of our
function $f(x)={p(x) \over g(x)}$ becomes zero for some value of $x$,
for example $f(x)={2 \over x-3}$ for $x=3$. In this case the formula
(\ref{limit2}) for limit cannot be used and other more careful
analysis is necessary. If using of calculator we substitute some numbers
into our function around point $3$ we will find the following:  if
$x$ becomes closer and closer to $3$ from the left,
e.g.  $x=3.1; 3.05; 3.01, 3.005; etc$  the
function value becomes larger and larger, while if $x$ becomes closer
and closer to $3$ from the right, e.g.  $x=2.9;
2.95; 2.99, 2.995; etc$  the function value is
negative and its absolute value also becomes larger and larger. We can
formally write it as $\mathop {\lim }\limits_{x \to 3^+ }  {2 \over x-3}=
+\infty$, while $\mathop {\lim }\limits_{x \to 3^- }  {2 \over x-3}= -\infty$. However, in a strict sense,
as there is no real number for which $f(x) $ approaches for  $x$ close  to $3$ (from either side) thus the limit here does not
exist.

We will use limits for drawing our functions and will see that limits
 at infinity give us horizontal asymptotes of our graphs, while blow
up of functions for some $x$  (as in the last example)  give
us  vertical asymptotes.

\subsection{Equations}

An equation is a mathematical relationship involving unknown
variables. These unknowns are usually expressed by  letters 'x', 'y',
however in biology we use many other letters 
(e.g. 'N', 'P'. 'T', 'V', etc.), which maybe somewhat related to the
name of the species they describe.  Solving  equations means  finding 
unknown(s) such that after substitution in the equation the left and
right hand sides will be equal to each other.  For example: equation
$2x-16=-10$ has a solution $x=3$, as $2*3-16=6-16=-10$.

The usual way to solve equations  which have
unknown variables in the first power only (linear equations), is to isolate the unknowns: 
$$x=[known \; numbers]$$
 We can  achieve that by using
the following rules of equation algebra: (1) we can multiply, or
divide  both sides of the equation by the same number, and (2) we can
move numbers/expressions from one to the other side of the equation,
by changing  their sign. The proof of these rules is
trivial. Indeed if two expressions  are the same $X=Y$, then if we
multiply (or divide) both of them by the same number $a$, they still
will be the same $aX=aY$. Similarly if $X=Y+a$, we can add $-a$ to the
both sides of the equation, which will not change the equality, but we
get: $X-a=Y+a-a$, or $X-a=Y$. Thus we see, that we were able to move
$a$ from the right hand side to the left hand side of our equation,
but it changed its sign as a result.

{\it Example(solve equation):} $4-2x=2-4x$. 

{\it Solution:}$ -2x+4x=2-4 \quad 2x=-2 \quad x=-1$

For a quadratic equation $ax^2+bx+c=0$ we use the 'abc' formula, which gives us the solutions as  $x_{1,2}={-b \pm  \sqrt {b^2-4ac} \over 2a}$.  

{\it Example (solve equation):} $x+1={4 \over 1+x}$.

{\it Solution:} $(1+x)(x+1)=4
\quad 1+x+x^2+x=4; \quad  x^2+2x+1-4=0; \quad x^2+2x-3=0$ from the 'abc'
formula $x_{1,2}={-2 \pm \sqrt{ 4-4*1*(-3)} \over 2}={-2 \pm \sqrt{16}
\over 2}={-2 \pm 4 \over 2} \quad x_1=1;x_2=-3$.

Equation may also contain parameters.  A parameter is an unknown
number (constant) that may have any value. It is different
from the unknown variable, as the parameter is just a
constant on which our solution depends.

{\it Example (solve equation):} $a{k \over P^2} -dP=0$, where  $P$ in unknown variable and
$a, k,d > 0$ are parameters.

{\it Solution:} By multiplying both sides by $P^2$ we get $ak -dP*P^2=0$, or $ak =dP*P^3$, or $P^3={ak
\over d}$, thus $P=\sqrt[3]{{ak \over d}}$. We see that the solution depends on 3 parameters and if someone
provides us with their values  we will be able
to find the solution by substituting  the parameter values into
the final formula.

\subsection{Systems of equations \label{SysEq}}

To  solve a system of two linear equations we  express one variable via the other  and substitute it into the other equation.

 {\it Example (solve the system of equations):} $
\left\{
\begin{array}{l}
2x +y =5\\ x+ y=3
\end{array}
\right.
$ 

{\it Solution:} From the second equation we find $ x=3-y$, so we substitute $x$ into the first equation:
$2(3-y)+y=5; \quad 6-2y+y=5; \quad -y=-1; \quad y=1$, now substitute this value to  $ x=3-y$ and find  $ x=3-1-2$, thus the solution is $x=2,y=1$.

Unfortunately, there are no  general rules  to solve a system of nonlinear
equations. The usual practical way is to start with a more simple
equation, try to obtain from it as much information  as possible and
then substitute it to the other equation. It is also very helpful to factor expressions in order  to simplify them.

 {\it Example (solve the system of equations):} $
\left\{
\begin{array}{l}
2n-2n^2-2np=0\\ np-2p=0
\end{array}
\right.
$ 

{\it Solution:} From the second equation by factoring  we find $
np-2p=p(n-2)=0$. The product is zero only if one of the multipliers is
zero, thus we have two possibilities $p=0$ or $n=2$. If we substitute
 $p=0$ into the first equation we find $2n-2n^2-0=0$,
or $2n(1-n)=0$, thus for $p=0$ we have two solutions $n=0$, or $n=1$; now substitute   $n=2$ into the first equation: $2*2-2*4-2*2*p=0$, $4-8-4p=0$, $-4p=4$, thus for $n=2$ we found $p=-1$. Overall,   we found the following three solutions of the given system  $(n=0,p=0),(n=1,p=0),(n=2,p=-1)$. 

Systems may also contain parameters.

 {\it Example (solve the system):} $
\left\{
\begin{array}{l}
an-an^2-bnp=0\\ np-kp=0
\end{array}
\right.$, where  $n,p$ are variables and $a,b,k >  0$, are the parameters.

{\it Solution:} We proceed similarly as in the previous case. From the
second equation: $np-kp=p(n-k)=0$, thus we have two cases $p=0$ or
$n=k$.  After substituting $p=0$ into the first equation we get:
$an-an^2-0=0$, $an(1-n)-0$, thus $n=0$, or $n=1$; after substituting
$n=k$ into the first equation we get: $ak-ak^2-bkp=0$, $ak(1-k)=bkp$,
$a(1-k)=bp$, thus $p={a(1-k) \over b}$. Therefore, we found three
solutions: $(n=0,p=0),(n=1,p=0),(n=k,p={a(1-k) \over b})$. It is easy
to see that if we substitute the parameter values $a=2,b=2,k=2$ to
these formulas we obtain the solution of the previous problem. Note
also, that for systems with parameters we need to be careful as not
all operations are allowed for arbitrary parameter values. In our
example in order to obtain the solution we had to make several
divisions by parameters $a,b$, and $k$. However we can
always do that as  the  parameters are positive numbers
($a,b,k > 0$) and thus they cannot be equal to zero.

Finally note, that we can solve systems of three and more equations  similarly, by subsequent substitutions from one equation to another, etc..
\section{Functions of one variable \label{sec1d}}
In science the relationships between quantities are
normally expressed using functions. The simplest type of functions are
functions of one variable. The function of one variable $f$ is a rule
that allows us to find the value of a variable (number) $f$ from a single
variable (number) $x$. We denote  it as $f(x)$. Below are examples of the most important functions:
\begin{enumerate}
\item[] power functions $x^a$, for example
\beq
f(x)=x^{1 \over 2}=\sqrt{x}; \;\;\;f(x)=x^{-2}={1 \over x^2}.
\eeq
\item[] polynomials, $ax^3+ ..+cx+d$,  for example:
\beq
\label{epar1}
f(x)=3x^3-2x^2+1
\eeq
\item[] rational functions $f(x)={p(x) \over g(x)}$:
\beq
f(x)={2-x \over x^2+1}
\eeq
\item[] trigonometric functions  $sin$, $cos$, $tan$: 
\beq
f(x)=2*sin(x);\;\; f(x)=cos(2x-1);\;\; f(x)=tan({x \over 2})
\eeq
\item[] exponential $a^x$  and  logarithmic function $^a log(x)$
\beq
f(x)=e^{x+1};\;\; f(x)=^{10}log(2x);\;\; f(x)=2^{2x}\;\; 
\eeq
\item[] etc.
\end{enumerate}
{\bf The derivative}  of a function $f(x)$  at point $x^*$ is given by the following limit:
\beq
\label{deriv}
f'(x^*)={df \over dx}=  \mathop {\lim }\limits_{x \to x^* }     {f(x)-f(x^*)  \over x-x^*}
\eeq

\begin{figure}[]
\centerline{
\psfig{type=pdf,ext=.pdf,read=.pdf,figure=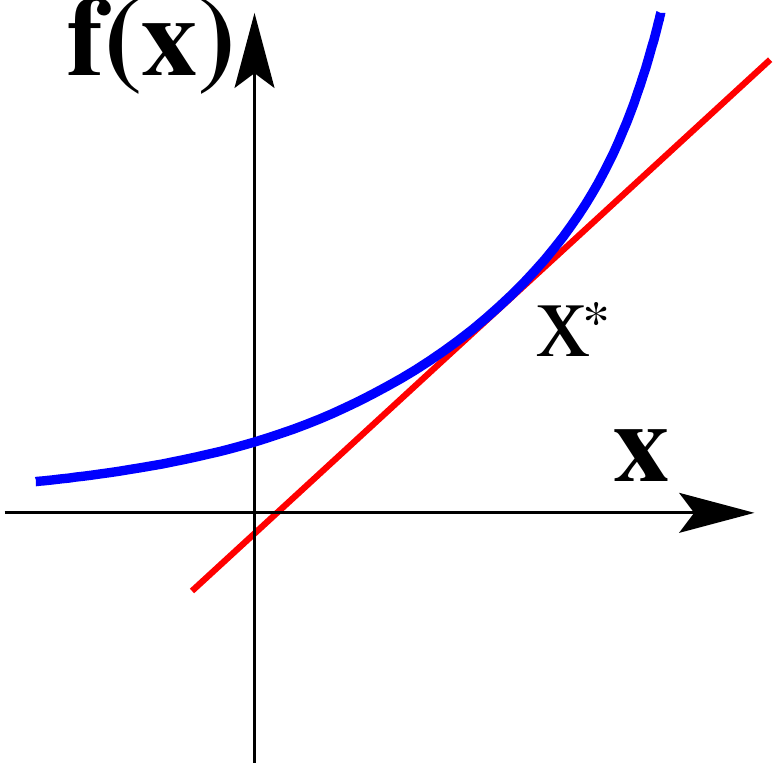,width=5.cm}
}
\caption{\label{gderiv}}
\end{figure}

The derivative $f^{'}(x^*)$ shows the rate of change of a function
$f(x)$ at a point $x^*$ and has many important applications:
\bit
\item If x(t) is the  distance traveled by a car as a function of time $t$, then $dx/dt$ gives the velocity of the car.
\item If $n(t)$ is the size of a population as the function of time, then $dn/dt$ gives the rate of growth of the population.
\eit
 Geometrically, the derivative $f^{'}(x^*)$ gives the slope of the tangent line to the graph of the function at the point $x^*$ (fig.\ref{gderiv}).

\noindent A graph of a line tangent to the function $f(x)$ at point $x^*$ (fig.\ref{gderiv}) is given by the following equation:
\beq
\label{diff}
y=  f(x^*) +  f'(x^*)*(x-x^*) 
\eeq
\noindent Equation (\ref{diff}) is also known as a  {\bf linear approximation} of function $f(x)$ at point $x^*$:

Let us check formula (\ref{diff}) by approximating the function
$y=2x^2+1$ at $x^*=1$.  We find: $f(x^*)=2*1^2+1=3,\; f'=4x,\;f'(1)=4$, hence
$f(x) \approx 3 + 4*(x-1)$.  At $x=1.1$ this approximate formula gives
$f(1.1) \approx 3+4*(1.1-1)=3.4$.  The exact value is
$f(1.1)=2*1.1^2+1=3.42$. So the error is just 0.6\%. However, if
$x=0,\; f(0) \approx 3+4*(0-1)=-1$ while the exact value is
$f(0)=1$. So we see, that the approximate formula works good if $x$ is
close to $x^*$ only.

\noindent {\bf Functions with parameters.} Functions may depend not only on variable(s)   but also on parameters.  We have already seen the following example of the function  $f(x)$ that depends on three parameters $a,b,c$:
\beq
\label{epar2}
f(x)=ax^2+bx+c
\eeq
Equation (\ref{epar2}) describes a general quadratic polynomial. If we
choose, for example $a=3,b=-2,c=1$ we will get the function given by
equation (\ref{epar1}). Studying functions with parameters allows us
to obtain results for whole classes of functions. We will frequently
use functions with parameters in our course. This is because
biological models usually depend on many (up to hundreds) parameters
and in many situations the exact values of these parameters are
unknown. One of practical difficulties in working with parameters is
that use of calculators is very limited, because calculators cannot do
calculations with unknown quantities. The most valuable methods to study functions with parameters 
are direct algebraic computations and analysis of the obtained
formulas. In this course we will widely use the graphical methods of
representation of function. Let us start with review of the basic
function graphs.

\section{Graphs of functions of one variable  \label{sec1dGraph}}

\noindent {\bf Example of graphs.} We usually represent functions using  graphs. To do that we plot
the value of the variable $x$ along the $x$-axis and the value of the function
$f(x)$ along the $y$-axis.  Let us start first by listing typical graph shapes that are important in this course.
\vskip 1pc
\begin{figure}[H]
\centerline{
\psfig{type=pdf,ext=.pdf,read=.pdf,figure=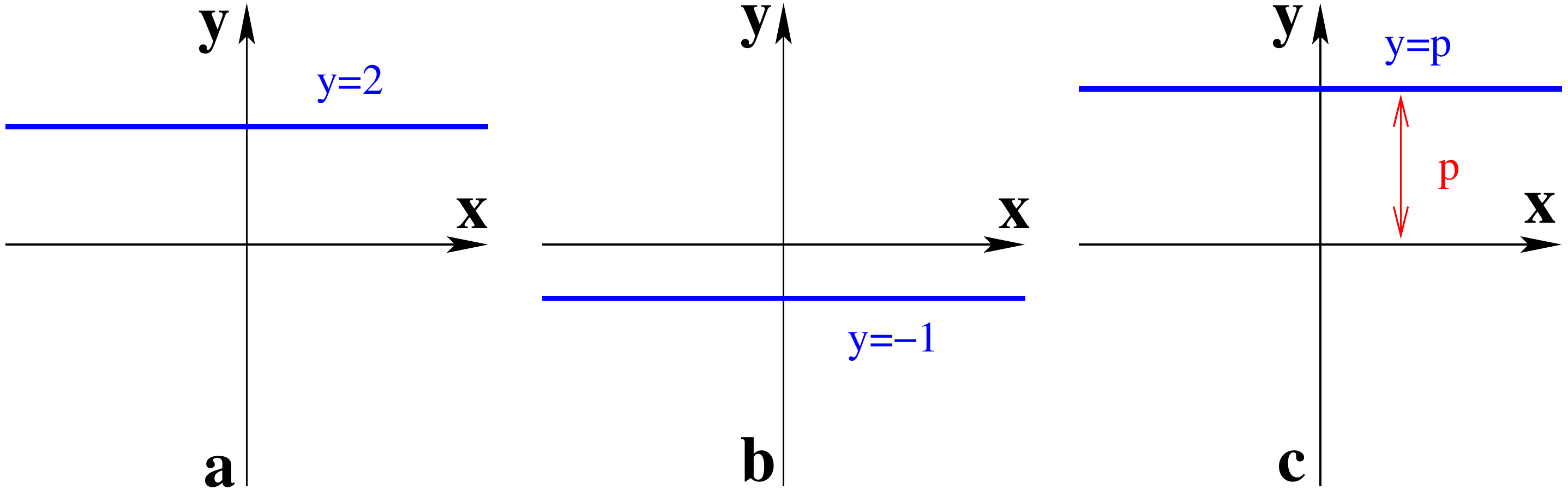,width=12.cm}
}
\caption{\label{g1d1}}
Equation $y=p$ produces a horizontal line at the level $p$ (fig.\ref{g1d1}). 
\end{figure}
\begin{figure}[H]
\centerline{
\psfig{type=pdf,ext=.pdf,read=.pdf,figure=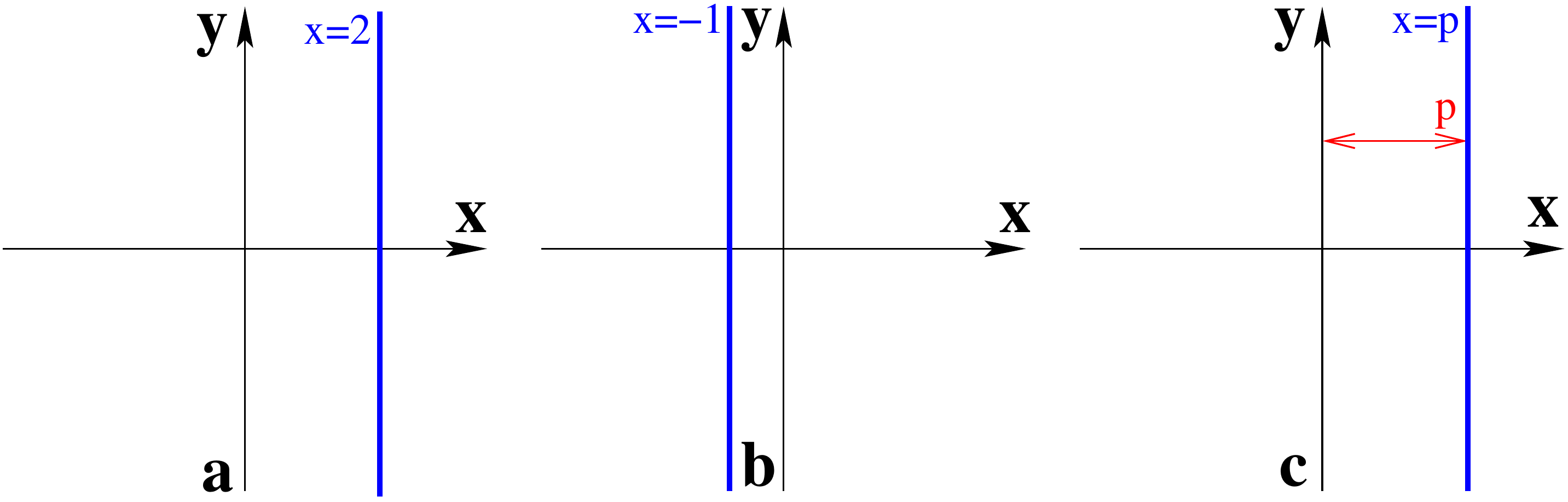,width=12.cm}
}
\caption{\label{g1d2}}
Equation $x=p$ produces a vertical line shifted by  $p$ from the $y$-axis (fig.\ref{g1d2})
\end{figure}
\begin{figure}[H]
\centerline{
\psfig{type=pdf,ext=.pdf,read=.pdf,figure=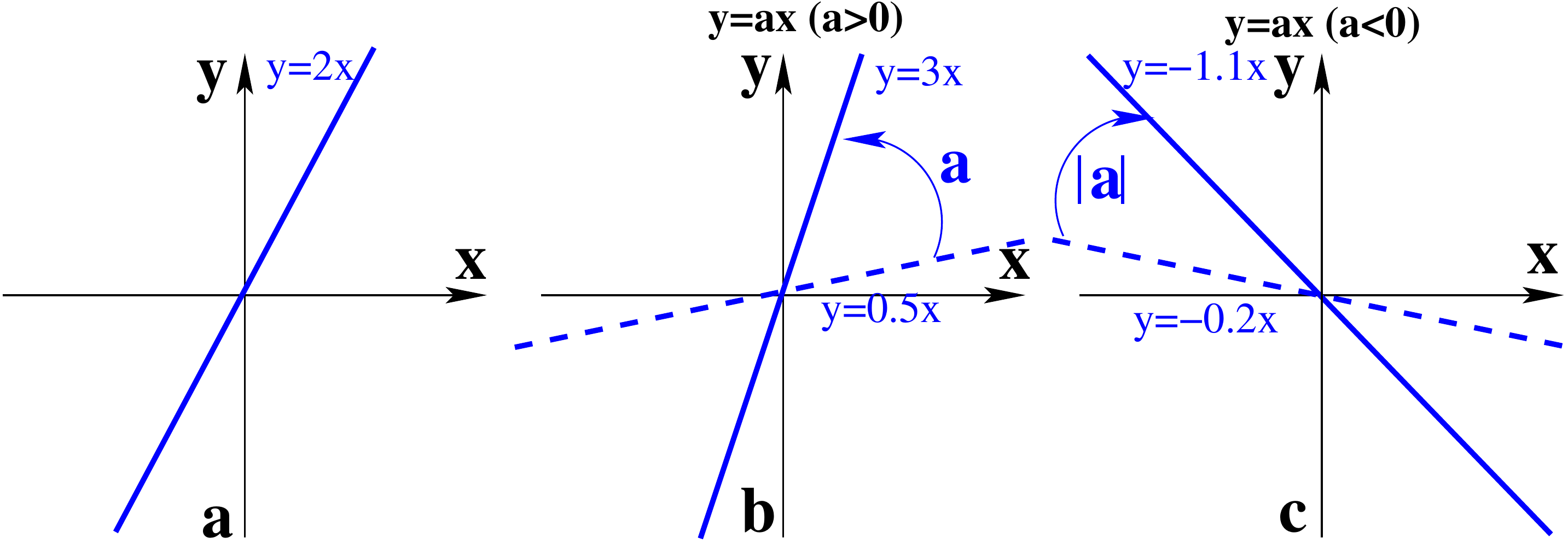,width=12.cm}
}
\caption{\label{g1d3}}Equation $y=ax+p$ (linear function) produces  a straight line
with the slope defined by the parameter $a$: the
larger   the absolute value of $a$, the steeper is the slope (fig.\ref{g1d3}).
\end{figure}
\begin{figure}[H]
\centerline{
\psfig{type=pdf,ext=.pdf,read=.pdf,figure=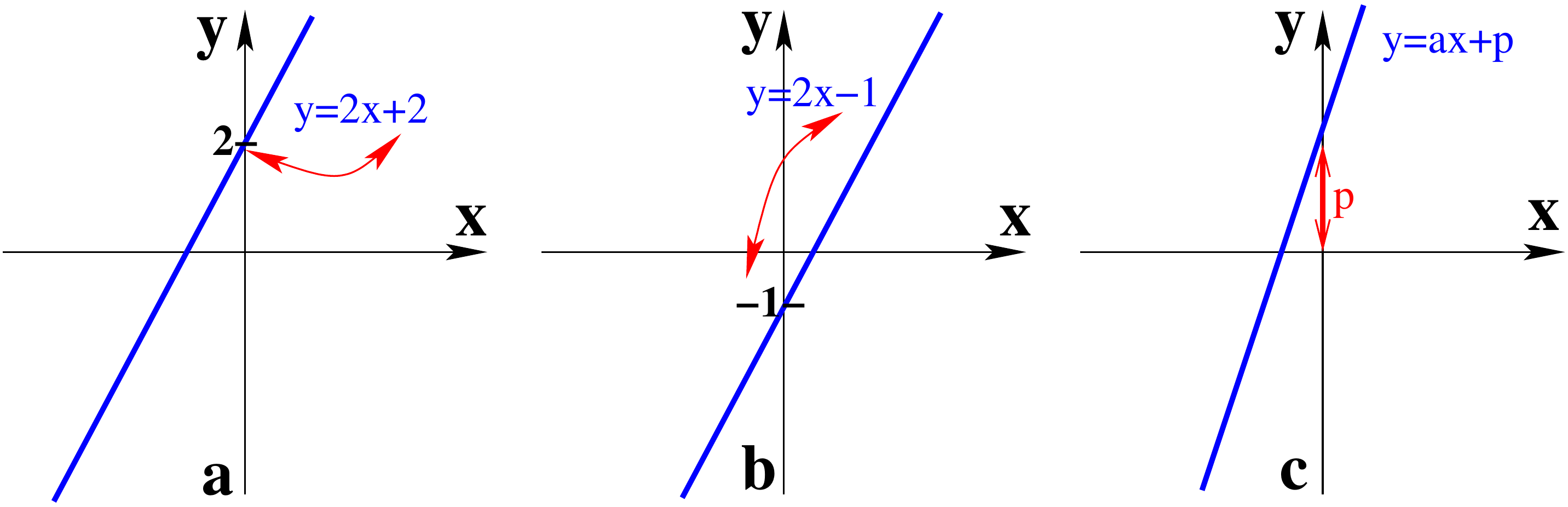,width=12.cm}
}
\caption{ \label{g1d4} }The
parameter $p$ in $y=ax+p$ accounts for the vertical shift of the graph fig.\ref{g1d3}.
\end{figure}

\begin{figure}[H]
\centerline{
\psfig{type=pdf,ext=.pdf,read=.pdf,figure=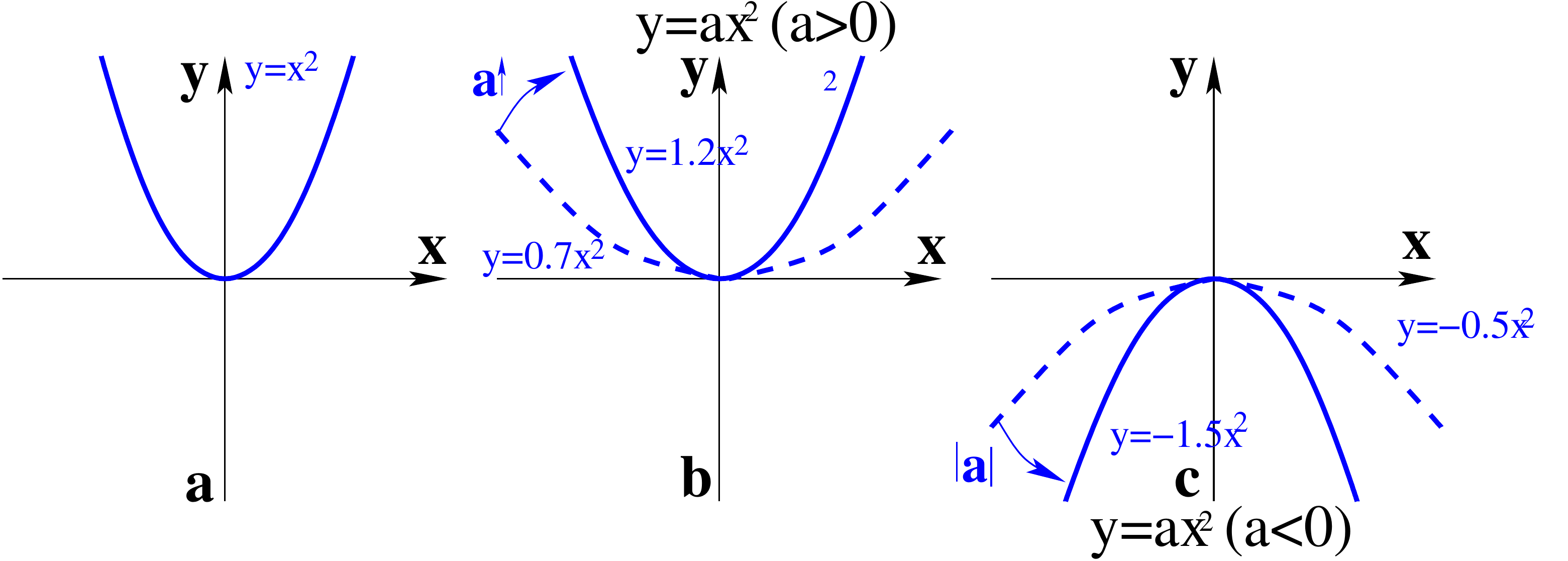,width=12.cm} 
}
\caption{  \label{g1d5} }
Equation $y=ax^2$ produces a parabola,  if $a>0$ the
parabola is opened upward (fig.a,b), and if $a<0$ the parabola is opened
downward (fig.c). The larger the absolute value of $a$ is, the steeper is  the parabola.
\end{figure}
\begin{figure}[H]
\centerline{
\psfig{type=pdf,ext=.pdf,read=.pdf,figure=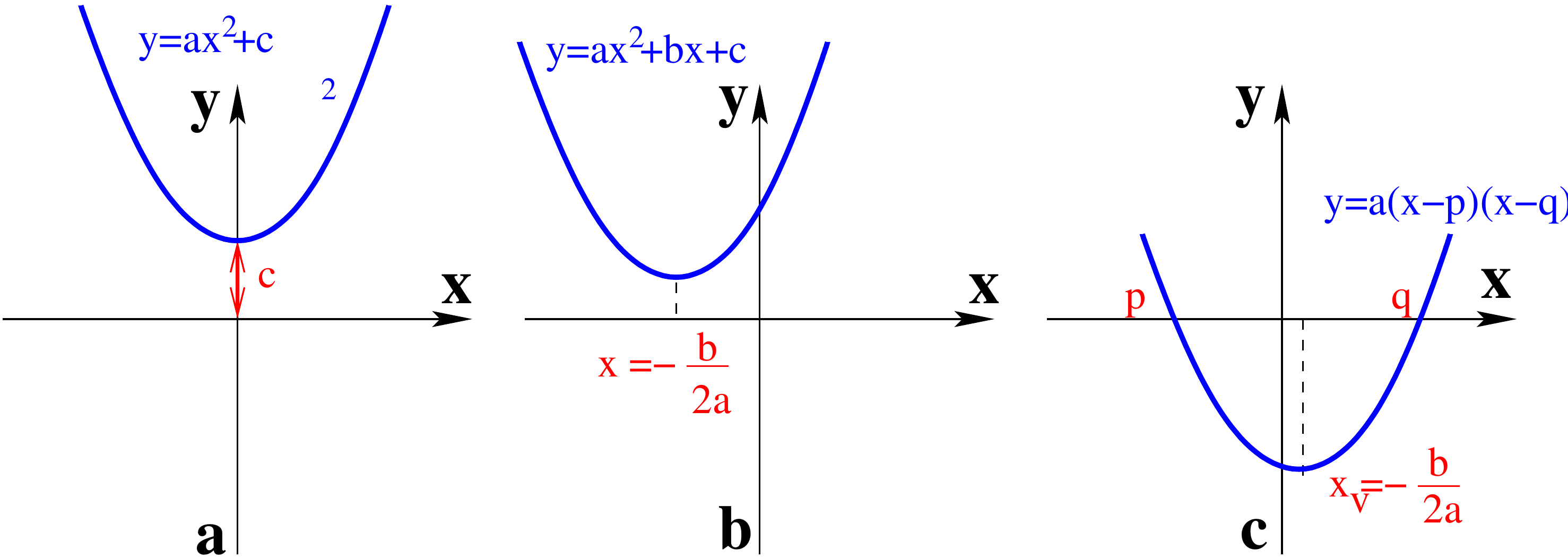,width=12.cm}
}
\caption{\label{g1d6}}
Equation $y=ax^2+bx+c$ also produces a parabola. 
Parameter $c$ (fig.\ref{g1d6}a) accounts
for the vertical shift of the graph. Parameter $b$ accounts for a horizontal shift of the parabola. It is possible to show   that the horizontal shift of the parabola is given by $-{b \over 2a}$ (fig.\ref{g1d6}b). We can calculate this shift by determining the location of the  vertex of the parabola which is a point of extremum (maximum or minimum) of the function. At this point  the derivative of the function to zero  $(ax^2+bx+c)'=2ax+b=0$, Thus the $x$ coordinate of the vertex is given by $x_v=-{b \over 2a}$, or in other words the (vertex of) parabola is shifted by $x_v=-{b \over 2a}$ from its central location in  (fig.\ref{g1d6}a. Note also, that a parabola may have up to two points of
intersection of the graph with the $x$-axis (zeros of the function).
They can be found from the 'abc' formula for roots of the equation
$ax^2+bx+c=0$, and if these roots ($p,q$) are known, the graph can easily be
depicted using them (fig.\ref{g1d6}c). Note,  that in this case the
vertex of the parabola is always located at the middle between these
two roots.
\end{figure}

\begin{figure}[H]
\centerline{
\psfig{type=pdf,ext=.pdf,read=.pdf,figure=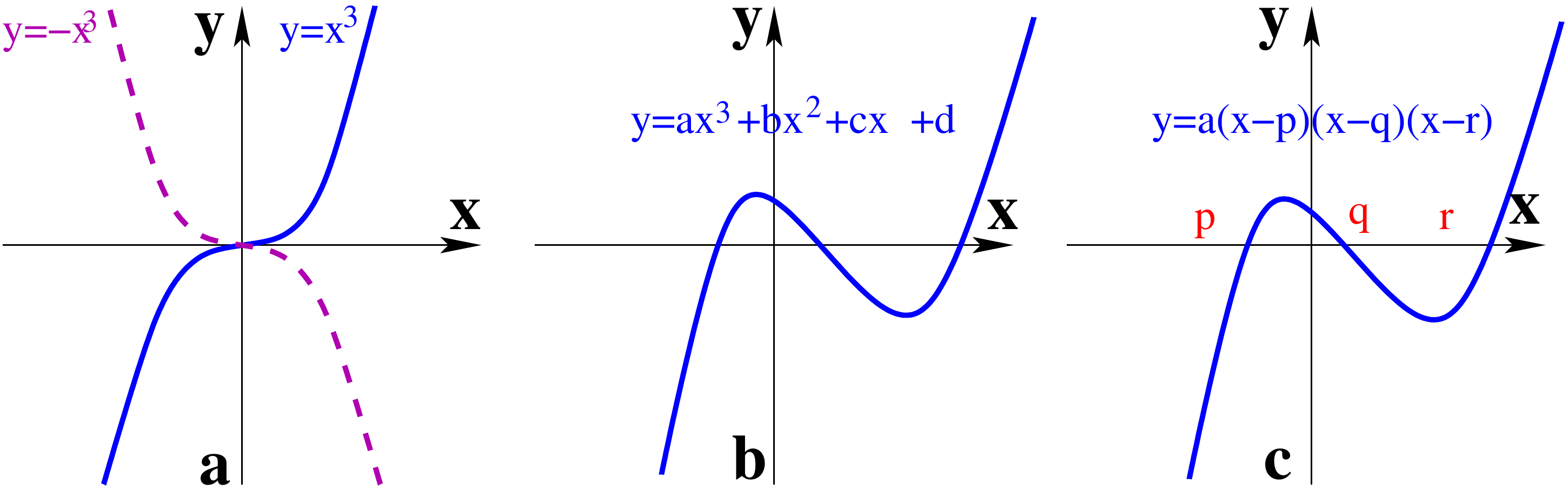,width=15.cm}
}
\caption{\label{g1d7}}
For a general cubic function $y=ax^3+bx^2+cx+d$ we have much more
possibilities and we will not discuss all of them here. The two basic
forms are given by the functions $y=x^3$ and $y=-x^3$ depicted in
Fig.\ref{g1d7}a. Important here is the asymptotic behavior of the
function at $x \rightarrow \pm \infty$. For $y=x^3$ we see that $y$
goes to $+\infty$ when $x$ increases and to $-\infty$ when $x$
decreases; for $y=-x^3$ we have the opposite situation.  A general graph
of $y=ax^3+bx^2+cx+d$ may have up to three zeros that can be found from
the solution of the equation $ax^3+bx^2+cx+d=0$, and up to two extrema
(fig.\ref{g1d7}b). The extrema are points where the derivative of the
function is zero,  which in this case results in the following quadratic
equation: $(ax^3+bx^2+cx+d)'=3ax^2+2bx+c=0$. If the zeros of the
function $(p,q,r)$ are known, the graph can  easily be drawn
as shown in Fig.\ref{g1d7}c.
\end{figure}

\begin{figure}[H]
\centerline{
\psfig{type=pdf,ext=.pdf,read=.pdf,figure=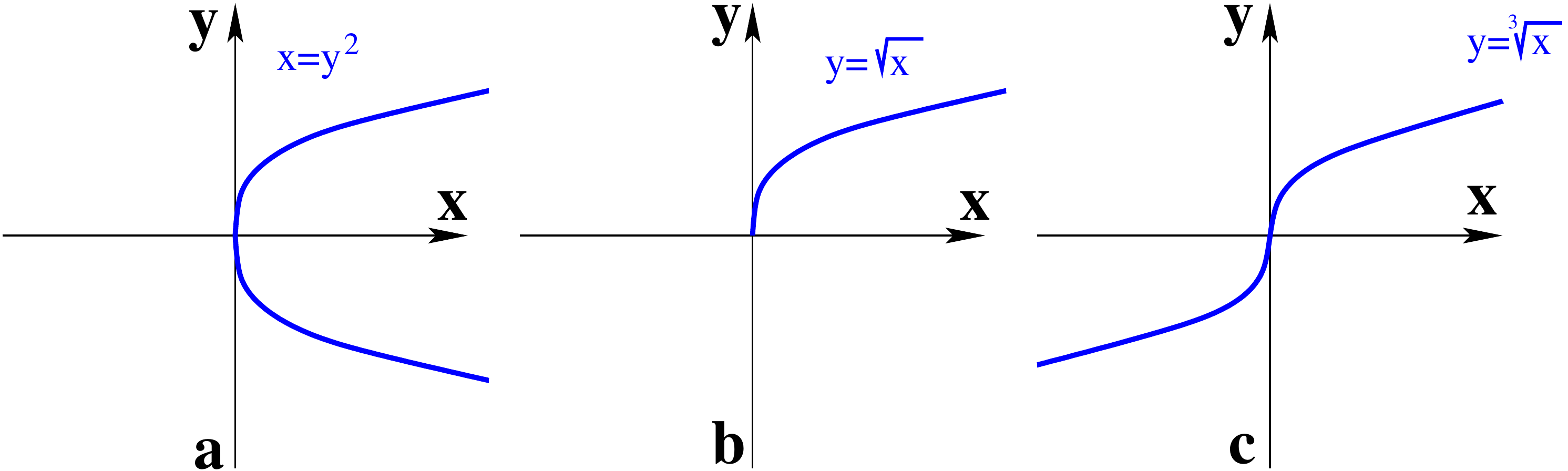,width=15.cm}
}
\caption{\label{g1d8}}
Three examples of graphs of the power function $x^a$ involving
fractional powers are shown in Fig.\ref{g1d8}. If $0<a<1$ than the
graph growth is slower than the function $y=x$ and is concave downward
(in the first quadrant). To draw graph $y=\sqrt{x}$ let us use the
graph of parabola $y=x^2$ discussed in Fig.{g1d5}a.  If in function
$y=x^2$ we switch the $x$ and $y$ we will get $x=y^2$, which is
equivalent to $y=\pm \sqrt{x}$.  The graph $x=y^2$ can be found by
switching the $x$ and the $y$-axis for the graph of the parabola
$y=x^2$ in Fig.\ref{g1d5}a and we get a curve depicted in
fig.\ref{g1d8}a in which the upper branch corresponds to
$y=\sqrt{x}$ (fig.\ref{g1d8}b) and the lower branch corresponds to
$y=-\sqrt{x}$. Similarly, the graph of the function $y=\sqrt[3]{x}$
(Fig.\ref{g1d8}c) can be found by a  $90^o$ rotation  of the graph of
the function $y=x^3$ from Fig.\ref{g1d7}a.
\end{figure}

\begin{figure}[H]
\centerline{
\psfig{type=pdf,ext=.pdf,read=.pdf,figure=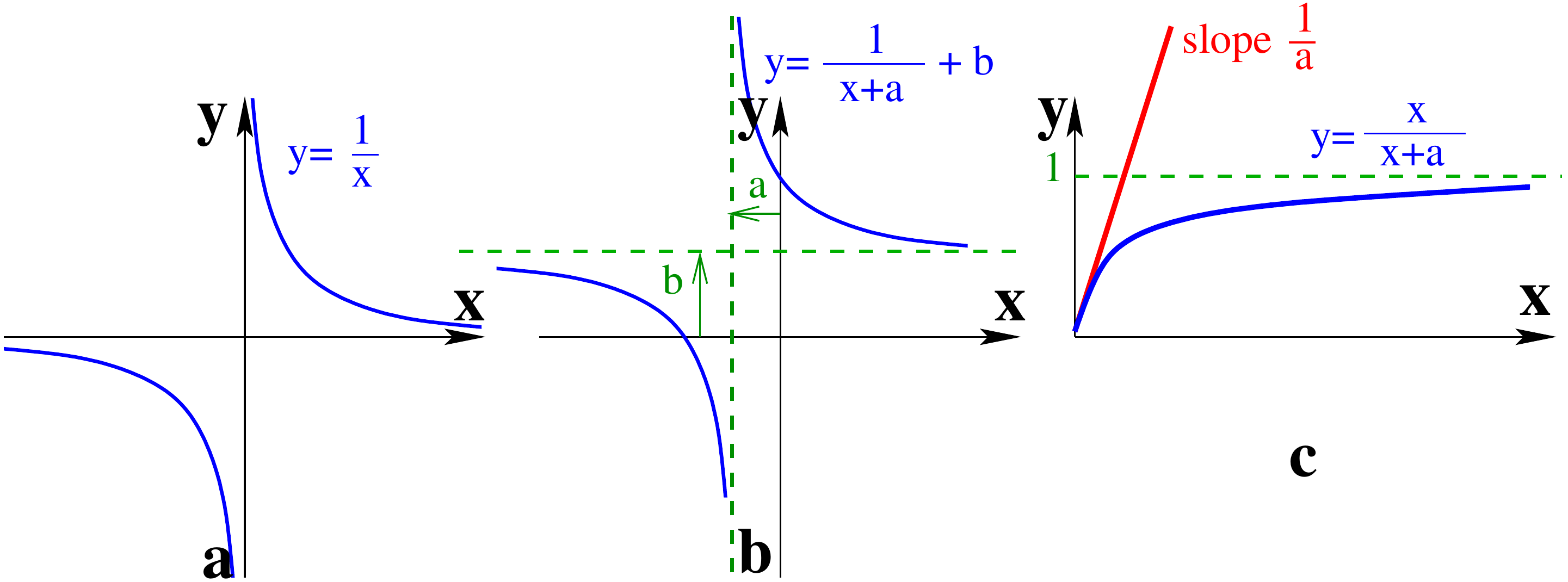,width=15.cm}
}
\caption{\label{g1d9}}
Rational functions ${p(x) \over q(x)}$ are very important in
theoretical biology. The graph of the function $y={1 \over x}$
(Fig.\ref{g1d9}a) has the vertical asymptote ($x=0$) and the
horizontal asymptote ($y=0$). The graph of function $y={1 \over
x+a}+b$ can be obtained by a shift of the graph $y={1 \over x}$ by $b$
units in the $y$ (vertical) direction and by $-a$ units in the $x$
(horizontal) direction. In this case the vertical asymptote ( $x$ at
which function goes to infinity) is $x=-a$, as at this point the
denominator in ${1 \over x+a}$ equals zero. The horizontal asymptote
of this graph is $y=b$, given by $\mathop {\lim }\limits_{x \to \infty
} {1 \over x+a}+b=b$. Another rational function $y={x \over x +a}$
occurs in the classical Michaelis-Menten kinetics.  Fig.\ref{g1d9}c
shows the graph of this function. Because for biological applications
$x$ and $a$ are always considered non-negative ($x \geq 0, a>0$), we
show the graph in the first quadrant only. We see that independent of
the value of the parameter $a$ the horizontal asymptote is always
located at $y=1$, as $\mathop {\lim }\limits_{x \to \infty } {x \over
x+a}=1$. The slope of this function at $x=0$ is given by the function
derivative $f'(x)=({x \over x +a})'={a
\over (x +a)^2}$ at $x=0$, which gives a slope of  $f'(0)={1 \over a}$.
\end{figure}

\begin{figure}[H]
\centerline{
\psfig{type=pdf,ext=.pdf,read=.pdf,figure=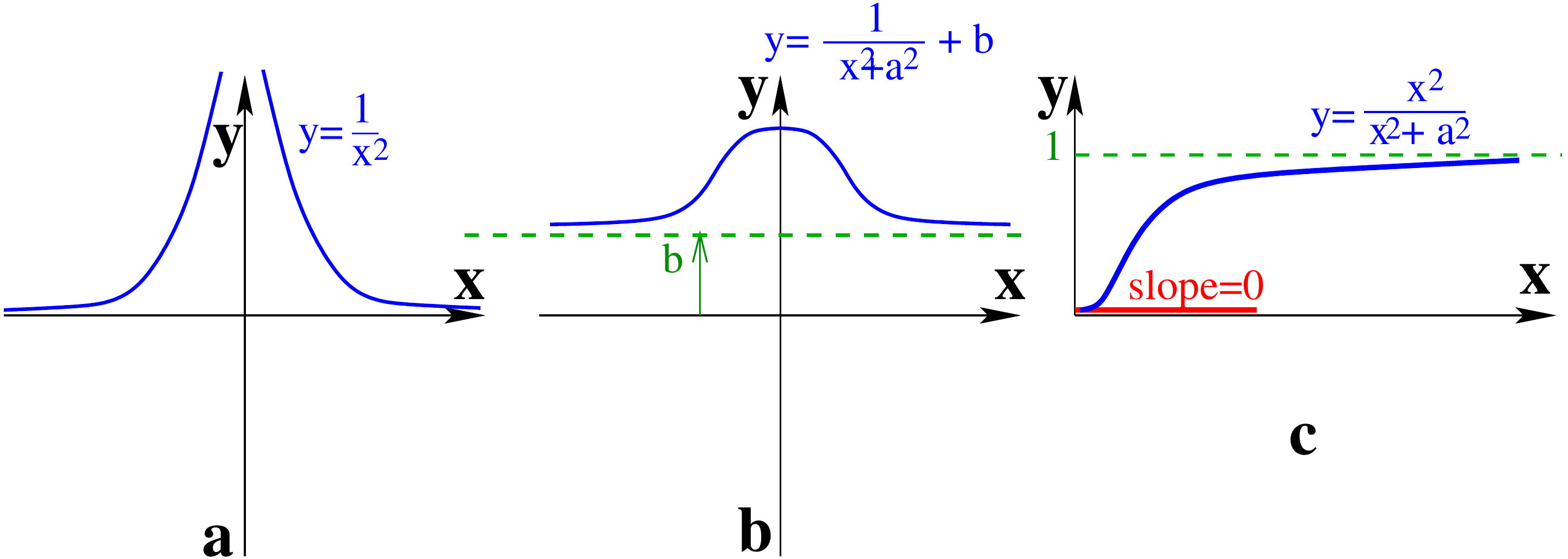,width=15.cm}
}
\caption{\label{g1d10}}
Graphs of similar functions involving a second power: $y={1 \over x^2}$
and $y={1
\over x^2+a^2}$, are shown in Fig.\ref{g1d10}a,b. We see that  function $y={1 \over x^2}$ has a graph similar to that of $y={1 \over x}$  but located in the first and second quadrants, rather than first and third.  One more difference is that  function $y={1 \over x^2+a^2}$ does not have a vertical asymptote, and always reaches a maximum at $x=0$. Function  $y={x^2 \over x^2+a^2}$ is an example of    famous  for its ecological applications Hill function  $y={x^n \over x^n+a^n}$ with $n=2$. Its graph (Fig.\ref{g1d10}c) has a horizontal asymptote at $y=1$ (similar to   $y={x \over x+a}$),  however, the rate of  growth of  $y={x^2 \over x^2+a^2}$  for small $x$ is slower than for   $y={x \over x+a}$:  the slope of the tangent line at $x=0$ here is $0$, which  can be found from the derivative of this function.
\end{figure}
\begin{figure}
\centerline{
\psfig{type=pdf,ext=.pdf,read=.pdf,figure=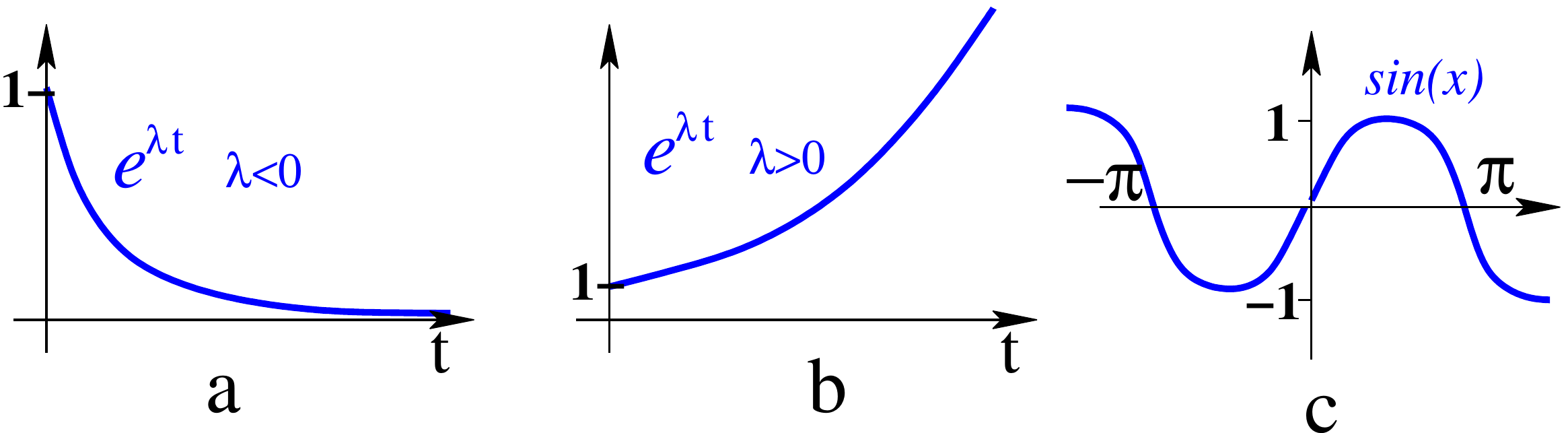,width=12.cm} 
}
\caption{\label{g1d12} }
Finally in fig.\ref{g1d12} we show graphs of two other functions that
are important in this course $e^{\lambda t}$ and $sin(x)$. Note that
if $t$ grows  the function $e^{\lambda t}$ approaches zero if $\lambda<0$ and
diverges to infinity if $\lambda>0$. The function $sin(x)$ oscillates with
a period of $2\pi$ between $-1$ and $+1$.\newline
\end{figure}
\noindent {\bf Tips on graphs}

\noindent Let us list  important rules that may help to plot
graphs of function  $y=f(x)$ with  parameters.

$\bullet$ The graph of the function $y=f(x)+p$ can be obtained by a
vertical shift by $p$ units of the graph of $y=f(x)$.

{\it Example:} In function ${ N \over N+b} +c$,  parameter $c$ just shifts the graph of  ${ N \over N+b}$ by $c$ units above. 

$\bullet$ Important points of the graph are points at which the graph
crosses the $y$-axis ($y$-intercept), given by $y=f(0)$, and points
where the graph crosses the $x$-axis (zeros of the function), given by
$f(x)=0$.  Note, that some graphs do not cross the $x$ or the $y$ axis
and thus do not have $y$-intercepts or zeros. For example graph of
function $f(x)={1 \over x}$ (Fig.\ref{g1d9}a) does not have finite
zeros or $y$-intercepts.

{\it Example:} For function  $f(N)={ N \over N+b} +c,\;b,c>0$, the $y$-intercept is 
${ 0 \over 0+b} +c=c$. Zeros can be found from  ${ N \over N+b} +c=0$, which gives $N  +c (N+b)=0$, or $N  +cN+cb=0$, or $N(1+c)=-cb$, thus zero is given by the formula  $N=-{cb \over 1+c}$, which is always valid as $c >0 $.

$\bullet$ Another important graph feature are asymptotes. To find a
{\it horizontal asymptote} we need to compute the $\mathop {\lim
}\limits_{x \to \infty } f(x)$. For functions without parameters, you
can try to compute this limit using calculator by filling in a large
numbers 10000, 20000, etc and looking if the function approaches some
constant value. For functions with parameter, you can try to fill in
some 'reasonable' parameter value and try to find similarly if the
asymptote exists, however the best way here is to find the limit using
our plan from section \ref{sec_limits}. {\it A vertical asymptote} is
usually a point where a denominator of a fraction is zero.  Not all
graphs have asymptotes, for example graph of function $f(x)=x^2 $ does
not have any vertical or horizontal asymptotes. However, even if the
asymptotes are absent it is still useful to understand behavior of the
functions at large $x$ and show it in the graph.

{\it Example:} For function  ${ N \over N+b} +c$ we can find $\mathop {\lim }\limits_{N \to \infty }$   as    ${ N \over N+b} +c ={ { N \over N} \over { N \over N} +{b \over N}} +c =
{ 1 \over 1 +{b \over N}} +c ={ 1 \over 1 +0} +c =1+c$, thus this graph has a horizontal asymptote $y=1+c$. The vertical asymptote here is at point where $N+b=0$, or  line $N=-b$.

$\bullet$ Several features of the graph can be found from the derivative of the function: 
a function grows if its derivative  $f'(x)>0$, decreases if
$f'(x)<0$ and has a local extremum (maximum or minimum) if  $f'(x^*)=0$.
We do not necessarily need to compute these feature for each graph, but it may be   helpful for some functions.

In many applications we will be interested in
points of intersection of graphs of two functions $f(x)$ and
$g(x)$. Because at the intersection point functions are equal to each
other, such points can be found from  the equation
$f(x)=g(x)$.

\begin{figure}[h]
\centerline{
\psfig{type=pdf,ext=.pdf,read=.pdf,figure=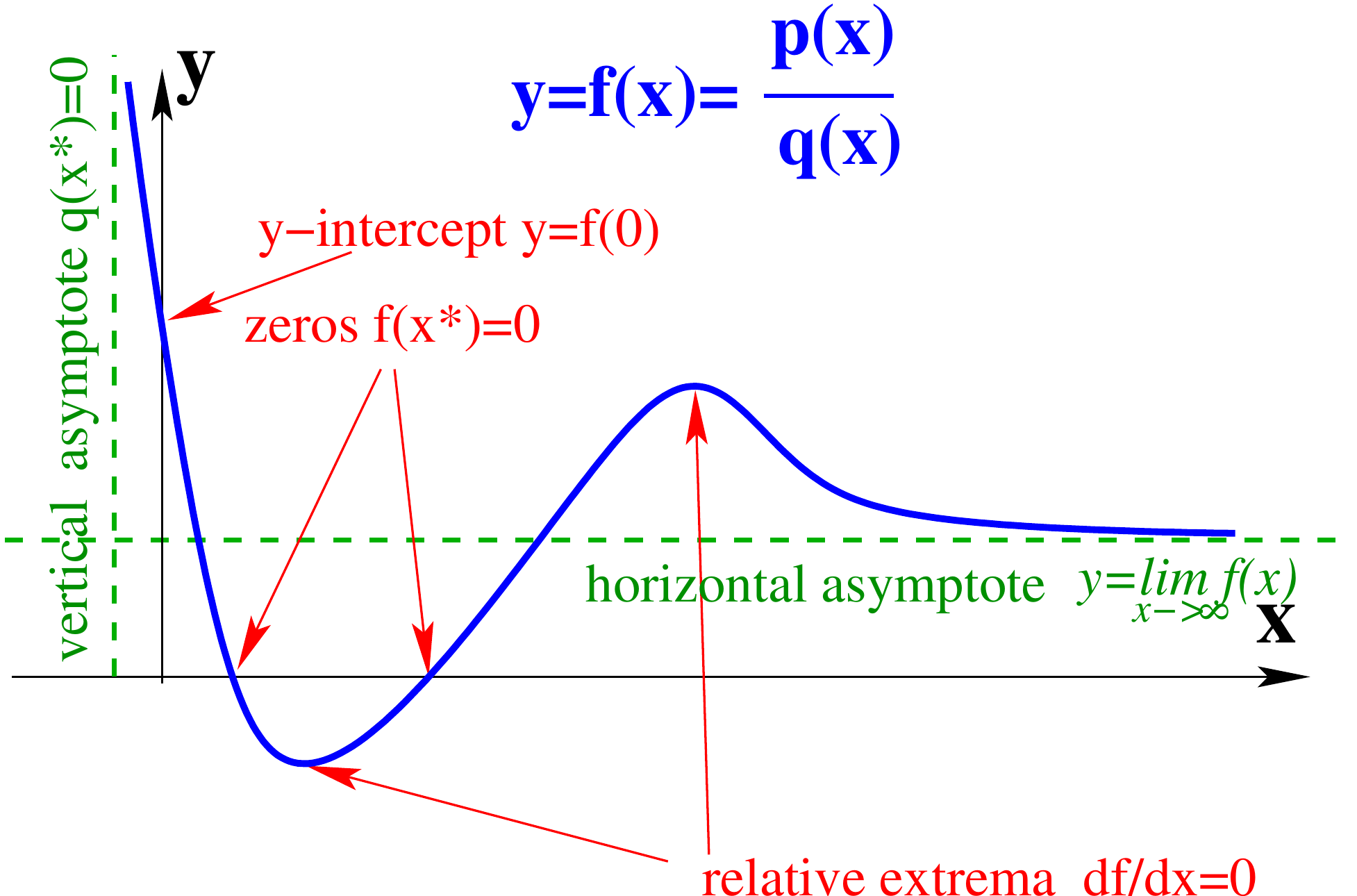,width=10.cm}
}
\caption{\label{g1d11}}
\end{figure}

The above mentioned tips are represented graphically in  fig.\ref{g1d11}.

\noindent Finally let us formulate the main rules for graphing functions with parameters.

\newpage
\noindent {\bf Plan for  graphing functions with parameters}
\ben
\item[1] Try to simplify the function and determine if it  belongs to a known  class of functions with graphs from  (Fig.\ref{g1d1}-\ref{g1d10}).

\item[2] {\sf Computer trail:}
\ben 
\item[(a)] Put  parameters   to 'reasonable' values and plot the graph  using a calculator. 
\item[(b)] Collect  qualitative information   such as :  number of zeros, existence of  vertical  and horizontal asymptotes. 

\item[(c)] Vary parameter values to see how this changes the shape of the graph.
\een
\item[3] {\sf Algebraic approach (note, not all steps may be  possible):}
\ben

\item[(a)] Find  $y$-intercept ($f(0)$), and  zeros of the function ($f(x)=0$).

\item[(b)] Find  horizontal asymptote from the limit $y=\mathop {\lim }\limits_{x \to \infty }    f(x)$ and vertical asymptote(s) (for rational function ${p(x) \over q(x) }$  they are at  the points where the denominator  becomes zero ($q(x)=0$))  ( fig.\ref{g1d11}). 

\item[(c)] Find other special points (e.g. maximum, minimum, etc), if they are important determinants of  the graph shape.

\item[(d)] Draw the graph and indicate  how  the graph shape changes  for  different parameter values.
\een
\een

{\bf Example} Plot the graph of the function $f(x)={ax \over x^2 +
c^2}\;\;x \geq 0 \;\;\; \; a>0\; c>0$. Find how the graph depends on the
parameters $a$ and $c$

{\bf Solution.}  

\ben
\item[1] We do not need to simplify the function.
 The function equation has some similarities with graph classes listed above, but does not coincide exactly with any of them.

\begin{figure}[h]
\centerline{
\psfig{type=pdf,ext=.pdf,read=.pdf,figure=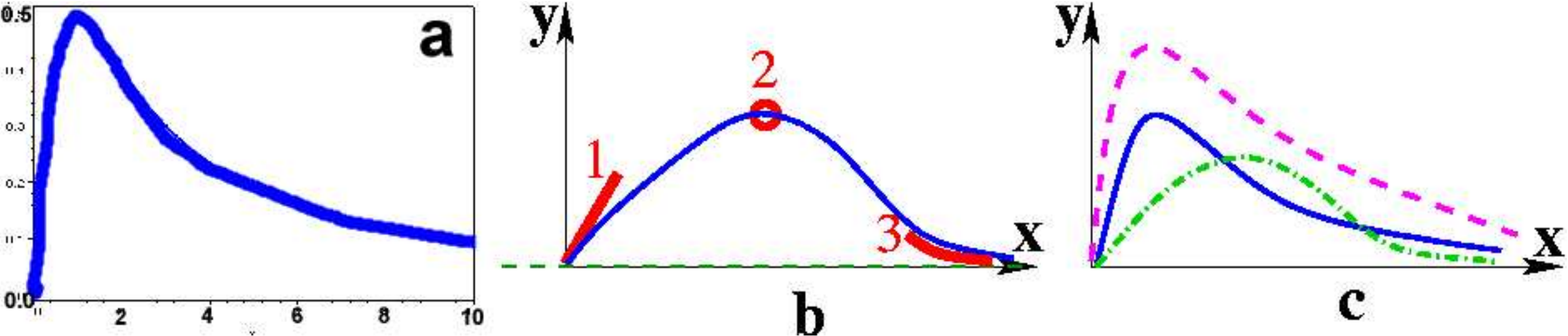,width=14.cm}
}
\caption{\label{exampleGr1}}
\end{figure}

\item[2a] Let us put $a=c=1$ and plot the graph using  calculator (fig.\ref{exampleGr1}a).

\item[2a] The graph (fig.\ref{exampleGr1}a) has the following characteristic features: the $y$-intercept here is $y=0$, we see one zero of the function $x=0$. If we fill in large values of $x$, we find that $f(1000)=0.00099$,$f(5000)=0.000199$ and  $f(10000)=0.00009$, thus we expect to have a horizontal asymptote $y=0$, we do not see any vertical asymptotes and function has an extremum point (maximum).
However, will these features persist for other parameter values? In
order to answer that let us perform an algebraic study.

\item[3a] $y$-intercept is $f(0)={a*0 \over 0^2 +
c^2}=0$. Zeros of the function are given by ${ax \over x^2 + c^2}=0$,
which has only one solution $x=0$, as $x^2 + c^2 \ne 0$ for all $x$.

\item[3b]  The function does not have vertical
asymptotes as the denominator cannot be  zero ($x^2 + c^2 > 0$
for all $x$ and $c>0$). To find the horizontal asymptote let us compute
$y=\mathop {\lim }\limits_{x \to \infty } {ax \over x^2 + c^2}=
\mathop {\lim }\limits_{x \to \infty } {{ a \over x^2} \over {x^2
\over x^2} + {c^2 \over x^2}} = {0 \over 1+0}= 0$, thus the horizontal
asymptote is the $x$-axis.

\item[3c] Important point here is the location of the maximum of the function. Let us find it.  For that let us
find the points where the derivative of the function is zero. The derivative of the function is 
$f'(x)={a*(x^2+c^2)-2x*ax \over (x^2 + c^2)^2}={ac^2-ax^2 \over (x^2 +
c^2)^2}=0$, thus the expression is  zero if $ac^2-ax^2=0$, or $x=\pm
c$. For $x \geq 0,c>0$ we have just one solution $x=c$.  The value of the
function at this extreme point is $f(c)={ac \over c^2 + c^2}={a
\over 2c}$,  thus the maximum is at    $(c,{a \over 2c})$.

\item[3d] Let us draw the  graph now. Because $f(0)=0$ the graph  always goes through the origin. Then the graph will reach the
maximum at $(c,{a \over 2c})$ (Fig.\ref{exampleGr1}a, symbol '1') and
then approaches the $x$ axis. If
we put all this information together we obtain a qualitative graph
shown in fig.\ref{exampleGr1}b. We see that it qualitatively coincides with the
calculator sketch. Now let us find how the graph shape depends on the  parameter values. We see that the
graph  has a bell-shape, with a single maximum at $(c,{a \over
2c})$. The $x$ location of this maximum depends on the parameter $c$
only, but the maximal value of the function increases if  $a$
increases. The solid and the dashed line in Fig.\ref{exampleGr1}c  illustrate how the graph shape changes if $a$ increases while we  keep $c$
constant. Alternatively, if we keep the $a$ value constant but increase the
value of the parameter $c$ (the solid and the dot-dashed line in
fig.\ref{exampleGr1}c), the $x$ location of the maximum shifts to the
right, and the  maximal value of the function ${a \over 2c}$ decreases.
\een

\section{Implicit function  graphs \label{secImplicitGr}}

As we know the relation between two variables $x$ and $y$ can be
expressed explicitly in terms of a function $y=f(x)$ that gives us the
value of $y$ if we know the value of $x$. It is also possible that the
relation between $x$ and $y$ is given implicitly as an equation.  Such
relations are called implicit functions, and their graphs are implicit
function graphs. One of the most effective methods to plot such graphs
is to try to solve that implicit equation and rewrite it as one or
several explicit functions. In some cases the relations between $x$
and $y$ can be plotted directly. Let us consider two examples:

{\bf Example:} Draw  a graph of the function(s) given by equation: $x^2+y^2=C^2$

\begin{figure}[h]
\centerline{
\psfig{type=pdf,ext=.pdf,read=.pdf,figure=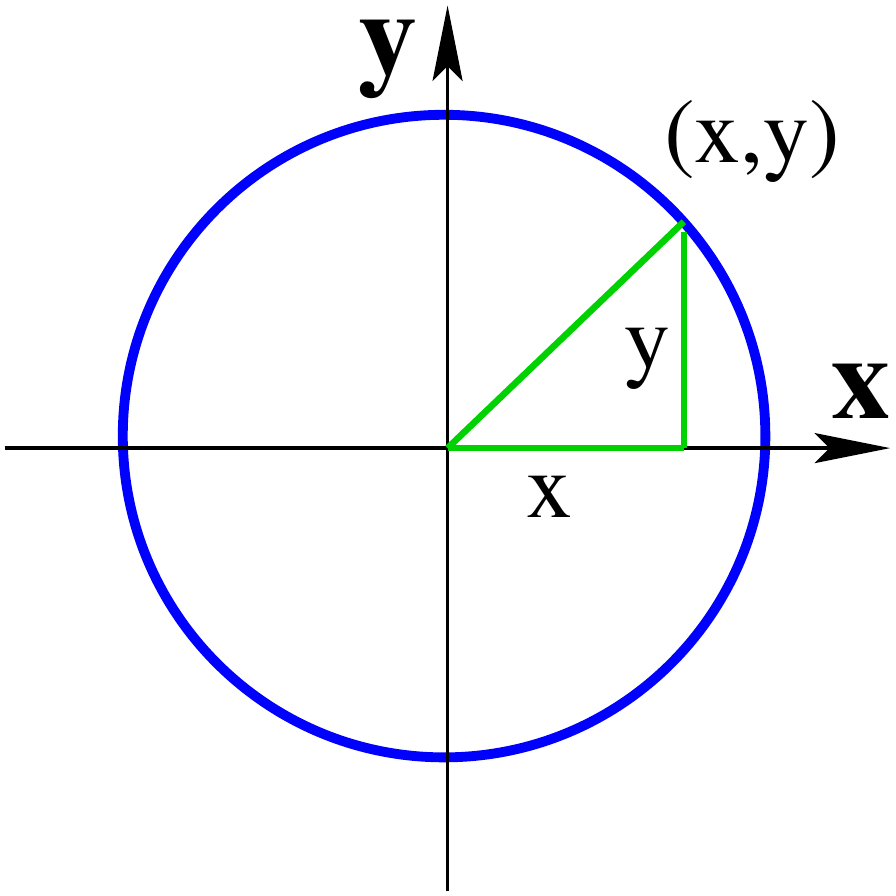,width=4.cm}
\psfig{type=pdf,ext=.pdf,read=.pdf,figure=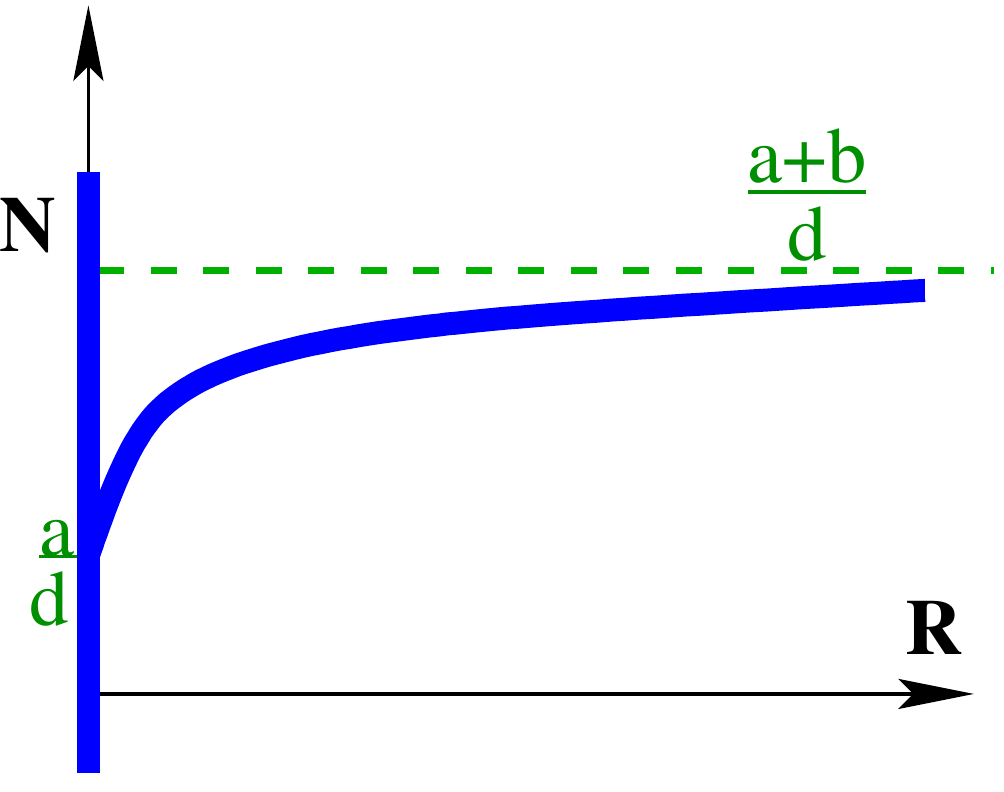,height=4.cm}
}
\caption{\label{circle}}
\end{figure}

{\bf Solution:} We can either rewrite it as two explicit functions  $y=\pm
\sqrt{C^2-x^2}$ and draw the two graphs given by this
equation. Alternatively, we can note that $x^2+y^2$ gives a square of
the distance from the point $(0,0)$ to the point $(x,y)$, thus
equation $x^2+y^2=C^2$ gives the points located at a distance $C$ from
the origin. That is a circle with a radius $C$ with the center at
$(0,0)$ (Fig.\ref{circle}a). We will use this graph later in our
course in chapter \ref{chap_2dlin} to plot fig.\ref{fig7.1}a.

{\bf Example:} Draw a graph of the function(s) given by equation: $aR+b{R^2
\over R+c}-dRN=0$, where $R,N \geq 0$ are the variables  and $a,b,c,d \geq 0$ are the parameters.

{\bf Solution:}
Let us factor the equation:
$$
aR+b{R^2\over R+c}-dRN=R(a+b{R\over R+c}-dN)= 0
$$
The product of two numbers is zero if one of these numbers is zero, therefore this equation is equivalent to:
\beqar
R=0\\
or,\\
a+b{R\over R+c}-dN=0
\eeqar
Graphing of $R=0$ is trivial. In order to graph $a+b{R\over R+c}-dN=0$ 
let us rewrite it as $dN=a+b{R\over R+c}$, or $N={a \over d}+{{b \over d}
R\over R+c}$.  The horizontal asymptote of this graph is: $N= \mathop {\lim }\limits_{R \to \infty }  {a \over d}+{{b \over d} R\over R+c}={a \over
d}+{b/d}={a+b \over d}$. The vertical asymptote occurs if the  denominator of the  fraction
is zero, i.e. at $R=-c$.  However, because $R \geq 0$ and $c>0$
this asymptote will be outside the range of our function. Additionally
note that this function is similar to the graph of Fig.\ref{g1d9}c, but it is shifted  upward by  ${a \over d}$. Thus the function graph here contains two branches  $R=0$ and $N={a \over d}+{b/d
R\over R+c}$  that are plotted  in Fig.\ref{circle}b

\section{Exercises}
\subsection*{Exercises for section \ref{secBasicAlgebra}}

\ben

\item Perform the indicated operations:
\ben
\item $(ax-2by)*(3y-4bx)+2b*(2ax^2+3y^2)-8xyb^2$
\item $ {6 \over r} -{5r \over 30r+5}$
\een

\item Find limits:
\ben
\item $ \mathop {\lim }\limits_{x \to \infty }  {a x + q \over c^2+x^2}  $
\item $ \mathop {\lim }\limits_{N \to \infty }  {a N^2 + q \over {b\over N}+ c^2+ dN^2}, \;d \ne 0 $
\een

\item Solve the equation for the specified variable:
\ben
\item find $r$ in: $3r + 2 - 5(r+1) = 6r+4$
\item find $x$ in: $x + {4 \over x} =4$
\item find $N$ in:  $(b-{N \over k})N=0$
\item find $N$ in:  $(b-d(1+{N \over k}))N=0,\;\;d \ne 0;\;k \ne 0$
\item find $N$ in: $(\frac{b}{1+N/h}-d)N=0, \;b \ne 0; $
\een

\item Solve the system of equations  for the specified variables:
\ben  
\item  find $x,y$ in: $\left\{
\begin{array}{l}
x -2y =-5\\ 2x+ y=10
\end{array}
\right.
$
\item  find $x,y$ in: $\left\{
\begin{array}{l}
 ax+by=0\\ cx+dy=-b
\end{array}
\right.
$
\item  find $x,y$ in: $\left\{
\begin{array}{l}
 x(1-2x)+xy=0\\ 4y-xy=0
\end{array}
\right.
$

\item  find $x,y$ in: $
\left\{
\begin{array}{l}
4x-xy-x^2=0 \\
 9y-3xy-y^2=0 
\end{array}
\right.     
$

\item  find $R,N$ in: $\left\{
\begin{array}{l}
 b(1-{R \over k} -d -aN)R=0\\ (R-\delta)N=0
\end{array}
\right. ,\;\;a,b,d,k,\delta  \ne 0;
$

\een

\subsection*{Exercises for section \ref{sec1d}}

\item Find the derivative of  $f(x)$:
\ben
\item $f(x)={1 \over x^3}$
\item $f(x)=e^{-5x};$
\item $y=(4x-x^2)*(2x+3)$
\item $y={x \over a^2-x^2}$
\een

\subsection*{Exercises for section \ref{sec1dGraph}}

\item Without plotting the function find  the following information  
about their graph: find the $y$-intercept and zeros; find horizontal
or/and vertical asymptotes (if they exist). (Proof of non-existence of  asymptotes  is not required). 
\ben

\item  function $y(n)$ given by $y={a n^2-b \over n^2+c^2}, \; \; a,b,c >  0$

\item  function $N(R)$ given by  $ N={r \over a} (h+R)(1 - {R \over K}), \;r,a,h,K \ne 0$

\een

\item Sketch graphs of the following functions:
\ben
\item $y=3-6x$

\item $y=x-3x^2$

\item $y={3x \over x + a}+4 $.  Find how the shape of the graph for $x,y \geq 0$ depends on the value of the parameter $a>0$.  

\een

\subsection*{Exercises for section \ref{secImplicitGr}}

\item  (a) Sketch qualitative graphs of the following implicit functions. (b) Find how  special points of these graphs (intercepts, zeros, asymptotes) depend on  parameters. (c) If graph contains several lines find their intersection points. {\it Note, all parameters represent  \underline{\bf positive} numbers.}
\ben 
\item  $x+3y^2=0$ on the $xy$ plane 
\item  $y^2+x^2=9$ on the $xy$ plane 
\item $xy=0$ on the $xy$ plane 
\item  $dN(a-P)=0$ on the $NP  $ plane 
\item $dN+{NP \over N+a}=0$ on the $NP $ plane 
\item $dR(b-R)={cRN \over R+a}$ on the $RN $ plane 
\item  $ bR(1-{R \over k} -dR) -aNR=0$ on the $RN $ plane 
\item $ aN+P(1-eP)+bP=0$ on the $NP $ plane 
\een

\subsection*{Additional exercises}
\item Perform the indicated operations:
\ben
\item $((x-2y)*(y-2x)+2y^2)* {1 \over x}$
\item ${a -2 b \over 2p} : {4b - 2a \over \sqrt{p}} $
\een

\item Solve the system of equations  for the specified variables:
\ben  

\item  find $A,B,C$ in: $\left\{
\begin{array}{l}
 rA(1-{A\over K}-AB)=0 \\ AB-dB-BC=0 \\ BC-fC=0
\end{array}
\right. \;\; r,K,d,f>0
$

\een

\item Find the derivative of  $f(x)$:
\ben
\item $f(x)=2^x$
\item  $f(x)=\sqrt{1 \over x^3}$
\item $f(g)=cos(x^2);$
\item $f=(cos(x))^2;$
\item $y=ax*e^{bx} \;a,b>0$
\item $y={x^2-5 \over  2x^2-3x}$
\item $y={a x^2 \over bx -c}, \; a,b,c>0$
\item $y={x \over 1+{x \over d}}, \;d>0$
\item $y={x^n \over x^n + a^n}, \;a>0$
\een

\item Solutions of differential equations
\ben 
\item
Show that function $N(t)=A*e^{3t}$, where $A$ is an arbitrary constant, is a solution  of  the differential equation: $\frac{d N(t)}{dt}=3N$. For that,  compute  derivative of this function and substitute this derivative and the function itself to the equation and show that the left hand side of the equation equals to the right hand side.

\item Show, using the same steps that  the function $N(t)=s(1-e^t)+Ae^{-t}$, where $a$ is an arbitrary constant and $s$ is a parameters, is a  solution of the differential equation $\frac{d N(t)}{dt}=s-N$
\een

\item Assume that $x(t)$ is an unknown function of $t$. For  $f(x)$ listed below find   the following derivatives: ${ df \over  dx}$ and ${df \over dt}$.
\ben
\item $f(x)=x^3$

\item $f(x)=e^{-ax}$

\item Find the expression for  ${df \over dt}$ for an arbitrary $f(x(t))$.

\een

\item Without plotting the function find  the following information  about their  graph:  find the $y$-intercept and zeros;  find  horizontal or/and  vertical asymptotes  (if they exist). 
\ben
\item function $y(x)$ given by $y={x -4\over x^2 -3x+2}$

\item  function $y(x)$ given by  $y= a:   {b \over x^3-c} \;\;a,b,c>0$

\een

\item Sketch graphs of the following functions:
\ben

\item $x=4e^{-3t}$

\item $y=x^2+2x-3$

\item $y={2 \over x+3}$

\item $y={bx^2 \over x^2 + a^2}+4$. Find how the shape of the graph depends on the value of the parameters $a>0$ and $b>0$.  

\item $y={bx^2 \over x^3 + c^3}$, Find how the shape of the graph depends on the value of the parameter $c>0$ and $b>0$.

\item $f(n)=rn*(1- {n \over k}) -h$,  $k \ne 0$ find for which values of the parameter $h>0$ the graph  touches the $n$-axis. {\it Tip: draw graph for $h=0$ and think about how $h$ affects this graph.}

\een

\item  Sketch qualitative graphs of the following pairs of implicit functions on the same graph.  Find  all intersection points.
\ben 
\item[(a)] $
aN-P(1+eP)-bP=0\quad and \quad  bP-c N=0 \quad parameters \;\; a,b,c,e >0
$

\item[(b)] $
rR(1-R/K)-\frac{NR}{h+R}=0 \quad and \quad 
\frac{NR}{h+R}-dN =0 \quad parameters\;\; r,K,h,d>0
$

\een

\een


\chapter{Selected topics of calculus  \label{chap_additional}}

In this chapter we introduce several new notions on calculus and
algebra which are important for our course.

\section {Complex numbers \label{sec_complex}}
Complex numbers were introduced for the solution of algebraic equations. It turns out
that in many cases we can not  find the solution of even very simple quadratic equations. Consider the general quadratic equation:
\beq
\label{eq7.1}
{\lambda}^2+B\lambda +C=0
\eeq
The roots of (\ref{eq7.1}) are given by the well known 'abc' formula:
\beq
 \lambda_{12}={-B \pm  \sqrt {B^2-4C} \over 2}={-B \pm  \sqrt {D} \over 2}
\eeq
where
\beq 
D=B^2-4C
\eeq
What happens with this equation if $D<0$? Does the equation have roots in this case? 

Complex numbers help to solve such kind of problems. The first step is to consider the equation
\beq
\label{eq7.4}
\lambda^2=-1
\eeq
Let us claim that (\ref{eq7.4}) has a  solution and denote it in the following way:
\beq
\lambda_{12}=\pm i
\eeq
where
\beq
i=\sqrt {-1}
\eeq
Here $i$ is the basic complex number which is similar to $'1'$ for
real numbers. Using it we can denote solutions of other similar
equations. For example if $$\lambda^2=-4, \;\lambda=\sqrt
{-1*4}=\sqrt{-1} \sqrt {4}=i*(\pm2)=\pm2i.$$ Similarly the equation
$\lambda^2=-a^2$, has solutions $\lambda=\pm ai$. Although we call
$ai$ a complex {\it number}, it is quite different from usual real
numbers. Using complex numbers $ai$ we cannot count how many
books we have in the library, for example. The only meaning of $i$ is
that  $i^2=-1$, and $ai$  is just a designation of a root
of the equation $\lambda^2=-a^2$.

Now we can solve equation (\ref{eq7.1}) for the case $D<0$.
If $D<0$, then $\sqrt {D}=i \sqrt{ -D}$ and
\beq
\label{eq7.7}
\lambda_{12}={-B \pm  i \sqrt {-D} \over 2} 
\eeq   

{\bf {Example}}. Solve the equation $\lambda^2+2\lambda+10=0$

{\bf {Solution}}.
\beq
\label{example1Ch2}
\lambda_{12}={-2 \pm  \sqrt {4-4*10} \over 2}=
{-2 \pm  \sqrt{ -36} \over 2}={-2 \pm  6i \over 2},
\eeq
 or 
$\lambda_{1}=-1+3i,\;\;\;\lambda_{2}=-1-3i. $

We see, that solution of this equation $\lambda_{1,2}$ has two parts, one part is just a real number '-1', which is the same for $\lambda_{1}$ and $\lambda_{2}$ and the other part,  is $i$ times another real number '3' which has opposite signs for $\lambda_{1}$ and $\lambda_{2}$. This is a general form of representation of complex number. Any complex number can be represented in the form:
\beq
z=\alpha+i\beta
\eeq
where $\alpha$ is called {\bf { the  real part}} of the complex number $z$, and $\beta$ is called  {\bf { the imaginary part}}  of  $z$. The notation for the real part is $Rez$ and  for the imaginary part is $Imz$. In our example $Re \lambda_1=-1; Im \lambda_1=3.$ and $Re \lambda_2=-1; Im \lambda_2=-3.$

We can work with complex numbers in the same way as with usual real numbers and expressions. The only thing which we need to remember, is that $i^2=-1$.

To add two complex numbers we need to add their real and imaginary parts. For example  
$$
z_1=3+10i,z_2=-5+4i,\;\; z_1+z_2=(3+10i)+(-5+4i)=3+10i+-5+4i=-2+14i.
$$

Similarly, multiplication by a real number results in  multiplication of the real and imaginary part by this number
$$
z_1=3+10i;\; 10z_1=10*(3+10i)=30+100i.
$$

Multiplication of two complex numbers is just an  exercise in multiplication of two expressions
$z_1=3+10i,z_2=-5+4i;\; z_1*z_2= (3+10i)*(-5+4i)= 3*(-5)+3*4i+10i*(-5)+10i*4i= -15+12i-50i+40i^2= -15-38i-40\; (as\; i^2=-1) = -55-38i.
$

Similarly 
$$
(z_1)^2=(3+10i)^2= 3^2+2*3*10i+(10i)^2= 9+60i+100i^2= 9+60i-100= -91+60i.
$$

Now we can check that $\lambda_{1}=-1+3i$ is a solution of the
equation in example (\ref{example1Ch2}). In fact:
$\lambda^2+2\lambda+10= (-1+3i)^2+2*(-1+3i)+10=
(-1)^2+2*(-1)*3i+(3i)^2-2+6i+10=1-6i-9-2+6i+10=(1-9-2+10)-6i+6i=0-0i=0$,
i.e. left hand side of this equation after substitution of
$\lambda_{1}=-1+3i$ equals zero and thus $\lambda_{1}=-1+3i$ is the
root of this equation.

One more definition. The number $z_2=a-ib$ is called {\bf the  complex conjugate} to the number $z_1=a+ib$  and is denoted as $\bar z_1=z_2=a-ib$. Complex conjugate numbers have the same real parts, but their imaginary parts have opposite signs.

Roots of a quadratic equation  with negative discriminant  $D<0$  are complex conjugate to each other. It follows from the formula (\ref{eq7.7}) 
\beq
\lambda_{1}={-B +   i \sqrt {-D} \over 2} \quad 
\lambda_{2}={-B -   i \sqrt {-D} \over 2}
\eeq
hence:
\beq
Re \lambda_{1}=Re \lambda_{2}={-B \over 2};
\quad Im  \lambda_{1}=  { \sqrt {-D} \over 2};\quad  Im \lambda_{2} =  -{ \sqrt {-D} \over 2}.
\eeq
Finally consider two more basic operations. If $z=a+ib$, then, $|z|=\sqrt {a^2 + b^2}$ is called {\bf {the absolute value}}, or  {\bf {modulus}} of $z$. Note, that 
${|z|}^2=z \bar{z}$, as $ (a+ib)*(a-ib)=a^2-(ib)^2=a^2+b^2$.

We use this trick to introduce division of two complex numbers
$$
{z_1 \over z_2}={z_1 \bar {z_2} \over z_2 \bar {z_2}} 
$$
So, to divide two complex numbers we multiply the numerator and the denominator  by a number which is the complex conjugate to the denominator, and we get the answer in the usual form.

{\bf {Example}}  $$
{1+3i \over 1-4i }= {1+3i \over 1-4i }* {1+4i \over 1+4i }=
 {(1+3i)(1+4i) \over 1^2 + 4^2  }={1 +3i+4i+12i^2  \over 17}= {-11+7i  \over 17}=
{-11  \over 17}+{7  \over 17}i \; 
$$

\section{Matrices \label{basicMatr0}}

From a very general point of view a matrix is a representation of data in the form of a rectangular table. An example of a matrix composed of numbers is given below:\beq
\label{matr1}
A=\left(\begin{array}{lcr}  1 & 4 & 5\\ 2 & 6 &  10 \end{array} \right)
\eeq
This matrix $A$ has two rows and three columns. We will call this
 a  matrix of the size $2 \times 3$. In general  matrix
size is defined as $number\;of\;rows \times number
\;of\;columns$. Even if you did not have  matrix algebra in school,
you probably know at least one matrix object, that is a vector. Indeed,
a vector is an object which is characterized by its
components: two numbers in two dimensions or  three numbers in three
dimensions.  In matrix algebra vectors can be represented in two
forms: as a column vector, i.e. as $n \times 1$ matrix (preferred representation), or as a row
vector, i.e. as $1 \times n$ matrix. For example a vector $V$ with the $x$-component $V_x=2$ and the $y$-component $V_y=1$ can be represented a column or as a row vector as:
$$
\vec{V}=\left(\begin{array}{c}  2\\ 1 \end{array} \right) \hskip 3pc or \hskip 3pc \vec{V}=\left(\begin{array}{lr} 2 & 1 \end{array} \right)
$$ 
Using matrices  we can perform the same operations on  large blocks of data simultaneously. For example, if we need to  multiply all 6 numbers of  matrix  $A$ in (\ref{matr1}) by  $4$, we can write it as $4A$ which will mean:
\beq
\label{matr2}
4A=4*\left(\begin{array}{lcr}  1 & 4 & 5\\ 2 & 6 &  10 \end{array} \right)
=\left(\begin{array}{lcr}  4*1 & 4*4 & 4*5\\ 4*2 & 4*6 &  4*10 \end{array} \right)=\left(\begin{array}{lcr}  4 & 16 & 20\\ 8 & 24 &  40 \end{array} \right)
\eeq
This operation is called {\it   multiplication of a matrix by a number}. For a general $2 \times 2$ matrix this can be written as:
$$
\lambda \left(\begin{array}{lr}  a & b\\ c &  d \end{array} \right)
=
\left(\begin{array}{lr}  \lambda a & \lambda b\\ \lambda c &  \lambda d \end{array} \right)
$$ Similarly, {\it addition of matrices} is adding the numbers that
have the same location. This operation is defined only for two
matrices of the same size:
\beq
\label{matr3}
A+B=\left(\begin{array}{lcr}  1 & 4 & 5\\ 2 & 6 &  10 \end{array} \right)
+
\left(\begin{array}{lcr}  2 & 1 & 4\\ 1 & 3 &  5 \end{array} \right)
=\left(\begin{array}{lcr}  1+2 & 4+1 & 5+1\\ 2+1 & 6+3 &10+  5 \end{array} \right)
=\left(\begin{array}{lcr}  3 & 5 & 6\\ 3 & 9 & 15 \end{array} \right)
\eeq
 For general $2 \times 2$ matrices it  can be written as:
\beq
\label{matr4}
\left(\begin{array}{lr}  a & b\\ c &  d \end{array} \right)
+
\left(\begin{array}{lr}  x & y\\ z &  w \end{array} \right)
=
\left(\begin{array}{lr}  a+x & b+y\\ c+z &  d+w \end{array} \right)
\eeq
Multiplication of matrices is not so trivial.  In general {\it matrix
multiplication} is defined as the products of the rows of the first
matrix with the columns of the second matrix.  Thus, to fund the
element in row $i$ and column $j$ of the resulting matrix we need to
multiply the $i$th row of the first matrix by the $j$th column of the second matrix. Thus we can multiply two
matrices $A*B$ only if the number of columns in matrix $A$ equals the
number of rows in matrix $B$.

For a product of two $2 \times 2$ matrices this  gives:
\beq
\label{matr8}
\left(\begin{array}{lr}  a & b\\ c &  d \end{array} \right) 
\left(\begin{array}{lr}  x & y\\ z &  w \end{array} \right)
=
\left(\begin{array}{lr}  ax+bz & ay +bw\\ cx+dz & cy+ dw \end{array} \right)
\eeq

From this it follows that multiplication of a matrix by a column
vector  is given by:
\beq
\label{matr5}
 \left(\begin{array}{lr}  a & b\\ c &  d \end{array} \right)
\left(\begin{array}{c}  v_x \\ v_y\end{array} \right)
=
\left(\begin{array}{c}  av_x + bv_y\\ c v_x + d v_y\end{array} \right)
\eeq

The last equation is useful for representation of linear systems as  can be seen  from the following example.  Assume we have a system of linear  equations:
\beq
\label{matr6}
\left\{
\begin{array}{l}
x -2y =-5\\ 2x+ y=10
\end{array}
\right.
\eeq
we can write the coefficients at $x$ and $y$  in the left hand side   as a square matrix:
$$
A= \left(\begin{array}{lr}  1 & -2\\ 2 &  1 \end{array} \right).
$$
We also have two numbers in the right hand side which we can write as a column vector:
$$
\vec{V}= \left(\begin{array}{c}  -5 \\ 10 \end{array} \right).
$$
Now if we write  $x$ and $y$ as a column vector:
$$
\vec{X}= \left(\begin{array}{c}  x \\ y \end{array} \right).
$$
we can represent  system (\ref{matr6}) using matrix multiplication (\ref{matr5}) as:
\beq
\label{matr7}
 A\vec{X}=\vec{V}; \; \; or \; \;  \left(\begin{array}{lr}  1 & -2\\ 2 &  1 \end{array} \right)
\left(\begin{array}{c}  x \\ y\end{array} \right)
=
\left(\begin{array}{c}  -5\\ 10\end{array} \right)
\eeq
Indeed, from (\ref{matr5}) we get $\left(\begin{array}{lr}  1 & -2\\ 2 &  1 \end{array} \right)
\left(\begin{array}{c}  x \\ y\end{array} \right)=\left(\begin{array}{c}  1*x  -2*y\\ 2*x  + 1*y \end{array} \right)$, that proves this result.

Another important matrix operation is the {\it determinant of a square
2x2 matrix}, which for the matrix $A=\left(\begin{array}{lr} a & b\\
c & d \end{array} \right)$ is defined as:
\beq
\label{matr9}
 det \left|\begin{array}{lr}  a & b\\
                         c &  d
\end{array} \right| =ad-cb
\eeq
The determinant of a matrix has many important applications in
algebra. For example using determinants it is possible to find solution of 
system of linear equations (e.g. system (\ref{matr6})) in the form of
so-called Cramer's rule, which was published by Gabriel Cramer as
early as in mid-18th century. Cramer's rule is briefly formulated in
exercise \ref{CramersRule} at the end of this chapter.

Now let us consider  one of the most important problems in  matrix algebra: the eigen value problem.

\section{Eigenvalues  and eigenvectors \label{basicEigen}}
Let us start with a definition:
\begin{D}
A  nonzero vector $\mathbf{v}$ and number $\lambda$ are called an eigen vector and an eigen value of a square matrix $A$ if they  satisfy    equation:
\beq
\label{evect}
A\mathbf{v}=\lambda\mathbf{v} 
\eeq
\end{D}

Eigen vectors are not unique, and it is easy to see that if we
multiply  it  by  an arbitrary constant $k$ we get another eigen vector corresponding to the same eigen value. Indeed by multiplying (\ref{evect}) by  $k$ we  get:
\beq
\label{evect1}
kA\mathbf{v}=k\lambda\mathbf{v} \quad or \quad A (k\mathbf{v})=\lambda (k\mathbf{v})
\eeq
therefore, we can say that   $k\mathbf{v}$ is also an eigen vector of (\ref{evect}) corresponding to eigen value $\lambda$.

For example, for  matrix $A=\left(\begin{array}{lr}  1 & 2\\
                         2 &  1
\end{array} \right)$,  number $\lambda=3$ and vector $\mathbf{v}=\left(\begin{array}{c}  1\\1 \end{array} \right)$ are an  eigen value and eigen vector as:

\beq
\label{exampleEV}
A\mathbf{v} =\left(\begin{array}{lr}  1 & 2\\  2 &  1 \end{array} \right) \left(\begin{array}{c}  1\\1 \end{array} \right) =
\left(\begin{array}{c}  1*1+2*1 \\  2*1 + 1*1 \end{array} \right)= \left(\begin{array}{c}  3\\3 \end{array} \right)=3\left(\begin{array}{c}  1\\1 \end{array} \right)=3\mathbf{v}
\eeq
If we multiply
$\mathbf{v}=\left(\begin{array}{c} 1\\1 \end{array} \right)$ by any
number, e.g. $2$, $25$, or etc., we will get new eigen vectors
$\mathbf{v}=\left(\begin{array}{c} 2\\2 \end{array}
\right)$,  $\mathbf{v}=\left(\begin{array}{c} 25\\25 \end{array}
\right)$ of this matrix for $\lambda=3$. You can check it in the same way as we did  in (\ref{exampleEV}) for a vector $\mathbf{v}=\left(\begin{array}{c} 1\\1 \end{array} \right)$.

Finding eigen values and eigen vectors is one of the most important
problems in applied mathematics. It arises in many biological
applications, such as population dynamics, biostatistics,
bioinformatics, image processing and many others.  In our course we
will apply it for the solution of systems of differential equations,  which
we will consider in chapter \ref{chap_2dlin}.

Let us consider how to  solve the eigen value problem for a 2x2 matrix
$A=\left(\begin{array}{lr} a & b\\ c & d \end{array} \right)$. For that  we need
to find $\lambda$ and $\left(\begin{array}{c} v_x\\v_y \end{array}
\right)$ satisfying:
\beq
\label{EVproblem}
A\mathbf{v}=\left(\begin{array}{lr} a & b\\ c & d \end{array} \right)
\left(\begin{array}{c} v_x\\v_y \end{array} \right) = \lambda
\left(\begin{array}{c} v_x\\v_y \end{array} \right).
\eeq
We can rewrite it as a system  of two equations with three unknowns $\lambda,v_x,v_y$:
\beq
\label{eveq02}
\left\{
\begin{array}{l}
 a*v_x+ b*v_y =\lambda v_x \\ c*v_x+d*v_y =\lambda v_y
\end{array}
\right.
\eeq
If we collect all unknowns at  the left hand side we will get the  following system:
\beq
\label{eveq2}
\left\{
\begin{array}{l}
 (a-\lambda)*v_x+ b*v_y =0\\ c*v_x+(d-\lambda)*v_y =0 
\end{array}
\right. \quad or \; in\;matrix\;form\; \left(\begin{array}{lr} a -\lambda & b\\ c & d -\lambda\end{array} \right)
\left(\begin{array}{c} v_x\\v_y \end{array} \right) = \left(\begin{array}{c} 0\\0\end{array} \right).
\eeq

This system always has a solution $v_x=v_y=0$, however it is not an
eigen vector, as in accordance with the definition the eigen vector
should be nonzero. In order to find non-zero solutions let us multiply
the first equation by $d-\lambda$, the second equation by $b$ and
subtract them. Multiplication gives:

\beq
\label{eveq22}
\left\{
\begin{array}{l}
(d-\lambda)*[ (a-\lambda)*v_x+ b*v_y]=0\\ b*[c*v_x+(d-\lambda)*v_y] =0 
\end{array}
\right.
\eeq
Subtraction of the equations results in:
$$
\begin{array}{l}
  (d-\lambda)*(a-\lambda)*v_x+ (d-\lambda)*b*v_y =0\\
-\\
b*c*v_x+b*(d-\lambda)*v_y =0 \\
gives\\
 (d-\lambda)*(a-\lambda)*v_x- b*c*v_x+ (d-\lambda)*b*v_y -b*(d-\lambda)*v_y =0\\ 
or
\end{array}
$$
\beq
\label{eveq1}
[(d-\lambda)*(a-\lambda)- b*c]*v_x=0 
\eeq
as $v_x \ne 0$ we get:
\beq
\label{charEQ}
(d-\lambda)*(a-\lambda)- b*c=\lambda^2 - (a+d) \lambda + (ad -cb)=0
\eeq

This is a quadratic equation with unknown $\lambda$ and for each
particular coefficients $a,b,c,d$ we can find two solutions:
$\lambda_1$ and $\lambda_2$ using the 'abc' formula. Thus we found that the eigenvalue problem for a $2x2$ matrix (\ref{EVproblem}) has  solutions for  the eigen values  $\lambda$. In general,  for a $nxn$ matrix  that the eigen value problem  has $n$ solutions for  $\lambda$.

Equation (\ref{charEQ}) is very important in our course and it
has a special name: {\bf characteristic equation}. In most of
the courses on mathematics this equation, however,  is written in a slightly different
matrix form. To derive it let us  recall the definition  of   the determinant of a matrix given in section \ref{basicMatr0}:\\
$det\left|\begin{array}{lr}  a & b\\
                         c &  d
\end{array} \right| =ad-cb.
$
 Similarly the determinant of   matrix $\left(\begin{array}{lr} a-\lambda & b\\ c & d-\lambda
\end{array} \right)$ is:
\beq
det
\left|\begin{array}{lr} a-\lambda & b\\ c & d-\lambda
\end{array} \right|=(a-\lambda)(d-\lambda)-bc 
\eeq
which coincides with the left hand side of characteristic equation (\ref{charEQ}) and thus  the characteristic  equation  can be rewritten
as:
\beq
\label{eigen_char}
det
\left|\begin{array}{lr} a-\lambda & b\\ c & d-\lambda
\end{array} \right|=0
\eeq

Let us use this approach  to find the eigen values of  matrix $A$ from example (\ref{exampleEV}).
We get  the following characteristic equation:
$$
\begin{array} {l}
 Det
\left|\begin{array}{lr} 1-\lambda & 2\\ 2 & 1-\lambda
\end{array} \right|=(1-\lambda)(1-\lambda)-2*2\\
=1-\lambda-\lambda+\lambda^2-4=\lambda^2-2\lambda-3=0
\end{array}
$$
From the 'abc' formula:
$$
 \lambda_{1,2}={2 \pm  \sqrt {4+12} \over 2}={2 \pm  \sqrt {16} \over 2}; \qquad \lambda_1=3 \quad \lambda_2=-1
$$
therefore we found two eigen values  $\lambda_1=3$ and  $\lambda_2=-1$.

Now, let us find eigen vectors.  For that let us substitute the found
eigen values to the original equation (\ref{eveq2}) and solve it for $v_x$ and $v_y$. Let us do it
first for a particular example (\ref{exampleEV}) for which we have
found eigen values $\lambda_1=3$ and $\lambda_2=-1$.  For eigen vector
corresponding to eigen value $\lambda_1=3$ we obtain:

\beq
\label{exev11}
\left\{
\begin{array}{l}
(1-3)v_x+2v_y=  0 \\ 2v_x+(1-3)v_y=0
\end{array}
\right. \quad or \quad 
\left\{
\begin{array}{l}
-2v_x+2v_y=  0 \\ 2v_x-2v_y=0
\end{array}
\right. \quad or \quad 
\left\{
\begin{array}{l}
-2v_x= -2v_y  \\ 2v_x=2v_y
\end{array}
\right. 
\eeq

\noindent Both equations give the same solution $v_x=v_y$. This means that if
$v_y=1$, then $v_x=1$ and a pair $\left( \begin{array}{c} 1 \\ 1
\end{array} \right)$  satisfies the system and  thus gives an eigen vector of
problem (\ref{exampleEV}).  We can also use any other value for
$v_y$. For example,  if we use $v_y=2$ then  $v_x$ will be $v_x=2$ and we get another eigen
vector $\left(
\begin{array}{c} 2 \\ 2 \end{array} \right)$, etc.
In general  any $v_x=k$, and  $v_y=k$  give an eigen vector. We can express it by  the following formula:
\beq
\label{tmp14}
 \left( \begin{array}{c} v_x
 \\  v_y  \end{array} \right)=
k\left( \begin{array}{c} 1
 \\  1  \end{array} \right)
\eeq
where $k$ is an arbitrary number. Formula (\ref{tmp14}) gives all
possible solutions of eq.(\ref{exev11}).  It also
illustrates  a general property of eigen vectors which we have
proven in (\ref{evect1}), that if we multiply an eigen vector by an
arbitrary number $k$ will get also an eigen vector of our matrix.
Using this property we can formulate an easy way to write a formula
for all eigen vectors. For that we take any found eigen vector and
multiply it by an arbitrary number $k$. Note, that if for problem
(\ref{exev11}) we use another found eigen vector $\left(
\begin{array}{c} 2 \\ 2 \end{array} \right)$, we can write an answer   as 
$\left( \begin{array}{c} v_x
 \\  v_y  \end{array} \right)=
k\left( \begin{array}{c} 2
 \\  2  \end{array} \right)$.
At the  first glance this formula is different from (\ref{tmp14}). However, it is easy to see that both formulas give  the same
result:  this is because $k$ in (\ref{tmp14}) is an arbitrary constant and any
vector given by the formula (\ref{tmp14}) with $\left( \begin{array}{c} 1 \\ 1
\end{array} \right)$ can be obtained using the formula  $k \left( \begin{array}{c} 2 \\ 2 \end{array} \right)$ for  another value of $k$. Thus the answer to our problem: to find eigen vectors of matrix (\ref{exampleEV})  for eigen value  $\lambda_1=3$, can be written as  $\left( \begin{array}{c} v_x
 \\  v_y  \end{array} \right)=
\left( \begin{array}{c} 1
 \\  1  \end{array} \right)$, or  $\left( \begin{array}{c} v_x
 \\  v_y  \end{array} \right)=
\left( \begin{array}{c} 2
 \\ 2 \end{array} \right)$, or etc. These vectors give particular
 solutions of this problem. We can also write a formula for all
 solutions as $\left( \begin{array}{c} v_x \\ v_y \end{array} \right)=
 k\left( \begin{array}{c} 1 \\ 1 \end{array} \right)$, or $\left(
 \begin{array}{c} v_x \\ v_y \end{array} \right)= k\left(
 \begin{array}{c} 2 \\ 2 \end{array} \right)$, or etc. As we discussed
 above all these answers will be correct and equivalent.

\noindent Similarly we find the eigen vector corresponding to the other eigen value  $\lambda_2=-1$:

\begin{enumerate}
\item {\it Substitution}:
\beq
\label{exev22}
\left\{
\begin{array}{l}
(1-(-1))v_x+2v_y= 0 \\ 2v_x+(1-(-1))v_y=0
\end{array}
\right. \quad or \quad
\left\{
\begin{array}{l}
2v_x+2v_y=  0 \\ 2v_x+2v_y=0
\end{array}
\right. \quad or \quad 
\left\{
\begin{array}{l}
2v_x=-2v_y\\ 2v_x=-2v_y
\end{array}
\right.
\eeq

\item {\it Relation between $v_x$ and $v_y$}: 
$$ 
v_x=-v_y
$$
\item  {\it Eigen vector}: use e.g. $v_y=1$, thus  $v_x=-1$

$$
\mathbf{v}= \left( \begin{array}{c} -1 \\ 1 \end{array} \right)
$$
The general form is $
\mathbf{v}=k\left( \begin{array}{c} -1 \\ 1 \end{array} \right)
$, where $k$ is an arbitrary number.

\end{enumerate}

Note, that in both cases in order to find eigen vectors we could use
the first equation only (see equations (\ref{exev11}) and
(\ref{exev22})), and the second equation in both cases did not provide
us any new information. It is not a coincidence, and this property is
the basis for the  following express
method for finding eigen vectors:

\subsection*{Express method for finding eigen vectors}

Let us derive a formula for finding the eigen vectors of a general system
(\ref{eveq02}). We assume that we have found eigen values $\lambda_1$
and $\lambda_2$ from the characteristic equation (\ref{eigen_char}). To
find the corresponding eigen vectors we need to substitute the found
eigen values into the matrix and solve the following system of linear
equations (\ref{eveq2}):
\beq
\label{eigen3}
\left\{
\begin{array}{l}
(a-\lambda_1)v_x+bv_y =0\\  cv_x+(d-\lambda_1)v_y=0
\end{array}
\right. \\
\eeq
It is easy to check  that if we use for $v_x$ and $v_y$ the values  $v_x=-b$ and $v_y=a-\lambda_1$ it gives the solution of the first equation:
\beqar
(a-\lambda_1)v_x+bv_y=(a-\lambda_1)(-b)+b(a-\lambda_1)=-b(a-\lambda_1)+b(a-\lambda_1)=0
\eeqar
If we substitute these expressions into the second equation we get:
\beqar
 cv_x+(d-\lambda_1)v_y=-cb+(d-\lambda_1)(a-\lambda_1)=0
\eeqar
To prove that this expression is also zero, note that
 $(d-\lambda_1)(a-\lambda_1)-cb$  is zero in accordance with the characteristic equation (\ref{charEQ}).
Therefore  $v_x=-b$ and $v_y=a-\lambda_1$ give a solution of (\ref{eigen3})  which is   an eigen vector corresponding to the eigen value
$\lambda_1$. Similarly we find the eigen vector corresponding to the the eigen value
$\lambda_2$.

However, this approach does not work if in  (\ref{eveq2}) both $b=0$ and $a-\lambda=0$. In this case we can use the second equation $cv_x+(d-\lambda_1)v_y=0$ and find an eigen vector as  $v_x=d-\lambda_1$ and $v_y=-c$. Indeed:
\beqar
cv_x+(d-\lambda_1)v_y=c(d-\lambda_1)+(d-\lambda_1)(-c)=0
\eeqar 
As in the previous case it is easy to show that this vector 
satisfies the other (first) equation as: $
(a-\lambda_1)v_x+bv_y=(a-\lambda_1)(d-\lambda_1)+b(-c)=0$ due to
(\ref{eigen3}).

 The final formulas are:
\beq
\label{eigen_ex}
{\bf v_1}=\left( \begin{array}{c} v_{1x}
 \\  v_{1y}  \end{array} \right)=
\left( \begin{array}{c} -b
 \\  a-\lambda_1  \end{array} \right) \quad
{\bf v_2}= \left( \begin{array}{c} v_{2x}
 \\  v_{2y}  \end{array} \right)=
\left( \begin{array}{c} -b
 \\  a-\lambda_2  \end{array} \right)
\eeq
or
\beq
\label{eigen_ex2}
{\bf v_1}=\left( \begin{array}{c} v_{1x}
 \\  v_{1y}  \end{array} \right)=
\left( \begin{array}{c} d-\lambda_1
 \\  -c  \end{array} \right) \quad
{\bf v_2}= \left( \begin{array}{c} v_{2x}
 \\  v_{2y}  \end{array} \right)=
\left( \begin{array}{c} d-\lambda_2
 \\  -c  \end{array} \right)
\eeq
where $a,b$ are the elements of the matrix $A= \left(\begin{array}{lr}  a & b\\
                         c &  d
\end{array} \right)$. 

Either  (\ref{eigen_ex}) or (\ref{eigen_ex2})  can be used to find eigen vectors. (Both answers will be valid.) If, however, one of the formulas gives a zero eigen vector, we should use the other one to obtain a non-zero vector.

Let us apply these  formulas for the system (\ref{exampleEV}) with  matrix     $A= \left(\begin{array}{lr}  1 & 2\\
                         2 &  1
\end{array} \right)$ and    eigen values $\lambda_1=3;\lambda_2=-1$.  The eigen vectors can be found from (\ref{eigen_ex}) as:
\beq
\label{expr1}
\lambda_1=3; \left( \begin{array}{c} v_{1x}
 \\  v_{1y}  \end{array} \right)=
\left( \begin{array}{c} -2
 \\  1-(3)  \end{array} \right)
=\left( \begin{array}{c} -2
 \\  -2 \end{array} \right)
 \quad
\lambda_2=-1;
\left( \begin{array}{c} v_{2x}
 \\  v_{2y}  \end{array} \right)=
\left( \begin{array}{c} -2
 \\  1-(-1)  \end{array} \right))=
\left( \begin{array}{c} -2
 \\  2  \end{array} \right)
\eeq

 and from (\ref{eigen_ex2}) as:
\beq
\label{expr2}
\lambda_1=3; \left( \begin{array}{c} v_{1x}
 \\  v_{1y}  \end{array} \right)=
\left( \begin{array}{c} 1-3
 \\  -2)  \end{array} \right)
=\left( \begin{array}{c} -2
 \\  -2 \end{array} \right)
 \quad
\lambda_2=-1;
\left( \begin{array}{c} v_{2x}
 \\  v_{2y}  \end{array} \right)=
\left( \begin{array}{c} 1-(-1)
 \\  -2  \end{array} \right))=
\left( \begin{array}{c} 2
 \\  -2  \end{array} \right)
\eeq

\noindent We see that the  vectors  differ from the  vectors found earlier, but it is easy to find that they are equivalent. For example,   if we multiply the first vector by $-{1 \over 2}$ we find $-{1 \over 2}\left( \begin{array}{c} -2 \\ -2 \end{array} \right)= \left( \begin{array}{c} 1 \\ 1 \end{array} \right)$, thus the same vector which we found earlier in (\ref{tmp14}).  We also see that formulas (\ref{expr1}) and (\ref{expr2}) give equivalent result. Indeed,   first vectors obtained form (\ref{expr1}) and (\ref{expr2}) are the same.  For  second vectors  note that: $-1\left( \begin{array}{c} -2 \\ 2 \end{array} \right)= \left( \begin{array}{c} 2 \\ -2 \end{array} \right).$

\section{Functions of two variables \label{sec_2DFunc}}

\begin{figure}[hhh]
\centerline{
\psfig{type=pdf,ext=.pdf,read=.pdf,figure=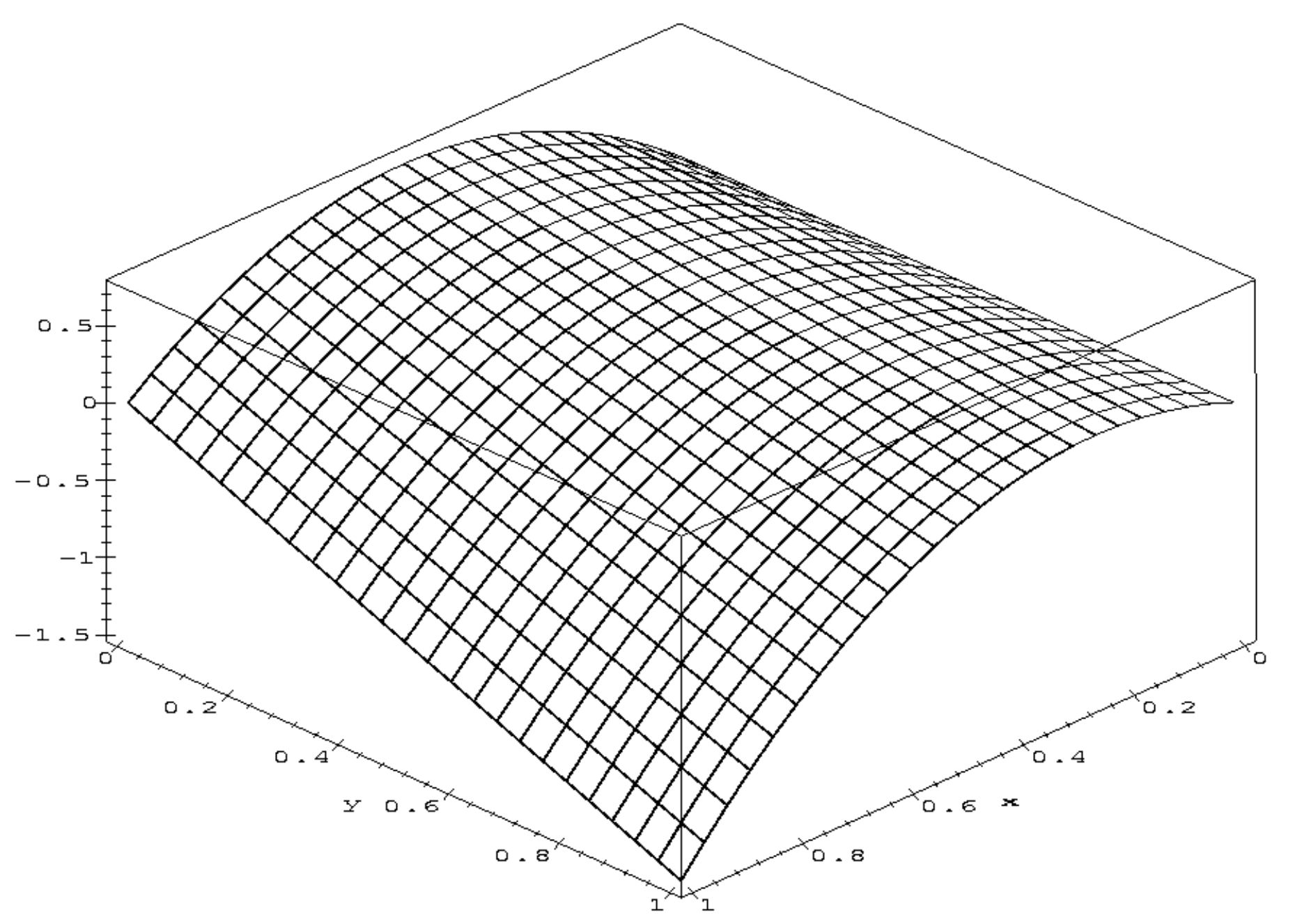,width=7.cm}
}
\caption{\label{fig3.1}}
\end{figure}

A function of two variables $f(x,y)$ describes  the rule of finding the
value of function $f$, if we know the values of the variables $x$ and
$y$.  For example, the area of a right-angled triangle with the sides
$x$, and $y$ is given by the following function of two variables:
$f(x,y)=xy/2$. Another example is the rate of growth of a prey
population in a typical ecological predator-prey model:
$f(x,y)=3x-3x^2-1.5xy$, where $x$ is the prey population and $y$ is
the predator population.  The graph of the function of one variable
$y=f(x)$ is a line on the $Oxy$-plane.  To sketch the graph of the
function of two variables $f(x,y)$, we must use a three dimensional
space $(x,y,z)$: the $Oxy$-plane for the values of the independent '
input' variables $x,y$, and the third axis $z$ for the function
'output' value $z=f(x,y)$. In such a representation the graph will be
a surface in a three dimensional space. Fig.\ref{fig3.1} shows a graph
of the function $f(x,y)=3x-3x^2-1.5xy$ plotted by a computer.

{\bf Derivatives}. The next step is the definition of the derivative
of $f(x,y)$. The main idea of finding the derivative of $f(x,y)$ is to
fix one variable at a constant value, say $x=x^*$. After that we will
get a function of one variable $y$ only ($f(x^*,y)$). Now, we can find
the derivative of $f(x^*,y)$, as the usual derivative of a function of
one variable $y$. For example, $f(x,y)=3x-3x^2-1.5xy$. Let us fix
$x=x^*=2$. We get the following function of one variable:
$f(2,y)=3*2-3*2^2-1.5*2y=-6-3y$. We can easily find the derivative
now: $df(2,y)/dy=d(-6-3y)/dy=-3$.

This type of derivative is called the {\bf partial derivative}  of $f(x,y)$
with respect to $y$ at $x=2$. We denote it as $$
\partial f / \partial y |_{x=2} =-3 
$$ 
We can find such a derivative at $x=3$, or at any other value of
$x$.  In fact for an arbitrary $x=x^*$,
$f(x^*,y)=3*x^*-3*{x^*}^2-1.5*x^**y$, and $$
\partial f / \partial y |_{x=x^*} = \partial (3*x^*-3*{x^*}^2-1.5*x^**y)/ \partial y =0-0-1.5*x^*
$$ Here $\partial (3*x^*)/ \partial y =0$ as we replaced $x$ by a
constant $x^*$ and the derivative of a constant is zero. Similarly,
$\partial (-3*{x^*}^2)/ \partial y =0$, and $\partial (-1.5*x^**y)/
\partial y = -1.5*x^*$, as $-1.5*x^*$ is a constant and the derivative
of $(ky)'=k$.  It is generally accepted to make all these
differentiations without explicitly  replacing $x$ by $x^*$. We just
should keep in mind, that for such a differentiation we treat $x$ as a
constant. Thus, to find the derivative of $f$ with respect to $y$ we
just write: 
$$
\partial f / \partial y  = \partial (3x-3{x}^2-1.5xy)/ \partial y =-1.5x.
$$ 
keeping in mind that $x$ is considered as a constant and not a variable during this
differentiation.

This expression is called  the  partial derivative of $f(x,y)$ with respect to $y$ and  is denoted  as  $\partial f / \partial y$.

Similarly, we can introduce a partial derivative of $f$ with respect to $x$:  $\partial f / \partial x$. To compute it, we fix $y$ (treat $y$ as a constant) and make the usual differentiations with respect to $x$. In our example it gives:
$$
\partial f / \partial x  = \partial (3x-3{x}^2-1.5xy)/ \partial x =3-3*2x-1.5y
$$
Here $\partial (3x)/ \partial x =3$, $  \partial (-3x^2 )/ \partial x =-3*2x$, and
 $\partial (-1.5xy)/ \partial x =-1.5y$ as  $y$ is fixed.

{\bf {Example}}.
Find  $\partial z/ \partial x$ and   $\partial z/ \partial y$ for $z=y^3\sin x$

{\bf {Solution}} $\partial z/ \partial x=y^3 \cos x$, as for
$\partial / \partial x $ we fix $y$, and $\partial (\sin x)/ \partial
x =\cos x$. Similarly, $\partial z/ \partial y= \partial (y^3\sin x) /\partial y=3y^2 \sin x $, as $x$ and hence $\sin x$ is treated as a constant.

\begin{figure}[hhh]
\centerline{
\psfig{type=pdf,ext=.pdf,read=.pdf,figure=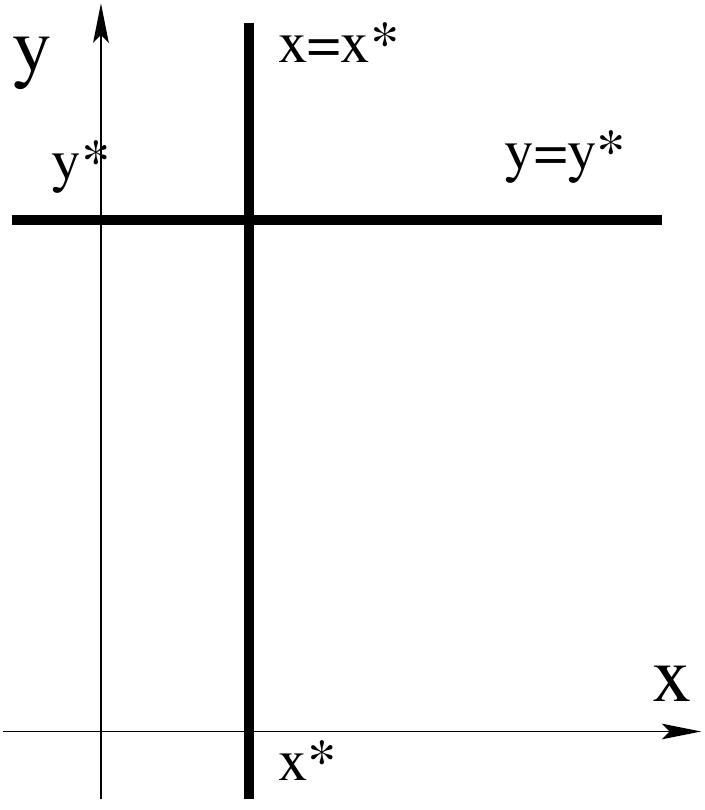,height=6.cm}
\psfig{type=pdf,ext=.pdf,read=.pdf,figure=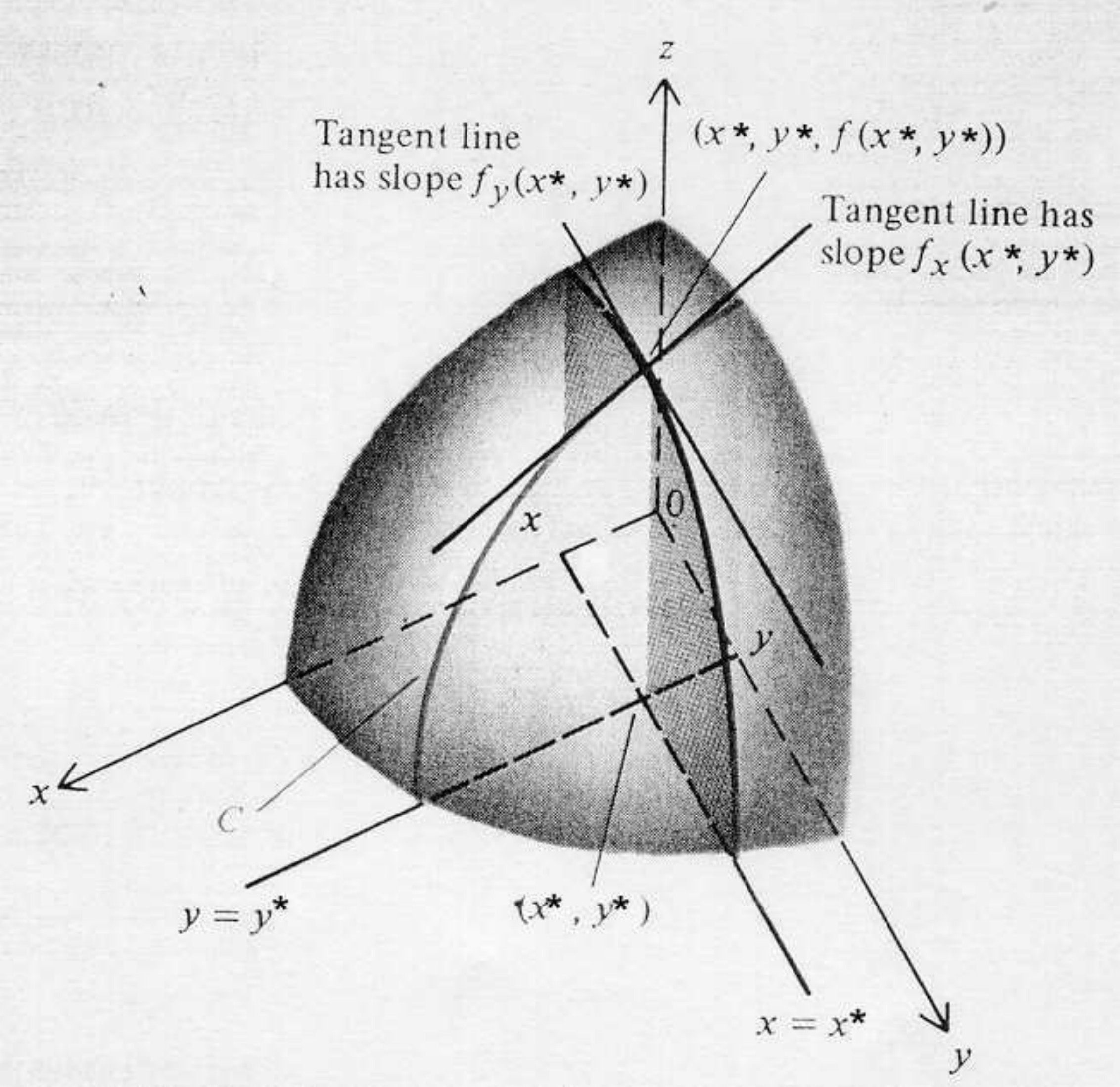,height=6.cm}
}
\caption{\label{partial}}
\end{figure}
{\bf {The geometrical representation}} of a partial derivative is
clear from fig.\ref{partial}. To compute $\partial f/ \partial x$ we
fix $y$, i.e. assume that $y$ has some value $y=y^*$. The condition
$y=y^*$ geometrically gives a horizontal line on the $Oxy$ plane
fig.\ref{partial}a, or a line parallel to the $x$-axis. In 3D this line
gives a curve on the 3D surface in   graph fig.\ref{partial}b, which is a
 1D function. The partial derivative with respect to $x$ for this particular $y^*$ will give us the   slope of the
tangent line to this 1D function. Thus (see fig.\ref{partial})  $\partial f/
\partial x$ gives the  slope of
the tangent line in the direction of the $x$-axis or the rate of change of $f(x,y)$ in the $x$ direction.  Similarly, computing
 $\partial f/ \partial y$ we
fix $x$, i.e. assume that $x$ has some value $x=x^*$. It gives us  a vertical line on the $Oxy$ plane
fig.\ref{partial}a, or a line parallel to the $y$-axis. Thus   $\partial f/
\partial y$ gives the slope of the 
tangent line in the direction of the $y$-axis, or the rate of change of
$f(x,y)$ in the $y$ direction.  If we consider $f(x,y)$ as a mountain $\partial f/ \partial x$ gives the slope of the mountain if we climb in the $x$-direction and  $\partial f/
\partial y$  gives the slope of the mountain if we climb in the $y$-direction.

Note, that in general at each point on
a surface we can draw a tangent line in any direction, and partial
derivatives $\partial f/ \partial x$  and $\partial f/
\partial y$ give the slopes  of two of these possible tangent lines. Note, that the slope of a tangent line any direction can be obtained as a combination of these two slopes. 

\begin{figure}[hhh]
\centerline{
\psfig{type=pdf,ext=.pdf,read=.pdf,figure=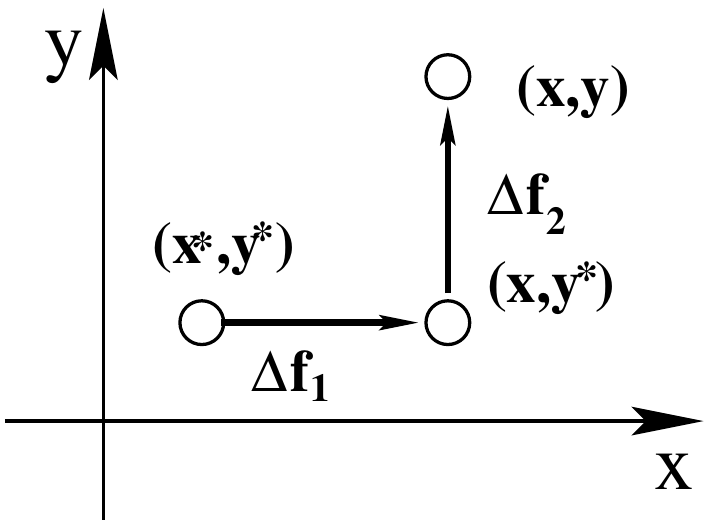,width=6.cm}
}
\caption{\label{fig3.2}}
\end{figure}

{\bf Linear approximation} Let us derive a  formula for approximating a  functions of two variables
$f(x,y)$. Let us assume  that  we know  $f(x,y)$ and its partial derivatives 
at some point $x^*,y^*$ and that we want to find the value of a function at the  close
point $x,y$ (fig.\ref{fig3.2}).
Let us move to the point $x,y$ in two steps. Let us first move
from the point $x^*,y^*$ to the point $x,y^*$, i.e. in the $x$-direction,
and then from $x,y^*$ to $x,y$, i.e. in the $y$-direction. Let
us apply the formula for approximation of a function of one variable in
formulation  (\ref{diff})  at each part of this
motion. Because on  the first part we move along the $x$
direction the change of the function ($\Delta f_1$) will be given as the 
product of the  rate of change of the function in the $x$ direction
($\partial f / \partial x$) times the distance between the points
($x-x^*$):
\beq
\label{3h1}
\Delta f_1=f(x,y^*)-f(x^*,y^*)=( \partial f / \partial x )(x-x^*)
\eeq
Similarly, on the second part of our motion, we move along the $y$-axis, and 
the change of the function here ($\Delta f_2$) will be given as the 
product of the  rate of change of the function in  the $y$ direction
($\partial f / \partial y$) times the distance between the points
($y-y^*$):
\beq
\label{3h2}
\Delta f_2=f(x,y)-f(x,y^*)=( \partial f / \partial y )(y-y^*)
\eeq

If we add equations (\ref{3h1}) and (\ref{3h2}) and solve it for
$f(x,y)$ we find the following formula which gives the approximation
for a function of two variables:
\beq
\label{2dapprox}
 f(x,y) \approx   f(x^*,y^*)+( \partial f / \partial x ) (x-x^*)+( \partial f / \partial y ) (y-y^*)
\eeq
This expression  is called a linear approximation, as  the independent variables $x,y$ are in the first power only, we do not have  terms  $x^2,y^2$, or $xy$, or etc.

{\bf {Example}} Find the linear approximation for the function $e^{x+2y}$ at the point $x=0,y=0$

{\bf Solution}. We use the formula (\ref{2dapprox}) with $f(x,y)=e^{x+2y}$ and  $x=0,y=0$.

$f(x,y)=e^{(0+0)}=1;$

$ \partial f / \partial x =e^{x+2y}$;  at  $x=0,y=0,$ $ \partial f / \partial x =e^{0+2*0}=1$

$ \partial f / \partial y =\partial (e^{x+2y}) / \partial y=e^{x+2y}*\partial ({x+2y}) / \partial y=2e^{x+2y}$, at  $x=0,y=0,$ $ \partial f / \partial y =2e^{0+2*0}=2$ 

Finally,
$  e^{x+2y} \approx 1+1*x+2*y$. 

At $x=0.1,y=0.1$ the approximate formula gives:
 $e^{x+2y} \approx 1+1*0.1+2*0.1=1.3$. The exact value of $ e^{x+2y}=e^{0.3}=1.3498$

\section{Exercises}
\subsection*{Exercises for section \ref{sec_complex}}
\ben

\item Find all roots of the given equations
\ben

\item $x^2+4x+5=0$

\item $x^2-5x+6=0$

\een
\subsection*{Exercises for section \ref{basicMatr0}}
\item Write the following linear systems in a matrix form $A\vec{X}=\vec{V}$. Find the determinant of matrix $A$. 
\ben  
\item  $\left\{
\begin{array}{l}
2x -4y =3\\ x+ y=1
\end{array}
\right.
$
\item $ \left\{
\begin{array}{l}
 ax+by=0\\ cx+dy=-b
\end{array}
\right.
$
\een
\subsection*{Exercises for section \ref{basicEigen}}
\item Find eigen values and eigen vectors of the following matrices:

\ben 
\item $\left(\begin{array}{lr}  -2 & 1\\
                         1 &  -2 \end{array} \right)$

\item $\left(\begin{array}{lr}  1 & 4\\
                         1 &  1 \end{array} \right)$

\item $\left(\begin{array}{lr}  -1 & 5\\
                         -1 &  3 \end{array} \right)$
\een

\subsection*{Exercises for section \ref{sec_2DFunc}}

\item Find partial derivatives  of these functions. After fining derivatives evaluate their  value  at the given point (if asked). 
\ben 
\item ${\partial z \over \partial x}$  and ${\partial z \over \partial y}$ for $z(x,y)=x^2+y^2-4; \quad$ at $\quad x=1;y=2$

\item ${\partial z \over \partial x}$ for $z(x,y)=x(25-x^2-y^2); \quad$at $\quad x=3;y=4$

\item ${\partial z \over \partial N}$ and ${\partial z \over \partial R}$  for $z(N,R)=N(bR-d)  \quad$ at $\quad R=0,N=0;$ and at $\quad R={d \over b},N=1:$

\item ${\partial z \over \partial P}$ and ${\partial z \over \partial M}$ for $
z(P,M)={a \over 1+P}-bM$.

\item ${\partial z \over \partial N}$ and ${\partial z \over \partial P}$ for $z(N,P)=aN-eN^2-bNP$

\item ${\partial z \over \partial M}$ and ${\partial z \over \partial A}$ for $z(M,A)=ML - \delta A - \frac{vMA}{h+A} $

\item ${\partial z \over \partial P_1}$ and ${\partial z \over \partial P_2}$  for $z(P_1,P_2)={-aP_2 \over h +P^2_1 + 2P_2} $

\item ${\partial z \over \partial N}$ and ${\partial z \over \partial T}$ for $z(N,T)=\frac{b^2N^2T}{1+cN+bTN^2} $

\een


\subsection*{Additional Exercises}

\item Perform  the indicated operations:
\ben

\item $\sqrt{3^2-90}$

\item $(-1+2i)+(4+7i)$

\item $(4+5i)*(7+2i)$

\item ${1 \over i}$

\een

\item Find all roots of the given equations
\ben
\item $x^2+121=0$

\item $x^2+2x+3=0$
\een

\item \label{CramersRule} Cramer's rule on an example. Cramer's rule for system of two linear  equations:
$$ \left\{
\begin{array}{l}
 ax+by=E\\ cx+dy=F
\end{array}
\right.
$$ allows us to find solutions from determinants of matrices. First we need to find a determinant of the main matrix $A$
$$
detA=det \left|\begin{array}{lr}  a & b\\
                         c &  d
\end{array} \right|.$$
 Then we need to find the determinant of a matrix formed by replacing the $x$-column values of the matrix $A$ with the answer-column values $\left(\begin{array}{c}  E \\ F \end{array} \right)$ as $det D_x=det \left|\begin{array}{lr}  E & b\\ F & d \end{array} \right|$  and similarly for the $y$-column: $det D_y= det \left|\begin{array}{lr} a & E\\ c & F \end{array} \right|$. The solutions of the system will be given by the ratios of these determinants as: 
$$x={det D_x \over detA}; \;\; y={det D_y \over detA }$$. 
\bit
\item Find solutions of the following system using the Cramer's rule:
\beq
\label{matr6_2}
\left\{
\begin{array}{l}
x +2y =5\\ 2x+ y=4
\end{array}
\right.
\eeq
\item Find also solution of (\ref{matr6_2}) by usual method. Show that Cramer's rule gives a correct result.
\eit
\item Find eigen values and eigen vectors of the following matrices:

\ben 
\item $\left(\begin{array}{lr}  -1 & 6\\
                         2 &  -2 \end{array} \right)$

\item $\left(\begin{array}{lr}  2 & 1\\
                         7 &  -4 \end{array} \right)$

\een

\item  Find a linear approximation for the function at the given point.
\ben 
\item $f(x,y)=x^2+y^2; \quad$ at $x=1,y=1$
\een
\een

\chapter{Differential equations of one variable}

Differential equations are equations that contain a derivative of an
unknown function. As we know derivative gives a velocity of a process
and differential equations occur when we describe various processes
via their velocities. Differential equations are widely used for
modeling in a variety of disciplines: in mathematics, physics,
chemistry, economics, engineering, medicine, life sciences,
etc. Development of methods of study of differential equations is the
main subject of this course.

In this chapter we will introduce differential equations, give
first definitions, show how to solve simple differential equations
analytically and consider a few examples. Then we will develop
qualitative methods for the analysis of differential equations of one
variable and  will  apply them for biological models.

\section{Differential equations of one variable and their solutions \label{sec_separation}}
\subsection{Definitions}

Let us construct a first differential equation. Consider a motion of a
car with a contact velocity, for example $v=10 m/sec$. If we denote
the distance traveled by the car at time $t$ as $l(t)$ we can write  $velocity=v={dl \over
dt}$, as  the velocity is the
derivative of the distance with respect to time. Thus we can write the following differential equations for this
process:
\beq
\label{e1.1}
dl/dt=10
\eeq

If the  car  travels  with an acceleration, then the velocity will linearly increase with time. If we assume that the acceleration is $a=1.2 m/sec^2$, then the velocity in the course of time will be given by  $v=1.2t$ and we get a differential equation as:

\beq
\label{e1.1_2}
dl/dt=1.2t
\eeq

Let us consider a biological example. If $N(t)$ is the population size
of a species at time $t$, then the rate of change of the population
size is:
\beq
\label{e1.2}
dN/dt=births-deaths
\eeq
Let us assume that the birth and the death terms are proportional to
$N$. This assumption is quite reasonable, as it means, that if, for
example,  we know the growth rate of a population of some insects on one
tree, then the total growth rate of the whole population of insects in
the forest will be proportional to the total number of trees. Thus
each term in equation (\ref{e1.2}) will be proportional to $N$ and we
get the following famous 'Malthus' equation for population dynamics:

\beq
\label{e1.4}
dN/dt=(\alpha-\beta)N=kN
\eeq
where $\alpha$ and $\beta$ are the rate constants for the birth and
death processes, and we see that $k > 0$ if $\alpha > \beta$, and
$k<0$, if $\alpha < \beta$.

Another model assumes that there is a maximum population size $K$  (called the carrying capacity) and that the the rate of growth of population depends on how close  the population is to this maximum size.  This yields the following differential equations:

\beq
\label{e1.4_2}
dN/dt=k(K-N)
\eeq

Now note, that in mathematical sense equations of the type (\ref{e1.1_2}) can be written as 
\beq
\label{nonauton}
{dx \over dt}=f(t)
\eeq
as the variable $t$ with respect to which we  differentiate function $x$ is also present in the right hand side of our equation.

Alternatively the equations of the type (\ref{e1.4}) and (\ref{e1.4_2}) can be written as:
\beq
{dx \over dt}=f(x) 
\eeq
as here  the variable $t$ is not present at the right hand side and we have only the unknown function $x$ (for  eqns(\ref{e1.4}),(\ref{e1.4_2}) $N$) there.

The latter equations  are the most important for us in this course and they have  a special name an 'autonomous differential equations':

\begin{D}
 Equation 
\beq
\label{1dgeneral}
{dx \over dt}=f(x) 
\eeq
 is called an autonomous differential equation
\end{D}

{\bf {Example}}

${dx \over dt}=2x-tan(x) \qquad autonomous$

${dx \over dt}=3 sin(x-t) \qquad non-autonomous$

${dx \over dt}=arcsin(x) \qquad autonomous$

\vskip 1pc
Before we find how to solve differential equations, let us discuss
from a very general point of view what kind of solutions can we expect
here. If we assume that a differential equation describes how the size
of a population will change in time, then we may think about two types
of problems. The first one it to find this size for a particular
population at each time moment. For that we obviously would need to
know the initial size of a population.  We can also ask a more general
question: to find the population size for an arbitrary initial
size. This solution will be called {\bf {the general solution}} of a
differential equation.  Because the general solution contains
information on solutions for arbitrary initial conditions, it normally
depends on an arbitrary constant. The differential equation with given
initial condition is called an initial value problem:

\begin{D}
  The problem $ {dx \over dt}=f(x), x(0)=x^* $ is called the initial value problem; Its solution is called the orbit or trajectory.
\end{D}
The initial value problem for most  $f(x)$ has a unique
solution.

Now let us consider how we can  solve differential equation.

\subsection{Solution of a differential equation }

e can easily solve an equations of type (\ref{nonauton}) using the
method of {\it separation of variables}. The main idea of this method
is to think about the derivative of an unknown function $x$ in ${dx
\over dt}=f(t)$ as a fraction $dx$ over $dt$. If we multiply  both
sides of this equation by $dt$ we will get: $$ dx=f(t)dt $$ If we
integrate  both sides of this equation and  get:
\beq
\label{separ1}
\begin{array}{l}
\int dx=\int f(t)dt\\
or\\
x=F(t)+C
\end{array}
\eeq
where $F(t)$ is an anti-derivative of $f$ and $C$ an arbitrary
constant. This is the general solution of differential equation
(\ref{nonauton}). We see that this solution  contains an arbitrary constant,
as expected from the discussion in the previous section.

Let us apply this method to a few examples.

For the simplest differential equation  $ {dx \over dt}=0$ we get:
\beq
\label{separ2}
\begin{array}{l}
{dx \over dt}=0\\
or\\
dx=0*dt\\
or\\
\int dx=\int 0 dt\\
or\\
x=C
\end{array}
\eeq
Therefore, the solution of this equation is  $x$ equals any constant $C$.

For the equation of motion of a car with a velocity $10 m/sec$ (\ref{e1.1}) we get:
\beq
\label{separ3}
\begin{array}{l}
{dx \over dt}=10\\
or\\
\int dx=\int 10 dt\\
or\\
x=10t+C
\end{array}
\eeq
Thus this solution shows us the  position of a car as a function of time $t$, and the  arbitrary constant $C$ here accounts for the  initial position of the car.

Finally, for  the motion from the rest with an acceleration $a=1.2 m/sec^2$   (\ref{e1.1_2}) we get:

\beq
\label{separ4}
\begin{array}{l}
{dx \over dt}=1.2t\\
or\\
\int dx=\int 1.2t dt\\
or\\
x=1.2{ t^2 \over 2}+C
\end{array}
\eeq
We obtained a formula which is well known to you from  your school physics and $C$  here also accounts for the  initial position of a car.

It turns out that we can  apply the method of separation of
variables also for an autonomous equation (\ref{1dgeneral}) (${dx \over
dt}=f(x)$). However, here we will need to separate variables as ${dx \over f(x)}=dt$ and as a result we will not  usually get the  explicit formula for $x(t)$. However, in many  cases we can do it after  some transformations.

Let us solve equation for population dynamics (\ref{e1.4}). 
\beq
\label{separ5}
\begin{array}{l}
{dN \over dt}=kN\\
or\\
\int {dN \over N}=\int k dt\\
or\\
ln(N)=kt+C
\end{array}
\eeq
To find $N$  note that equation $ln(x)=a$ has a solution $x=e^a$, thus we find $
 N=e^{kt+C}=e^C*e^{kt}$ and if we denote $e^C=A$ we will get:
\beq
\label{e1.5}
N(t)=Ae^{kt},
\eeq where $A$ is an arbitrary constant.

Let us apply it for the following initial value problem with  $k=4$ and  the initial population size of  $10$:
\beq
\label{e1.7}
dN/dt=4N \qquad N(0)=10
\eeq
The general solution here is given by $N=A*e^{4t}$. To find the solution of the initial value problem we note that at $t=0$ the population size was $N(0)=10$, i.e. we can write: $N(0)=Ae^{4*0}=10$, or we find that $A=10.$ Hence, we have the following solution of the initial value problem (\ref{e1.7}): $N=10e^{4t}$.

Similarly, we can solve the general initial value problem (\ref{e1.4})
for an arbitrary  initial size of $N(0)$ and find:
\beq
\label{e1.8}
N=N(0)e^{kt}
\eeq
thus $A$ here  gives the initial size of the population.

Finally let us solve equation (\ref{e1.4_2}), for specific parameter values $K=20,k=4$:
\beq
\label{separ6}
\begin{array}{l}
{dN \over dt}=4(20-N)=-4(N-20)\\
or\\
\int {dN \over N-20}=-\int 4 dt\\
or\\
ln(N-20)=-4t+C\\
or\\
N-20=Ae^{-4t}\\
or\\
N=20+Ae^{-4t}
\end{array}
\eeq
This is a general solution. Let us find a particular solution,
corresponding to the initial size of the population $N(0)=10$, for
example. For that we find: $N(0)=20+Ae^{-4*0}=10$, or $A=10-20=-10$, thus
 the population dynamics in the course of time will be given: $N=20-10e^{-4t}$.

\begin{figure}[hhh]
\centerline{
\psfig{type=pdf,ext=.pdf,read=.pdf,figure=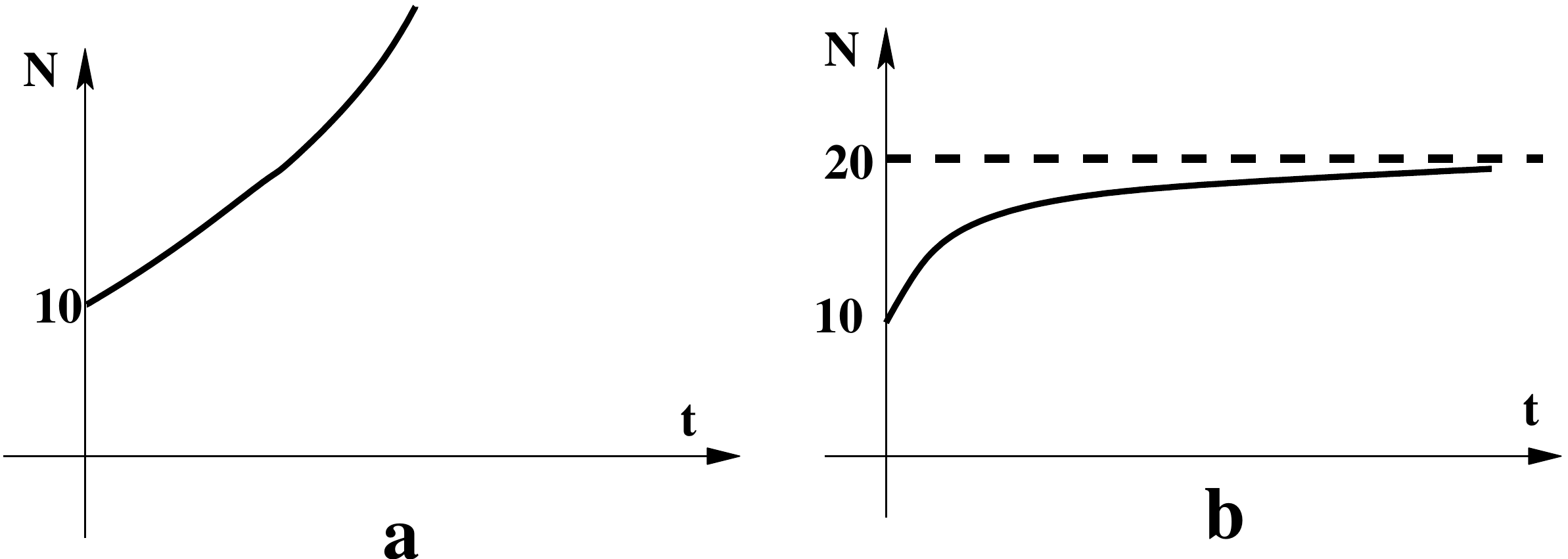,width=10.cm}
}
\caption{\label{figSep1} a- Size of population obtained form the solution of the initial values problem (\ref{e1.4}) with $k=4,N(0)=10$; b-the same for equation (\ref{e1.4_2}) with $K=20, k=4,N(0)=10$; }
\end{figure}

The  solutions of equation (\ref{e1.4}) and equation (\ref{e1.4_2}) are shown in fig.\ref{figSep1}.

We see  that the size of the population for
eq.(\ref{e1.4}) goes to infinity (fig.\ref{figSep1}a).  Quantitatively  the size of
population increases in $e
\approx 2.73$ times each ${1 \over 4}=0.25$ seconds. Indeed,  from
the particular solution $N=10e^{4t}$ we find
$N(0)=10;\;N(0.25)=10*e;\;N(0.5)=10*e^2;\; etc. $ In general, for
arbitrary $k$ in (\ref{e1.4}) this {\it characteristic time} of change
is given by $\tau={1 \over k}$, which follows from equation
(\ref{e1.8}). For equation (\ref{e1.4_2}), the population approaches
its carrying capacity value of $K=20$ (see fig.\ref{figSep1}b) and the
characteristic time is also determined by the value of $k$ in the
following sense: the difference between the current population size
and its stationary value decreases in $e$ times over the time period
$\tau={1 \over 4}$. This follows from solution
(\ref{separ6}), which gives for this difference
$N-20=Ae^{-4t}$. Similar expression for an arbitrary value of $k$ is
given by $N-K=Ae^{-kt}$, which gives for the characteristic time $\tau={1 \over k}$.

This concludes our analytical study of differential equations. In the
next chapter we will formulate another method of analysis of 
differential equations that does not require  direct integration
of these equations.

\section{Qualitative methods of analysis of differential equations of one variable \label{sec_1Dqualitative}}
In this section we will consider a general non-linear differential equation
${dx \over dt}=f(x)$ and develop an effective method  for  qualitative analysis of  this equation without finding solutions analytically. In the next chapters this method will be extended to the systems of two differential equations.

\subsection{Phase portrait \label{sec_1Dphase}}
Let us start with equation (\ref{e1.40}) which we considered in
section \ref{sec_separation}.
\beq
\label{e1.40}
dN/dt=kN 
\eeq

We found that $N=10e^{4t}$ is a solution for this equation for $k=4$ and the
initial population size of $N(0)=10$ (see (\ref{e1.7})). This solution was represented graphically in fig.\ref{figSep1}.

\begin{figure}[h]
\centerline{
\psfig{type=pdf,ext=.pdf,read=.pdf,figure=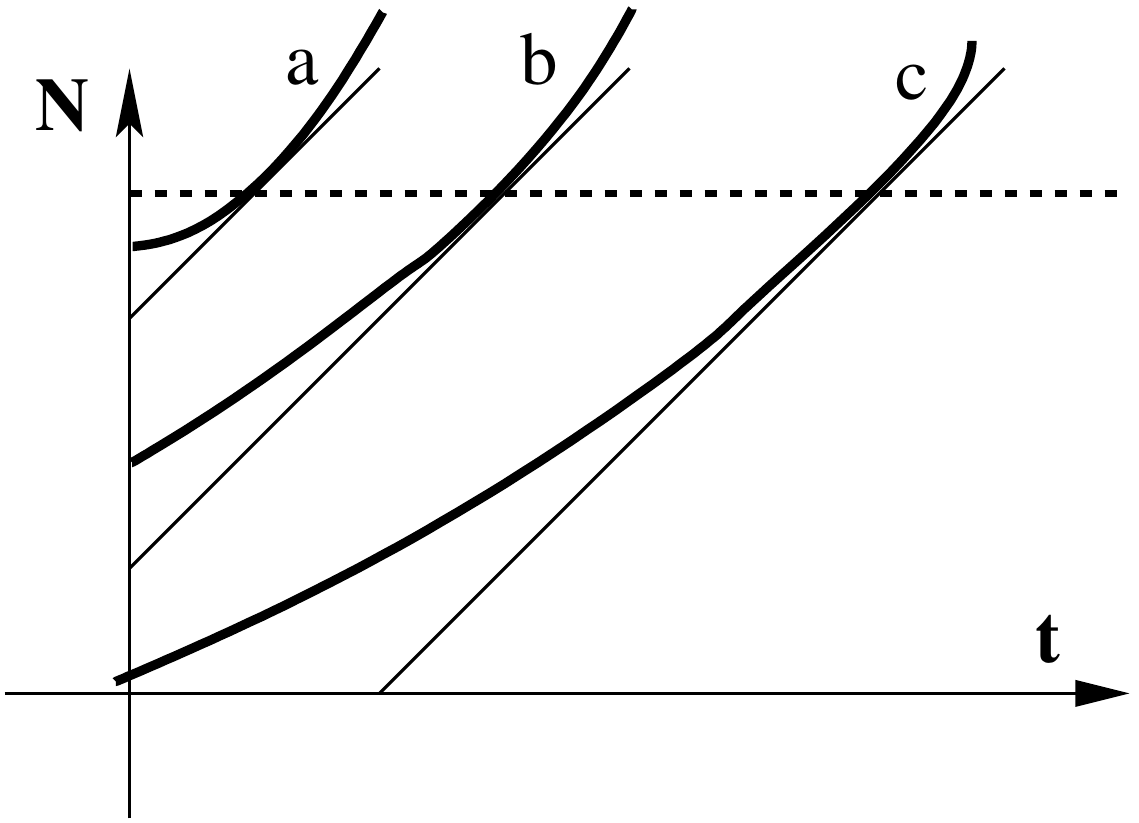,width=5.cm}
}
\caption{\label{fig1.2} Three solutions of equation (\ref{e1.40}) together with tangent lines for these solutions at $N=N^*$.  }
\end{figure}

If the initial size of the population was different, for example
$N(0)=5,$ or $N(0)= 3$, $N(0)= 0.1$ etc., we get other solutions of
equation (\ref{e1.40}) and if we represent these solutions graphically
we will obtain the following curves ($a,b,c$) shown in
fig\ref{fig1.2}. Let us analyze them. An important characteristic of
any line is its slope. It turns out that we can easily find the slope
for the solutions of (\ref{e1.40}): $slope=dN/dt=4N$. We see that the
slope depends only on $N$ (the size of the population) and does not
depend on other factors, for example on the initial conditions. For
example, if $N=3$ the slope of the line representing solution at point
$N=3$ is 4*3=12 for any initial condition.  Geometrically this means
if we graph several solutions (as in fig\ref{fig1.2})
and determine slopes of these functions for given $N$ (at points of
intersection of the dotted line $N=N^*$  with the solution
curves) we will find that all slopes are the same ($slope= 4*N^*$).

We can use this information and represent a qualitative picture of
solutions of (\ref{e1.40}) using only one $N$-axis. For that let us
denote the slope of the solution on the $N$-axis for each value of $N$
(see numbers 4,8,12,16, etc.). We  see, however, that this  is not very helpful for representation of the solution of our equation.
\begin{figure}[h]
\centerline{
\psfig{type=pdf,ext=.pdf,read=.pdf,figure=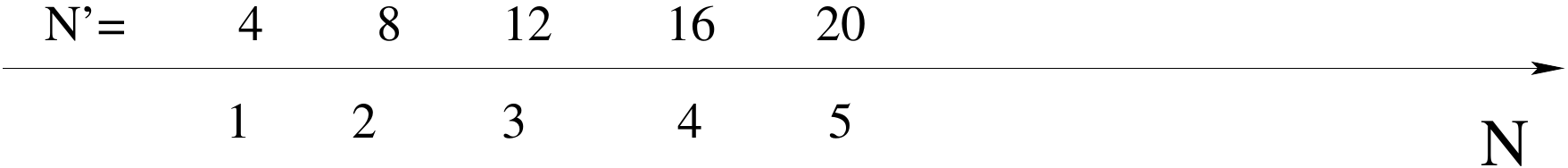,width=12cm}
}
\end{figure}
To improve it let us think about biological interpretation of the
slope of the curve in fig\ref{fig1.2}. The slope of the curve gives us
the rate of change of the function and because fig\ref{fig1.2} shows
how the size of population depends on time, the slope values
(represented above the $N$-axis) show the growth rate of population at
given $N$.  The most important qualitative aspect of the dynamics here
is the gowth of  population. We can represent it  in the following qualitative way:
\begin{figure}[h]
\centerline{
\psfig{type=pdf,ext=.pdf,read=.pdf,figure=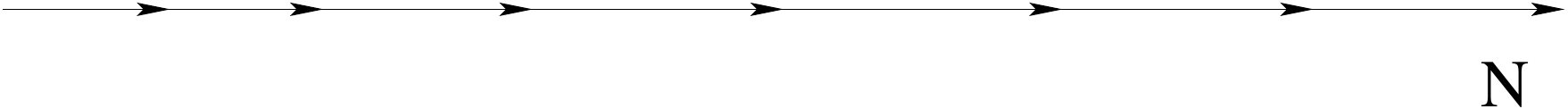,width=12cm}
}
\end{figure}

i.e. we show the growth of the population by an arrow which is  directed to the right. Of course, this is  not a complete description of our system, but it gives a good idea about  the behavior of our system. It shows that if we start at some initial value of $N^*$, then $N$ will grow and the size of the  population will be continuously increasing. Note  that to obtain this result we have only used the direction of the arrows in the figure.

As we will see in the next section, such representation can be easily
obtained for any autonomous differential equation $({dx \over
dt}=f(x))$ from the graph of the function $f(x)$ at the right hand
size. Such representation is called the phase portrait:
\begin{D}
 The collection of all possible orbits of a differential equation together with the direction arrows is called the phase portrait.
\end{D}
\subsection{Equilibria, stability, global plan}
Let us consider two differential equations that arise in population ecology:
\beq
\label{e1.9}
dx/dt=4x
\eeq
and 
\beq
\label{e1.10}
dx/dt=240-0.01x
\eeq
\begin{figure}[hhh]
\centerline{
\psfig{type=pdf,ext=.pdf,read=.pdf,figure=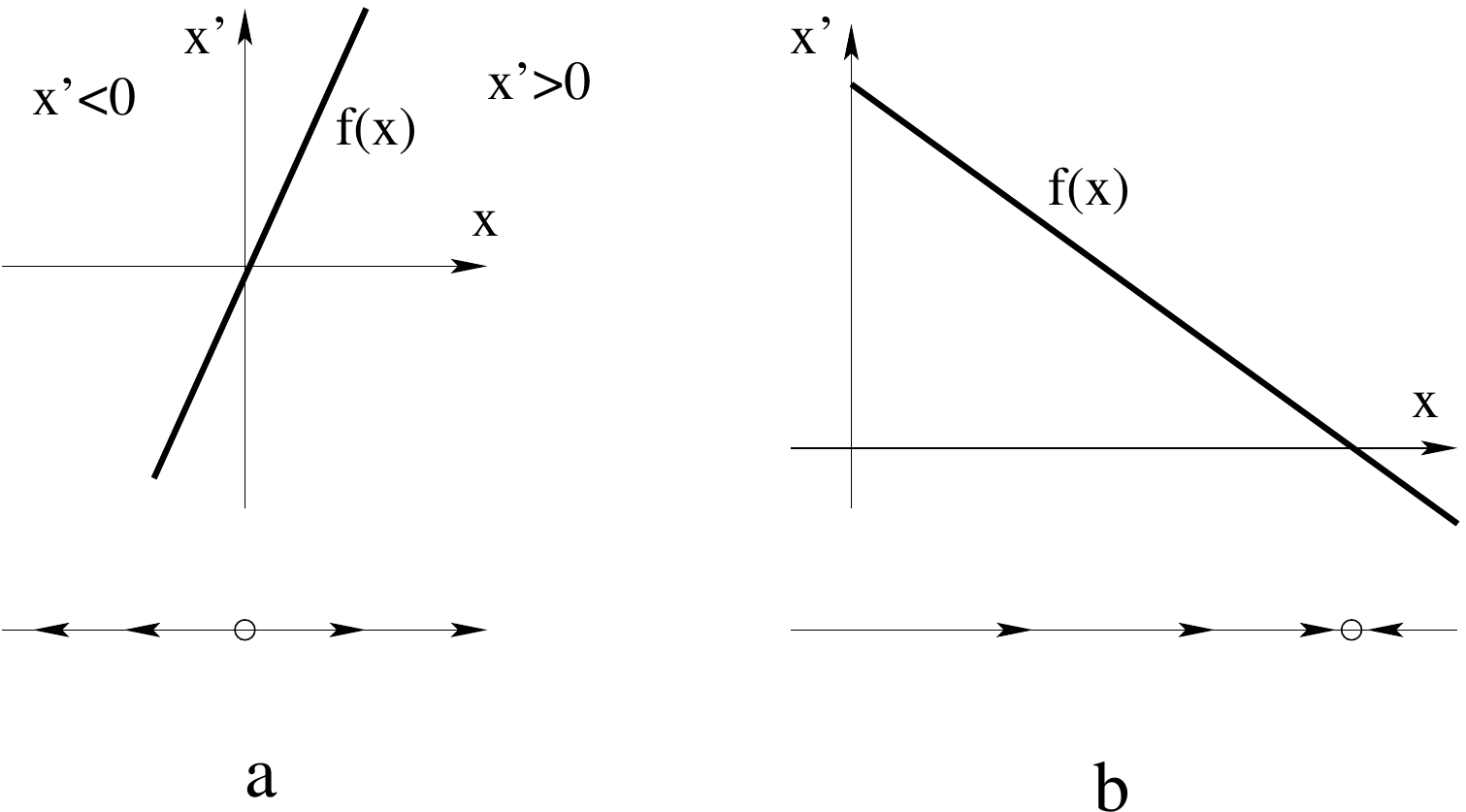,width=11.0cm}
}
\caption{\label{fig1.3}}
\end{figure}
Let us sketch phase portraits for (\ref{e1.9}) and (\ref{e1.10}). We
can do it without finding a solution. In general, to sketch a phase
portrait of an equation $({dx \over dt}=f(x))$ we need to draw
$\rightarrow$ or $\leftarrow$ arrows on the $x$-axis. The
$\rightarrow$ arrow means growth of $x$, i.e. ${dx \over dt}>0$. The
$\leftarrow$ arrow means decreasing of $x$, or ${dx \over
dt}<0$. Because ${dx \over dt}=f(x)$ we will  graph function
$f(x)$ and then assign  $\rightarrow$ to that regions where the graph
 is above the $x$-axis and $\leftarrow$ to that regions
where the graph is below the $x$-axis. For equation (\ref{e1.9})
$dx/dt=4x$, the graph of $f(x)=4x$ is shown at the top panel of
fig.\ref{fig1.3}a and we draw the right arrow $\rightarrow$ for $x>0$
, and the left arrow $\leftarrow$ for $x<0$ (fig.\ref{fig1.3}a bottom).

The interesting point here is $x=0$.  Here $dx/dt=4x=0$, thus the rate
of change of $x$ here is zero and we cannot assign any direction for
the arrow at this point. However, the dynamics of our system here is
trivial: $x$ do not change in the course of time, or $x(t) = 0$ for
all $t$. This means, that if the initial size of the
population was zero, it will be zero forever. Such points of a phase
portrait are called equilibria. They occur at points where the rate of
change is zero (${dx \over dt}=0$) For equation $({dx \over dt}=f(x))$ equilibria occur if $f(x)=0$, which is also used as a definition of an equilibrium. 

\begin{D}
A point $x^*$ is called an equilibrium point of $dx/dt=f(x)$, if $f(x^{*})=0$
\end{D}

Finally the phase portrait of eq.(\ref{e1.9}) in fig.\ref{fig1.3}a
gives the following dynamics of $x$: if the initial value of $x$ is to
the right or to the left left of the equilibrium point $x=0$, it will
go to plus or minus infinity respectively.

Let us study equation (\ref{e1.10}). Again, our plan is:
$f(x)\; graph \rightarrow  phase$ $portrait$ (fig.\ref{fig1.3}b).
Here we have an equilibrium point $x=24000$  which is the root of the function $240-0.01x$, and the arrows (flow) for this case are  shown in fig.\ref{fig1.3}b. So, the dynamics of solutions of our equation will be as in the following figure:

\begin{figure}[hhh]
\centerline{
\psfig{type=pdf,ext=.pdf,read=.pdf,figure=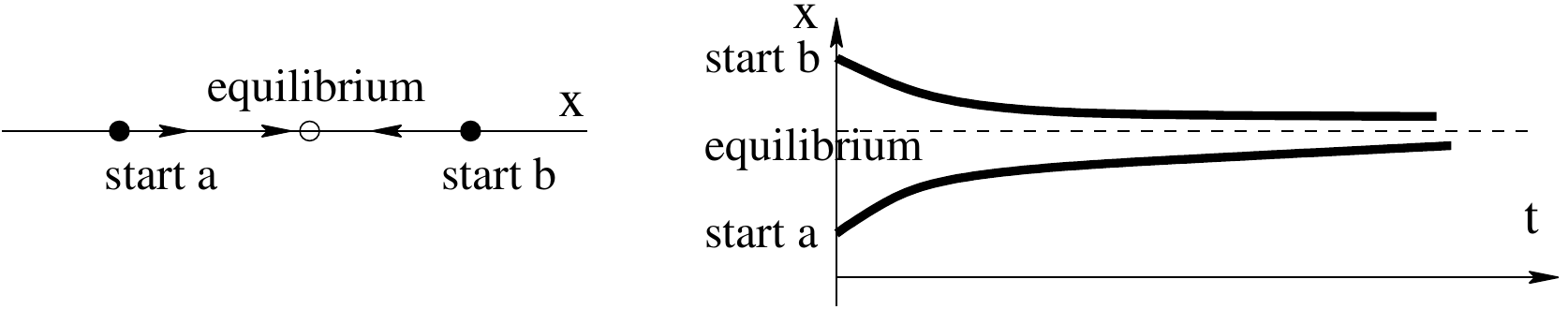,width=12.5cm}
}
\caption{\label{fig1.4}}
\end{figure}

i.e. for any initial condition, $x$ will eventually approach the
equilibrium point.

If we compare the equilibria of equations (\ref{e1.9}) and
(\ref{e1.10}), we can see that they are different. The variable $x$
diverges from the equilibrium point of equation (\ref{e1.9}). Such
equilibrium points are called non-stable equilibria. On the contrary,
the variable $x$ converges to the equilibrium point of equation
(\ref{e1.10}). Such equilibria points of are called stable equilibria
or attractors.

Now we can formulate  a general plan  for finding the phase portrait of $dx/dt=f(x)$.

{\bf {Global plan}}.
\ben 
\item  Sketch  the graph of $f(x)$. 

\item Draw the phase portrait. For that transform the points where $f(x)=0$ to equilibria points,  the regions where  $f(x)>0$ to right headed arrows ($\rightarrow$), and  the regions where  $f(x)< 0$ to  left  headed arrows  ($\leftarrow$). This gives the overall phase portrait.
\een

Let us apply it  to the  logistic equation for population growth
\beq
\label{e1.11}
dn/dt=rn(1-n/k) \;\;\;n \geq 0
\eeq
This equation describes growth of a population in a medium with limited resources. We can  study (\ref{e1.11}) for arbitrary values of parameters $r,k$. However for simplicity let us fix $r=2$ and $ k=3$. 
The equation becomes
\beq
\label{e1.12}
dn/dt=2n(1-n/3)=(2/3)*n*(3-n)
\eeq
Let us find the phase portrait  and the dynamics of the solutions of   (\ref{e1.12}).
 First we use the global plan. 
\ben
\item The right hand side of our equation is $(2/3)*n*(3-n)=2n-(2/3)n^2$. The graph of this function is a parabola, opened  below with the roots $n=0;n=3$.
\begin{figure}[hhh]
\centerline{
\psfig{type=pdf,ext=.pdf,read=.pdf,figure=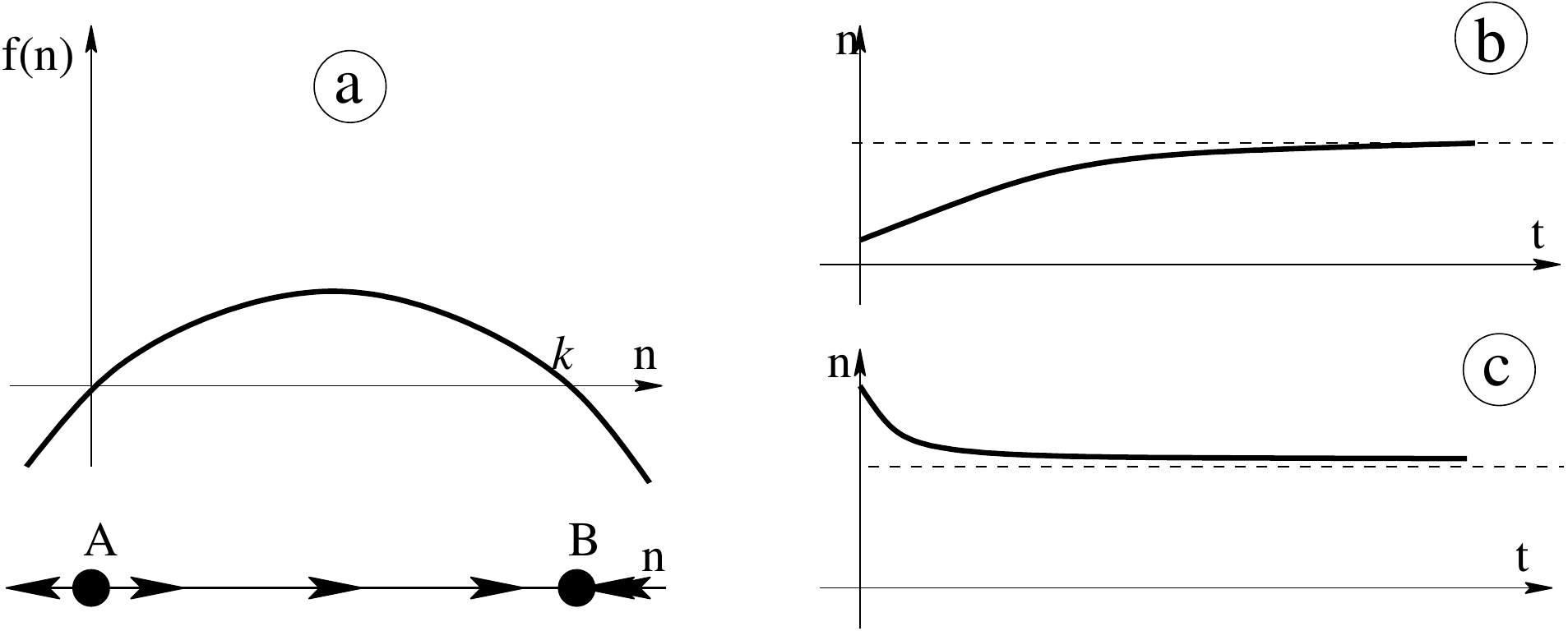,width=12.5cm}
}
\caption{\label{fig1.6}}
\end{figure}

\item The construction of the phase portrait is clear from fig.\ref{fig1.6}a.

\een
The dynamics of the system for different initial conditions is shown in 
 fig. \ref{fig1.6}b,  for an    initial size of the population  $0<n<3$, and in fig.\ref{fig1.6}c for  $n>3$. We see that in the course of time the size of the population  becomes $n=3$, i.e. the stable equilibrium point  $n=3$ is the only attractor of our system.

\begin{figure}[hhh]
\centerline{
\psfig{type=pdf,ext=.pdf,read=.pdf,figure=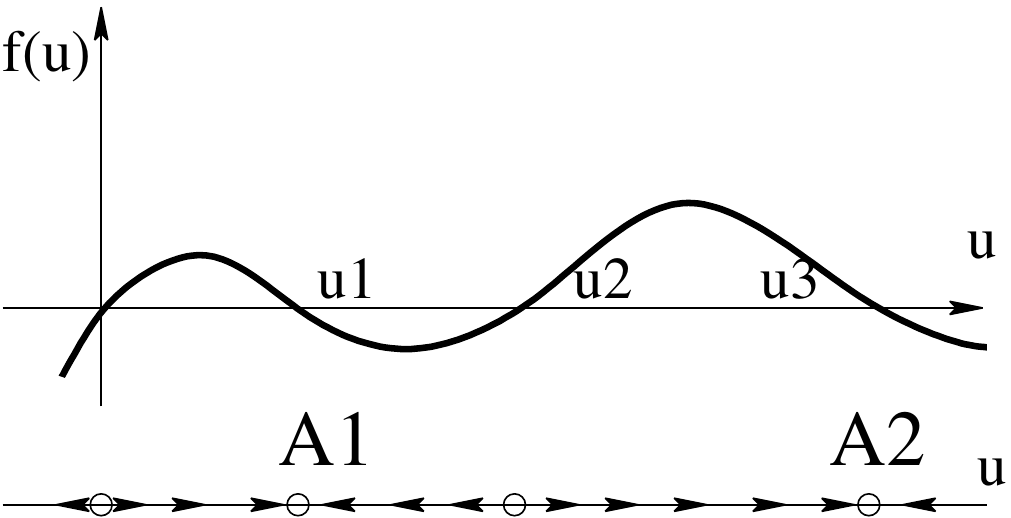,width=8cm}
}
\caption{\label{fig1.8}}
\end{figure}

Sometimes,  differential equations  have several stable equilibria (attractors). For example, the model for the spruce bud-worm population (\ref{e1.13}) has the following  phase portrait (fig.\ref{fig1.8}).
\beq
\label{e1.13}
du/dt=f(u)
\eeq

We see  that there are two attractors: A1 and A2 which correspond to bud-worm populations  of different size. We see that if the initial size of the population is $0<u_0<u_2$, then the population eventually reaches A1; if  $u_2<u_0<\infty$, then  population eventually reaches A2. These intervals are called basins of attraction.
\begin{D}
The basin of attraction of a stable equilibrium point $x^{*}$ is the set of values of $x$ such that, if  $x$ is  initially somewhere in that set, it will subsequently move  to the equilibrium point  $x^{*}$.
\end{D}
In the case of fig\ref{fig1.8}, the basin of attraction of the equilibrium $u_1$ (A1) is the interval  $0<u_0<u_2$;  the basin of attraction of the equilibrium $u_3$ (A2) is the interval  $u_2<u_0<\infty$. 
It is very important to know  basins of attraction of a system in order to predict its  behavior.
\section{Systems with parameters. Bifurcations.\label{secBif}}

One of the aims of modeling in biology is to predict the behavior of a
system for different conditions.  In that case we differential
equations will contain parameters. Let us consider two examples. The
first is a general linear equation with one parameter $k$

\beq
\label{e1.140}
dn/dt=k*n
\eeq

\begin{figure}[hhh]
\centerline{
\psfig{type=pdf,ext=.pdf,read=.pdf,figure=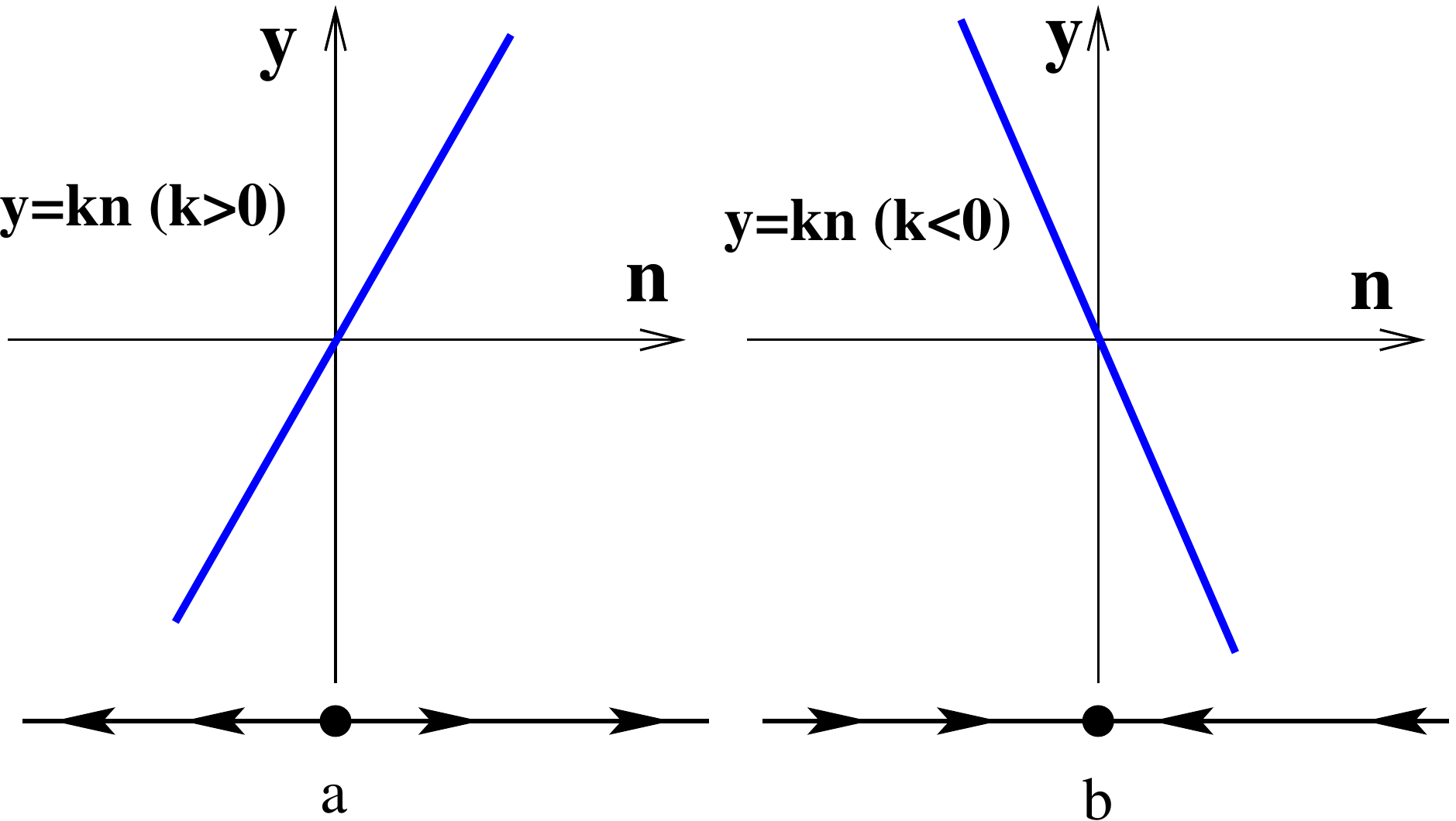,width=10cm}
}
\caption{\label{fig1.ex4}}
\end{figure}

If we draw the graph of the right hand side function $y=kn$ for different values of the parameter $k$ we find that we can have two possibilities  depending on the sign  of the parameter $k$ (fig.\ref{fig1.ex4}): a non-stable  equilibrium point at $n=0$ for $k>0$ (fig.\ref{fig1.ex4}a) and a stable equilibrium point at $n=0$ for $k<0$ (fig.\ref{fig1.ex4}b).

\vskip 1pc

\noindent Another example of an equation with parameters is the logistic equation
for population growth (\ref{e1.11}). This equation depends on two
parameters $r$ and $k$, where $r$ accounts for the growth rate and $k$
accounts for the carrying capacity.  Let us consider a slight
modification of equation (\ref{e1.11}) for a population which is
subject to harvesting at a constant rate $h$:
\beq
\label{e1.14}
dn/dt=r*n*(1-n/k)-h
\eeq
where $h$ is the  harvesting rate (an extra parameter).

\begin{figure}[hhh]
\centerline{
\psfig{type=pdf,ext=.pdf,read=.pdf,figure=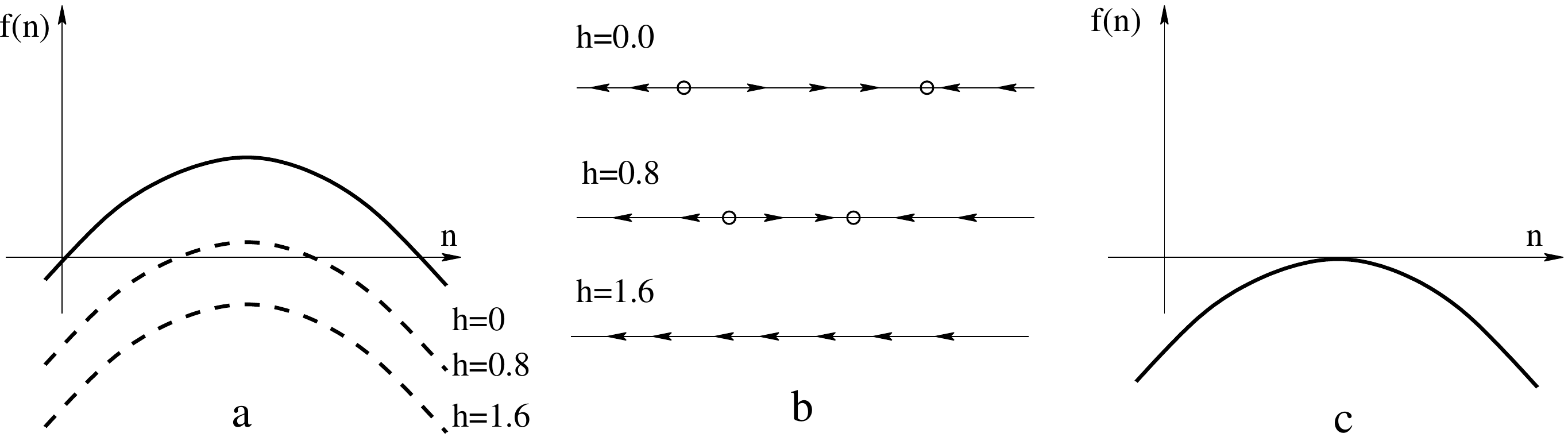,width=16cm}
}
\caption{\label{fig1.9}}
\end{figure}

Let us fix the parameters $r=2$ and $k=3$  (the same values as in equation (\ref{e1.12})), and study only the effect of varying the harvesting  parameter $h$ on, the dynamics of the population.
\beq
\label{e1.15}
dn/dt=2*n*(1-n/3)-h
\eeq
 When $h=0$ equation  (\ref{e1.15})  coincides with equation   (\ref{e1.12}) which was studied in fig.\ref{fig1.6}.  Now assume that the harvesting $h$  is not zero.
Let us plot graphs of $2*n*(1-n/3)-h$ for $h=0;h=0.8;h=1.6$ (fig.\ref{fig1.9}a). 
The phase portraits for $h=0.0;h=0.8;h=1.6$ are shown in fig\ref{fig1.9}b.
We see that at $h=0.8$ the behavior of the system is qualitatively
similar to the behavior of the system without harvesting ($h=0.0$): the
population eventually approaches the stable non-zero
equilibrium. However the final size of the population in this case is
slightly smaller than for the population without harvesting
(fig.\ref{fig1.6}a) At $h=1.6$ the situation is different. We do not
have a stable and non-stable equilibrium anymore. The flow is always directed
to the left and the size of the population decreases. In this simple
model this means the extinction of the population. The important
question here is, what is the maximal possible harvesting rate at
which the population still survives. From the previous analysis it is
clear that the critical harvesting is reached when the parabola
$2*n*(1-n/3)-h$ touches the $n$-axis (fig.\ref{fig1.9}c). To find this
critical value we note, that $h$ just shifts the parabola $2n(1-n/3)$
downward. Therefore, the situation of (fig.\ref{fig1.9}c) occurs, when the
shift equals the maximum of the parabola $2*n*(1-n/3)$. To find the
maximal value we find a point where the derivative $df/dn=0$
\beq
\label{e1.16}
\begin{array}{l}
df/dn=2/3(3-n)-(2/3)*n*3=2-4n/3=0\\
n_{max}=3/2; \quad f(n_{max})=(2/3)*(3/2)*(3-3/2)=3/2
\end{array}
\eeq
So the maximal value of $(2/3)*n*(3-n)$ equals $3/2$ and therefore the maximal harvesting is $h=3/2$

If $h>3/2$ there are no equilibria and the population will go
extinct. If $h<3/2$ there is a stable and non-stable equilibrium and
the population will go to the stable equilibrium. At $h=3/2$ we are at
a boundary between these two qualitatively different cases. Such a
qualitative change in system behavior is called a {\bf bifurcation}.
Bifurcations are studied in a special section of mathematics: theory
of dynamical systems.

\section{Exercises}
\subsection*{Exercises for section \ref{sec_separation}}
\ben
\item Assume that a population grows in accordance with the following equation:
$$
{dn \over dt}=1.5 n
$$
If the initial size of the population was $n(0)=30$, find what will be size at time $t=4$. Find at what time the population will double its initial size.

\item A bacterial population doubles its size each  20 minutes. The growth  of this population $N$ satisfies the differential equation ${dN \over dt}=k N$. Find the value of $k$ in $sec^{-1}$.

\subsection*{Exercises for section \ref{sec_1Dqualitative}}

\item Study the listed differential equations by answering  the following questions:
\bit

\item Draw the phase portrait.

\item How many equilibria do we have here? At which $x$?

\item For each equilibrium tell whether is it stable or non-stable

\item What will be the final value of $x$ if  $t \rightarrow \infty$.
(e.g. $x$ converges to equilibrium, or $x$ goes to infinity, etc.)

\item  List  attractor/attractors and determine  their  basin/basins
 of attraction. 
\eit
\ben
\item 
$
{dx \over dt} = -15+8x-x^2
$

\item $
{dx \over dt} = -4+ 5x -x^2 
$

\item $
{dx \over dt} = -x(x^2+x-6)
$

\item $
{dx \over dt} = 8x-x^3
$

\item $
{dx \over dt} = f(x)
$
with the following graphs of   $f(x)$:
\begin{figure}[hhh]
\centerline{
\psfig{type=pdf,ext=.pdf,read=.pdf,figure=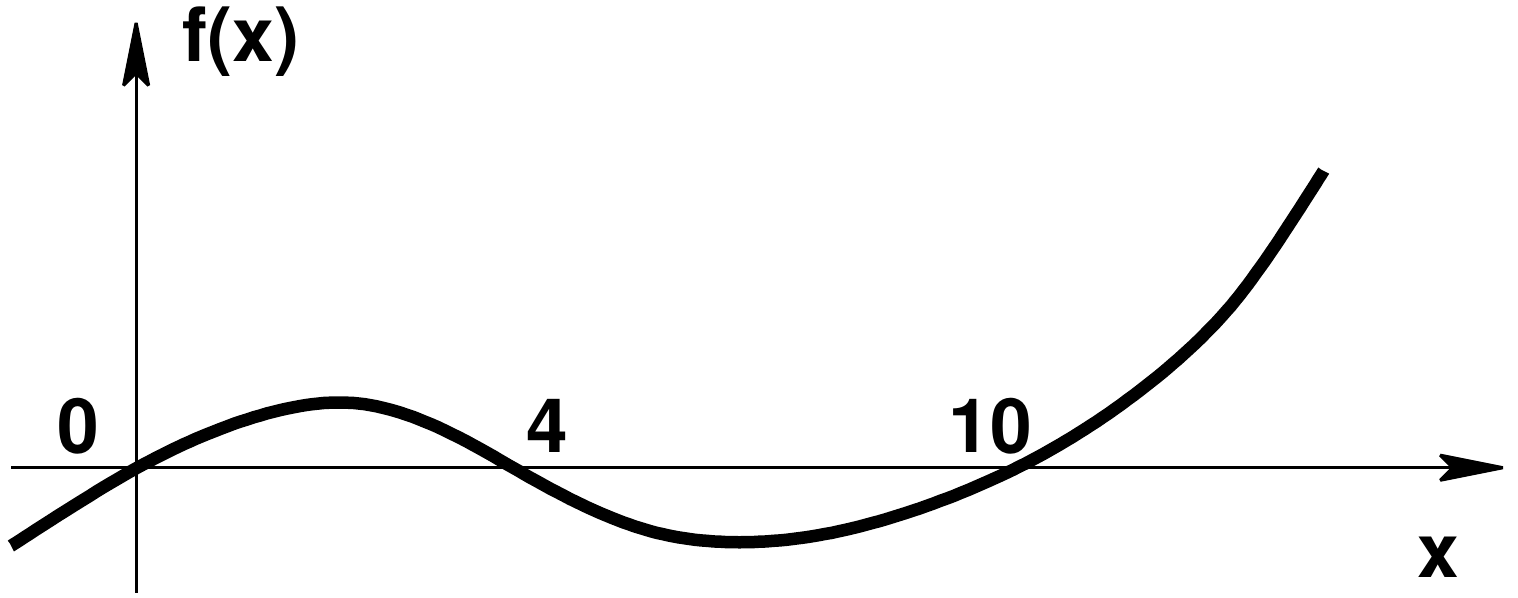,width=5cm}
\psfig{type=pdf,ext=.pdf,read=.pdf,figure=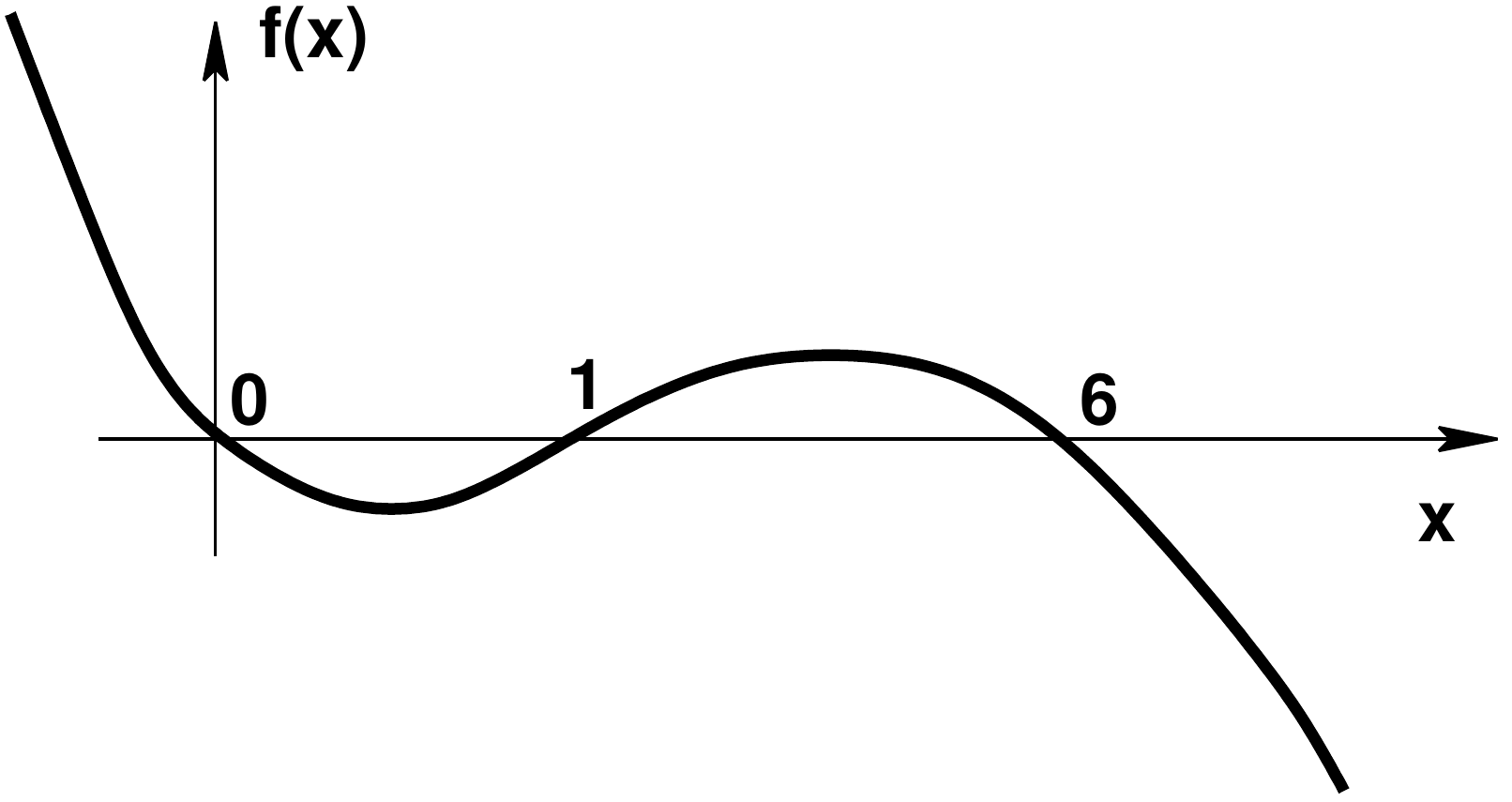,width=5cm}
}
\end{figure}

\item $
{dx \over dt}={3x^2 \over 2+x^2}-x
$  (this equation (Adler 1996) describes the dynamics of population of  species which cannot bread successfully when numbers are too small or too large)

\een

\item
The following equation describes the production of a gene product with
concentration $x$: 
$$ 
dx/dt=-x(x-0.2)(x-1)+s 
$$ 
Here $s$ is the
parameter accounting for a chemical which produces the gene
product. The initial state of the system was: $x=0$ and the value of
$s=0$.  At some moment of time the value of $s$ was slowly increased
from $s=0$ to $s=s_{max}$ and then slowly decreased back to the value $s=0$.
\ben

\item The  value of function $-x(x-0.2)(x-1)$ at its   local minimum is $-0.009$. {\it (Optional) Show this   from function derivative, or using calculator. } 

\item What will be the value of the concentration of the gene product $x$ at the
end of the described  process for $s_{max}=0.005$?

\item The same for $s_{max}=0.02$?
 
\item Show that there is a  critical value of $s_{max}$ that separates different outcomes of this process. Find this critical value of $s_{max}$.

\een

\item Consider a model population with logistic growth which is subject to harvesting at a constant rate $h$
\beq
\label{eharv1}
dn/dt=r*n*(1-n/k)-h
\eeq 

Find the maximal yield $h$.

\subsection*{Additional Exercises}

\item Assume that  the growth of a mass of an animal  can  given by the following differential equation:
$$ {dW \over dt}=400-0.3W $$ 
where W is the weight in grams and $t$
time in weeks. 
\ben
\item Find the solution for the initial $W(0)=10$. 
\item At what
time the mass will reach  half of the saturated value. 
\item If we
assume that the linear size of the animal is proportional to the cubic
root of the mass, find at what time the object will reach half of its
saturated linear size.
\een

\item The dynamics of the ionic channels in the famous Hodgkin-Huxley model for a nerve  cell is described by the following type of equations:
$$
{dm \over dt}= \alpha (1-m)-\beta m
$$
where $m$ is a gating variable and $\alpha,\beta$ are the parameters.

Find the steady state values of the gating variable $m$ and the characteristic time of approaching this steady state.

\item Consider a model where the harvesting $(hn)$  is proportional to the size of the population:

\beq
\label{eharv2}
dn/dt=r*n*(1-n/k)-hn
\eeq
Find the maximal yield.

\item Compare the harvesting strategies (\ref{eharv1}) and (\ref{eharv2}). Which strategy is better. Why?

\een

\chapter{System of two linear differential equations \label{chap_2dlin}}

Many biological systems are described by several differential
equations. One of the most simple types of such systems is  a system of two linear differential equations that on a general
form can be written as:

\beq
\label{2dlin2}
\left\{
\begin{array}{l}
{dx \over dt}=  a x+ b y \\ {dy \over dt} =   c x+ d y
\end{array}
\right.
\eeq
here $x(t)$ and $y(t)$ are unknown functions of time $t$, and $a,b,c,d$ are constants (parameters).

System (\ref{2dlin2}) by itself has many practically important
applications, for example compartmental models in biology and
pharmacology, electrical circuits in physics, models in economics,
etc.  System (\ref{2dlin2}) will also be very important for study
so-called nonlinear system of differential equations which is widely
used in theoretical biology and will be considered in chapter
\ref{chap2dnonlin}.

This chapter we will introduce main definitions for linear systems
(phase portrait and equilibria points) we will derive a formula for
general solution of this system and classify possible solutions of this
system and their phase portraits. These results will be used later to
study models of biological processes.

\section{Phase portraits and equilibria \label{sec_2dlin} }

\begin{figure}[hhh]
\centerline{
\psfig{type=pdf,ext=.pdf,read=.pdf,figure=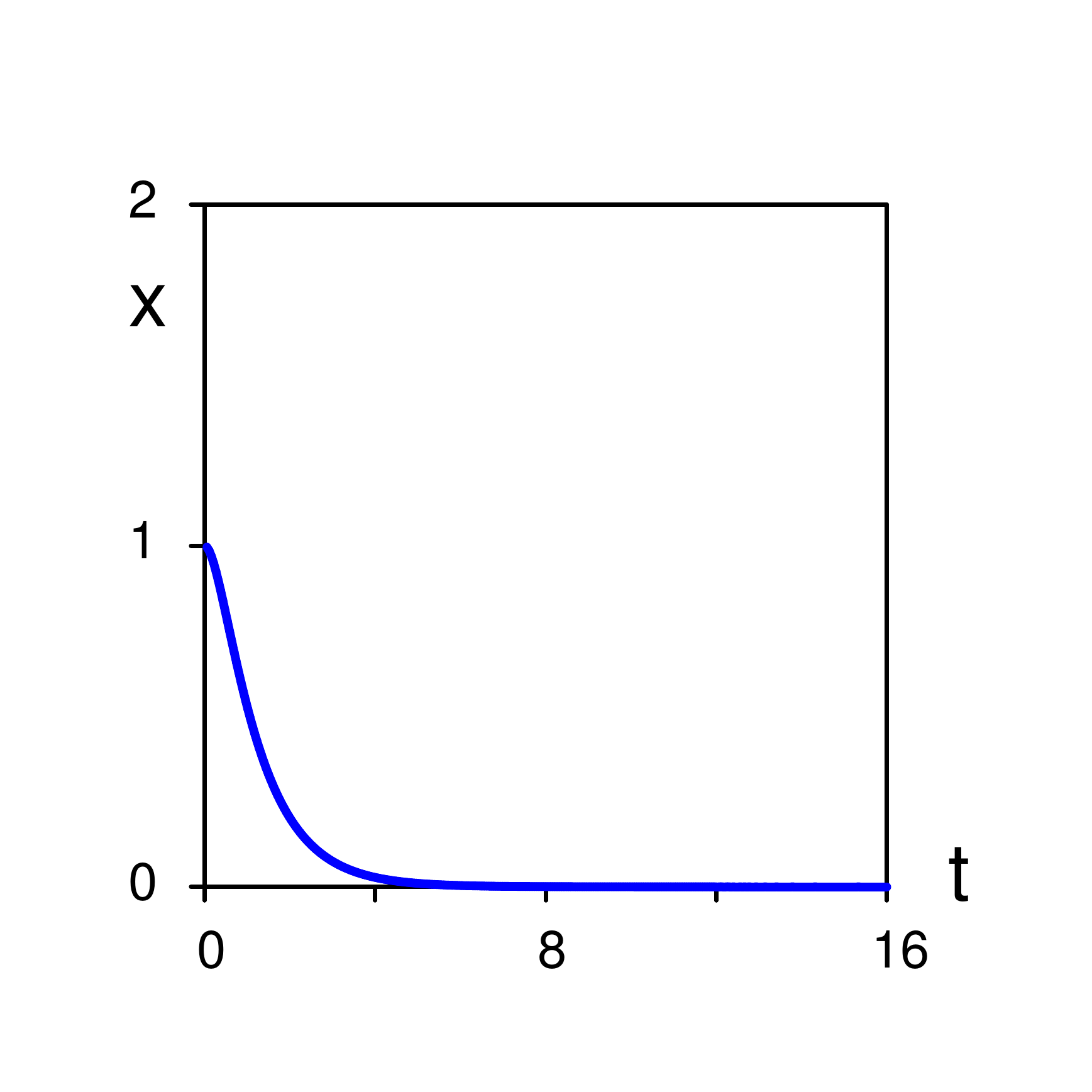,width=5cm}
\psfig{type=pdf,ext=.pdf,read=.pdf,figure=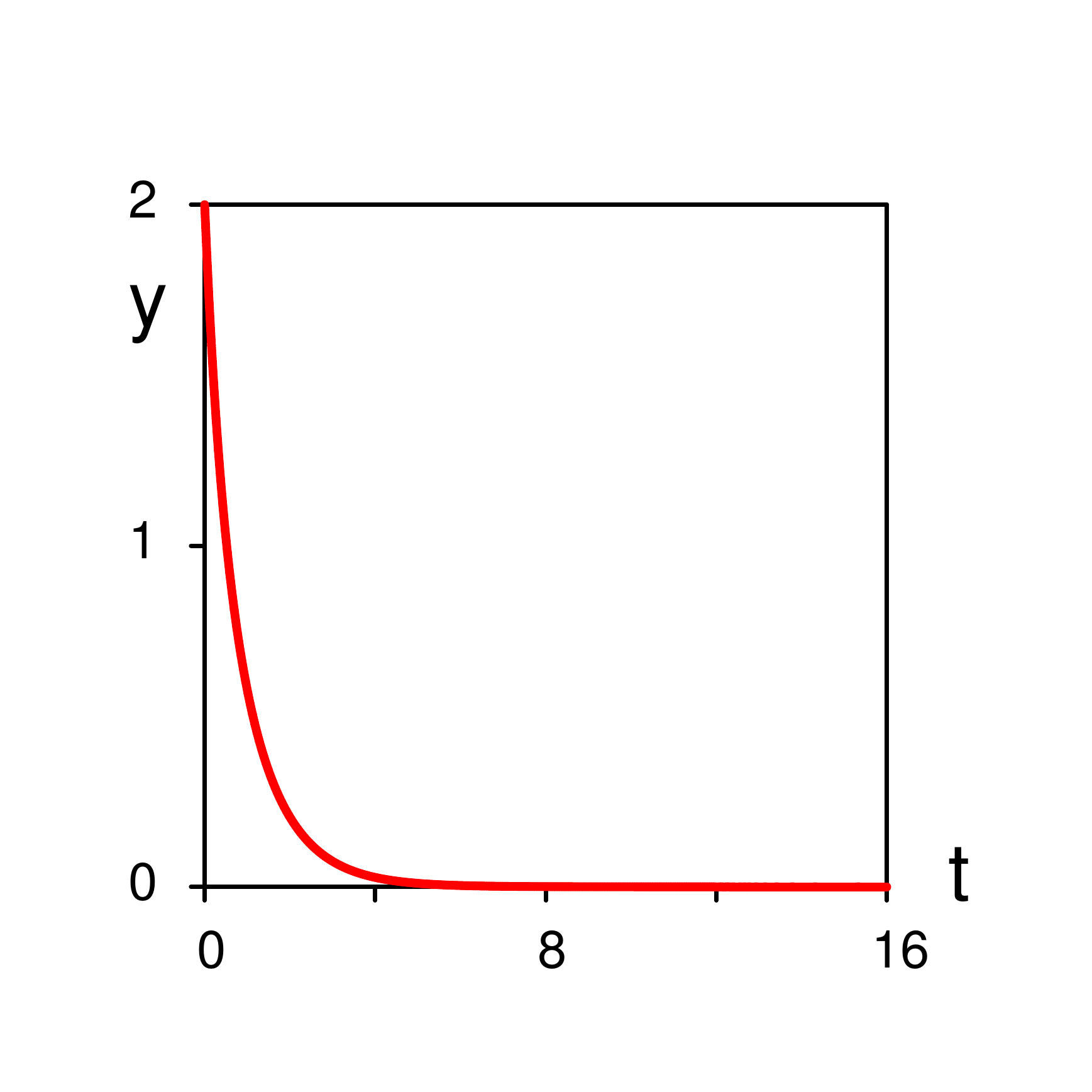,width=5cm}
}
\caption{\label{2dlin_fig1}}
\end{figure}

\begin{figure}[hbt]
\centerline{
\psfig{type=pdf,ext=.pdf,read=.pdf,figure=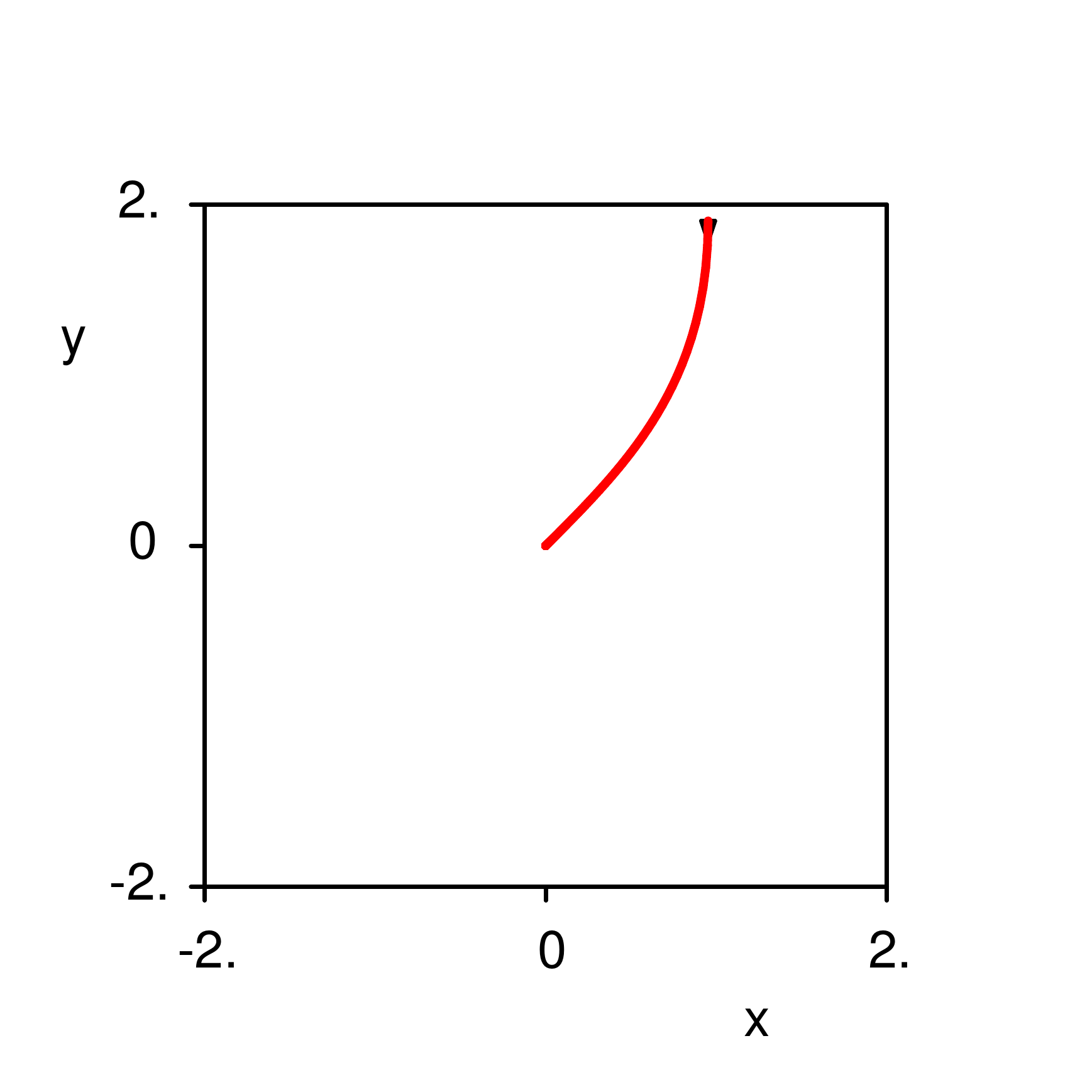,width=6.cm}
\psfig{type=pdf,ext=.pdf,read=.pdf,figure=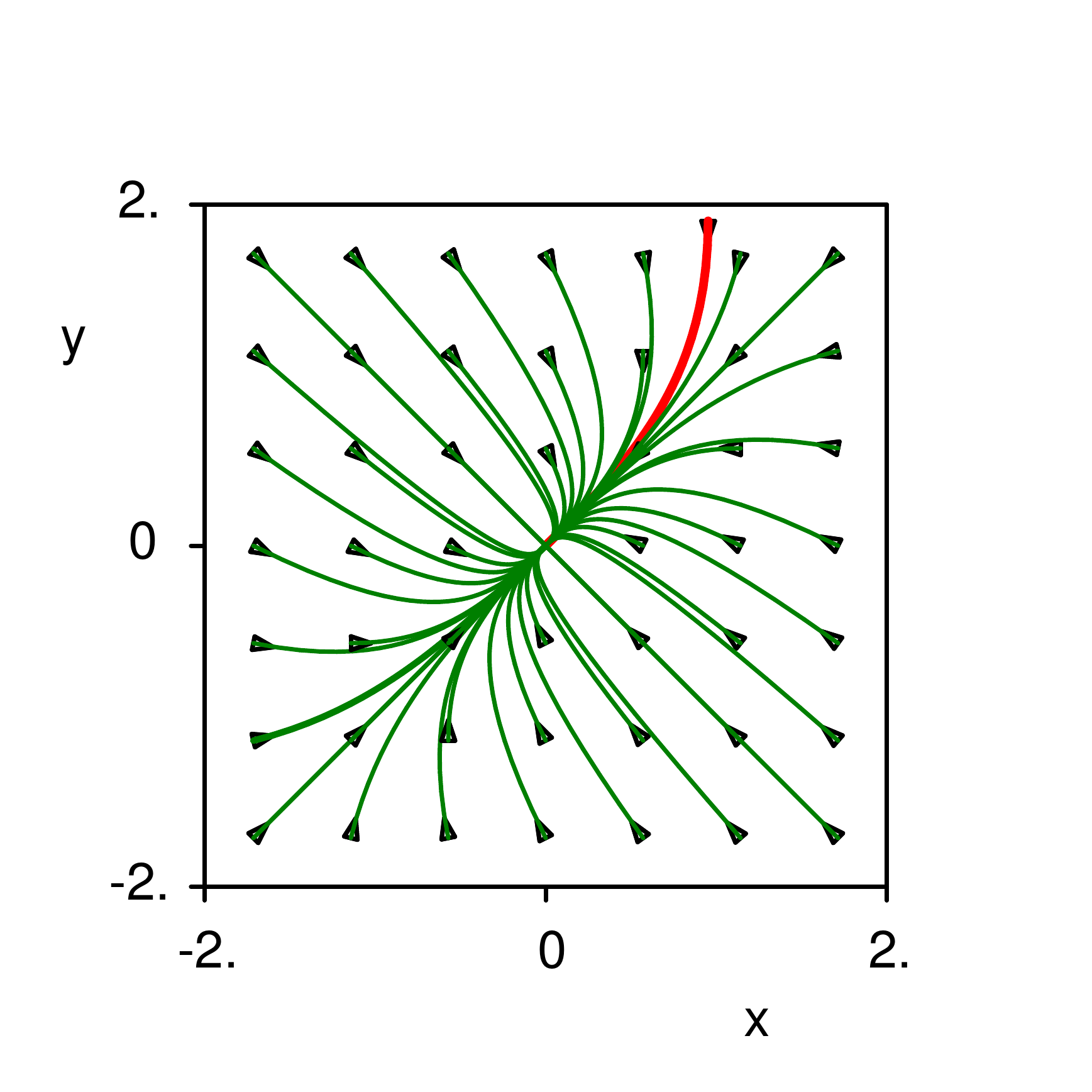,width=6.cm}
}
\caption{ \label{2dlin_fig2}} Phase portrait of system (\ref{2dlin20}): (${dx \over dt}=  -2x+  y; {dy \over dt} =   x -2 y$) found numerically.
\end{figure}

Let us consider an example of  system (\ref{2dlin2}) with particular values for the parameters $a,b,c,d$:

\beq
\label{2dlin20}
\left\{
\begin{array}{l}
{dx \over dt}=  -2x+  y \\ {dy \over dt} =   x -2 y
\end{array}
\right.
\eeq
\noindent Let us first solve this system on a computer. For that we need to
choose  initial values for $x$ and $y$ and let the computer 
find their dynamics in the course of time. Solutions for $x(0)=1,y(0)=2$ are shown in
fig.\ref{2dlin_fig1}. We see, that in the course of time, $x$ and $y$ approach the  stationary values $x=0.,y=0.$
Let us represent this solution graphically. For a differential equation
with one variable $({dx \over dt}=f(x)$) we presented the solutions in
terms of a one-dimensional phase portrait using the $x$-axis. For
system with two variables, we need to use two axes to represent the
dynamics. Let us consider a two dimensional coordinate system $Oxy$
with the $x$-axis for the variable $x$ and the $y$-axis for the
variable y. Such a coordinate system is called {\bf a phase space}.
Let us represent the trajectory from fig.\ref{2dlin_fig1} on the
$Oxy$-plane. The initial sizes of the populations were
$x(0)=1,y(0)=2$, thus we put this point ($1,2$) on the
$Oxy$-plane. At the next moment of time we get other values for $x$
and $y$ and we also put them on the $Oxy$-plane and  the $x$ and the
$y$ coordinates of the next point, etc. Finally, we will get the line shown
in fig.\ref{2dlin_fig2}a that starts at ($2,1$) and ends at 
($0,0$).  To show how $x$ and $y$ change in the course of time we draw
an arrow as in fig.\ref{2dlin_fig2}a. This trajectory is the first
element of the phase portrait. If we start trajectories from many 
different initial conditions we will get the complete phase portrait
of system (\ref{2dlin20}) (fig.\ref{2dlin_fig2}b). Each trajectory
represents a certain type of dynamics of $x(t),y(t)$, which can be
easily shown on time plots similar to fig.\ref{2dlin_fig1}. The phase
portrait in fig.\ref{2dlin_fig2}b give us the overall qualitative dynamics of
our system: the variables $x$ and $y$ approach $0,0$ from all possible initial
conditions.

The main aim of our course is to develop a procedure of drawing a
phase portrait of a general system of two differential equations
without using a computer which will allow us  to study models of biological processes. In the 1D case the phase portrait consisted of  two
main elements:  equilibria points and flows (trajectories) between
them.  Similar elements also compose the phase portrait
of a system of two differential equations. Let us
start with the first main element and define   equilibria of the system.
\vskip 1pc
\noindent In the 1D case  equilibria were  points where our system is
stationary: placed at an equilibrium point the system will stay there
forever. Mathematically equilibria for eq. ${dx \over dt} =f(x)$
(\ref{1dgeneral}) were determined as the points where ${ dx \over dt}
=0$, i.e. where $f(x)=0$. In the 2D case it is 
required  that at the equilibrium point both variables $x$ and $y$ do not change their values, i.e. both ${dx \over dt}=0$ and ${dy \over dt}=0$. For system   (\ref{2dlin20}) these conditions give the following the system of two algebraic equations for finding equilibria:

t
\beq
\label{eqlin}
\left\{
\begin{array}{l}
{dx \over dt} =0=-2x+  y \\ {dy \over dt} =0=x -2 y
\end{array}
\right.
\begin{array}{l}
\Rightarrow
\end{array}
\left\{
\begin{array}{l}
y=2x \\ x-2y=0
\end{array}
\right.
\begin{array}{l}
\Rightarrow
\end{array}
\left\{
\begin{array}{l}
y=2x \\ x-2*2x=0
\end{array}
\right.
\begin{array}{l}
\Rightarrow
\end{array}
\left\{
\begin{array}{l}
y=2x \\ -3x=0
\end{array}
\right.
\eeq
that have the only solution $x=0,y=0$. Therefore system (\ref{2dlin20}) has an equilibrium point $(0,0)$. As we  see   in fig.\ref{2dlin_fig2}b this equilibrium  is an attractor for  all trajectories.

\vskip 1pc
\noindent For a general linear system (\ref{2dlin2}) the equilibria will be
given by:
\beq
\label{eqlin2}
\left\{
\begin{array}{l}
ax+ b y=0 \\ cx +d y=0
\end{array}
\right.
\eeq
This system always has a solution $x=0,y=0$ and thus the general
linear system (\ref{2dlin2}) always has an equilibrium at the point
$x=0,y=0$. In the next sections we will find out how to sketch a phase
portrait of (\ref{2dlin2}) around this equilibrium.
Our plan will be the following. We will first find the general
analytical solution of this system and then will use it to draw the
phase portraits.

\section{General solution of linear system \label{sec_2dlin2}}
Consider a general system of two differential equations with constant coefficients:
\beq
\label{2dlin201}
\left\{
\begin{array}{l}
{dx \over dt}=  a x+ b y \\ {dy \over dt} =   c x+ d y
\end{array}
\right.
\eeq
 The general  solution of (\ref{2dlin201}) can be written in the following form
\beq
\label{2dlin_sol}
\left( \begin{array}{c} 
x \\  y \end{array} \right) = 
C_1 \left( \begin{array}{c} v_{1x}
 \\  v_{1y}  \end{array} \right)e^{\lambda_1 t } 
+ C_2 \left( \begin{array}{c} v_{2x}
 \\  v_{2y}  \end{array} \right)e^{\lambda_2 t }
\eeq
where $\lambda_1,\lambda_2$ are eigen values of the matrix 
$A=\left(\begin{array}{lr}  a & b\\
                         c &  d
\end{array} \right)$, and 
$\left( \begin{array}{c} v_{1x}
 \\  v_{1y}  \end{array} \right)$,
$\left( \begin{array}{c} v_{2x}
 \\  v_{2y}  \end{array} \right)$
are the corresponding eigen vectors.

We will not derive the formula (\ref{2dlin_sol}), we will show the
main ideas behind the real derivation. For that let us consider first
the one dimensional case, and then extend it to a two dimensional
system.

The easiest way to find solutions of system (\ref{2dlin201})  is by  the method of substitution. Let us illustrate this method on example of 1D analog of system (\ref{2dlin201}) which is 1D linear differential equation:

\beq
\label{1dlin}
{dx \over dt}=ax
\eeq

We can easily solve (\ref{1dlin}) using the direct method of
separation of variables and subsequent integration. However, let us
find the solution using another method, the method of
substitution. The main idea of this method is to look for a solution
in some known class of functions which should be chosen in advance
from some preliminary analysis. It was found that for linear systems
this class is the class of exponential functions $Ce^{\lambda
t}$. Important questions such as: how was this class found and is this
class unique etc, will not be discussed. Our aim here will be
an illustration of the main components of the solution rather than
comprehensive analysis of linear systems, which is a large special
section of mathematics. Once the class of functions is given (in our case the class of exponential functions), we need to
check under which circumstances it will satisfy the equations we are
solving. Thus we will look for a solution of (\ref{1dlin}) of the form
$x=Ce^{\lambda t}$, where $C$ and $\lambda$ are unknown
coefficients. The main idea of the method of substitution is to find
these unknown coefficients for a particular system. Let us substitute
$x=Ce^{\lambda t}$ into (\ref{1dlin}). We find: 
 ${dx \over dt}=(Ce^{\lambda t})' =\lambda C e^{\lambda t} $ , or:   

$$
\lambda C e^{\lambda t}=a C e^{\lambda t}
$$
We can cancel $e^{\lambda t}$ and $C$, and we get:
\beq
\lambda =a 
\eeq
 Hence we found the following the solutions of (\ref{1dlin}):
\beq
x= C e^{a t} 
\eeq 
where $C$ is an arbitrary constant. 

Now, let us use the same approach for the two dimensional system (\ref{2dlin201}). It turns out that the class of functions in two dimensions will be the same as in  one dimension, but because we have two variables, we need to introduce different constants for $x$ and $y$, so our substitution will be
\beq
\label{2dsubs}
x=C_x e^{\lambda t};y=C_y e^{\lambda t}\quad or \quad 
\left( \begin{array}{c} 
x \\  y \end{array} \right) = 
 \left( \begin{array}{c} C_x
 \\  C_y  \end{array} \right)e^{\lambda t } 
\eeq 
where $C_x,C_y,\lambda$ are unknown coefficients.
Let us make this substitution for a particular example:
\beq
\label{2dlinour}
\left( \begin{array}{c} 
{dx \over dt} \\  {dy \over dt}  \end{array} \right) = 
\left(\begin{array}{lr}  1 & 4\\
                         1 &  1
\end{array} \right)  
\left( \begin{array}{c} x
 \\  y  \end{array} \right)
\eeq
Substitution for ${dx \over dt}=\lambda C_x e^{\lambda t},{dy \over dt}=\lambda C_y e^{\lambda t},x=C_x e^{\lambda t},y=C_y e^{\lambda t}$ gives:
\beq
\left\{
\begin{array}{l}
 \lambda C_x e^{\lambda t}=C_x e^{\lambda t}+4C_y e^{\lambda t} \\
 \lambda C_y e^{\lambda t}=C_x e^{\lambda t}+C_y e^{\lambda t}
\end{array}
\right.
\eeq
we can cancel $e^{\lambda t}$, (but not $C_x,C_y$ as in one dimensional case), and get:
\beq
\left\{
\begin{array}{l}
 \lambda C_x =C_x +4C_y  \\ 
\lambda C_y =C_x +C_y
\end{array}
\right.
\eeq
or in the matrix form:
\beq
\label{eigen}
\lambda \left( \begin{array}{c} 
C_x \\  C_y  \end{array} \right) = 
\left(\begin{array}{lr}  1 & 4\\
                         1 &  1
\end{array} \right)  
\left( \begin{array}{c} C_x
 \\  C_y  \end{array} \right)
\eeq
We see that to find the solution of (\ref{2dlinour}) we need to solve the problem (\ref{eigen}), which  is the eigen value problem considered in section \ref{basicEigen}.  To solve it  we  find eigen values from the characteristic equation (\ref{eigen_char}):
 $ Det \left|\begin{array}{lr} 1-\lambda & 4\\ 1 & 1-\lambda
\end{array} \right|=(1-\lambda)*(1-\lambda) -4=\lambda^2-2\lambda-3=0 \;\;  \lambda_{12}={2 \pm  \sqrt {2^2+3*4} \over 2}=1 \pm  {1 \over 2} \sqrt {16}=1 \pm 2$. Hence $\lambda_1=-1,\lambda_2=3$

We use the express  formula (\ref{eigen_ex}) for eigen vectors and get:
\beq
\label{eig35_0}
\lambda_1=-1; \left( \begin{array}{c} v_{1x}
 \\  v_{1y}  \end{array} \right)=
\left( \begin{array}{c} -4
 \\  1-(-1)  \end{array} \right)
=\left( \begin{array}{c} -4
 \\  2 \end{array} \right)
 \quad
\lambda_2=3;
\left( \begin{array}{c} v_{2x}
 \\  v_{2y}  \end{array} \right)=
\left( \begin{array}{c} -4
 \\  -2  \end{array} \right).
\eeq

Formula (\ref{eig35_0}) gives just one eigen vector for each eigen value. We also know (see formula (\ref{tmp14} in section \ref{basicEigen}) that all eigen vectors can be found by multiplication of this eigen vector by an arbitrary constant. If we denote by  $C_1$ the constant for the first eigen vector and by $C_2$ the constant for the second eigen vector we will get the following solution  of  the eigen value problem  (\ref{eigen}):
\beq
\label{eig35}
\lambda_1=-1; \left( \begin{array}{c} v_{1x}
 \\  v_{1y}  \end{array} \right)
=C_1\left( \begin{array}{c} -4
 \\  2 \end{array} \right)
 \quad
\lambda_2=3;
\left( \begin{array}{c} v_{2x}
 \\  v_{2y}  \end{array} \right)=
C_2\left( \begin{array}{c} -4
 \\  -2  \end{array} \right).
\eeq

\noindent If we substitute these eigen vectors into the formula (\ref{2dsubs}) we find the following solutions of (\ref{2dlinour})
\beq
\label{sol2dlinour}
\left( \begin{array}{c} 
x \\  y \end{array} \right) = 
C_1 \left( \begin{array}{c} -4
 \\  2  \end{array} \right)e^{-1* t }
\quad
\left( \begin{array}{c} 
x \\  y \end{array} \right) = 
C_2 \left( \begin{array}{c} -4
 \\  -2 \end{array} \right)e^{3* t }
\eeq
\noindent Now let us prove the following property of system (\ref{2dlin201}): \\
$\bullet$  If $x_1,y_1$ and $x_2,y_2$ are two solutions of (\ref{2dlin201}), then $x_1+x_2,y_1+y_2$ is also a solution of  (\ref{2dlin201}).\\
{\it Proof:} As $x_1,y_1$ and $x_2,y_2$ are the solution this means that they satisfy (\ref{2dlin201}), i.e.
\beq
\left\{
\begin{array}{l}
{dx_1 \over dt}=  a x_1+ b y_1 \\ {dy_1 \over dt} =   c x_1+ d y_1
\end{array}
\right.
\qquad
\left\{
\begin{array}{l}
{dx_2 \over dt}=  a x_2+ b y_2 \\ {dy_2 \over dt} =   c x_2+ d y_2
\end{array}
\right.
\eeq
If we  add equations for ${dx_1 \over dt}$ and ${dx_2 \over dt}$ we get:
${dx_1 \over dt}+{dx_2 \over dt}=  a x_1+ b y_1+ a x_2+ b y_2$, which can be re-written as: ${d (x_1+x_2) \over dt}=  a (x_1+x_2)+ b(y_1+y_2)$. If we do the same for  equations for ${dy_1 \over dt}$ and ${dy_2 \over dt}$ we will finally get:
\beq
\left\{
\begin{array}{l}
{d (x_1+x_2) \over dt}=  a (x_1+x_2)+ b(y_1+y_2)\\ {d(y_1+y_2) \over dt} =  c (x_1+x_2)+ d(y_1+y_2)
\end{array}
\right.
\eeq
which explicitly shows  that  $x_1+x_2,y_1+y_2$ is a solution of  (\ref{2dlin201}).

\noindent If we apply this result for two found solutions (\ref{sol2dlinour}) of  (\ref{2dlinour}) we can conclude that  the sum of these two solutions is also a solution of (\ref{2dlinour}):
\beq
\label{2dlinour_sol}
\left( \begin{array}{c} 
x \\  y \end{array} \right) = 
C_1 \left( \begin{array}{c} -4
 \\  2  \end{array} \right)e^{-1* t }
+ 
C_2 \left( \begin{array}{c} -4
 \\  -2 \end{array} \right)e^{3* t }
\eeq

We proved formula (\ref{2dlin_sol}) for a particular system.  If we
apply the same steps for a general system (\ref{2dlin201}) we will get
the general solution in the form (\ref{2dlin_sol}).

\noindent So, we solved system (\ref{2dlin201}). In the next sections we will find out how to draw its  phase portraits.

\section{Real eigen values. Saddle, node. \label{realEV}}

As we know, the general solution of (\ref{2dlin2}) is given by
(\ref{2dlin_sol}). We can use this formula  to sketch a phase portrait of this
system. It turns out, that we can have several different 
types of equilibria depending on the eigen values
$\lambda_1,\lambda_2$.  As we know $\lambda_1,\lambda_2$ are the roots
of the characteristic equation (\ref{eigen_char}), which is a general
quadratic equation. Therefore, the roots can be real or complex
numbers. In this section we consider the case of real roots.

So, assume that the eigen values $\lambda_1$ and $\lambda_2$ are real. This yields  the following three cases:
\ben
\item Eigen values have different signs ($\lambda_1<0;\lambda_2>0$, or $\lambda_1>0;\lambda_2<0$).
\item Both eigen values are positive ($\lambda_1>0;\lambda_2>0$)
\item Both eigen values are negative ($\lambda_1<0;\lambda_2<0$) 
\een
Note, that we do not consider the case when $\lambda=0$. This situation is quite rare   and is not considered  in this course.
\subsection{Saddle; $\lambda_1<0;\lambda_2>0$, or $\lambda_1>0;\lambda_2<0$}
System (\ref{2dlinour}), which we solved in section \ref{sec_2dlin2}, had  eigen values 
$\lambda_1=-1,\lambda_2=3$.  Let us draw its  phase portrait.
The general solution of this system is given by (\ref{2dlinour_sol})
\beq
\left( \begin{array}{c} 
x \\  y \end{array} \right) = 
C_1 \left( \begin{array}{c} -4
 \\  2  \end{array} \right)e^{-1* t }
+ 
C_2 \left( \begin{array}{c} -4
 \\  -2 \end{array} \right)e^{3* t }.
\eeq
Because $C_1,C_2$ are arbitrary constants let us consider  three simple cases, in which 
one of these constants, or both of them are zero.

1) If $C_1=0,C_2=0$, then $x=0,y=0$ and do not depend on time. The trajectory  is just one point (0,0), which is the equilibrium point of the system (\ref{2dlinour}).

\begin{figure}[h]
\centerline{
\psfig{type=pdf,ext=.pdf,read=.pdf,figure=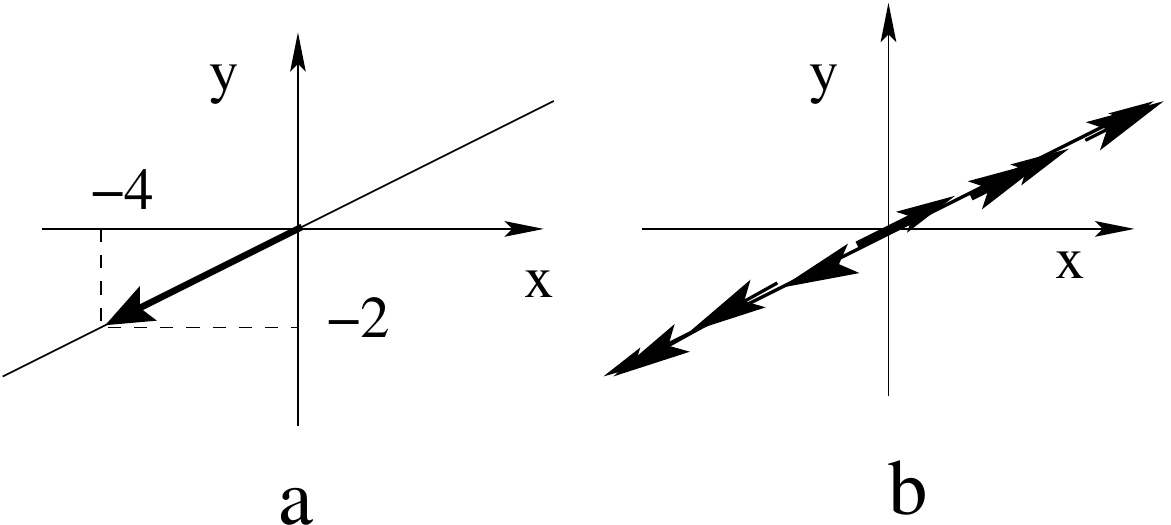,width=11cm}
}
\caption{\label{fig6.1}}
\end{figure}

\begin{figure}[h]
\centerline{
\psfig{type=pdf,ext=.pdf,read=.pdf,figure=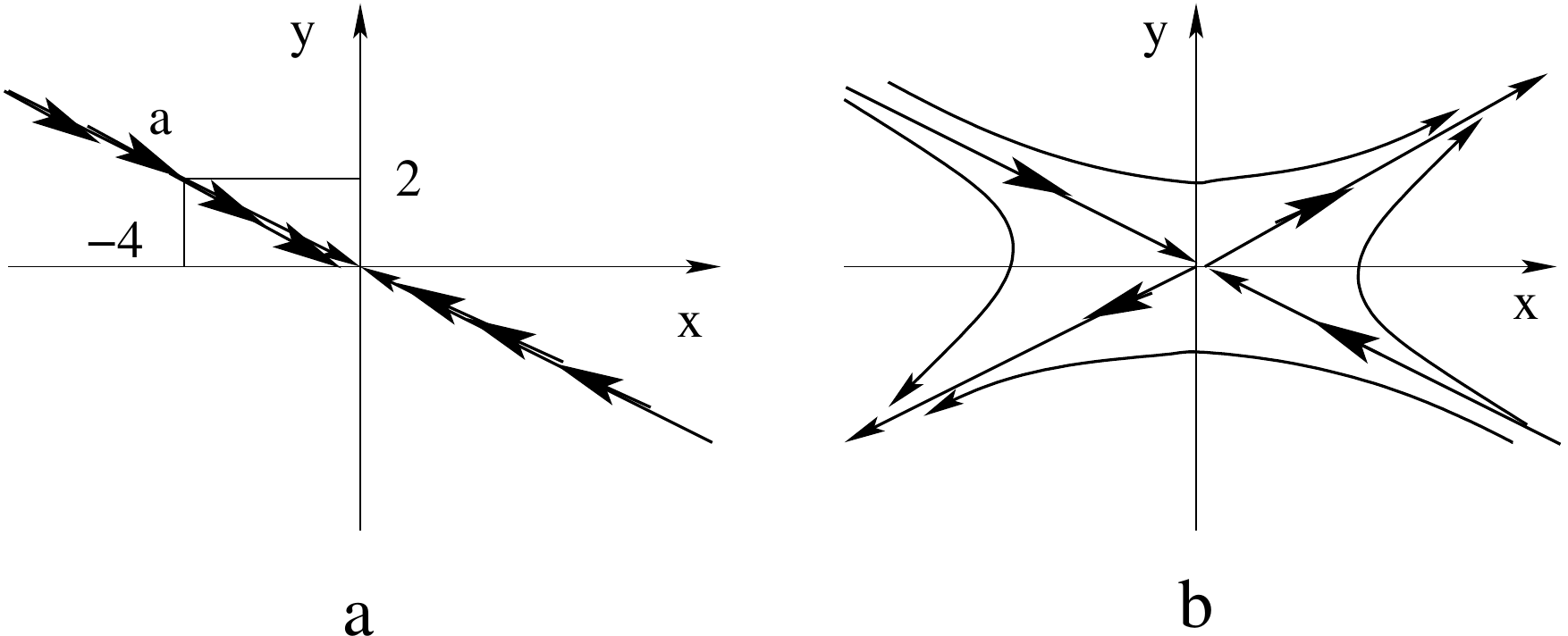,width=11.cm}
}
\caption{\label{fig6.2}}
\end{figure}

2) If $C_1=0,C_2=$ {\it any number}, then
$\left( \begin{array}{c} 
x \\  y \end{array} \right) = 
C_2 \left( \begin{array}{c} 2
 \\  1  \end{array} \right)e^{3* t } $. Because $e^{3* t }$ can change from 1 (at $t=0$) to any infinitely  large number and $C_2$ is an arbitrary positive or negative number, this expression can be rewritten as:
\beq
\label{tmp20}
\left( \begin{array}{c} 
x \\  y \end{array} \right) = 
C_2 \left( \begin{array}{c} -4
 \\  -2  \end{array} \right)e^{3* t } = \left( \begin{array}{c} -4
 \\  -2  \end{array} \right)K={\bf V_2}K
\eeq
where $K$ is an arbitrary number from $-\infty < K < \infty$ and ${\bf
V_2}$ is a vector $ \left( \begin{array}{c} -4 \\ -2 \end{array}
\right)$. Thus  expression (\ref{tmp20}) means   multiplying of the 
eigen vector ${\bf V_2}$ by an arbitrary number $K$. In general, if we
multiply a vector by a positive number $K$ we get a vector with the
same direction but the length will be increased $K$ times. If we
multiply the vector by a negative number, the direction of the vector
will be changed to the opposite and the length will be changed by a
factor $|K|$. Because in (\ref{tmp20}) $K$ assumes all values from
$-\infty < K < \infty$, this will give a straight line along this
vector $ \left( \begin{array}{c} 2 \\ 1 \end{array}
\right)$ (fig.\ref{fig6.1}a). To complete drawing  the
trajectory we need to show  arrows indicating the motion of a point along
the trajectory in the course of time. Because time dynamics is given
by $e^{3t}$,  $|K|$ in (\ref{tmp20}) will grow in the
course of time, i.e. it will become either more  positive, or
more  negative depending on its initial sign. Geometrically
this means that a point will move apart from the origin of the $Oxy$
coordinate system and we will get a picture as in fig.\ref{fig6.1}b.
3) The third case is $C_1= \;${\it any number}, $C_2=0$. The solution in this case is
\beq
\left( \begin{array}{c} 
x \\  y \end{array} \right) = 
C_1 \left( \begin{array}{c} -4
 \\  2  \end{array} \right)e^{-t }.
\eeq
As in previous case we conclude, that all trajectories in this case will be located on a straight line along the vector  $ \left( \begin{array}{c} -4 \\  2  \end{array} \right)$ and we just need to show the direction of flow along this line.
In this case, time dynamics is given by the function $e^{-t }$, which approaches zero when $t$ goes to infinity. Therefore, independent of initial conditions (independent of the value of $C_1$) we will approach the point $x=0,y=0$, and  the arrows will have the following direction (fig.\ref{fig6.2}).
Finally let us draw the phase portrait for arbitrary $C_1$ and
$C_2$ (fig.\ref{fig6.2}b). Let us consider one
trajectory for which $C_1 \neq 0;C_2 \neq 0$, for example: $C_1=0.1, \;
C_2=0.1$. The solution in this case is given by:
\beq
\label{aa2}
\left( \begin{array}{c} 
x \\  y \end{array} \right) = 
0.1 \left( \begin{array}{c} -4
 \\  2  \end{array} \right)e^{-t }
+ 
0.1 \left( \begin{array}{c} -4
 \\  -2  \end{array} \right)e^{3t }
\eeq
This trajectory starts as  point $x=-0.4 + 0.2=-0.2; \;
y=0.2-0.1=0.1$. In the course of time $e^{-t }$
will become smaller and smaller, while $e^{3t }$ will grow. So, the
first term in (\ref{aa2}) $0.1 \left( \begin{array}{c} -4 \\ 2 \end{array}
\right)e^{-t }$ will be small compared to the second term $0.1 \left(
\begin{array}{c} -4 \\ -2 \end{array} \right)e^{3t }$, and the dynamics
at large $t$ will be described by an approximate formula : $\left(
\begin{array}{c} x \\ y \end{array} \right) \approx 0.1 \left(
\begin{array}{c} -4 \\ -2 \end{array} \right)e^{3t }$, thus the trajectory   will approach the line (\ref{tmp20}) presented in fig.\ref{fig6.1}b. There
will be similar behavior for any other trajectories: independently on
staring points they will approach  line of fig.\ref{fig6.1}b
from various directions.  The qualitative picture will be as in
fig.\ref{fig6.2}b.

Such a phase portrait is called a {\bf {saddle point}}. It has the
following important features: (1) There is an equilibrium point at
$x=0,y=0$. (2) There are two lines associated with eigen vectors of our
system (fig.\ref{fig6.1}b,fig.\ref{fig6.2}a). These lines are called
manifolds. The manifolds in fig.\ref{fig6.1}b and fig.\ref{fig6.2}a
are different. If we follow the trajectory along the manifold in
fig.\ref{fig6.1}b the distance to the equilibrium {\it {increases}}
(see fig.\ref{fig6.2}). On the contrary, if we follow the trajectory
along the manifold in fig.\ref{fig6.2}a the distance to the
equilibrium {\it {decreases}}. The manifold from fig.\ref{fig6.1}b is
called {\bf {a non-stable manifold}}. The manifold from
fig.\ref{fig6.2}a is called {\bf {a stable manifold}}.

\begin{figure}[h]
\centerline{
\psfig{type=pdf,ext=.pdf,read=.pdf,figure=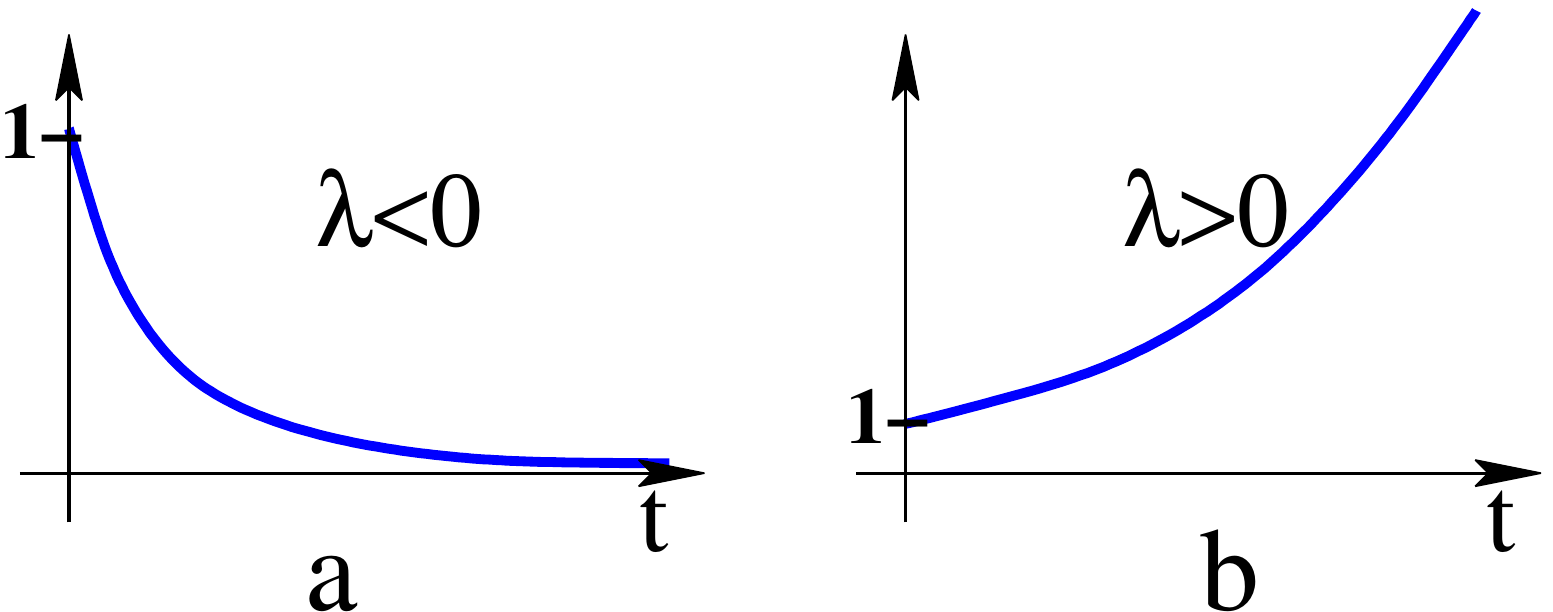,width=10.cm}
}
\caption{\label{fig6.3}}
\end{figure}

Conclusions of this study can be easily generalized. If we consider an expression:
\beq
\label{manifold}
\left( \begin{array}{c} 
x \\  y \end{array} \right) = 
C \left( \begin{array}{c} v_{x}
 \\  v_{y}  \end{array} \right)e^{\lambda t },
\eeq 
it will obviously determine a manifold (straight line) along the
vector $\left( \begin{array}{c} v_{x} \\ v_{y} \end{array} \right)$
and the stability of this manifold will be determined by $e^{\lambda t
}$.  There are two main types of behavior of the function $e^{\lambda
t }$ (fig.\ref{fig6.3}). If $ \lambda<0, e^{\lambda t }$ approaches
zero, when $t$ increases.  If $ \lambda>0, e^{\lambda t }$ grows to
infinity with increasing $t$. Hence, if $ \lambda<0$,
eq.(\ref{manifold}) will determine a stable manifold: $x,y$ will
approach $0$ in the course of time (as in fig.\ref{fig6.2}a). If $
\lambda>0$, then $x,y$ will diverge to infinity and we will get a
non-stable manifold.
\bec
The equation (\ref{manifold}) on a phase portrait gives a manifold in
the form of a straight line. This line goes through the origin and is
directed along the vector $ \left( \begin{array}{c} v_{x} \\ v_{y}
\end{array} \right) $. This manifold is stable if $ \lambda<0$ and
non-stable if $\lambda>0$.
\eec
Finally the formal definition of a saddle point:
\bec
 Fig.\ref{fig6.2} shows the phase portrait of a saddle point. It occurs close to equilibrium, at which eigen
 values of the system are real and have different signs, i.e.  $
 \lambda_1<0; \lambda_2>0$, or $ \lambda_1>0; \lambda_2<0$. The phase
 portrait of a saddle point has two manifolds directed along the eigen
 vectors. One manifold is stable
 (corresponding to the negative eigen value of the system). The other
 manifold is non-stable (corresponding to the positive eigen value of
 the system).
\eec
\subsection{Non-stable node;  $\lambda_1>0;\lambda_2>0$ \label{nnode_sec}}
Let us draw the phase portrait for the case when eigen values are real and are both positive. The general solution of the system is given by (\ref{2dlin_sol}):
\beq
\label{node1}
\left( \begin{array}{c} 
x \\  y \end{array} \right) = 
C_1 \left( \begin{array}{c} v_{1x}
 \\  v_{1y}  \end{array} \right)e^{\lambda_1 t } 
+ C_2 \left( \begin{array}{c} v_{2x}
 \\  v_{2y}  \end{array} \right)e^{\lambda_2 t }
\eeq

\begin{figure}[h]
\centerline{
\psfig{type=pdf,ext=.pdf,read=.pdf,figure=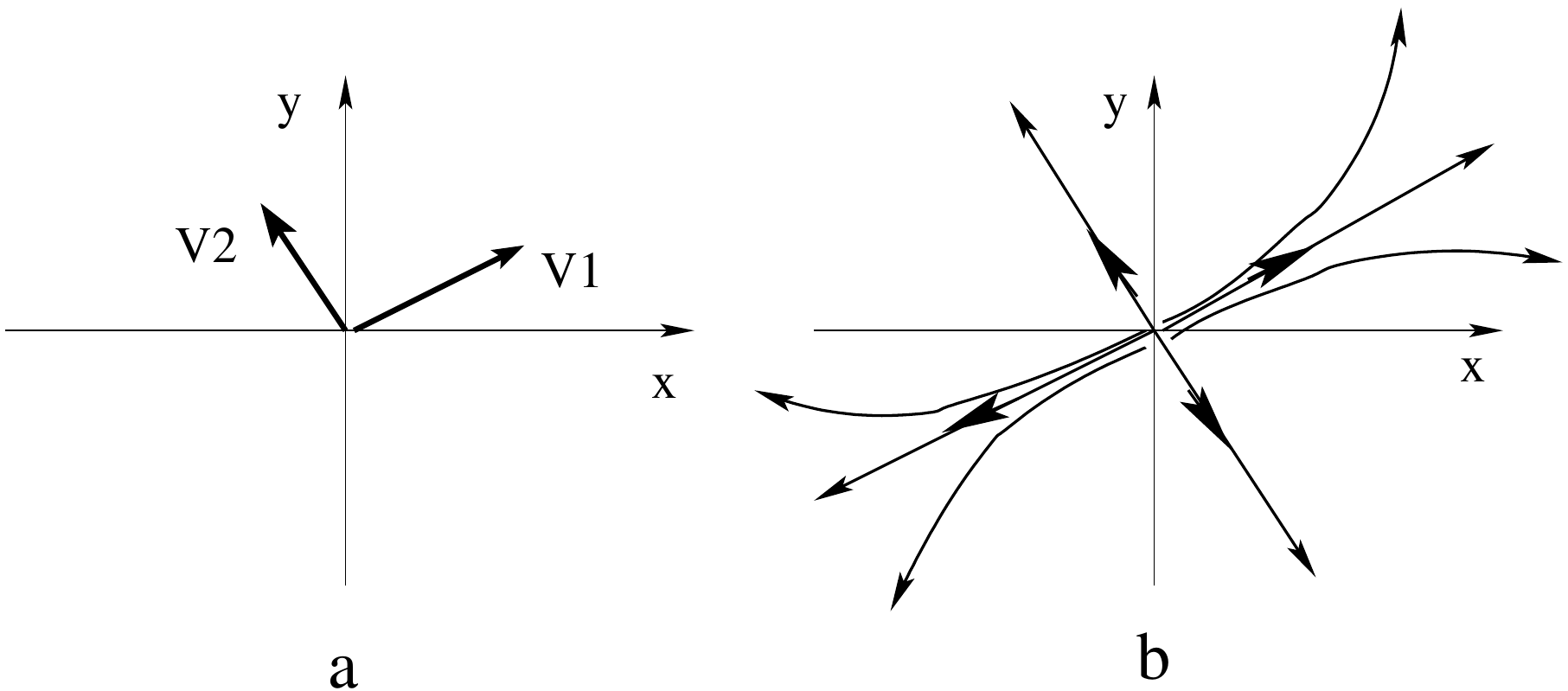,width=12.5cm}
}
\caption{\label{fig6.4}}
\end{figure}

From the previous analysis we immediately conclude, that the phase portrait in this case has the equilibrium point at $(0,0)$ and two unstable manifolds along the vectors $ \left( \begin{array}{c} v_{1x}
 \\  v_{1y}  \end{array} \right)$ and
 $\left( \begin{array}{c} v_{2x}
 \\  v_{2y}  \end{array} \right)$ (fig.\ref{fig6.4}a). Let us put that on the  graph (fig.\ref{fig6.4}b).

To complete the picture we need to add several trajectories which
start between the manifolds. For such trajectories $C_1 \neq 0\;C_2
\neq 0$ and both terms in (\ref{node1}) will diverge to plus or minus
infinity. So, we get trajectories as in fig.\ref{fig6.4}b. Such an equilibrium
is called {\bf {a non-stable node}}.
\bec
 If the eigen values of system (\ref{2dlin2}) are real and both positive ($\lambda_1 >0,\lambda_2>0$) we have an equilibrium point
 called a non-stable node. To draw a phase portrait at this
 equilibrium we need to show two non-stable manifolds along the eigen
 vectors of system (\ref{2dlin2}) and add several diverging
 trajectories between the manifolds.
\eec
\subsection{Stable node;  $\lambda_1<0;\lambda_2<0$ \label{snode_sec}}

\begin{figure}[h]
\centerline{
\psfig{type=pdf,ext=.pdf,read=.pdf,figure=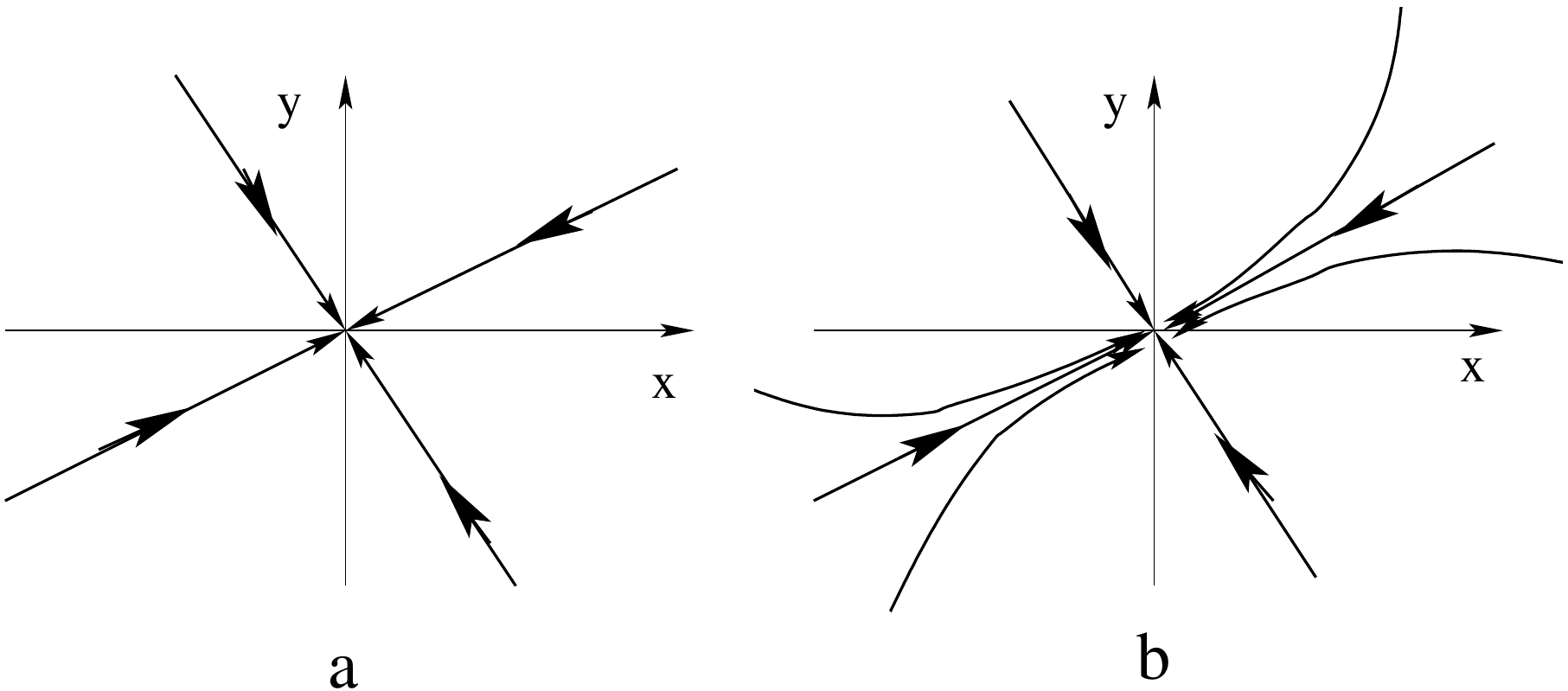,width=12.5cm}
}
\caption{\label{fig6.5}}
\end{figure}

The general solution in this case has the same form  (\ref{node1}). The phase portrait will be similar to fig.\ref{fig6.4}, but because   $\lambda_1<0;\lambda_2<0$ both manifolds will be stable.
So we get a picture fig.\ref{fig6.5}a

If the trajectory starts between the manifolds  $(C_1 \neq 0\;C_2 \neq 0)$ it will also approach equilibrium as both terms in  (\ref{node1}) will converge to 0, because  $\lambda_1<0;\lambda_2<0$ (fig.\ref{fig6.5}b). Such an equilibrium  is called {\bf {a stable node}}.
\bec
If the eigen values of system  (\ref{2dlin2}) are real and both negative ($\lambda_1 <0,\lambda_2<0$) we have an equilibrium point called  {\bf {stable node}}. To draw a phase portrait at this equilibrium we need to show two stable manifolds along the eigen vectors of system  (\ref{2dlin2}) and add several trajectories converging to the equilibrium (0,0).

\eec

\section{Phase portraits for complex eigen values: spiral, center \label{complexEV}}

In the previous section we have studied the case when the roots of the
characteristic equation (\ref{eigen_char}) are real, and found three
possible types of phase portrait (equilibria): saddle, stable node and
unstable node. Here we will study the case when the roots of the
characteristic equation (\ref{eigen_char}) are complex.

\subsection {General ideas on equilibria with complex eigenvalues}
From section \ref{sec_2dlin2} we know that in order to find the type of
phase portrait of linear system we need to solve the  characteristic
equation (\ref{eigen_char}) which is a general quadratic equation, and
because the parameters $a,b,c,d$ of the linear system are arbitrary the coefficients of the
characteristic  equation are also arbitrary. Therefore, it may happen that the discriminant of this quadratic equation will be negative and we will have  complex eigen values. From equation (\ref{eq7.7}) we know that these eigen values will be given by 
\beq
\label{complex3}
\lambda_{1}={-B +   i \sqrt {-D} \over 2} \quad 
\lambda_{2}={-B -   i \sqrt {-D} \over 2}
\eeq

or if we denote the real part of these complex numbers as :
$
\alpha={-B \over 2}$ and the imaginary part as $\beta={ \sqrt {-D} \over 2}
$ we can rewrite (\ref{complex3}) as

\beq
\label{complex4}
\lambda_{1,2}= \alpha \pm i \beta
\eeq

Which type of dynamics do  we expect here.  If
we drop for a while the imaginary part $i\beta$ in (\ref{complex4}), we
get that $\lambda_{1}=\lambda_{2}=
\alpha$, i.e. both eigen values will be the same and hence they will have the same sign.  Which type of
equilibria do we have for two real eigen values which have the same sign? 
We can have either a non-stable node if $\lambda_1
>0,\lambda_2>0$ , or a stable node, if
$\lambda_1 <0,\lambda_2<0$ (see section \ref{sec_2dlin2}). In the case
of a non-stable node all the trajectories diverge from the equilibrium
(fig.\ref{fig6.4}b), while in the case of a stable node all the
trajectories converge to it (fig.\ref{fig6.5}b). It turns out that 
we will get  similar behavior for complex eigen values: if $Re
\lambda_{1,2} >0 $ we will get  divergence of the trajectories, if $Re
\lambda_{1,2} < 0 $ we will get  convergence of the trajectories to the
equilibrium. However, the picture will be slightly different from
fig.\ref{fig6.4}b and fig.\ref{fig6.5}b, as the imaginary part of the eigen values $Im
\lambda_{1,2}$ will add {\it {rotation}} to the trajectories.

The complete derivation of the formula for the general real solution of
system (\ref{2dlin2}) with complex eigen values is given in the
appendix at the end of this chapter for extra reading. Here we will
just demonstrate that $Im \lambda_{1,2}$ adds rotation.

\subsection{Center, spiral \label{sec_center}}
Let illustrate that the imaginary part of the complex eigen value of
the characteristic equation results in a rotation of the trajectories.
For that let us consider a system:
\beq
\label{complex6}
\left( \begin{array}{c} 
{dx \over dt} \\  {dy \over dt}  \end{array} \right) = \left(\begin{array}{lr}  0 & 2\\
                         -2 &  0
\end{array} \right)  \left( \begin{array}{c} x
 \\  y  \end{array} \right) \qquad or \qquad \left\{
\begin{array}{l}
{dx \over dt}=2y\\ {dy \over dt}=-2x
\end{array}
\right.
\eeq
In this case the eigen values are given by:
$ Det \left|\begin{array}{lr}  -\lambda & 2\\
                         -2 &  -\lambda
\end{array} \right|=\lambda^2 +4=0 \;\;  \lambda_{1,2}= \pm 2i$. 
Thus we have  purely  imaginary eigen values. It turns out that we can
draw a phase portrait of  this system using the following trick. Let us multiply the
first equation in (\ref{complex6}) by $x$, the second equation by $y$
and let us add them together. We get

\beq
\label{complex66}
\begin{array}{l}
x{dx \over dt}=x2y\\ + \\ y{dy \over dt}=-y2x\\
gives\\
x{dx \over dt}+y{dy \over dt}=2xy-2xy=0
\end{array}
\eeq
Now note that  
$$
x{dx \over dt}= {1 \over
2}{dx^2 \over dt}
$$
 (just check this  by applying the chain rule for
${dx^2 \over dt}$), and similarly 
$$
y{dy \over dt}={1 \over 2}{dy^2 \over dt}
$$ 
hence eq.(\ref{complex66}) can be rewritten as:

$$
\begin{array}{l}
{1 \over 2}{dx^2 \over dt}+{1 \over 2}{dy^2 \over dt}=0\\
\\
{d(x^2+y^2) \over dt}=0
\end{array}
$$

We know that if the derivative of the function $f$ is zero (${df \over
dt}=f'(x)=0$) then the function $f$ is a constant, thus the above
equation implies:
\beq
\label{complex7}
x^2+y^2=Const
\eeq

Expression (\ref{complex7}) gives a so-called first integral of our
system: combination of variables which are preserved in time. It is
not equivalent to the solution of our system, but using it we can draw
the phase portrait of system (\ref{complex6}). 

\begin{figure}[hhh]
\centerline{
\psfig{type=pdf,ext=.pdf,read=.pdf,figure=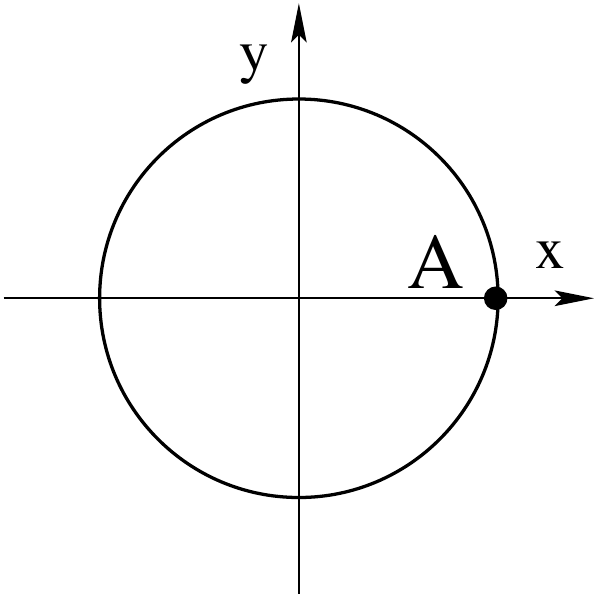,height=5cm}
\psfig{type=pdf,ext=.pdf,read=.pdf,figure=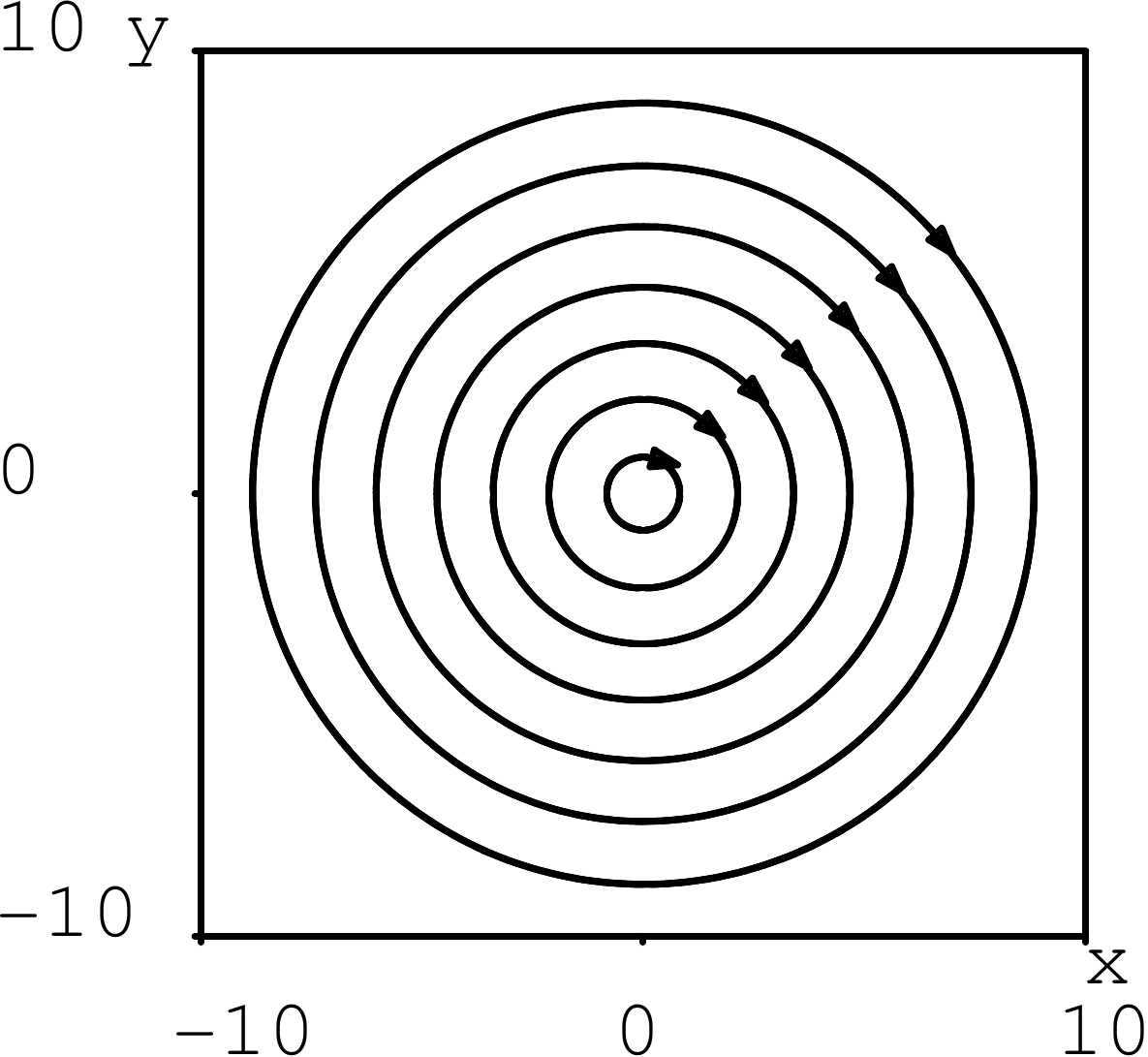,height=5cm}
\psfig{type=pdf,ext=.pdf,read=.pdf,figure=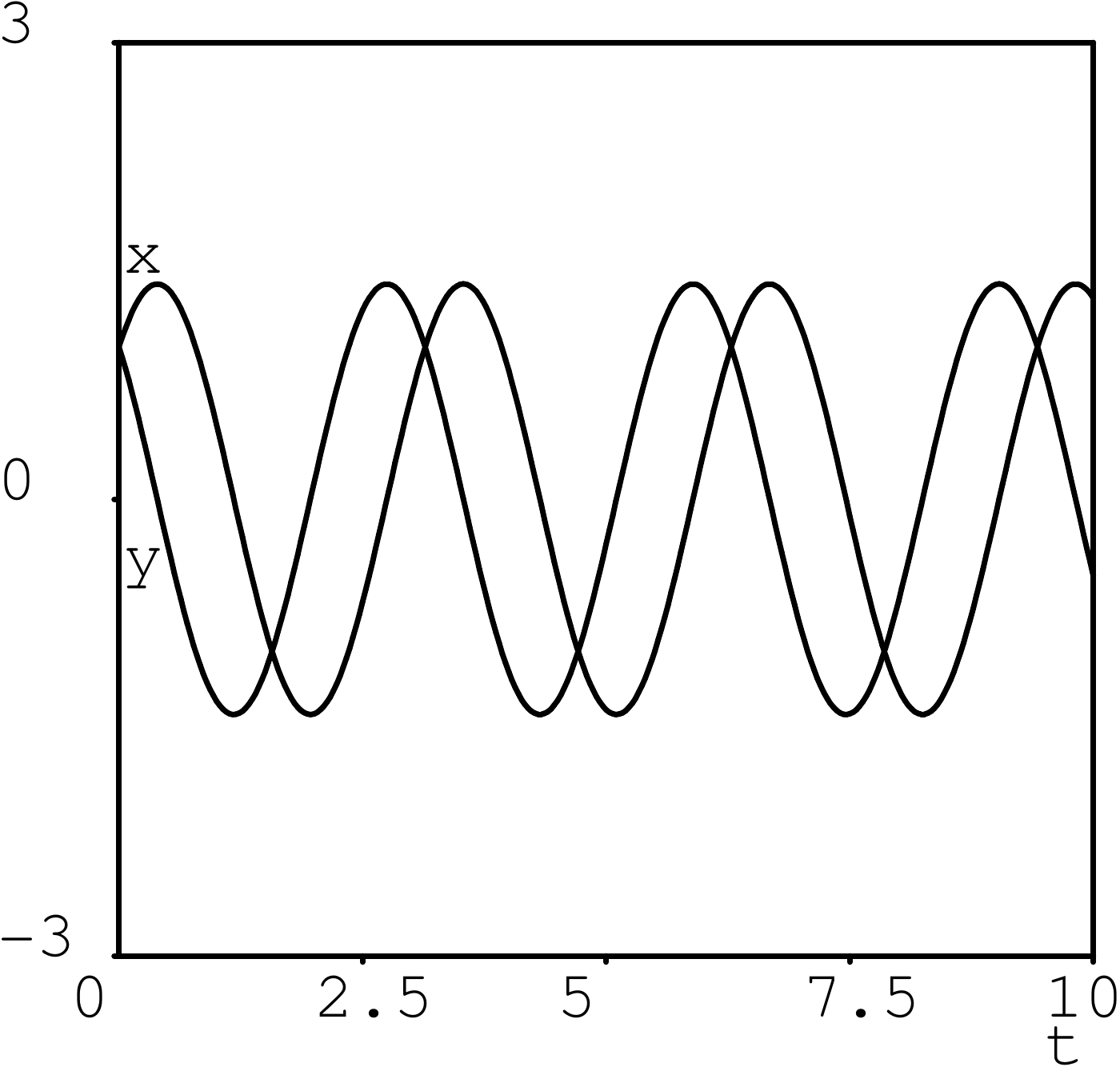,height=5cm}
}
\caption{\label{fig7.1}}
\end{figure}

Because $Const$ in (\ref{complex7}) is an arbitrary positive number,
let us denote it as $Const=A^2$, where $A$ is just another arbitrary
constant. Thus we will get the  equation $x^2+y^2=A^2$, which represents
a graph of a circle with  radius $A$ and with the
center at the origin (fig.\ref{fig7.1}a).

To find the dynamics on this circle let us find ${dy \over dt}$ at a
point $x=A,y=0$ (the bold point in fig.\ref{fig7.1}a). From
eq.(\ref{complex6} ${dy \over dt}=-2x=-2A<0$, i.e. the coordinate $y$
decreases at this point in the course of time. Thus we conclude that
motions of trajectory along the circle will be in the clockwise
direction. If we draw a series of such circles for different $A$ we
will get the phase portrait as shown in fig.\ref{fig7.1}b. The
dynamics of the $x$ and $y$ variables can be found by considering
motion along trajectories of the system: motion here is a rotation of
a point along the circle. During this rotation the variable $x$
changes periodically (between the values $+A$ and $-A$), and thus we
will see periodic oscillations. The same is valid for the variable
$y$. An example of such dynamics is represented in
fig.\ref{fig7.1}c. The equilibrium point $(0,0)$ in such system is
called {\bf {a center}}. The phase portrait in Fig.\ref{fig7.1} is a
set of circles, which is a consequence of the symmetry of system
(\ref{complex6}). In a general case, if $\lambda_{1,2}=\pm i
\beta$, we will also get a center point, similar to that in
fig.\ref{fig7.1}, but instead of circles we can get a series of
embedded  ellipses. The dynamics of the
variables will always be oscillations.

\begin{figure}[hhh]
\centerline{
\psfig{type=pdf,ext=.pdf,read=.pdf,figure=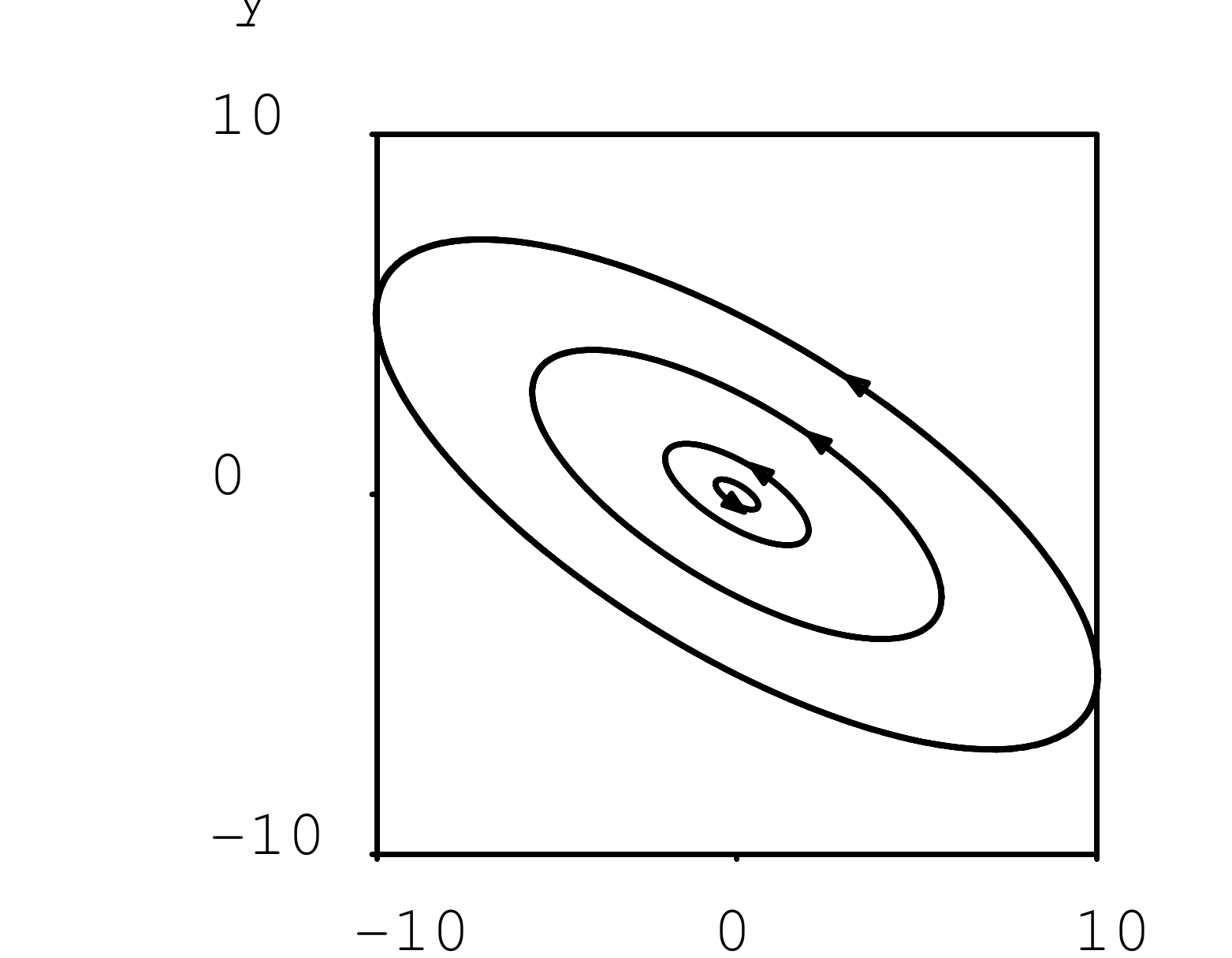,width=6.5cm}
\psfig{type=pdf,ext=.pdf,read=.pdf,figure=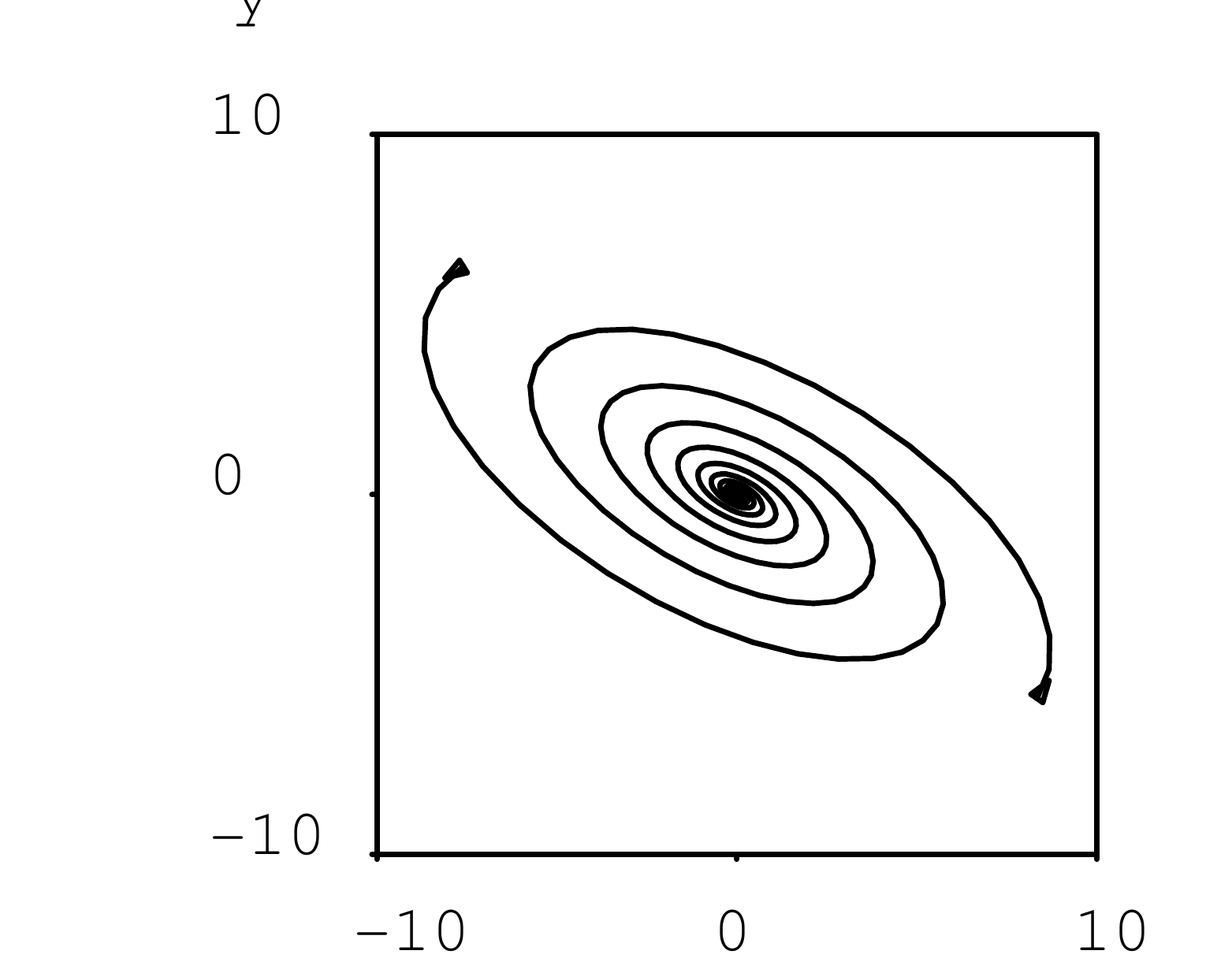,width=6.5cm}
\psfig{type=pdf,ext=.pdf,read=.pdf,figure=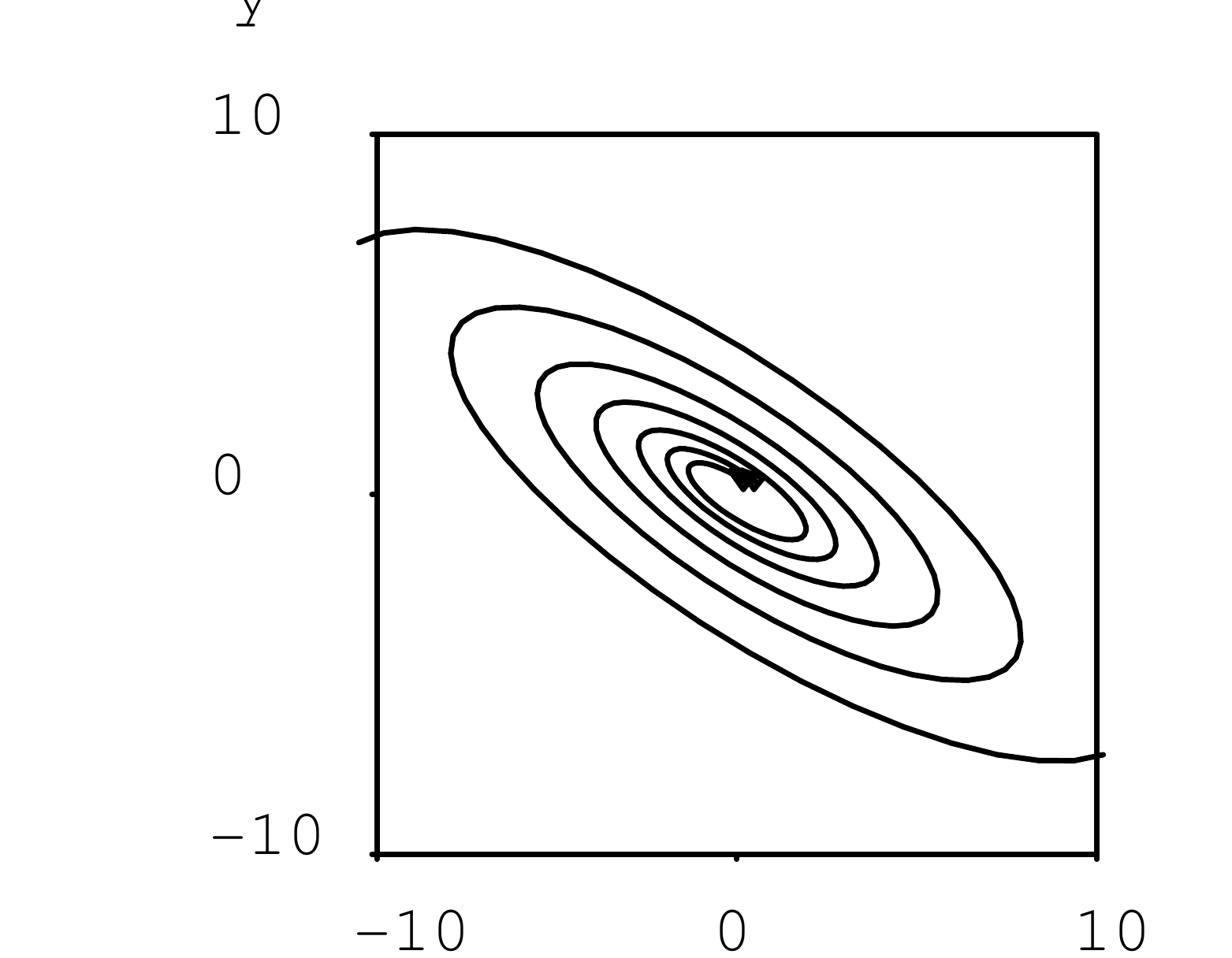,width=6.5cm}
}
\caption{\label{fig7.3}}
\centerline{
\psfig{type=pdf,ext=.pdf,read=.pdf,figure=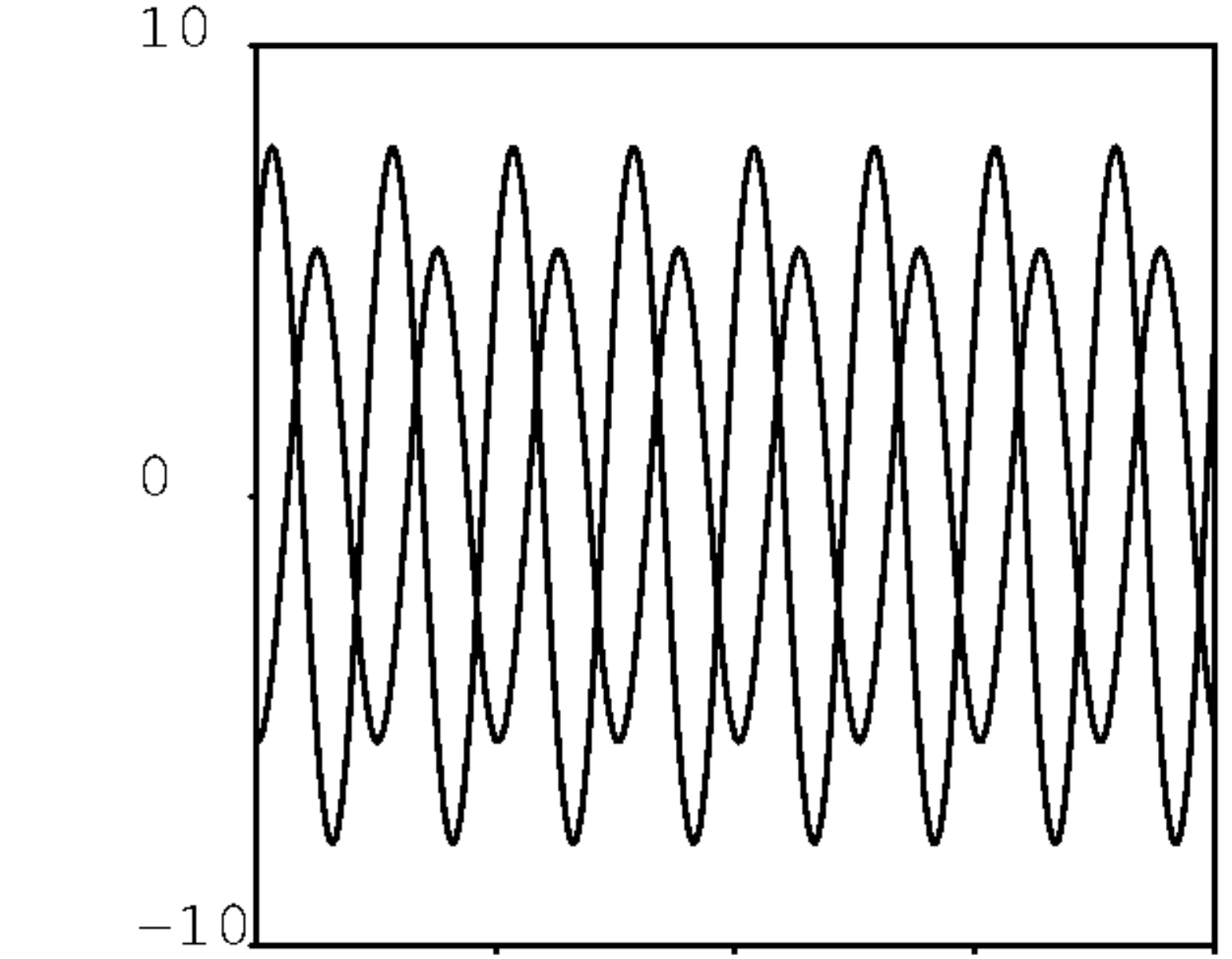,width=6.cm}
\psfig{type=pdf,ext=.pdf,read=.pdf,figure=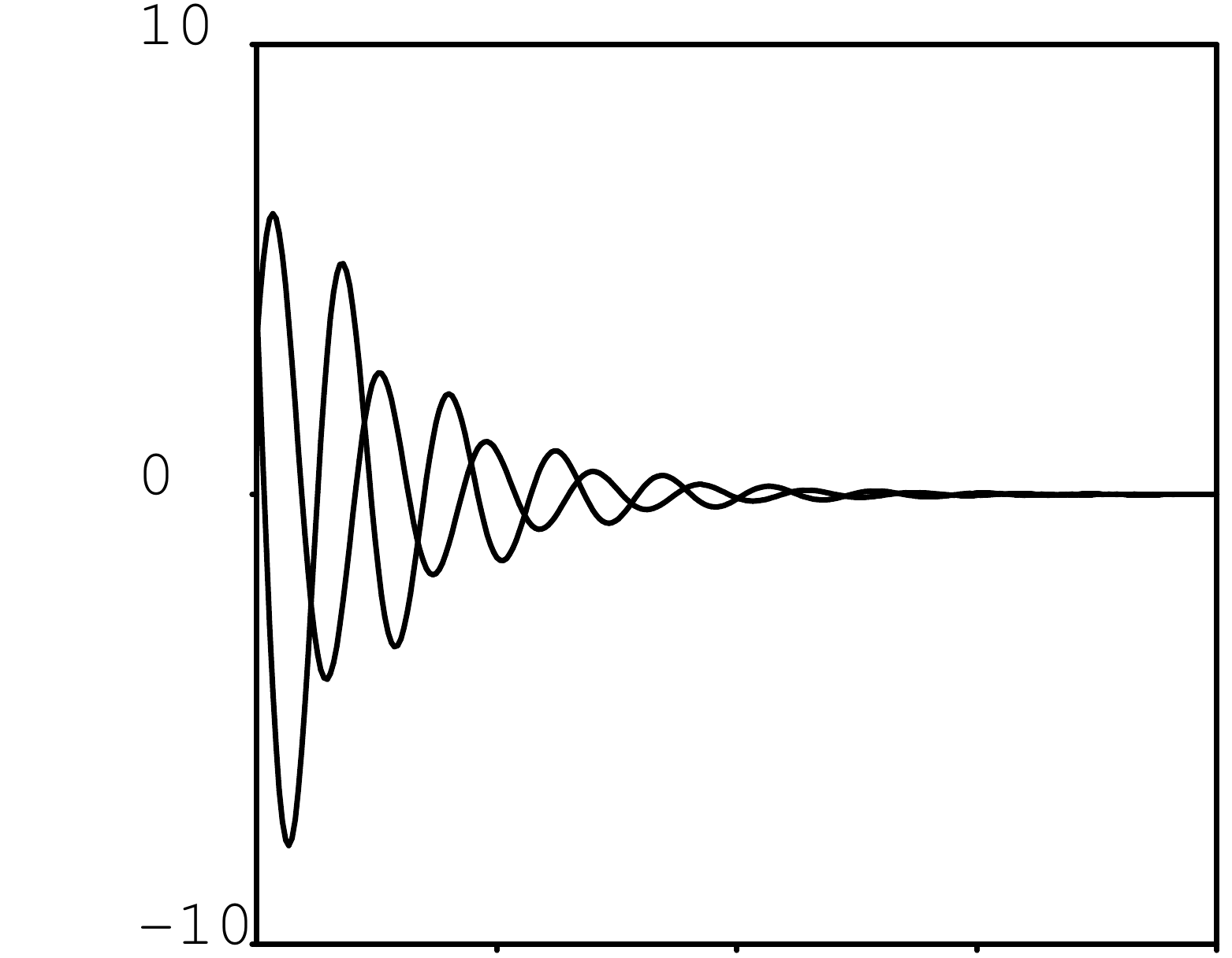,width=6.cm}
\psfig{type=pdf,ext=.pdf,read=.pdf,figure=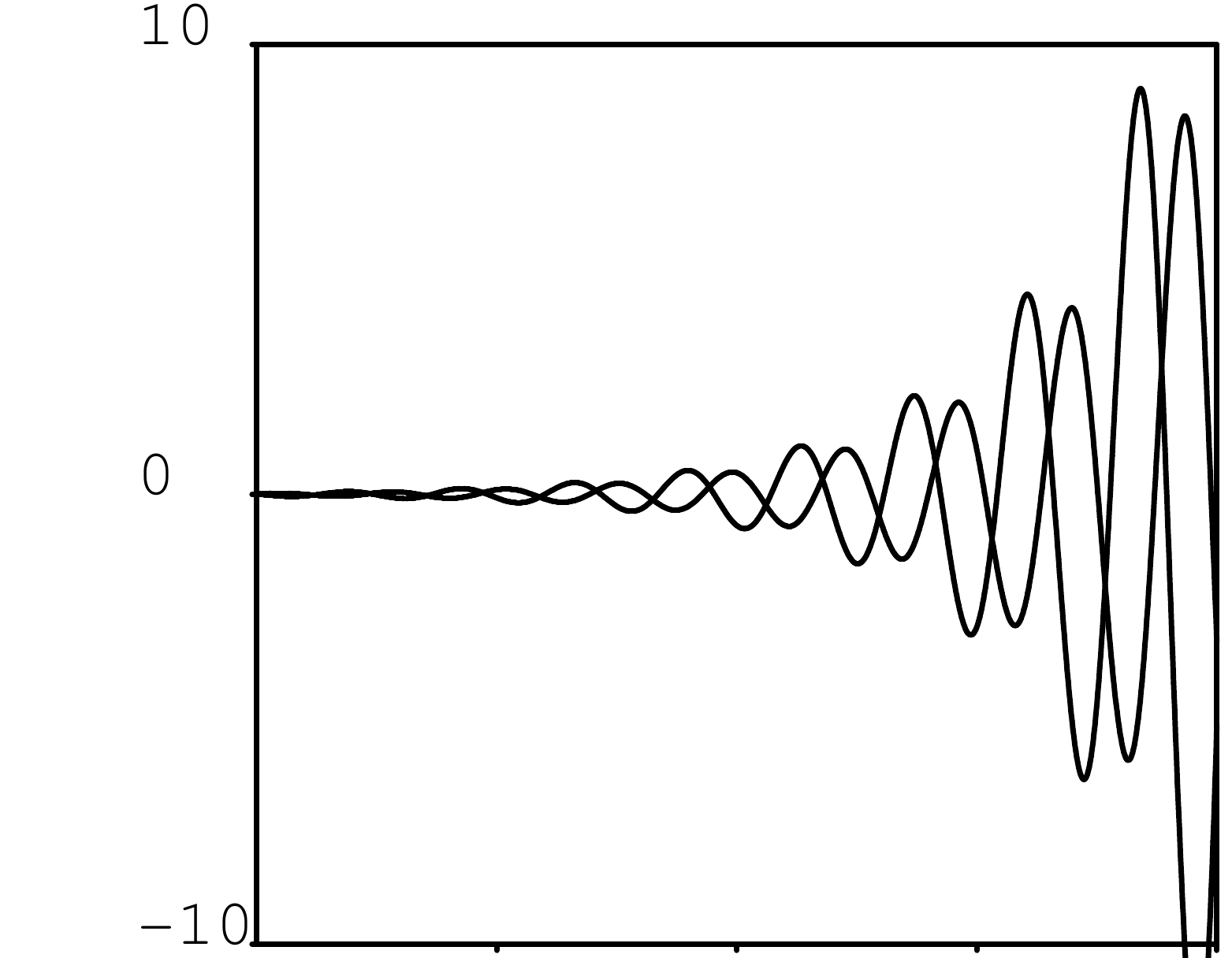,width=6.cm}
}
\caption{The dynamics of variables $x$ and $y$ for the phase portraits shown in fig.\ref{fig7.3}.
The left picture corresponds to the equilibrium point center, the middle picture corresponds to the equilibrium point stable spiral and the right picture corresponds to the equilibrium point non-stable spiral. \label{fig7.4}}

\end{figure}

\bec
If the eigen values of system (\ref{2dlin2}) are
 $\lambda_{1,2}=\pm i \beta$,  we have an equilibrium point called  a center. 
The dynamics of variables $x,y$ are oscillations. The phase portrait is a set of embedded ellipses.
\eec
A computer generated phase portrait of the system ${dx \over dt}=-x-2y; {dy \over dt}=x+y$ is shown in fig.\ref{fig7.3}a. The eigen values in this case are $\lambda_{1,2}= \pm i 0.1$. The time-plot for the $x$ and $y$ variables is shown in fig.\ref{fig7.4}(left).

\bec
The imaginary part of the eigen values results in the  rotation of  trajectories   on a phase portrait.
\eec
Now let us consider the next two cases: $\lambda_{1,2}= \alpha \pm i \beta$.

As we discussed in section \ref{realEV}, $Re \lambda_{1,2}$ determines
convergence or divergence of the trajectories to the equilibrium. In
the case $\lambda_{1,2}= \alpha \pm i\beta; \alpha<0 $ we expect that
the real part of the eigenvalue will give a behavior similar to the
case of both negative real eigenvalues, hence we expect the
convergence of trajectories to the equilibrium, as for a stable
node. In addition to this, as we saw in the previous section, the
imaginary part of the eigen values causes the rotation of the
trajectory.  If we add these two processes together we will get
convergence to the equilibrium with rotation, hence trajectories will
have the form of spirals.  A computer generated phase portrait for
this case is shown in fig.\ref{fig7.3}b.  The system is ${dx \over
dt}=-x-2y; {dy \over dt}=x+0.7 y$. The eigen values in this case are
$\lambda_{1,2}= -0.15 \pm i 0.13$. The time-plot for the $x$ and $y$
variables is shown in fig.\ref{fig7.4}(middle). The dynamics of the
system are oscillations with gradually decreasing amplitude.

\bec
 If the eigen values of system (\ref{2dlin2}) are
 $\lambda_{1,2}= \alpha \pm i\beta; \alpha<0 $,
  we have an equilibrium point called  {\bf {a stable spiral}}, fig.\ref{fig7.3}b. 
\eec

The last case occurs if $\lambda_{1,2}= \alpha \pm i\beta; \alpha>0 $.
This case is similar to the previous one. The only difference is that
because $Re \lambda_{1,2}= \alpha>0$, the real part gives motion
equivalent to a non-stable node, or divergence of trajectories from
the equilibrium. So, together with rotation from the imaginary part of
$\lambda$ we get the following phase portrait (fig.\ref{fig7.3}c).
This phase portrait is generated by computer for the system ${dx \over
dt}=-x-2y; {dy \over dt}=x+1.2y$. The eigen values in this case are
$\lambda_{1,2}= 0.1 \pm i 0.89$.  The time-plot for the $x$ and $y$
variables is shown in fig.\ref{fig7.4}(right). The dynamics of the
system are oscillations with gradually increasing amplitude.

\bec
If the eigen values of system (\ref{2dlin2}) are
 $\lambda_{1,2}= \alpha \pm i\beta; \alpha>0 $,
  we have an equilibrium point called  {\bf {a non-stable spiral}}, fig.\ref{fig7.3}c. 
\eec
We have found all possible types of equilibria which can occur in 2D systems: saddle, non-stable node, stable node, center, non-stable spiral and stable spiral. The next  question which we will discuss here is the stability of these equilibria.

\section{Stability of equilibrium}
\begin{figure}[h]
\centerline{
\psfig{type=pdf,ext=.pdf,read=.pdf,figure=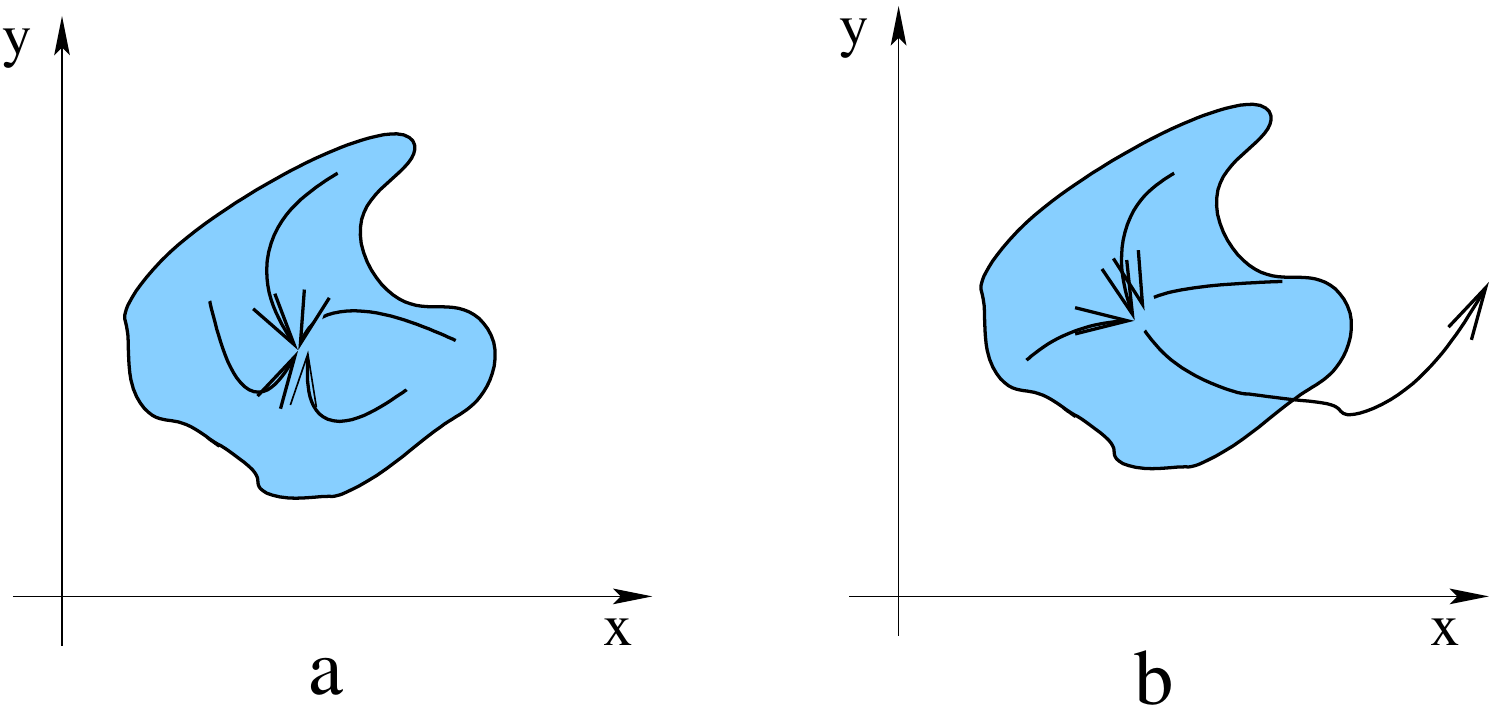,width=12.5cm}
}
\caption{\label{fig8.1}}
\end{figure}

We will call an equilibrium point {\bf {stable}}, if there is a neighborhood of this equilibrium, such that all trajectories which start in  this neighborhood will converge to the equilibrium (Fig.\ref{fig8.1}a). 
We will call the equilibrium point {\bf {non-stable}}, if there   is at least one diverging trajectory from the close neighborhood of this equilibrium (fig.\ref{fig8.1}b).

If we analyze the stability of  the 6 types of equilibria studied in the previous section we find the following:

{\bf Stable equilibria}
\ben
\item Stable node  $\lambda_1<0;\lambda_2<0$ real

\item Stable spiral $\lambda_{1,2}= \alpha \pm i\beta; \alpha<0 $
\een
 {\bf Non-stable equilibria}
\ben
\item Non-stable node  $\lambda_1>0;\lambda_2>0$  real

\item Non-stable spiral $\lambda_{1,2}= \alpha \pm i\beta; \alpha>0 $

\item Saddle point   $\lambda_1<0;\lambda_2>0$; or  $\lambda_1>0;\lambda_2<0$  real
\een

In case of spirals and nodes the stability and non-stability is
obvious.  In case of a saddle point we have a converging trajectory,
however the existence of the diverging trajectories (fig.\ref{fig6.2}b)
implies that this equilibrium point is non-stable. The last case,
$\lambda_{1,2}= \pm i\beta$ (center point), is non conclusive.
Trajectories do not converge and do not diverge from the
equilibrium. Usually this case is treated as neutrally stable. All
these cases can be formulated in the following theorem.
\begin{theorem}
If all eigenvalues of the linear system (\ref{2dlin2}) have negative real parts, then the equilibrium point $x=0,y=0$ is stable.
\end{theorem}
It is easy to see, that this theorem includes all listed stable
equilibria. It obviously works for a stable spiral, but it also works
for a stable node, because any real number can be considered as a
complex number with imaginary part equal to zero. For example: $-3$
can be represented as $z=-3=-3+i0$, and $Re z=-3; Im z=0$.

\section{Exercises}

\subsection*{Exercises for section \ref{sec_2dlin2}}

\ben

\item
Find the general solution of the following systems of ordinary differential equations:
\ben
\item 
$\left( \begin{array}{c} 
{dx \over dt} \\  {dy \over dt}  \end{array} \right) = 
\left(\begin{array}{lr}  -2 & 1\\
                         1 &  -2
\end{array} \right)  
\left( \begin{array}{c} x
 \\  y  \end{array} \right)$

\item $\left( \begin{array}{c} 
{dx \over dt} \\  {dy \over dt}  \end{array} \right) = 
\left(\begin{array}{lr}  3 & -1\\
                         -2 &  4
\end{array} \right)  
\left( \begin{array}{c} x
 \\  y  \end{array} \right)$
\een
\item Find the solution for the following initial value problem:

$$\left( \begin{array}{c} 
 {dx \over dt} \\  {dy \over dt}  \end{array} \right) = 
\left(\begin{array}{lr}  1 & -2\\
                         5 &  8
\end{array} \right)  
\left( \begin{array}{c} x
 \\  y  \end{array} \right) \quad 
\left( \begin{array}{c} x(0)
 \\  y(0)  \end{array} \right)=\left( \begin{array}{c} 3
 \\  -3  \end{array} \right)$$
(Hint:  See  an example of solution  in section \ref{solutionIVP}).

\item Two different concentrations of a solution are separated by a membrane through which the solute
can diffuse. The rate at which the solute diffuses is proportional to the difference in concentrations between two solutions.
The differential equations governing the process are:
$$
\left\{
\begin{array}{l}
dC_1/dt=-{k \over V_1}(C_1-C_2)\\
  dC_2/dt={k \over V_2}(C_1-C_2)
\end{array}
\right.
$$
where $C_1$ and $C_2$ are the two concentrations, $V_1$ and $V_2$ are the volumes of the respective compartments, and $k$ is a constant
of proportionality. If $V_1=20liters$, $V_2=5 liters$, and $k=0.2 \; liters/min$ and if initially $C_1=3 \; moles/liter$ and $C_2=0$, find
$C_1$ and $C_2$ as functions of time.

\subsection*{Exercises for sections \ref{realEV} and \ref{complexEV}}

\item  Find eigen values and eigen vectors (for real eigen values only) of the following systems. Determine equilibrium type and sketch  phase portraits. For real eigen values show non-stable, stable manifolds and several trajectories between the manifolds. 

\ben
\item 
$
\left\{
\begin{array}{l}
{dx \over dt}=   x+ 4 y \\ {dy \over dt} =   2 x+ 3 y
\end{array}
\right. 
$

\item
$
\left\{
\begin{array}{l}
{dx \over dt}=   5x-y \\ {dy \over dt} =   3x+ y
\end{array}
\right. 
$

\item
$
\left\{
\begin{array}{l}
{dx \over dt}=   3x-5 y \\ {dy \over dt} =   x-y
\end{array}
\right. 
$

\item 
$
\left\{
\begin{array}{l}
{dx \over dt}=   -2x+ y \\ {dy \over dt} =   x-2y
\end{array}
\right. 
$

\item 
$
\left\{
\begin{array}{l}
{dx \over dt}=   -2y \\ {dy \over dt} =   x-2y
\end{array}
\right. 
$

\item 
$
\left\{
\begin{array}{l}
{dx \over dt}=   -x-y \\ {dy \over dt} =   2x+y
\end{array}
\right. 
$

\item 
$
\left\{
\begin{array}{l}
{dx \over dt}=   -2x-y \\ {dy \over dt} =   3x+2y
\end{array}
\right. 
$

\een

\item Study the following linear system with a parameter $a$:

$$
\left\{
\begin{array}{l}
{dx \over dt}=   -2x-ay \\ {dy \over dt} =   3x-y
\end{array}
\right. 
$$

 Find  the types of equilibrium  which are possible for different values of $ -\infty < a < \infty $. Give the parameter region for each equilibrium and draw  qualitative phase portraits. For which parameter values is the equilibrium stable?

\item For which values of parameters $a$ and $b$ does the following system has periodic oscillations (i.e. a center equilibrium point):
$$
\left\{
\begin{array}{l}
{dx \over dt}=   -ax+y \\ {dy \over dt} =   (2a-3)x-by
\end{array}
\right. 
$$

\begin{figure}[H]
\centerline{
\psfig{type=pdf,ext=.pdf,read=.pdf,figure=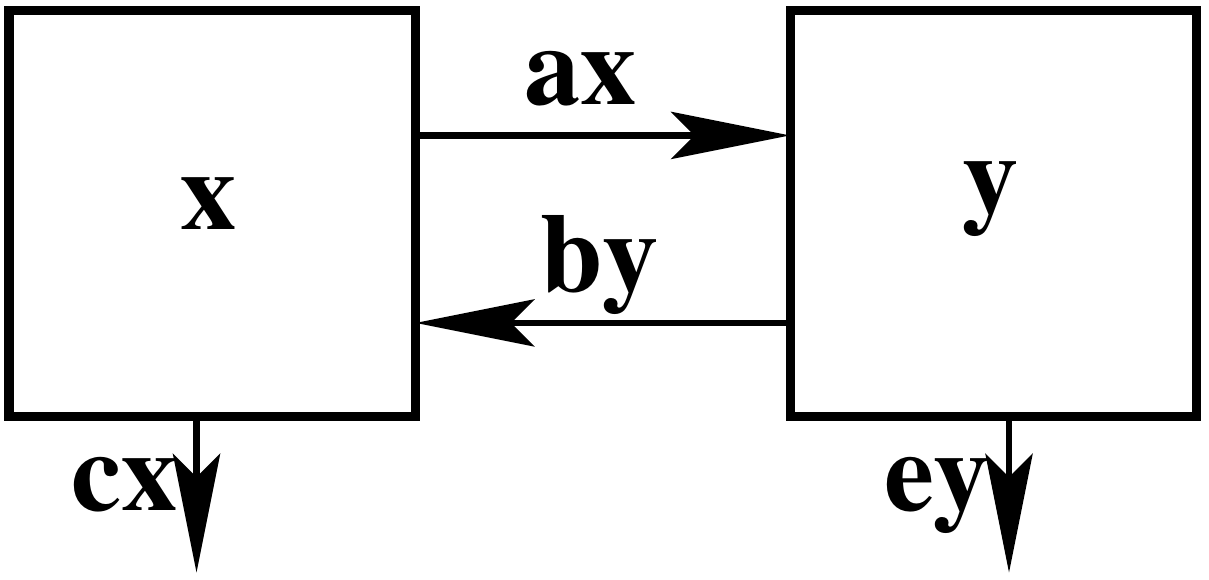,width=5.cm}
}
\caption{\label{figCOMP1}}
\end{figure}

\item Compartmental models play an important role in different parts of population biology, pharmacology and biochemistry. They describe the interaction between  several processes, which  may be interactions of populations,  chemical reactions,   etc. A two compartment model is schematically shown in fig.\ref{figCOMP1}. It represents two interacting species $x$ and $y$. The concentration of the species $x$ can be changed either due to a transition to species $y$ with the rate given by $ax$, or $x$ can die with the rate $cx$. Similar transitions exist also for the species $y$. The rates of these processes are specified in the figure.

\ben
\item Derive a system of  two differential equations for  species $x$ and $y$.

\item If $a=0.5,b=2,c=4.5,e=3$ find equilibrium type. Is it stable?

\item Determine the stability of the equilibrium in a general case when  $a>0,b>0,c>0,e>0$
\een

\een

\newpage{\small
\section{Additional concepts (appendix) \label{append}}

\subsection{General solution for complex eigen values}

\subsubsection {Notes about general solution of a system with complex eigenvalues \label{notes_sec}}
It turns out that if $\lambda_{1,2}= \alpha \pm i \beta$  formula
(\ref{2dlin_sol}) is still valid, so the solution can be represented
in the form:
\beq
\left( \begin{array}{c} 
x \\  y \end{array} \right) = 
C_1 \left( \begin{array}{c} v_{1x}
 \\  v_{1y}  \end{array} \right)e^{\lambda_1 t } 
+ C_2 \left( \begin{array}{c} v_{2x}
 \\  v_{2y}  \end{array} \right)e^{\lambda_2 t }
\eeq
but as $\lambda_1=\alpha+i\beta; \lambda_2=\alpha-i\beta$ we get
\beq
\label{complex5}
\left( \begin{array}{c} 
x \\  y \end{array} \right) = 
C_1 \left( \begin{array}{c} v_{1x}
 \\  v_{1y}  \end{array} \right)e^{(\alpha+i\beta) t } 
+ C_2 \left( \begin{array}{c} v_{2x}
 \\  v_{2y}  \end{array} \right)e^{(\alpha-i\beta) t }
\eeq
This expression gives the correct solutions of
(\ref{2dlin201}). However, the form of the solution is not good. First,
because $\lambda_{1,2}$ are complex,  the eigenvectors $\left(
\begin{array}{c} v_{1x} \\ v_{1y} \end{array} \right)$ and $\left(
\begin{array}{c} v_{2x} \\ v_{2y} \end{array} \right)$ will also be
complex. Next, we should consider $C_1,C_2$ as general complex
constants. Therefore (\ref{complex5}) is a quite  complicated
expression which gives $x$ and $y$ as complex valued functions of time
$t$. Mathematically it is correct and if we substitute
(\ref{complex5}) into the original equation (\ref{2dlin201}) we get an
equality. However, we need to draw a phase portrait of
(\ref{2dlin201}) on the $Oxy$-plane, i.e. we only need the {\it
{real}} solutions of our system. Such real solutions are  present
in the general expression (\ref{complex5}), i.e. they are a subset of
all the possible solutions. However, extracting them from (\ref{complex5})
is not a simple task.  We will highlight the main idea behind this
derivation below. To derive the general expressions we will need the
following Euler formula.

\subsubsection*{Euler formula}
The Euler formula gives a representation of  $e^{i\beta t } $ 
in terms of trigonometric functions. It is quite unexpected:

\beq
\label{euler1}
e^{i\beta t } =\cos {\beta t} + i \sin {\beta t}
\eeq
or in another representation:
\beq
\label{euler}
e^{i \phi } =\cos {\phi } + i \sin {\phi }
\eeq
When you see this formula for the first time it looks quite crazy. We
know that $\sin \phi $ and $\cos \phi$ come from the simple geometry
of triangles, $i=\sqrt{-1}$ and $e$ is a special exponential
function. Why are these functions connected together in such a simple
way (\ref{euler})?

To prove this formula, one should use Taylor series. However, here I will
present another simpler derivation  of (\ref{euler}) on the
basis of differential equations.

At the beginning of this chapter, in order to find a solution of
(\ref{2dlin201}), we first considered a one dimensional differential
equation ${dx \over dt}= ax$, and we found its solution
$Ce^{at}$. Consider the following initial value problem for this
equation:
\beq
\label{tmp30}
{dx \over dt}=ax \qquad x(0)=1
\eeq

This initial value problem has the unique solution $x(t)=e^{at}$. So
we say, that the solution of (\ref{tmp30}) is $e^{at}$. But we can
also say it vice versa: we can define the exponential function
$e^{at}$ as the function which satisfies the initial value problem
(\ref{tmp30}). For example, if we give to a person just this equation
and a computer, he will be able to solve it and to draw the graph of
$e^{at}$, even without knowledge about exponential functions.  The
advantage of such a definition is that it can be easily extended to
complex numbers. So, let us define $e^{it}$ as a function which
satisfies the initial value problem (\ref{tmp30}) with $a=i$

\beq
\label{tmp31}
{dx \over dt}=ix \qquad x(0)=1
\eeq

In other words: $e^{it}$ must be the function $x(t)$, such that
$x(0)=1$, and the derivative of this function ${dx \over dt}(t)$ must
be equal to this function times $i$, i.e.  ${dx \over
dt}(t)=ix(t)$. Let us find an expression which satisfies these
conditions.  It turns out that it will be $x(t)=\cos {t } + i \sin {t
}$. Let us check it. The first condition is satisfied:

\beqar
x(0)=\cos {0 } + i \sin {0 }=1+0i=1,
\eeqar
To check the second condition we write:
\beqar
{dx \over dt}(t)=(\cos {t } + i \sin {t })'=\cos' {t } + i \sin' {t } =\\
-\sin {t } + i \cos {t }
\eeqar
if we replace $-1$ by $i^2$ we get:
\beqar
{dx \over dt}(t)= -\sin {t } + i \cos {t }=i^2\sin {t }+i \cos {t }=\\
i(\cos {t} + i \sin {t })=ix(t), 
\eeqar
i.e. the second condition is also satisfied.  So the function $x(t)=\cos {t } + i
\sin {t }$ gives the solution of  (\ref{tmp31}), hence it is the same as $e^{it}$ or
$e^{it}=\cos {t } + i \sin {t }$ and we get the Euler formula
(\ref{euler}). The formula (\ref{euler1}) is just the formula (\ref{euler}) in
which instead of $\phi$ the letters $\beta t$ are used. To find
$e^{(\alpha+i\beta) t }$ we write:
\beq
e^{(\alpha+i\beta) t }=e^{\alpha t} e^{i\beta t }=e^{\alpha t}(\cos {\beta t} + i \sin {\beta t})
\eeq
\subsubsection*{General solution}

Now let us find a solution of a  system with imaginary eigenvalues. As we know the general solution is given by the formula (\ref{complex5}) and because of the Euler formula we can rewrite it in the following way:
\beq
\label{tmp32}
\left( \begin{array}{c} 
x \\  y \end{array} \right) = 
C_1 \left( \begin{array}{c} v_{1x}
 \\  v_{1y}  \end{array} \right)e^{\alpha t}(\cos {\beta t} + i \sin {\beta t})
+ C_2 \left( \begin{array}{c} v_{2x}
 \\  v_{2y}  \end{array} \right)e^{\alpha t}(\cos {\beta t} - i \sin {\beta t}),
\eeq
where $C_1,C_2$ are arbitrary complex constants and $\left(
\begin{array}{c} v_{1x} \\ v_{1y} \end{array} \right)$ and $\left(
\begin{array}{c} v_{2x} \\ v_{2y} \end{array} \right)$ are complex
eigen vectors. Now we should find a real part of this complicated
expression and get a general real solution of our system in this case.
A general solution of a system of two differential equations should
depend on two arbitrary constants as the initial value of each of the
variables can be arbitrary. We will use this fact to find the general
solution.  Our idea is instead of extracting all real solutions from
(\ref{tmp32}) we will find just two real solutions. By multiplying
them by two arbitrary constants we will get a general solution.

We will be able to find these two solutions from the first term of
(\ref{tmp32}):
\beq
\label{tmp33}
Y_1=\left( \begin{array}{c} v_{1x}
 \\  v_{1y}  \end{array} \right)e^{\alpha t}(\cos {\beta t} + i \sin {\beta t})=
{\bf v_1}e^{\alpha t}(\cos {\beta t} + i \sin {\beta t})
\eeq
Let us extract real and an imaginary parts of this term. If we use the
formula (\ref{eigen_ex}) for the eigen value $\lambda=\alpha + i \beta$ we find
the eigen vector ${\bf v_1}$: It has the real and imaginary parts:
\beq
{\bf v_1}=\left( \begin{array}{c} v_{1x}
 \\  v_{1y}  \end{array} \right)=
\left( \begin{array}{c} -b
 \\  a-\lambda_1  \end{array} \right)=\left( \begin{array}{c} -b
 \\  a-\alpha -i\beta  \end{array} \right)=\left( \begin{array}{c} -b
 \\  a-\alpha  \end{array} \right)+i \left( \begin{array}{c} 0
 \\   -\beta  \end{array} \right)={\bf v_{r}}+i{\bf v_{i}}
\eeq
 The vector ${\bf v_1}$ has the real part ${\bf v_{r}}$ and imaginary
 part ${\bf v_{i}}$.  So, the term (\ref{tmp33}) can be written as:
\beqar
Y_1=({\bf v_{r}}+i{\bf v_{i}})e^{\alpha t}(\cos {\beta t} + i \sin
{\beta t})=\\ e^{\alpha t}({\bf v_{r}}\cos {\beta t}+i{\bf v_{r}} \sin
{\beta t}+i{\bf v_{i}}\cos {\beta t} -{\bf v_{i}}\sin {\beta t})\\
=e^{\alpha t}({\bf v_{r}}\cos {\beta t}-{\bf v_{i}}\sin {\beta t})+
i e^{\alpha t}({\bf v_{r}} \sin {\beta t}+{\bf v_{i}}\cos {\beta t})
\eeqar
If we denote:
\beq
\label{complex_y}
\begin{array}{l}
{\bf y_1}=e^{\alpha t}({\bf v_{r}}\cos {\beta t}-{\bf v_{i}}\sin
{\beta t})\\ {\bf y_2}=e^{\alpha t}({\bf v_{r}} \sin {\beta t}+{\bf
v_{i}}\cos {\beta t})
\end{array}
\eeq
the term (\ref{tmp33}) can be rewritten as $$Y_1={\bf y_1}+i{\bf y_2}.
$$ Let us prove that both ${\bf y_1}$ and ${\bf y_2}$ are the real
solutions of (\ref{2dlin201}). For that we will use the fact that
(\ref{tmp32}) gives a solution of (\ref{2dlin201}) and hence $Y_1$ is a
complex solution of (\ref{2dlin201}) as it is a part of (\ref{tmp32}).

System (\ref{2dlin201}) in a matrix form can be written as:
\beq
\label{tmp34}
{\bf {dX \over dt}} =A {\bf X}
\eeq
As $Y_1$ is a solution, it satisfies (\ref{tmp34}):

\beqar
 {dY_1 \over dt} =A Y_1\\ {\bf {dy_1 \over dt }} + i{\bf {d y_2 \over
 dt}} = A ({\bf y_1}+i{\bf y_2})\\ {\bf {d y_1 \over dt}} + i{\bf
 {dy_2 \over dt}}=A {\bf y_1} +iA{\bf y_2}
\eeqar
Equating the real and imaginary parts yields
\beqar
{\bf {d y_1 \over dt }} =A {\bf y_1} \qquad {\bf {d y_2 \over
dt}}=A{\bf y_2}
\eeqar
Hence ${\bf y_1}$ and ${\bf y_2}$ are real solutions of
(\ref{2dlin201}).  Finally, because ${\bf y_1}$ and ${\bf y_2}$ are real
solutions of (\ref{2dlin201}) the general solution is given by the
formula:
\beq
\label{complex_sol}
\left( \begin{array}{c} 
x \\ y \end{array} \right) = C_1{\bf y_1}+C_2{\bf y_2}
\eeq
where $C_1$ and $C_2$ are arbitrary constants and ${\bf y_1}$ and
${\bf y_2}$ are given by (\ref{complex_y}).

{\bf {Example}} Find the general solution of the following system.
\beq
\label{complex_ex}
\left( \begin{array}{c} 
{dx \over dt} \\  {dy \over dt}  \end{array} \right) = \left(\begin{array}{lr}  0 & 2\\
                         -2 &  0
\end{array} \right)  \left( \begin{array}{c} x
 \\  y  \end{array} \right) \qquad or \qquad \left\{
\begin{array}{l}
{dx \over dt}=2y\\ {dy \over dt}=-2x
\end{array}
\right.
\eeq

{\bf {Solution}}. The eigen values are given by: $ Det \left|\begin{array}{lr}
  -\lambda & 2\\
      -2 &  -\lambda \end{array} \right|={\lambda}^2+4=0,\; \lambda_{1,2}=\pm \sqrt {-4}= 
\pm2i$. So the eigen vector ${\bf v_1}$ corresponding to the eigen value $\lambda=2i$ is:
\beq
{\bf v_1}=\left( \begin{array}{c} v_{1x}
 \\  v_{1y}  \end{array} \right)=
\left( \begin{array}{c} -2
 \\  0-2i  \end{array} \right)=\left( \begin{array}{c} -2
 \\  0  \end{array} \right)+i \left( \begin{array}{c} 0
 \\   -2  \end{array} \right)={\bf v_{r}}+i{\bf v_{i}}
\eeq
So,
\beqar
{\bf y_1} = e^0( \left( \begin{array}{c} -2
 \\  0  \end{array} \right) \cos{2 t}- \left( \begin{array}{c} 0
 \\  -2  \end{array} \right) \sin{2 t})=\left( \begin{array}{c} -2\cos{2 t}
 \\  2  \sin{2 t}  \end{array} \right) \\
{\bf y_2} = e^0( \left( \begin{array}{c} -2
 \\  0  \end{array} \right) \sin{2 t}+ \left( \begin{array}{c} 0
 \\  -2  \end{array} \right) \cos{2 t})=\left( \begin{array}{c} -2\sin{2 t}
 \\  -2  \cos{2 t}  \end{array} \right) 
\eeqar

Therefore the solution from the formulas (\ref{complex_sol}) is:
\beq
\left( \begin{array}{c} x
 \\  y  \end{array} \right) =
C_1\left( \begin{array}{c} -2\cos{2 t}
 \\  2  \sin{2 t}  \end{array} \right) +
C_2\left( \begin{array}{c} -2\sin{2 t}
 \\  -2  \cos{2 t}  \end{array} \right) 
\eeq
Or if we denote $-2C_1=A$ and $-2C_2=B$ we get:

\beq
\label{tmp36}
\left\{
\begin{array}{l}
x=A\cos 2t +B \sin 2t \\
y=-A \sin 2t+B\cos 2t 
\end{array}
\right.  
\eeq
Let us check that (\ref{tmp36}) does gives a solution of (\ref{complex_ex}).
Substitution of (\ref{tmp36}) into equation (\ref{complex_ex}) yields:
\beqar
\left\{
\begin{array}{l}
(A \cos 2t +B \sin 2t)'=2*(-A \sin 2t +B \cos 2t)\\
(-A \sin 2t +B \cos 2t)'=-2*(A \cos 2t +B \sin 2t)
\end{array}
\right. \quad or \\
\left\{
\begin{array}{l}
-2A \sin 2t +2B \cos 2t=2*(-A \sin 2t +B \cos 2t)\\
-2A \cos 2t -2B \sin 2t=-2*(A \cos 2t +B \sin 2t)
\end{array}
\right.
\eeqar
So,  (\ref{tmp36}) is a solution of  (\ref{complex_ex}).

Finally note that formula (\ref{complex7}) which we got for the same equation (\ref{complex_ex}) is of course valid for functions given in  (\ref{tmp36}).
For that you need to find $x^2+y^2$ with $x$ and $y$ given by (\ref{tmp36}).
A direct computation will give us
$$
x^2+y^2=A^2+B^2=Const
$$
i.e. the same result as in (\ref{complex5}). (Note that in order to get the final result you need to apply several times the well-known formula from trigonometry $sin^2(\alpha)+cos^2(\alpha)=1$).
} 

\chapter{System of two non-linear differential equations  \label{chap2dnonlin}}
\section{Introduction and first definitions \label{sec2dnonlin}}

\subsection{Phase portrait}

After analyzing linear systems let us consider a general non-linear
system which can be written in the following general form:
\beq
\label{2dgeneral}
\left\{
\begin{array}{l}
{dx \over dt}=f(x,y) \\ {dy \over dt}=g(x,y)
\end{array}
\right.
\eeq

\begin{figure}[hbt]
\centerline{
\psfig{type=pdf,ext=.pdf,read=.pdf,figure=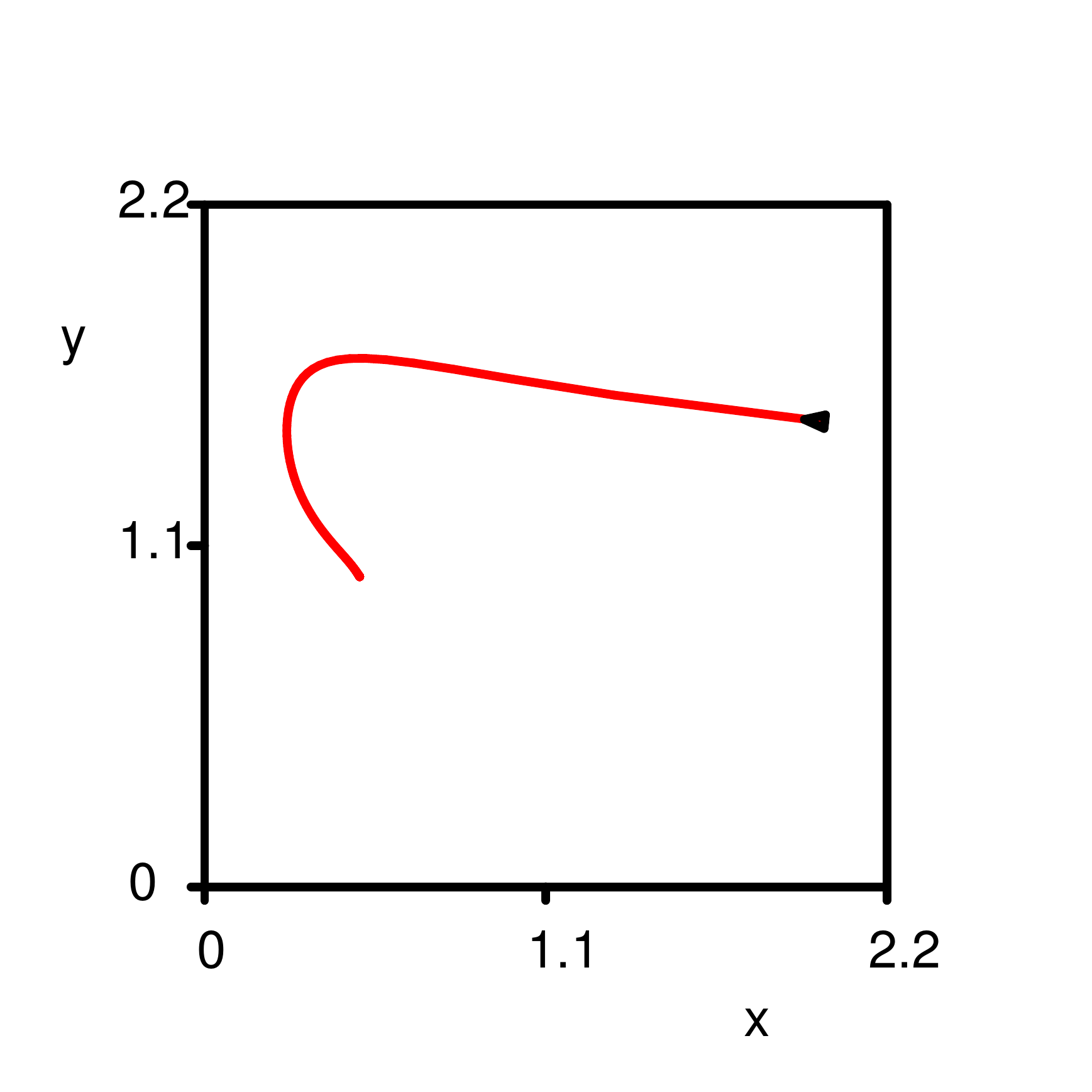,width=6.cm}
\psfig{type=pdf,ext=.pdf,read=.pdf,figure=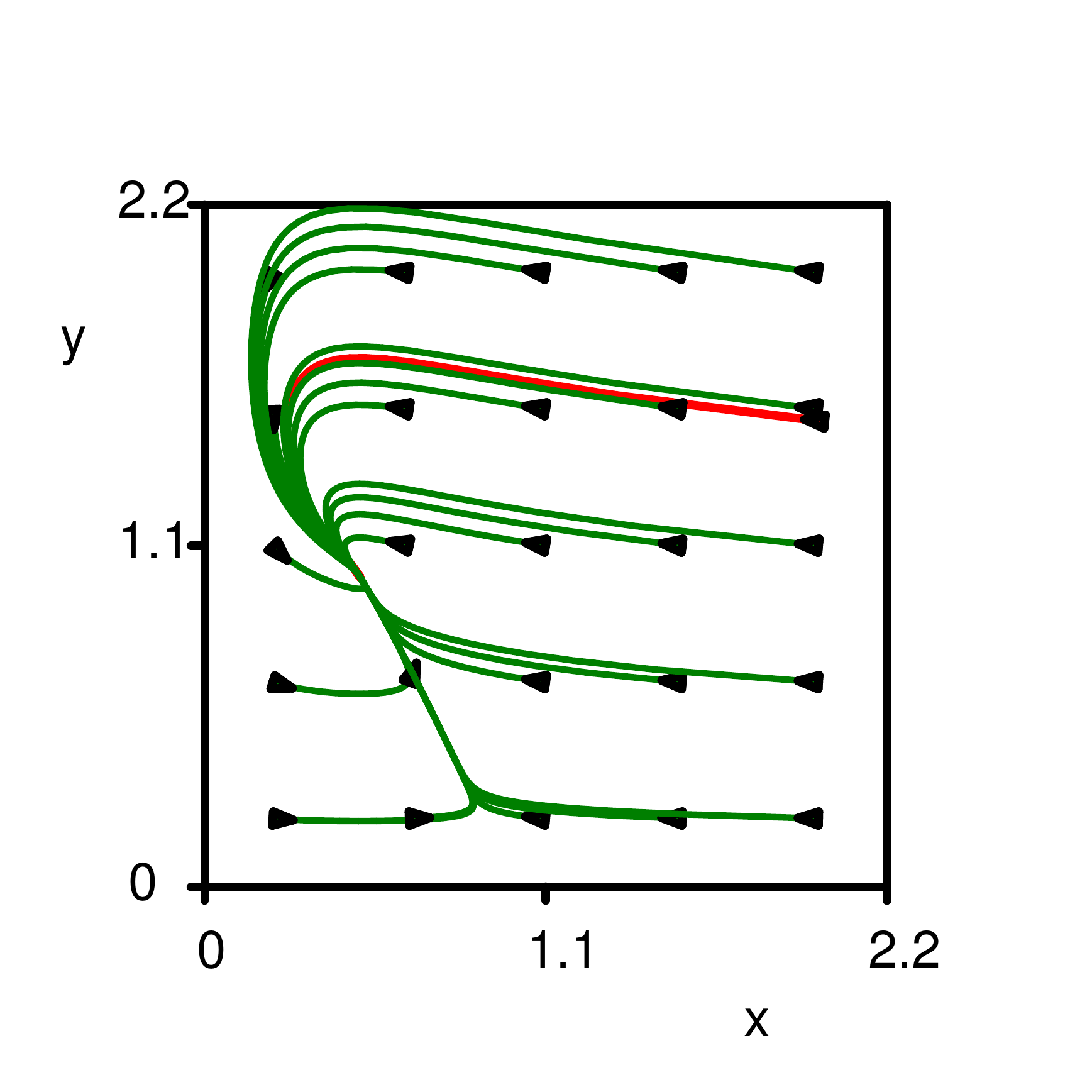,width=6.cm}
}
\caption{ A trajectory starting at point $x(0)=2,y(0)=1.5$ (a) and the complete phase portrait of system (\ref{PPour}) generated by a computer \label{2d2_fig}}
\end{figure}

Many biological systems are described by such systems. One of the
classical examples of ecological models (the predator-prey model) can
be derived as follows.  Let us consider the prey population $x$ with a
logistic growth given by eq.(\ref{e1.11}): ${dx \over dt}=rx(1-x/k)$ , which
interacts with the predator $y$ and let us assume that the effect of
the predator on the prey population is given by the term $-bxy$. Then,
if we assume that the growth of the predator population is
proportional to the predator prey interaction $cxy$ and that the death
rate of the predator is given by $-dy$, we will get the following
system of differential equations:

\beq
\label{PP1}
\left\{
\begin{array}{l}
{dx \over dt}=rx(1-x/k)-bxy\\ 
  {dy \over dt}=cxy-dy
\end{array}
\right.
\eeq
Formally  system (\ref{PP1}) describes the predator-prey
interactions with competition in the prey population. It has several
parameters, which account for the specific properties of the
populations. Let us study it for $r=3,k=1,b=1.5,c=0.5,d=0.25$:
\beq
\label{PPour}
\left\{
\begin{array}{l}
{dx \over dt}=3x(1-x)-1.5xy\\ 
 {dy \over dt}=0.5xy-0.25y
\end{array}
\right.
\eeq
If, as in section \ref{sec_2dlin}, we  solve this system on a
computer, we will get the following phase portrait
(fig.\ref{2d2_fig}b).  We see, that in the course of time all
trajectories approach the stationary values of $x=0.5;y=1$.

The main aim of our course is to develop the procedure of drawing a
phase portrait of a general system of two non-linear differential
equations without using a computer. We expect that as for 1D
differential equation (section \ref{sec_1Dphase}) and for 2D linear
systems (section \ref{sec_2dlin}) the  phase portrait should
include two main elements: equilibria points and flows (trajectories)
between them.  Let us define first equilibria points of a general
non-linear system (\ref{2dgeneral}).
\subsection{Equilibria \label{equilibrium_sec}}
In the 1D case and for 2D linear systems the equilibria were the
points where our system is stationary: placed at equilibrium point the
system will stay there forever. Therefore, for 1D equation ${dx \over dt}=f(x)$
equilibria were determined as the points where ${dx \over dt}=0$, i.e. where
$f(x)=0$. For 2D linear system (section \ref{sec_2dlin}) we required
that both variables $x$ and $y$ are stationary at equilibria points,
i.e. both ${dx \over dt}=0$ and ${dy \over dt}=0$. For a general
non-linear system (\ref{2dgeneral}) these conditions yield the
following definition of equilibria:
\begin{D}
A point $(x^*,y^*)$ is called an equilibrium point of a system (\ref{2dgeneral}) if 
\beq
\label{equilibrium}
f(x^*,y^*)=0, \quad g(x^*,y^*)=0 
\eeq
\end{D}
Equilibria in two dimensions are also stationary points, i.e.  if
system is placed at  the equilibrium  it will stay there
forever. Thus this trajectory will contain  just one point.

{\bf {Example}}. Find the equilibria of the system (\ref{PPour}): 

{\bf {Solution}}
To find the equilibria we need to solve a system of algebraic equations (\ref{equilibrium}) which in our case becomes:
\beq
\label{ex_eq}
\left\{
\begin{array}{l}
3x(1-x)-1.5xy=0\\ 
0.5xy-0.25y=0
\end{array}
\right.
\eeq
From the second equation we find $y(0.5x-0.25)=0$, which can be either when $y=0$ or when $x=0.5$.
Substitution of $y=0$ to the first equation yields $3x(1-x)-0=0$. This equation has two solutions $x=0$ and $x=1$. Substitution of the other case $x=0.5$ to the first equation gives $3*0.5*(1-0.5)-1.5*0.5y=0$, or $y=1$. Thus we have found three equilibria points: $(0,0),(1,0)$ and $(0.5,1)$.

We see in fig.\ref{2d2_fig} that point $(0.5,1)$ is indeed an
important attractor of our system which determines the final state of
the  populations. The other two points  are not
apparent in fig.\ref{2d2_fig}, however, as we will see later they also
account for important changes of trajectories of our system.

Thus we have defined  equilibria for  2D systems. Our next step is to
understand what is the 2D analog of flows,  which on the 1D phase portrait were represented by the 
'$\rightarrow$' or $\leftarrow$ arrows.

\subsection{Vector field \label{secVecField}}

In 1D flows,  visualizations of the direction of change of the
variable $x$  were given via the sign of its derivative ${dx \over dt}$. In 2D,  both
variables can change and the rate of their change is given by the
derivatives ${dx \over dt}$ and ${dy \over dt}$. In 1D we were able to find the
direction of flow at any point $x$ from the right hand side function
of the equation ${dx \over dt}=f(x)$.  Similarly in 2D we can find ${dx \over dt}$ and
${dy \over dt}$ at any point $(x,y)$ from the right hand sides of system
(\ref{2dgeneral}) ( functions $f(x,y)$ and $(g(x,y)$). For example for
system (\ref{PPour}) at a point $x=1,y=1$ we find
${dx \over dt}=f(x,y)=3x-3x^2-1.5xy=3-3-1.5=-1.5$, and ${dy \over dt}=
g(x,y)=0.5xy-0.25y.=0.5-0.25=0.25.$ However, what do these two numbers
show?  They tell us that if the size of the prey population $x=1$ and
the size of the predator population is $y=1$, then the prey population
decreases with the rate of $-1.5$ and the predator population grows
with the rate of $0.25$. On the phase plane $x,y$ this will result in
a shift of a point representing populations from point (1,1) (point
$A$ in fig.\ref{2d3_fig}a) to some point $B$ which is to the left and
upward from point $A$. Let us make it more quantitative. We know that
the rate of change of $x$ in our case is $1.5/0.25$ times larger than the
rate of change of $y$. This determines the direction of shift of point
$B$ relative to point $A$. The easiest way to represent it is to draw
from point $(1,1)$ a horizontal arrow heading to the left with the
length of $1.5$ and a vertical arrow heading upward with the length of
$0.25$. The direction of the overall shift will be given by the
resultant vector of these two vectors fig.\ref{2d3_fig}b. The
resultant vector will give us the direction tangent to the trajectory
which goes through the given point. We can generalize this result as:
\begin{figure}[hhh]
\centerline{
\psfig{type=pdf,ext=.pdf,read=.pdf,figure=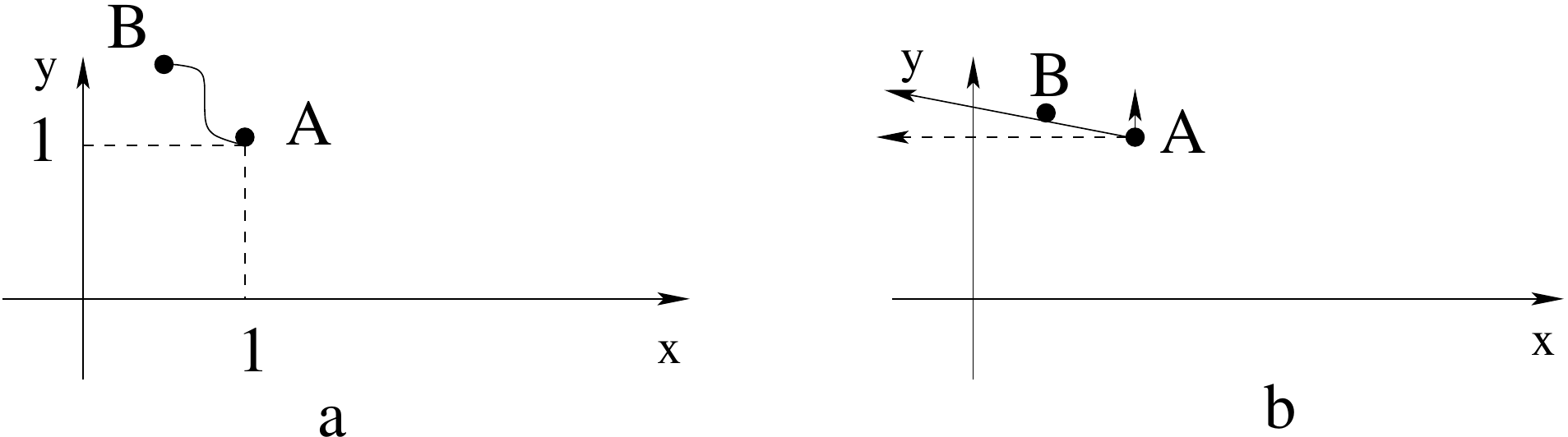,width=12.5cm}
}
\caption{\label{2d3_fig}}
\end{figure}

\bec
 At any point $(x,y)$ of a phase space for  system (\ref{2dgeneral}),
 we can define the vector $\vec{v}$ with the components
 $(f(x,y),g(x,y))$. Such vectors will be tangent to the trajectories
 of our system. We can find this vector field without a solution of
 our system, just from the right hand sides of our system.
\eec

\begin{figure}[hhh]
\centerline{
\psfig{type=pdf,ext=.pdf,read=.pdf,figure=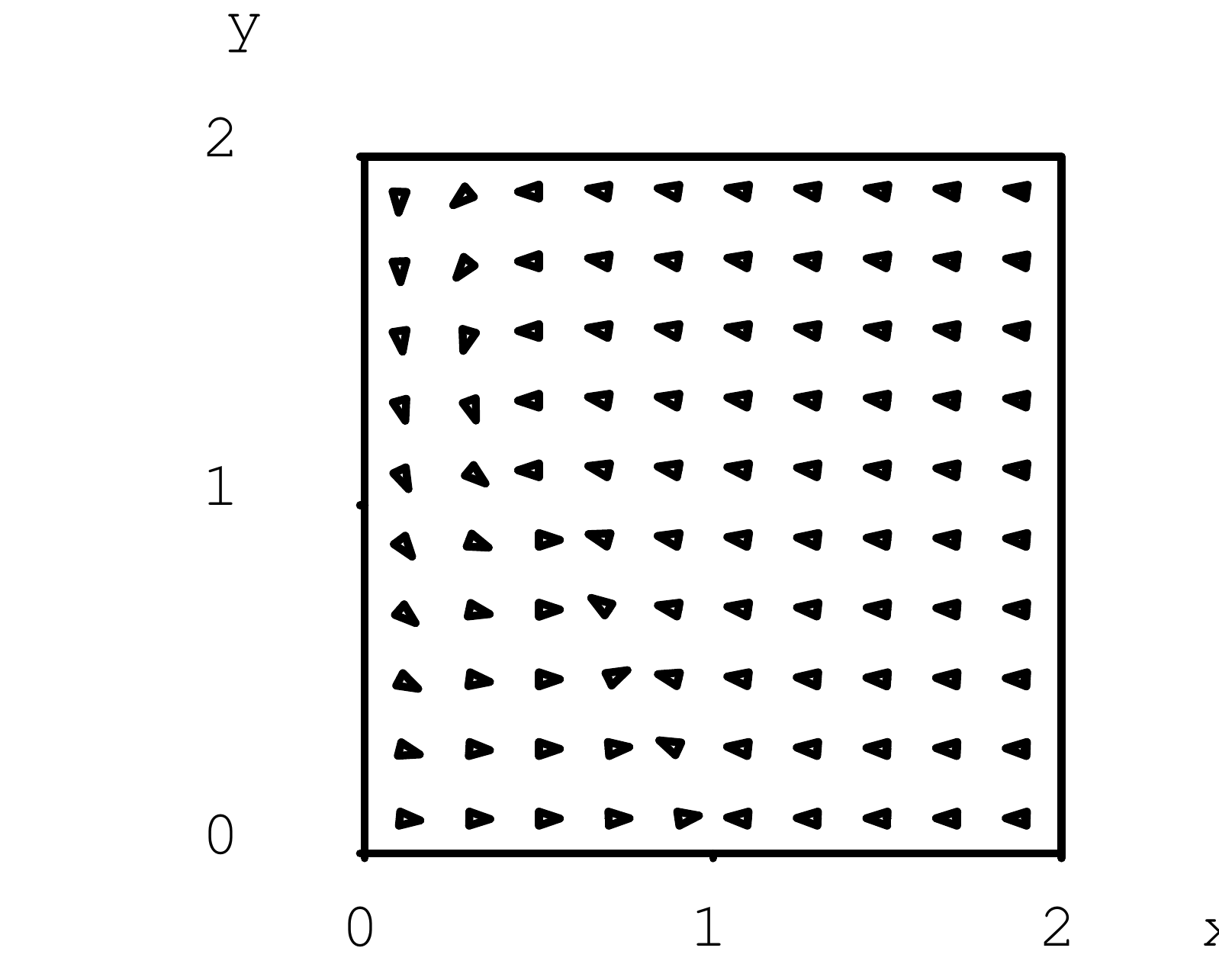,width=7.cm}
}
\caption{\label{2d4_fig}}
\end{figure}

 Note, that the length of the vector in fig.\ref{2d3_fig} is not important, as we are interested in the direction, only. If we
 apply the same procedure at  many points and represent the directions by
 shorter vectors we will get the following vector field of the  system
 (\ref{PPour}) (fig.\ref{2d4_fig}).
We see that although the vector field describes qualitatively the
direction of trajectories in the phase portrait of our system, it does
not give us information on convergence/divergence of trajectories, thus we cannot determine the  attractors of our system, which  is crucial for our
study. However, as we will see in chapter \ref{chap_graph} we can
elaborate the  vector field based methods and in many cases
will be able to obtain   convergence/divergence information using the 
so-called graphical Jacobian approach.

But in order to derive this approach we need to understand how to
find attractors of a general non-linear system (\ref{2dgeneral}) using
an analytical approach. For that we need to establish the relation of 
the general non-linear system (\ref{2dgeneral}) and the general linear system
(\ref{2dlin2}) that we studied in chapter \ref{chap_2dlin}.

\section{Linearization of a system: Jacobian \label{secJac}}

Consider a general system of two differential equations:
\beq
\label{2dgen}
\left\{
\begin{array}{l}
{dx \over dt}=f(x,y) \\ {dy \over dt}=g(x,y)
\end{array}
\right.
\eeq 

In this section we will show that close to equilibrium point the phase
portrait of this general non-linear system (\ref{2dgen}) can be found
from the solution of the linear system (\ref{2dlin2}) that we studied in
chapter \ref{chap_2dlin}. As a consequence we will get six possible
types of equilibria, that are:
saddle, non-stable and stable node, non-stable and stable spiral, and
a center. 

Our main tool here will be a formula (\ref{2dapprox_2}) for the 
approximation of function of two variables $f(x,y)$ around point
$(x^*,y^*)$ which we derived in section
\ref{sec_2DFunc}:
\beq
\label{2dapprox_2}
 f(x,y) \approx   f(x^*,y^*)+( \partial f / \partial x ) (x-x^*)+( \partial f / \partial y ) (y-y^*)
\eeq
where $\partial f / \partial x$ and $\partial f / \partial y$  are the values of the partial derivatives at the point $(x^*,y^*)$, i.e. they are just numbers.

We will apply this formula to approximate  the functions $f(x,y)$
and $g(x,y)$ of our system (\ref{2dgen}) and later solve the
approximated system and find the phase portraits close to the
equilibrium.

Let us start the derivation. Assume that system (\ref{2dgen}) has an
equilibrium point at ($x^*,y^*$). This means (see (\ref{equilibrium}))
that:
\beq
\left\{
\begin{array}{l}
f(x^*,y^*)=0 \\ g(x^*,y^*)=0
\end{array}
\right.
\eeq
Let us approximate $f(x,y)$ close to the equilibrium using the formula (\ref{2dapprox_2}):
$$
 f(x,y) \approx   f(x^*,y^*)+( \partial f / \partial x ) (x-x^*)+( \partial f / \partial y ) (y-y^*)
$$
As we assumed  ($x^*,y^*$) is an equilibrium, i.e. $f(x^*,y^*)=0$ and we get 
\beq
\label{2tmp1}
f(x,y) \approx   ( \partial f / \partial x ) (x-x^*)+( \partial f / \partial y ) (y-y^*)
\eeq
A similar approach for $g(x,y)$ yields:
\beq
\label{2tmp2}
 g(x,y) \approx   ( \partial g / \partial x ) (x-x^*)+( \partial g / \partial y ) (y-y^*)
\eeq

\begin{figure}[hhh]
\centerline{
\psfig{type=pdf,ext=.pdf,read=.pdf,figure=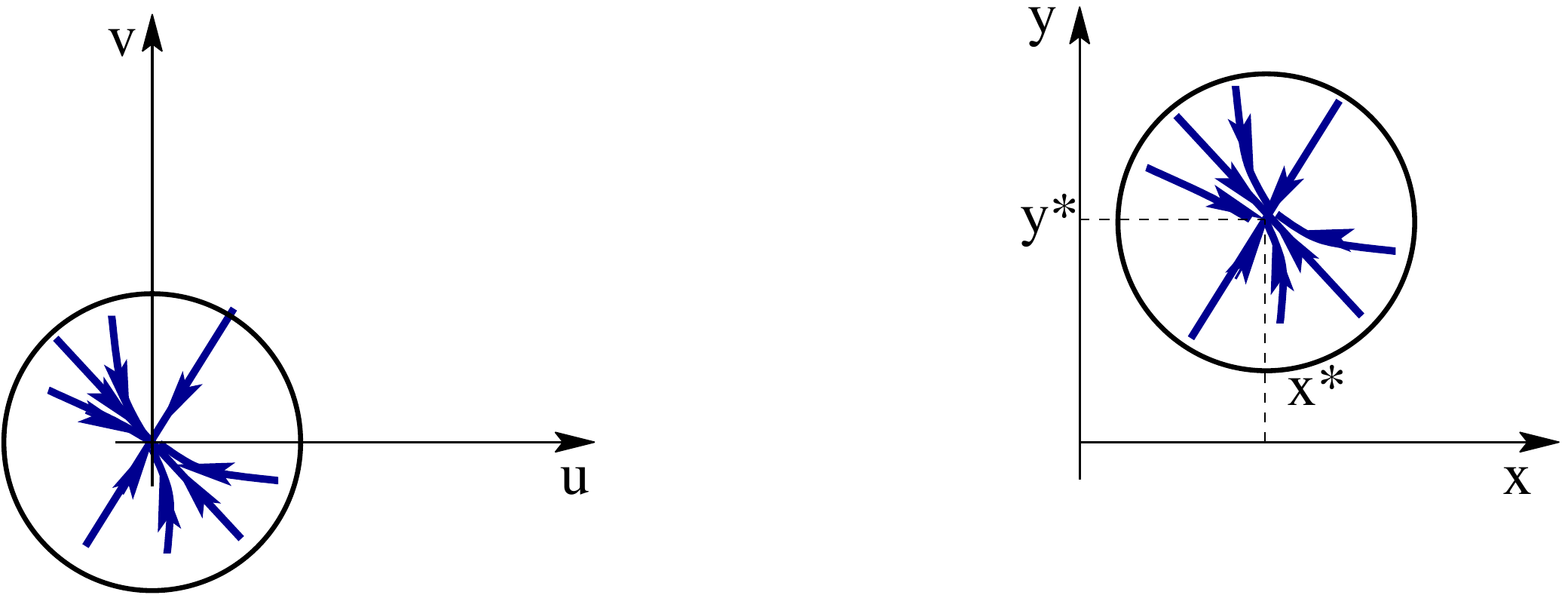,width=12cm}
}
\caption{\label{fig4.3}}
\end{figure}

If we replace the right hand sides of (\ref{2dgen}) by their
approximations (\ref{2tmp1}), (\ref{2tmp2}), we get the following
system:
\beq
\label{eq4.5}
\left\{
\begin{array}{l}
{dx \over dt}=  ( \partial f / \partial x ) (x-x^*)+( \partial f / \partial y ) (y-y^*) \\ {dy \over dt} =   ( \partial g / \partial x ) (x-x^*)+( \partial g / \partial y ) (y-y^*)
\end{array}
\right.
\eeq
The system (\ref{eq4.5}) is simpler than the original system
(\ref{2dgen}), as the partial derivatives in (\ref{eq4.5}) are {\it
{constants}} (numbers, as they are evaluated {\it {at the equilibrium
point}} $x^*,y^*$). So we can rewrite (\ref{eq4.5}) as :
\beq
\label{eq4.6}
\left\{
\begin{array}{l}
{dx \over dt}=  a (x-x^*)+ b (y-y^*) \\ {dy \over dt} =   c (x-x^*)+ d (y-y^*)
\end{array}
\right.
\eeq
where $a= \partial f / \partial x;b= \partial f / \partial y ; c=
\partial g / \partial x,d=\partial g / \partial y$. We can simplify
(\ref{eq4.6}) even more. For that let us introduce new variables:
\beq
\label{eq4.7}
u=x-x^* \qquad v=y-y^*
\eeq
were $u,v$ are new unknown functions of $t$. If we  substitute them
into the right hand side of 
(\ref{eq4.6}),we get:
\beq
\label{eq4.8}
\left\{
\begin{array}{l}
{dx \over dt}=  a u+ b v \\ {dy \over dt} =   c u+ dv
\end{array}
\right.
\eeq
In order to substitute $u$ and $v$ into the left hand side of
(\ref{eq4.8}), note that $u(t)=x(t)-x^*$, i.e. ${du \over dt}={dx
\over dt}-0$ (here ${dx^* \over dt}=0$ because $x^*$ is a
constant). Similarly, ${dv \over dt}={dy \over dt}$. After replacing
${dx \over dt}$ by ${ du \over dt}$ and ${dy \over dt}$ by ${ dv \over
dt}$ in (\ref{eq4.8}) we get:
\beq
\label{eq4.9}
\left\{
\begin{array}{l}
{du \over dt}=  a u+ b v \\ {dv \over dt} =   c u+ dv.
\end{array}
\right.
\eeq
System (\ref{eq4.9}) coincides with a general linear system
(\ref{2dlin2}) that we studied in chapter \ref{chap_2dlin} and for
which we found six possible types of solutions resulting  in six
possible phase portraits: saddle, non-stable and stable node,
non-stable and stable spiral, and a center. However, how can we use
the results of study of system (\ref{eq4.9}) for the study of the original
system (\ref{2dgen})?  In order to derive (\ref{eq4.9}) we made two
steps: (1) we used formula (\ref{2dapprox_2}) for function
approximation; (2) we changed variables $x,y$ to $u,v$. Let us analyze
each of these two steps. 
(1) As we discussed earlier, formula (\ref{2dapprox_2}) gives a good
approximation of the function $f(x,y)$ only if $x,y$ is  close to the
point of approximation $x^*,y^*$. So we can use the linear system 
 (\ref{eq4.9}) for approximating the solutions of the non-linear system 
(\ref{2dgen}) only  close to the equilibrium point $(x^*,y^*)$.
(2) about change of variables. Equation (\ref{eq4.7}) gives the
variables $u,v$ via the variables $x,y$. However, we can also solve the 
equations and find how $x,y$ will be expressed in terms of $u,v$:
\beq
\label{eq4.7a}
x=u+x^* \qquad y=v+y^*
\eeq
Using this expression we can draw the trajectories for our original
variables $x(t),y(t)$ if we know the trajectories $u(t),v(t)$ of the
linearized system (\ref{eq4.9}).  Indeed, as the $x$ coordinate of the
trajectory equals $u$ plus number $x^*$, and the $y$ coordinate equals
$v$ plus number $y^*$, the only thing what we need to do is just to
draw the trajectory $u(t),v(t)$ and shift its $x$-coordinate by $x^*$
and the $y$-coordinate by $y^*$ units. As we know in a general linear
system (\ref{eq4.9}) the phase portrait is centered around a point
$(0,0)$, therefore all what we would need to do in order to draw the
phase portrait of the non-linear system (\ref{2dgen}) close to its
equilibrium $(x^*,y^*)$ is just to shift the phase portrait of the
linear system (\ref{eq4.9}) to a location of equilibrium in non-linear
system (\ref{2dgen}). If, for example, the linearized system
(\ref{eq4.9}) will have a stable node type phase portrait, then a
non-linear system (\ref{2dgen}) will have the same stable node point
but around an equilibrium $(x^*,y^*)$ as shown in fig.\ref{fig4.3}.

\bec
The phase portrait of liner system (\ref{eq4.9}) close to the origin
($u=0,v=0$) is similar to the phase portrait of non-linear system
(\ref{2dgen}) close the to equilibrium point $(x^*,y^*)$. To draw the
phase portrait of non-linear system (\ref{2dgen}) close to
equilibrium, we need to shift  a phase portrait of linear system (\ref{eq4.9}) from the origin  the equilibrium point $(x^*,y^*)$.
\eec

To find the linearized system (\ref{eq4.9}) we need to find the
equilibrium point ($x^*,y^*$) and compute the following four numbers:
the values of derivatives of right hand sides of our system at this
equilibrium: 
$$ a= \partial f / \partial x \;\; b= \partial f /
\partial y \;\; c=
\partial g / \partial x \;\; d=\partial g / \partial y $$ So, 
system (\ref{eq4.9}) can be written as:

\beq
\left\{
\begin{array}{l}
{du \over dt}=    {\partial f \over  \partial x} u+  { \partial f  \over  \partial y } v \\ {dv \over dt} =    { \partial g   \over  \partial x} u+ {\partial g  \over  \partial y }  v.
\end{array}
\right.
\eeq

 From coefficients of this system we can construct a matrix $J$  that is
 called {\bf {the Jacobian}}
\beq
\label{jac}
J=\left(\begin{array}{lr}  {\partial f \over  \partial x} & { \partial f  \over  \partial y }\\
                         { \partial g   \over  \partial x} &  {\partial g  \over  \partial y } 
\end{array} \right)
\eeq

 {\bf {Example}} Find qualitative phase portrait of system
 (\ref{PPour}) close to a nontrivial equilibrium point (i.e. at the
 equilibrium where $x\neq 0,y \neq 0$).  $$
\left\{
\begin{array}{l}
{dx \over dt}=3x(1-x)-1.5xy\\ 
 {dy \over dt}=0.5xy-0.25y
\end{array}
\right.
$$

{\bf {Solution}} In section \ref{equilibrium_sec} we found that this
system has three equilibria $(0,0)$, $(1,0)$, and $(0.5,1)$. The
nontrivial equilibrium is $(0.5,1)$. To find the Jacobian of our system
we compute the partial derivatives (\ref{jac})  and evaluate them at the
equilibrium point. In our case $f(x,y)=3x(1-x)-1.5xy;\; g(x,y)=0.5xy-0.25y$.

$ {\partial f /  \partial x}=3-6x-1.5y$ at  point $(0.5,1)$ this derivative equals $ {\partial f /  \partial x}=3-3-1.5=-1.5$. Similarly: $ {\partial f /  \partial y}=-0.75; \; {\partial g /  \partial x}=0.5;\; {\partial g /  \partial y}=0$, hence the linearization of our system at  point $(0.5,1)$ is
$$
\left\{
\begin{array}{l}
{du \over dt}=  -1.5u- 0.75v \\ {dv \over dt} =   0.5u
\end{array}
\right.
$$
or the Jacobian is:
$$
J=\left(\begin{array}{lr}  -1.5 & -0.75\\
                         0.5 &  0
\end{array} \right).
$$

In order to find a phase portrait of this linear system we need to find eigen values of the Jacobian matrix from the following characteristic equation (\ref{eigen_char}):
\beqar
 Det
\left|\begin{array}{lr} -1.5-\lambda & -0.75\\ 0.5 & 0-\lambda
\end{array} \right|=(-1.5-\lambda)(-\lambda)+0.5*0.75\\
=\lambda^2+1.5\lambda+0.375=0
\eeqar
Using 'abc' formula we find that:

\beq
 \lambda_{12}={-1.5 \pm  \sqrt {2.25-1.5} \over 2}={-1.5 \pm  \sqrt{0.87} \over 2}
\eeq
or, $\lambda_1=-1.36$ and  $\lambda_1=-0.138$. Because both eigen values are real and negative we will have a stable node and a phase portrait qualitatively similar to that in fig.\ref{fig4.3}b.

\section{Determinant-trace  method for finding the type of  equilibrium \label{dettr_sec}}

In this section we derive a simple method for finding signs of eigen
values and the type of equilibrium of the linear system
(\ref{2dlin2}).  The results will also be applicable for a general
non-linear system (\ref{2dgen}), because, as discussed in the previous
section, the equilibrium type of a non-linear system can be found from
its linearization.

The eigen values of system (\ref{2dlin2}) are the roots of the characteristic equation (\ref{eigen_char}):
\beq
 Det \left|\begin{array}{lr}  a-\lambda & b\\
                         c &  d-\lambda
\end{array} \right|=0,
\eeq
or:

$ Det \left|\begin{array}{lr}  a-\lambda & b\\
                         c &  d-\lambda
\end{array} \right|= (a-\lambda)(d-\lambda)-cb={\lambda}^2-\lambda(a+d)+ad-cb=0$

Let us express the last equation in a slightly different form using the following definition:
\begin{D}
The trace of the  matrix 
$\left(\begin{array}{lr}  a & b\\
                         c &  d
\end{array} \right) $ is $trA=a+d$. 

The determinant of the matrix $A$ is $det A=ad-cb$
\end{D}
So we can rewrite the characteristic equation using the definition of the  trace and the determinant as follows:
\beq
\label{simple2}
{\lambda}^2-trA \lambda +det A=0
\eeq

We see that although the original system depends on four parameters
$a,b,c,d$ the characteristic equation depends only on two parameters
$trA$ and $det A$, thus if we know the determinant and the trace of
our system we can find the eigen values and the type of the
equilibrium of the system. Indeed, from (\ref{simple2}) we can easily
find that the eigen values are:
\beq
\label{simple3}
 \lambda_{1,2}={trA \pm  \sqrt {D} \over 2}
\quad where \quad  D=(trA)^2-4detA
\eeq
 Roots of the equation (\ref{simple3}) are as the roots of
 any quadratic equation,  connected in the following way to the coefficients of the equation:
\beq
\label{simple4}
 \lambda_1+\lambda_2 =trA
\eeq
\beq
\label{simple5}
\lambda_1 * \lambda_2 = det A
\eeq
To prove the properties (\ref{simple4}) and (\ref{simple5}),  just note that if $\lambda_1$ and $\lambda_2$ are the roots of a quadratic (\ref{simple2}), then it can be written as: ${\lambda}^2-trA \lambda +det A=(\lambda- \lambda_1)((\lambda- \lambda_2)$. The direct computation yields:
$$
\begin{array}{l}
 {\lambda}^2-trA \lambda +det A=(\lambda- \lambda_1)((\lambda- \lambda_2)\\
or\\
 {\lambda}^2-trA \lambda +det A=\lambda^2- \lambda_1 \lambda- \lambda_2 \lambda+ \lambda_1 \lambda_2\\
or\\
 {\lambda}^2-trA \lambda +det A=\lambda^2- (\lambda_1+\lambda_2) \lambda+ \lambda_1 \lambda_2
\end{array}
 $$
If we now compare  the left and the right hand sides of the last equation we will get both properties  (\ref{simple4}) and (\ref{simple5}).

Let us start classification.
\ben
\item
If $detA<0$, then $D=(trA)^2-4detA >0$,  so we have real roots. From (\ref{simple5}) we conclude that their product is negative, i.e. roots have different signs, i.e.
$\lambda_1<0,\lambda_2>0$, or  $\lambda_1>0,\lambda_2<0$ and we have a {\it saddle point}.
\item If $detA> 0$, then $D=(trA)^2-4detA$  can be negative as well 
as positive. This means the roots can be real, or complex. 
Let us consider the case of real roots first, i.e.
\beq
\label{simple6}
D=(trA)^2-4detA  \geq 0
\eeq
If (\ref{simple6}) holds,  the roots are real. Next, let us use the property $\lambda_1 * \lambda_2 = det A$. In the case of  $detA> 0$,  the product of the roots is positive, i.e.  the  roots have the same sign. They can be  both positive, or both negative. The sign of the roots can be found from the trace of the matrix ($ \lambda_1+\lambda_2 =trA$). When  $trA>0$, then  $\lambda_1>0$, and $\lambda_2>0$ and we have a {\it non-stable node}. When   $trA<0$,     $\lambda_1<0$ and $\lambda_2<0$ and we have a stable node. Let us formulate it as a separate case:

\item If $detA> 0$, $D>0$ and  $trA<0$ the equilibrium is a {\it stable node}.

\begin{figure}[hhh]
\centerline{
\psfig{type=pdf,ext=.pdf,read=.pdf,figure=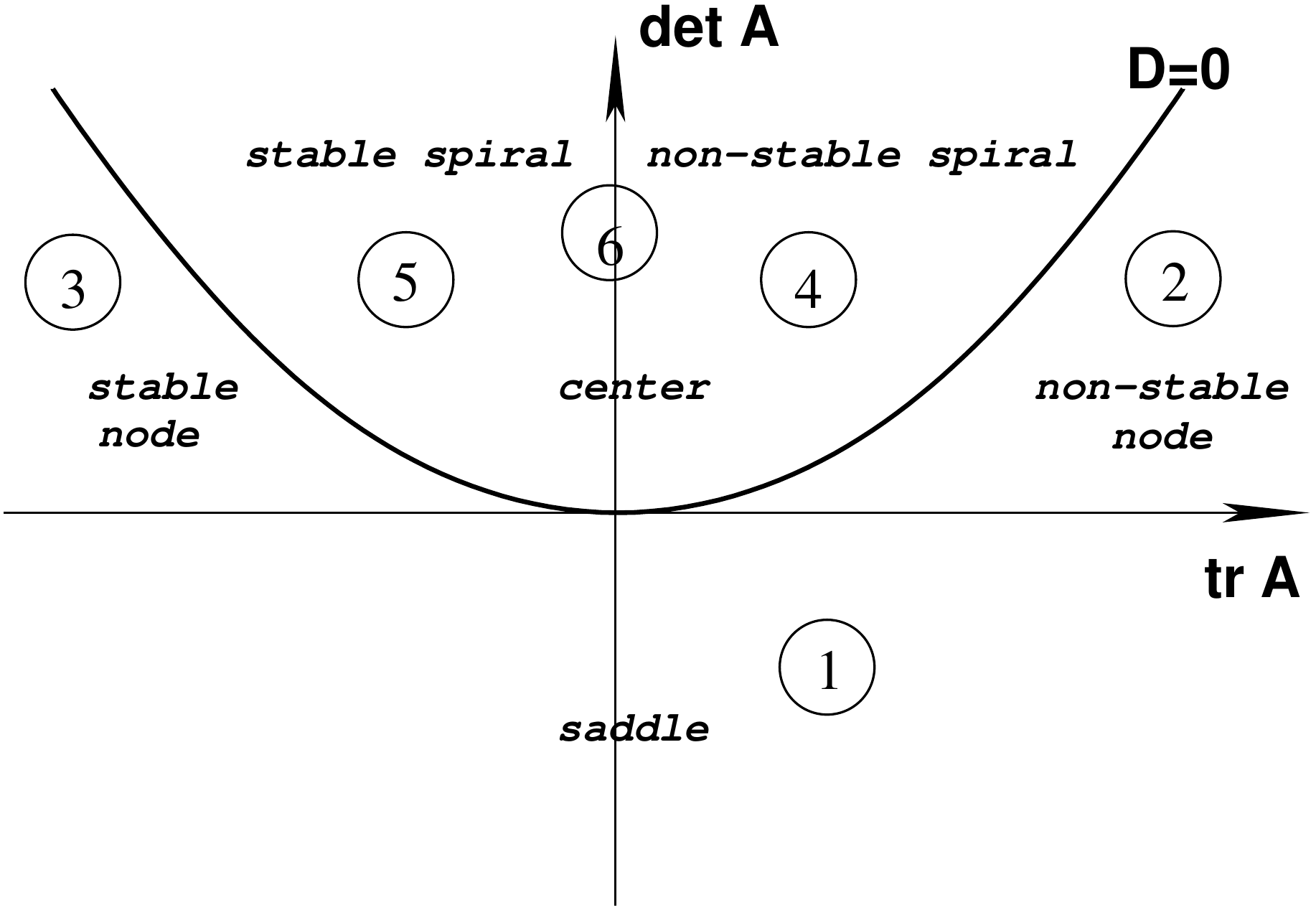,width=11cm}
}
\caption{\label{fig8.2}}
\end{figure}

Let us put this information into a  graph (fig.\ref{fig8.2}). On this
graph let us use $trA$ as the $x$-axis and $detA$ as the
$y$-axis. Case 1 of a saddle point then corresponds to the lower half
plane (region 1). The line (\ref{simple6}), which separates real and complex
roots,  is the parabola given by $detA=(trA)^2/4$, or $y=x^2/4$. Real roots are
below this line. Therefore, in region 2, where $trA>0$, we have case 2 of  
 a non-stable node.  In  region 3  $trA<0$ and we have  case 3 of a stable
node.

\item If  $detA> 0$, and  $D<0$,  we have complex roots. In accordance with (\ref{simple4}) and (\ref{complex3}) they are:
$$
 \lambda_{1,2}={trA \over 2} \pm i{ \sqrt {-D} \over 2}
$$
Hence, $Re\lambda_{1,2}={trA / 2}$. From this we immediately see  that if
$trA>0$,  we have a {\it non-stable spiral} (region  4).

\item If  $detA> 0$,   $D <0$ and $trA<0$,  we have a {\it stable spiral} (region  5).

\item The last case of a {\it center point} appears when  $Re\lambda_{1,2}={trA / 2}=0$, or when  $trA=0$ and  $detA> 0$. In our graph is it  the upper part of the $detA$ axis (region  6)
\een

Let us apply this method for several examples:

{\bf {Example}}. Find the  equilibrium type for the following systems:

a)$\left( \begin{array}{c} 
{dx \over dt} \\  {dy \over dt}  \end{array} \right) = \; \left(\begin{array}{lr}  1 & 2\\
                         3 &  4
\end{array} \right)  \left( \begin{array}{c} x
 \\  y  \end{array} \right); \quad$\hskip 1pc
b)$\left( \begin{array}{c} 
{dx \over dt} \\  {dy \over dt}  \end{array} \right) = \; \left(\begin{array}{lr}  4 & 1\\
                         1 &  2
\end{array} \right)  \left( \begin{array}{c} x
 \\  y  \end{array} \right); \quad$\\
c)$\left( \begin{array}{c} 
{dx \over dt} \\  {dy \over dt}  \end{array} \right) = \left(\begin{array}{lr}  -2 & 3\\
                         -2 &  1
\end{array} \right)  \left( \begin{array}{c} x
 \\  y  \end{array} \right); \quad$
d)$\left( \begin{array}{c} 
{dx \over dt} \\  {dy \over dt}  \end{array} \right) = \left(\begin{array}{lr}  1 & -1\\
                         2 &  -1
\end{array} \right)  \left( \begin{array}{c} x
 \\  y  \end{array} \right); \quad$

{\bf {Solution}} Our plan is to find $detA, trA$ for the corresponding matrices and make a conclusion.

a) The matrix is $A= \left(\begin{array}{lr}  1 & 2\\
                         3 &  4
\end{array} \right);\; trA=1+4=5;\;detA=1*4-2*3=-2$. Therefore  we have case   1, hence the system has a saddle point.

b) The corresponding matrix is $A= \left(\begin{array}{lr}  4 & 1\\
                         1 &  2
\end{array} \right);\; trA=4+2=6;\;detA=4*2-1*1=7$. As $detA>0$ we need to check that  the roots are real. We have $  D=(trA)^2-4detA=6*6-4*7=36-28 >0$, hence the roots are real and we have  case   2, hence the system has a non-stable node.

c) $trA=-2+1=-1;\; detA=-2+6=4;\;D=(trA)^2-4detA=1-4*4=-15<0$, we have complex roots and we are in  region   5 and have a stable spiral

d) $trA=1-1=0;\;detA=-1+2=1$, so we have  case   6 and we have a center point, or oscillation in our system.
\section{Exercises \label{exCh5} }
\subsection*{Exercises sec.\ref{sec2dnonlin}}
\ben
\item \label{ProbSec5} (A) Find equilibria of the  following systems  $\left\{
\begin{array}{l}
{dx \over dt} =f(x,y)\\ 
{dy \over dt}=g(x,y)
\end{array}
\right.$ (see definition in section \ref{equilibrium_sec}).\\ (B) Find the following partial derivatives at each equilibrium point ($ {\partial f \over  \partial x}, { \partial f  \over  \partial y }, { \partial g   \over  \partial x},  {\partial g  \over  \partial y } $).\\
(a) $\left\{
\begin{array}{l}
{dx \over dt} =-4y\\ 
{dy \over dt}=4x-x^2-0.5y
\end{array}
\right.$
(b) $\left\{
\begin{array}{l}
{dx \over dt} =9x + y^2\\ 
{dy \over dt}=x-y
\end{array}
\right.$

(c) $\left\{
\begin{array}{l}
{dx \over dt} =2x-xy\\ 
{dy \over dt}=-y +y^2x
\end{array}
\right.$
(d)
 $\left\{
\begin{array}{l}
{dx \over dt} =(1-x-3y)x\\ 
{dy \over dt}=(1-2x-2y)y
\end{array}
\right.$

(Hint: See an example of solution  in section \ref {solutionsCh5.1}).

\item Find equilibria of the following  Lotka Volterra model  with competition in the prey population. Determine for which parameter values
all  equilibria are non-negative. 
$$
\left\{
\begin{array}{l}
dN/dt=aN-eN^2-bNP \\ 
  dP/dt=cNP-dP \;\;\;\;\;\;\;a,b,c,d,e >0
\end{array}
\right.
$$

\item Protein synthesis depends on DNA transcription ($a$) making mRNA molecules ($M$) and translation ($c$) of mRNA
into proteins ($P$). Some proteins inhibit the transcription of their own mRNA ($\frac{1}{1+P}$). mRNA and proteins are degraded at rates b and d, respectively. This process gives the following of two differential equations.
Find equilibria of this model.
$$
\left\{
\begin{array}{l}
dM/dt={a \over 1+P}-bM \;\; P,M \geq 0 \\ 
  dP/dt=cM-dP \;\;\;\;\;\;\;a,b,c,d >0
\end{array}
\right.
$$

\item 
Mathematical epidemiology also makes use of simple ODE models. One of
these models describes the number of susceptible individuals ($S$) and
infected individuals ($I$). Individuals are born at rate ($B$), and
die at rate ($\mu$). Susceptible individuals can become infected when
they come into contact with infected individuals ($-\beta S I$). Once
infected, an individual has a certain death rate ($\alpha$); this may
be different from the death rate of non-infected individuals. This
process can therefore be modeled by the following equations:
$
\left\{
\begin{array}{l}
dS/dt=B-\beta SI-\mu S\\ 
  dI/dt=\beta SI-\alpha I \;\;\;\;\;\;\;B,\alpha,\beta,\mu>0
\end{array}
\right.
$\\
Find equilibria of this model.  Determine for which parameter values
all equilibria are non-negative.

\subsection*{Exercises sec. \ref{secJac} and \ref{dettr_sec} }

\item  Find the type of equilibria using the $det$-$tr$ method. Determine the stability of  the equilibrium.\\
(a) $\left\{
\begin{array}{l}
{dx \over dt} =3x+y\\ 
{dy \over dt}=-20x+6y
\end{array}
\right.$
(b) $\left\{
\begin{array}{l}
{dx \over dt} =2x+y\\ 
{dy \over dt}=2x-10y
\end{array}
\right.$
(c) $\left\{
\begin{array}{l}
{dx \over dt} =2x+y\\ 
{dy \over dt}=5x-2y
\end{array}
\right.$
(d) $\left\{
\begin{array}{l}
{dx \over dt} =x+10y\\ 
{dy \over dt}=-10x-y
\end{array}
\right.$

\item \label{ProbAlgae} Consider the following model for the algae population:  $\left\{
\begin{array}{l}
{dx \over dt} =2x(1-y) \;\;x \geq 0;\\ 
{dy \over dt}=2-y-x^2 \;\; y \geq 0.
\end{array}
\right.$ 
\ben 
\item Find equilibria

\item Find the general expression for the Jacobian of this system

\item Determine type of each equilibrium using $det$-$tr$ method

\item Draw qualitative local phase portraits around each equilibrium point 
\een

\item
 Consider the following biological model:

\begin{equation}
\label{8ex}
\left\{
\begin{array}{l}
de/dt=b*e-e^3-g\\
  dg/dt=e-g \;\;\;b \geq 0
\end{array}
\right.  
\end{equation}

\bit
\item For which values of parameter $b$  the system has only one equilibrium?

\item Determine  stability and type of this equilibrium in found   parameter range (in which    system  (\ref{8ex})   has only
one equilibrium).
\eit

\item \label{ProbLV} Study the Lotka-Volterra model for a predator-prey system:
$
\left\{
\begin{array}{l}
{dN \over dt} =aN -bNP\\ 
{dP \over dt}=cNP-dP
\end{array}
\right.
$\\ 
here $N$ denotes the prey population, $P$ denotes the predator population and $a>0,b>0,c>0,d>0$ are  parameters
\ben
\item Find  the nontrivial equilibrium of the system (i.e. an equilibrium where $N\neq 0,P \neq 0$).

\item Find the linearization of the  system at this point (i.e. the Jacobian matrix)

\item Determine the type of the equilibrium

\item Sketch the phase portrait around this equilibrium. Which kind of dynamics do we expect here?
\een


\subsection*{Additional exercises }
\item Consider the system: ${dx \over dt}=x+4y+e^x-1; {dy \over dt}=-y-y*e^x$

\ben 
\item Check that ($0,0$) is an equilibrium point of the system

\item Find the general expression for the Jacobian of this system

\item Find  the Jacobian at the point ($0,0$)

\item Write the linearization of the system close to the equilibrium ($0,0$)
\een

\item Find the equilibria of the following systems. 
Compute the Jacobian at  the equilibria points. Determine type of each equilibrium using $det$-$tr$ method. Draw qualitative local phase portraits around each equilibrium point.\\
(a)$\left\{
\begin{array}{l}
{dx \over dt} =y^2-3x+2\\ 
{dy \over dt}=x^2-y^2
\end{array}
\right.$
(b) $\left\{
\begin{array}{l}
{dx \over dt} =y\\ 
{dy \over dt}=-x+x^3
\end{array}
\right.$

\item Revisit solutions of problems \ref{ProbSec5}. Use found equilibria  and partial derivatives at these equilibria points to determine type of each equilibrium using $det$-$tr$ method. Draw qualitative local phase portraits around each equilibrium point.

\item Consider a modification of the Lotka-Volterra model, which includes competition in the prey population ($-eN^2$):
$
\left\{
\begin{array}{l}
{dN \over dt}=aN-eN^2-bNP\\ 
 { dP \over dt}=cNP-dP
\end{array}
\right.
$,
where the parameters $a,b,c,d,e>0$ and the variables $N \geq 0$, $P \geq 0$.
\ben
\item Find all equilibria of this system.

\item Compute the Jacobian at each equilibrium point.

\item At which parameter values do we have a non-trivial equilibrium (i.e. an equilibrium at
which both $N$ and $P$ are positive). Find stability of this equilibrium.

\een

\item 
Consider the following model for  cardiac tissue:

\begin{equation}
\label{9ex}
\left\{
\begin{array}{l}
de/dt=-e(e-a)(e-1)-g \qquad 0<a<1\\
  dg/dt=\varepsilon e \qquad \varepsilon >0
\end{array}
\right.  
\end{equation}
Here the variable $e$ accounts for the transmembrane potential,
the variable $g$ accounts for the refractory period and $a,\varepsilon$ are the parameters.

\begin{figure}[hhh]
\centerline{
\psfig{type=pdf,ext=.pdf,read=.pdf,figure=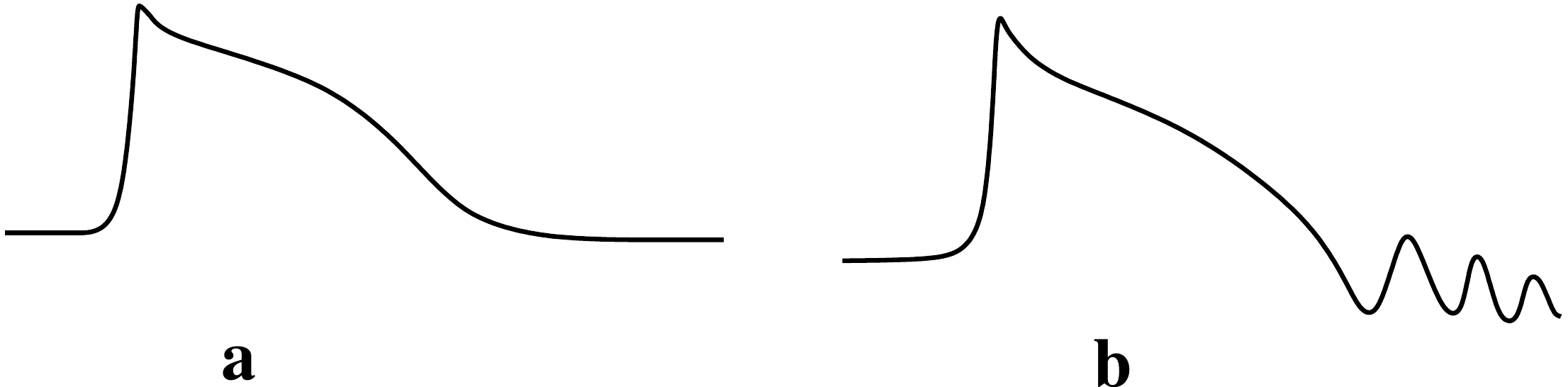,width=7.cm}
}
\caption{\label{fig8.7} The monotonous (a) and oscillatory recovery (b) in an excitable medium }
\end{figure}

The shape of the action potential in cardiac tissue  is an important characteristic 
of myocardium.  If the recovery of the transmembrane potential shows oscillation as in fig.\ref{fig8.7}b.
it can cause  dangerous cardiac arrhythmias. 
From a mathematical point of view the oscillations in fig\ref{fig8.7}b occur when system (\ref{9ex})
has an  equilibrium point which is a stable spiral. Monotonous recovery occurs
when this  equilibrium is a stable node. 

Determine for which parameter values  we will have 
situation fig.\ref{fig8.7}a and for which parameter values  we will have 
 situation fig.\ref{fig8.7}b.

\een

\chapter{Graphical  methods to study systems of differential equations \label{chap_graph}}
In the first section of this  chapter we will start from the vector field of a general
non-linear system introduced in section \ref{secVecField},  and find how
we can approximate the vector field by the so-called null-cline method,
without using a computer. Then, in section
\ref{gjac_sec},  we will show
that the null-clines can be used not only for vector field
approximation but also for determining the type of an equilibrium point
without explicit computation of the Jacobian of the system.

\section{Null-clines \label{secNull}}
\begin{figure}[H]
\centerline{
\psfig{type=pdf,ext=.pdf,read=.pdf,figure=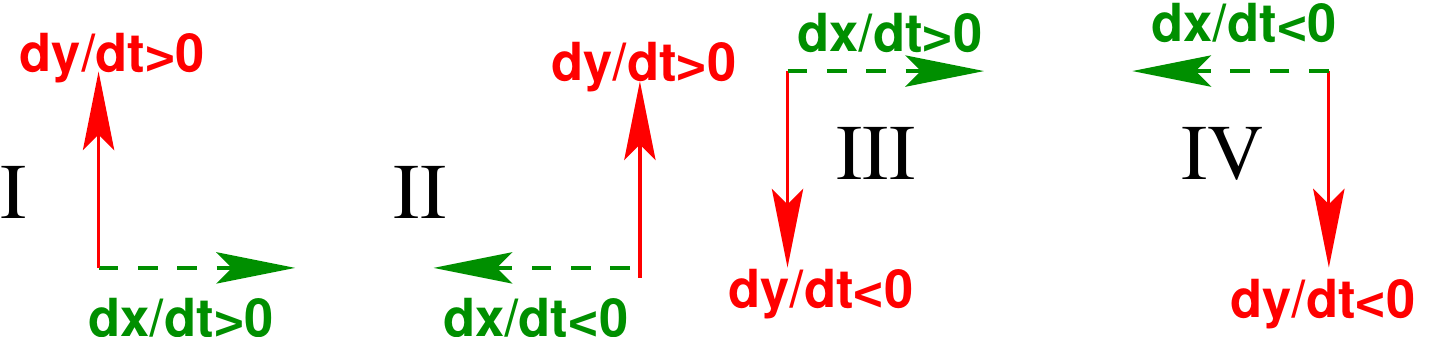,width=10.cm}
}
\caption{\label{2d5_fig}}
\end{figure}

We introduced the vector field in section \ref{secVecField}, where we showed that for a 
general 2D system
\beq
\label{2dgeneral_3}
\left\{
\begin{array}{l}
{dx \over dt}=f(x,y) \\ {dy \over dt}=g(x,y)
\end{array}
\right.
\eeq
the direction of the trajectories on the phase portrait will be along
the vector $\vec{v}$ with the components $(f(x,y),g(x,y))$. Thus, to
draw the vector field of a particular system we need to evaluate
the values of the functions $(f(x,y),g(x,y))$ in many points which  usually
requires using  a computer. In this section we will develop a
so-called method of null-clines, which will allow us to sketch a
qualitative picture of the vector field analytically.  The main idea
here is similar to what we did in the 1D case, where we have represented
the derivative ${dx \over dt}$ by arrows of two types: '$\rightarrow$' 
for ${dx \over dt}>0$ and '$\leftarrow$'  for ${dx \over
dt}<0$. In 2D the vector field has two components $V=(v_x,v_y)=({dx
\over dt},{dy \over dt})$. Each of these components ${dx \over dt},{dy
\over dt}$ can be positive, or negative. Therefore, we can have the following
four cases shown fig.\ref{2d5_fig}. (Note, that we use different line types there that will be important in the future).

\begin{figure}[hhh]
\centerline{
\psfig{type=pdf,ext=.pdf,read=.pdf,figure=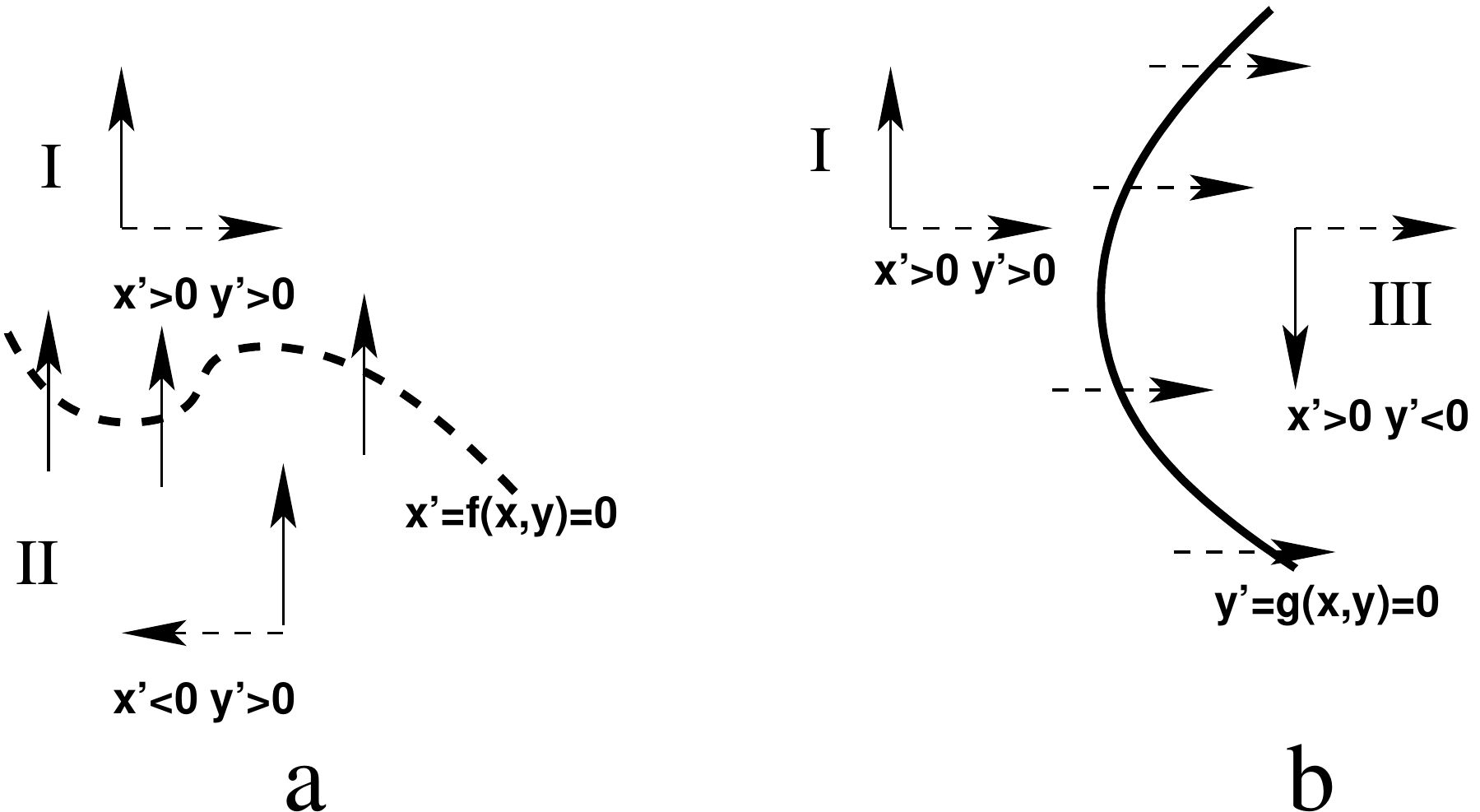,width=10.cm}
}
\caption{\label{2d6_fig}}
\end{figure}

The method of null-clines represents the vector field using these
four cases. 
The main idea behind this method is the following.
If we compare cases I and II we see that they differ by the sign of the ${dx \over dt}$ derivative:
for case I ${dx \over dt}>0$ and for case II  ${dx \over dt}<0$. Therefore these cases are separated by the boundary
 where ${dx \over dt}=0$ (fig.\ref{2d6_fig}a). We know that  for system (\ref{2dgeneral_3}) ${dx \over dt}=f(x,y)$, therefore
the boundary between cases I and II is given by the condition 

\beq
\label{cond1}
f(x,y)=0
\eeq

Geometrically equation (\ref{cond1}) gives a graph of one or more
lines in the $Oxy$-plane (see section \ref{secImplicitGr}).
Note that at this line the horizontal component of the vector field is zero, therefore the vector
is {\it vertical}.

Similarly, the transition from case I to case III occurs when ${dy
\over dt} =0$ (fig.\ref{2d6_fig}b). As we know for the system (\ref{2dgeneral_3})
${dy \over dt}=g(x,y)$, hence the separation line in this case is
given by the equation:
\beq
\label{cond2}
g(x,y)=0
\eeq
and the direction of the vector field at this line is {\it horizontal}.

Equation (\ref{cond2}) will give us not only the boundary between
cases I and III but also a boundary between cases II and IV, as the
transition between cases II and IV also occurs when ${dy
\over dt} =0$.  In
general equations (\ref{cond1}),(\ref{cond2}) give two (or more) lines
on the $Oxy$-plane, which separate the plane into several regions with
different directions of vectors (cases I-IV).

These lines are called null-clines.

\begin{D}
The  $x$-null-cline (or ${dx \over dt}=0$ null-cline) is the set of points satisfying the condition  $f(x,y)=0$. The  $y$-null-cline (or ${dy \over dt}=0$ null-cline) is the set of points satisfying the condition  $g(x,y)=0$.
\end{D}

To use the method of null-clines it is useful to note that at the
$x$-null-cline the $x$-component of the vector changes its sign,  and at
the $y$-null-cline the $y$-component of the vector changes its sign.
To use this rule effectively we will always denote the vector field
components using lines of different types: the horizontal component as a
dashed line and the vertical component as a solid line.  Let us use
these ideas and formulate a  plan for finding the vector field using
null-clines.

{\bf Plan of null-cline analysis for system (\ref{2dgeneral_3})}

We assume that on the $Oxy$-plane the $x$-axis is the horizontal axis  and the $y$-axis is the vertical axis.

\ben
\item Draw  ${dx \over dt}=0$ null-clines from the equation $f(x,y)=0$ using  a dashed line and   ${dy \over dt}=0$ null-clines from  the equation $g(x,y)=0$ using  a solid line.

\item  Choose a point in one of the regions in the $x,y$ plane and find
  the $x$ and the $y$ -components of the vector field. Use the dashed
  line for the $x$ component and the solid line for the $y$ component.
  For finding the directions use the following rule: if $f(x,y) >0$ the
  $x$ component is directed as '$\rightarrow$', if $f(x,y) <0$ it is
  directed as '$\leftarrow$'; if $g(x,y) >0$ the $y$-component is directed
  as $\uparrow$, $g(x,y) <0$ it is directed as $\downarrow$.

\item Find the vector field in the adjacent regions  using the following rule:
\ben
\item change the direction of the dashed component of the vector
field 
 if in order to get to
the  adjacent region you cross the dashed null-cline

\item change the direction of the solid component of the vector
field 
 if in order to get to
the  adjacent region you cross the solid null-cline

\item show the direction of the vector field on the null-clines.

\een

 \een

Note, that instead of dashed and solid lines you can use  lines of different colors.
Then the last step of this plan would  be: change  the direction of the component of the same color as the color of
the null-cline which we cross to get to the adjacent region.

Note, that although this plan works good in most of cases,  there are some situations 
when components of the vector field do not change their  sign at the corresponding 
null-cline. These are special so-called degenerate cases (exceptions). We will not consider them
in our course.

{\bf {Example}}
Find the vector field of the following system  using null-clines.
\beq
\label{PPour_2}
\left\{
\begin{array}{l}
{dx \over dt}=3x(1-x)-1.5xy\\ 
 {dy \over dt}=0.5xy-0.25y
\end{array}
\right.
\eeq 

\begin{figure}[hhh]
\centerline{
\psfig{type=pdf,ext=.pdf,read=.pdf,figure=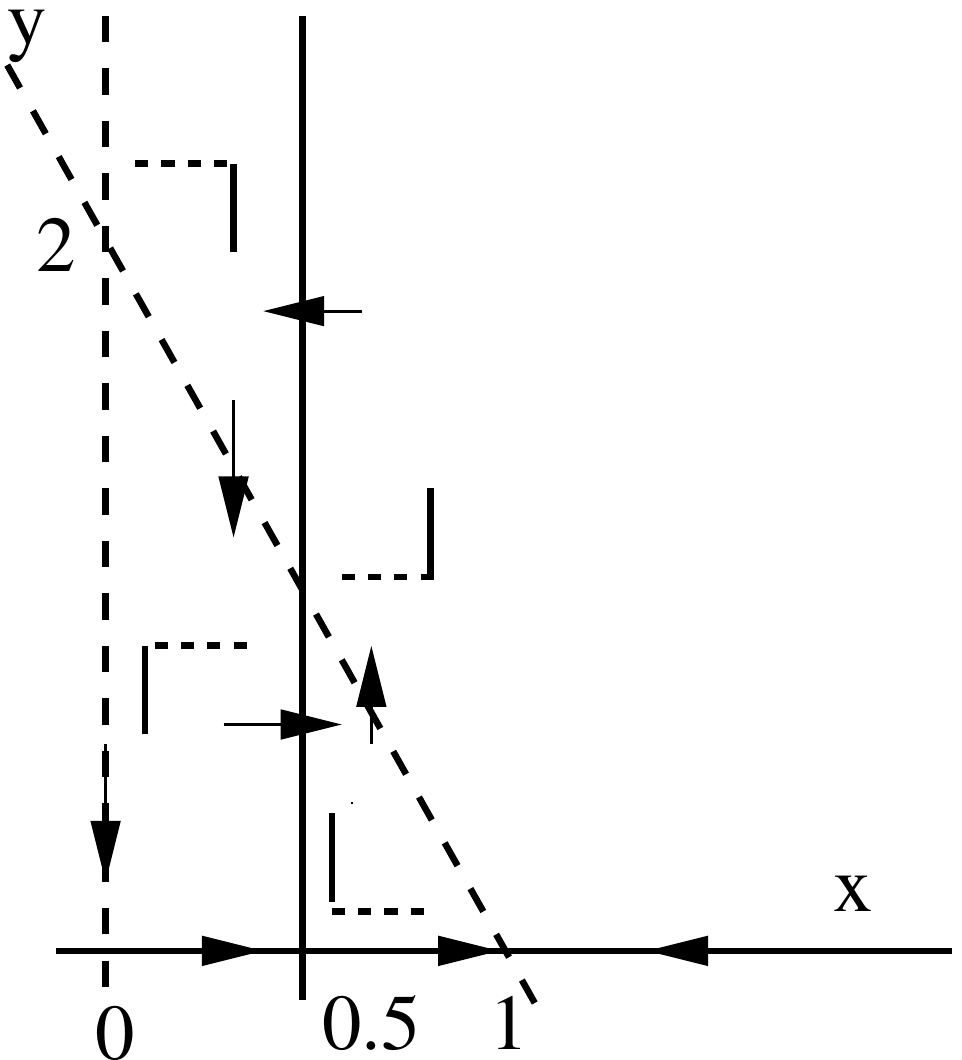,width=6.cm}
}
\caption{\label{2d7_fig}}
\end{figure}

{\bf {Solution}}.
We follow our plan as follows
\ben
\item  The ${dx \over dt}=0$ null-clines are given by  the equation $f(x,y)=0$, i.e. $3x(1-x)-1.5xy=x(3-3x-1.5y)=0$. This equation  has two solutions: $x=0$ (the vertical line which coincides with the $y$-axis) and $y=2-2x$  which is a straight line with the negative slope $-2$ which goes through the point $(2,0)$.  The graphs  are  shown using  dashed lines
in fig.\ref{2d7_fig}.  The ${dy \over dt}=0$ null-clines are given by the equation
$g(x,y)=0$, i.e. $0.5xy-0.25y=y(0.5x-0.25)=0$, which also has two
solutions: $y=0$ (horizontal  line which coincides with the $x$-axis) and $x=0.5$ (vertical line through the point $x=0.5$). Graphs are  shown by  solid lines in
fig.\ref{2d7_fig}.

\item We find that at point $(2,2)$ ${dx \over dt}=3*2*(1-2)-1.5*2*2=-12<0$ and ${dy \over dt}=0.5*2*2-0.25*2=1>0$, hence the direction of the dashed arrow is to the left '$\leftarrow$' and of the solid arrow is
upward $\uparrow$.

\item Now we complete the picture. For example,   to get into the region to the left from point ($2,2$) 
we cross the solid line, thus we change the direction of the solid  component here.
 Similarly for the other regions. We get the picture as in
 fig.\ref{2d7_fig}. Finally we show the vector field on the
 null-clines.

\een

We see that the vector field in fig.\ref{2d7_fig} is a good
approximation of the flow in fig.\ref{2d4_fig}.  We also see that
attractor ($0.5,1$) is a special point in fig.\ref{2d7_fig}: the point
of intersection of the $x$ and $y$ null-clines. This is not a
coincidence. If we compare the conditions for finding the equilibrium
point (\ref{equilibrium}) and equations for null clines (\ref{cond1}),
(\ref{cond2}), we see that the first equation for finding equilibria
$f(x,y)=0$ is also the equation for $x$ null-cline and the second
equation for finding the equilibrium point $g(x,y)=0$, is the equation
for the $y$-null-cline. Thus, the  solution of system (\ref{equilibrium}),
which gives the points satisfying both equations, gives the points
which belong to both null-clines, i.e. the points of intersection of
the null-clines. So we found that:
\bec
Equilibria are the points of intersection of the $x$ and $y$-null-clines.
\eec

Note that this definition applies points of intersection of {\it
different} null-clines only. For example intersection of null clines
of the same type in fig.\ref{2d7_fig} at points (0.5,0) and (0,2) do
not give equilibria of system (\ref{PPour}), while intersection of the different null-clines at points (0,0),(1,0), (0.5,1) give the equilibria of this system. 

Our next step will be to study how can we apply the null-clines for
finding the types of equilibrium.

\section{Graphical Jacobian \label{gjac_sec}}

\begin{figure}[hhh]
\centerline{
\psfig{type=pdf,ext=.pdf,read=.pdf,figure=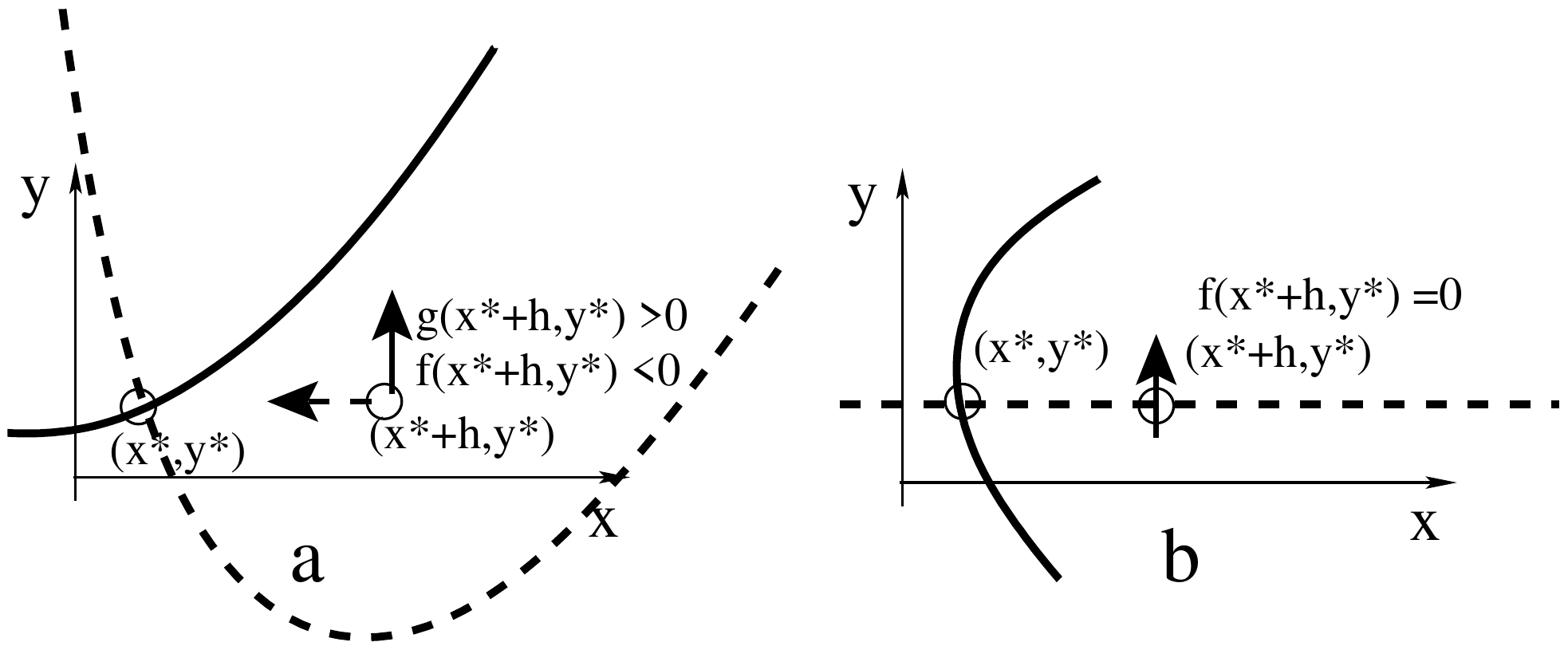,width=10cm}
\psfig{type=pdf,ext=.pdf,read=.pdf,figure=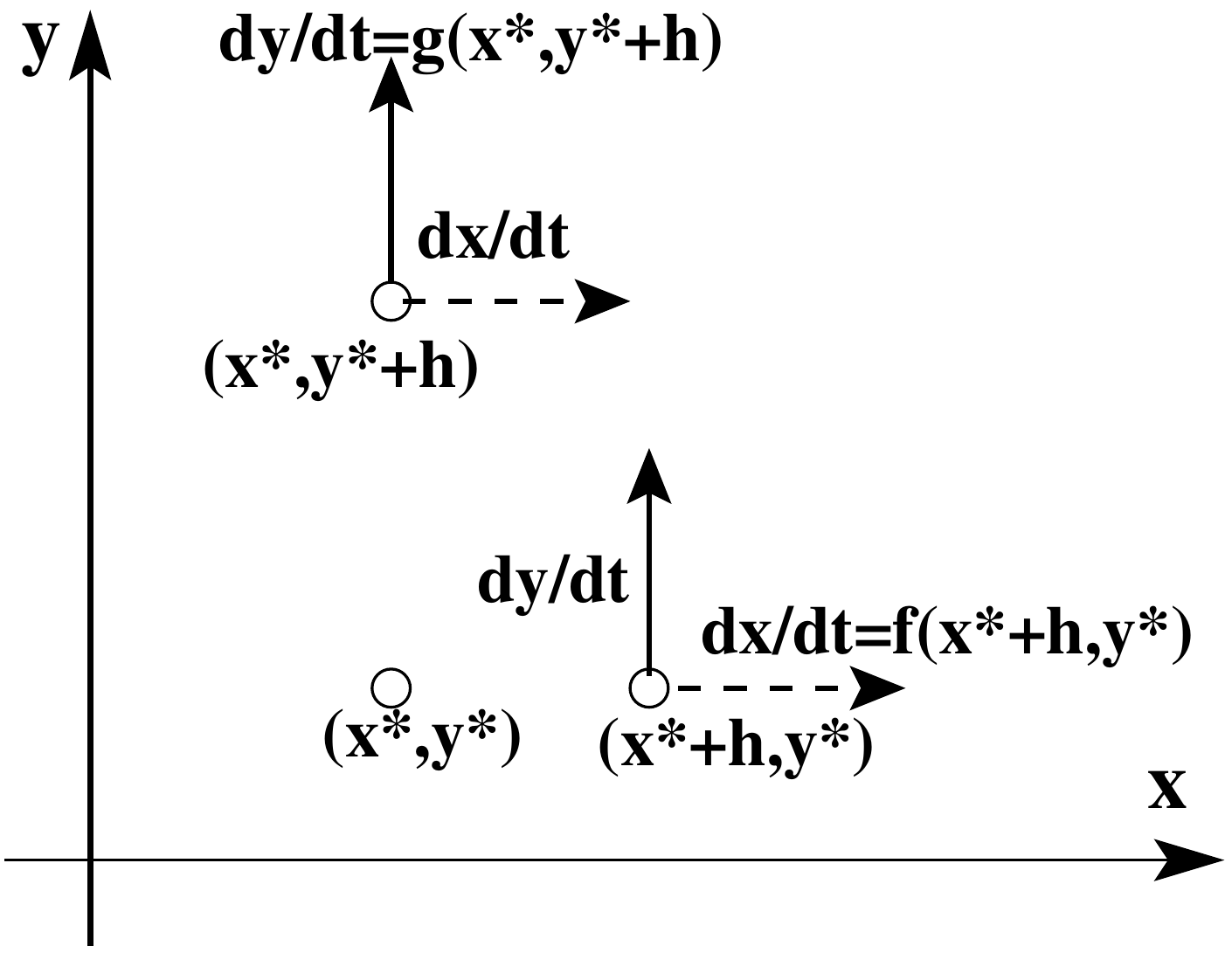,width=5cm}
}
\caption{\label{graph1_fig}}
\end{figure}

The method which we present here allows, in a number of cases,
to find the type of an equilibrium from the null-clines of the system,
i.e., even without computation of partial derivatives of the Jacobian
in the equilibrium. The main idea behind this method  can be seen from
the scheme in Fig.\ref{graph1_fig}. Let us consider an equilibrium
point $(x^*,y^*)$ and two close points: one located to the right, with
the coordinates $(x^*+h,y^*)$, the other upward with the coordinates
$(x^*,y^*+h)$. Because we assume that $(x^*,y^*)$ is an equilibrium
point, $f(x^*,y^*)=$ and  $g(x^*,y^*)=0$ (see
(\ref{equilibrium})). We can approximate the partial derivative
${\partial f\over \partial x}$ at $(x^*,y^*)$ using a formula similar
to (\ref{deriv}) as : $ {\partial f\over \partial x} \approx
{f(x^*+h,y^*)-f(x^*,y^*) \over h}$, but because $f(x^*,y^*)=0$, we get
$ {\partial f\over \partial x}
\approx {f(x^*+h,y^*) \over h}$. If we apply the same approach for all
derivatives constituting the Jacobian at the equilibrium point
$(x^*,y^*)$ we get:
\beq
\label{gjac}
J=\left(\begin{array}{lr}  {\partial f \over  \partial x} \approx  {f(x^*+h,y^*) \over h} & { \partial f  \over  \partial y }  \approx  {f(x^*,y^*+h) \over h}\\
                         { \partial g   \over  \partial x} \approx  {g(x^*+h,y^*) \over h} &  {\partial g  \over  \partial y }  \approx  {g(x^*,y^*+h) \over h}
\end{array} \right)
\eeq
This approximation will be better if points $(x^*+h,y^*)$ and
$(x^*,y^*+h)$ are closer to the equilibrium point, i.e., if $h$ is
small. It turns out that in many cases the exact values of the
derivative will not be important for us and we will be able to find
the equilibrium type from just the {\it sign} of the components of the
Jacobian. From (\ref{gjac}) it is clear that the sign of the Jacobian
components is the same as the sign of the functions at the appropriate
points, e.g.  the sign of ${\partial f \over \partial x}$ is the same
as the sign of $f(x^*+h,y^*)$, etc. Let us now recall, that the sign
of the functions $f(x,y),g(x,y)$ is represented on  the vector field  of
our system.  In fact, '$\rightarrow$' means that ${dx \over dt}>0$ and
occurs at the points where $f(x,y)>0$, the '$\uparrow$' means that
${dy \over dt}>0$ and occurs at the points where $g(x,y)>0$ (see
fig.\ref{2d6_fig}), etc. Thus from the vector field  we can easily determine
the sign of the components of the Jacobian matrix. For example, the
negative direction of the $x$-component in fig.\ref{graph1_fig}a will
mean that $f(x^*+h,y^*)<0$ and hence  ${\partial f \over \partial x}<0$,
the positive direction of the $y$-component in fig.\ref{graph1_fig}a
will mean that $g(x^*+h,y^*)>0$ and thus ${\partial g \over \partial
x}>0$. Similarly, in Fig.\ref{graph1_fig}b the vector field at point $(x^*+h,y^*)$ is vertical, thus ${\partial g \over
\partial x}> 0$ and ${\partial f \over
\partial x}=0$.

We will formulate this result  as the following conclusion:
\bec
The sign of the $x$ and $y$ vector field components to the
right from the equilibrium point give the sign of  ${\partial f
\over \partial x}$ and ${\partial g \over \partial x}$ components of the Jacobian. The
sign of the $x$ and $y$ vector field components upward from the
equilibrium point give the sign of ${\partial f \over \partial y}$ and
${\partial g \over \partial y}$ components of the Jacobian matrix
$J=\left(\begin{array}{lr} {\partial f \over \partial x} & {
\partial f \over \partial y }\\ { \partial g \over \partial x} &
{\partial g \over \partial y }
\end{array} \right)$ at this equilibrium point (fig.\ref{graph1_fig}c).
\eec

Three notes on the application of this rule:
\begin{itemize}
\item The testing points
must be {\it exactly horizontal} $(x^*+h,y^*)$  and  {\it exactly vertical }
$(x^*,y^*+h)$  with respect to the equilibrium point.

\item Testing points should be as close as possible to the equilibrium and should never cross a null cline when going from the equilibrium. (Putting the testing point to the right  from the dashed line  in fig.\ref{graph1_fig}a will be wrong).

\item If a null cline is exactly horizontal or exactly vertical then one of the  corresponding derivatives will be zero. ( In fig.\ref{graph1_fig}b,${\partial f \over \partial x}=0;{\partial g \over \partial x}>0$ ).

\begin{figure}[hhh]
\centerline{
\psfig{type=pdf,ext=.pdf,read=.pdf,figure=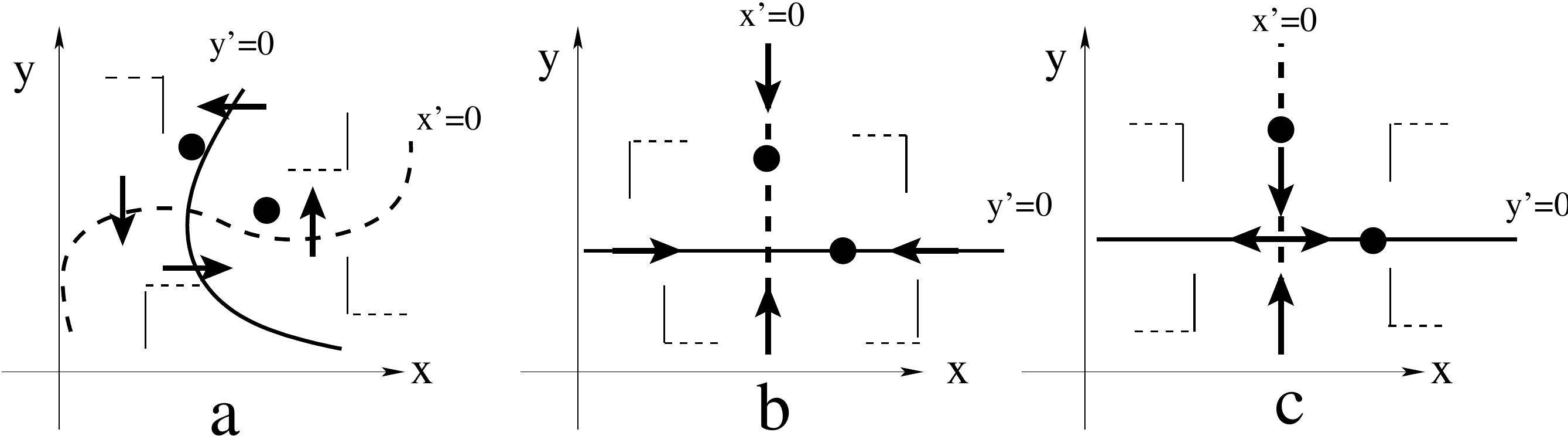,width=15cm}
}
\caption{\label{graph3_fig}}
\end{figure}

\end{itemize}

\begin{figure}[hhh]
\centerline{
\psfig{type=pdf,ext=.pdf,read=.pdf,figure=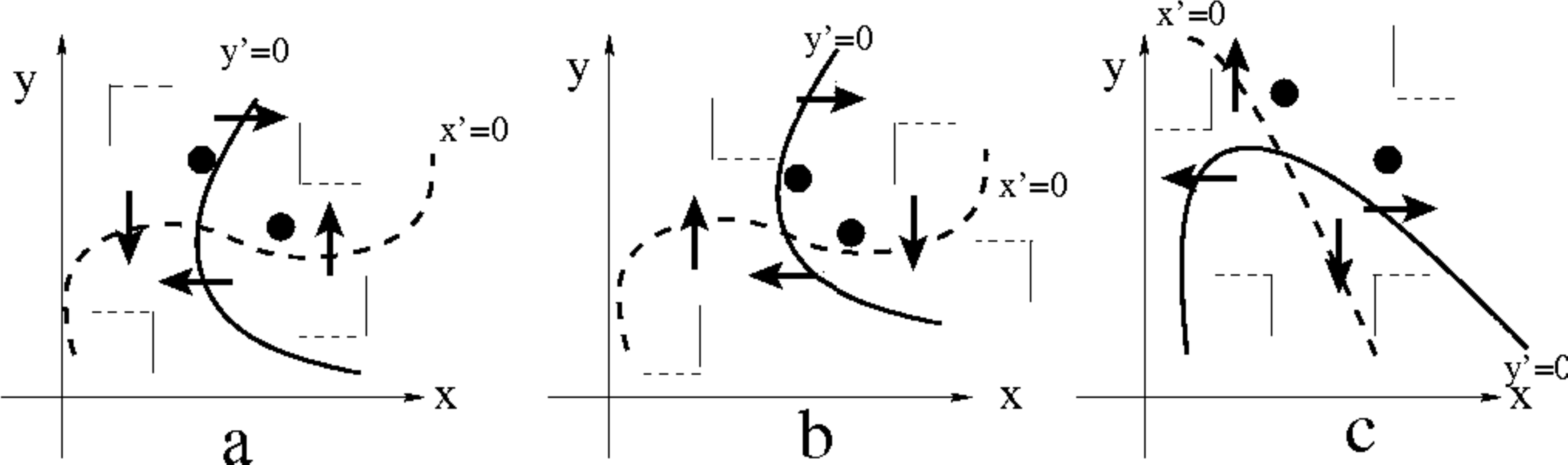,width=15cm}
}
\caption{\label{graph4_fig}}
\end{figure}

{\bf Examples.} Let us consider several examples of the application of
this graphical approach for null-clines presented in
Fig.\ref{graph3_fig} and Fig.\ref{graph4_fig}. On all these figures the
testing points are marked by filled circles.

From fig.\ref{graph3_fig}a we find  the following components of the Jacobian:
$J=\left(\begin{array}{lr}  -\alpha  & -\beta\\
                          +\gamma &  -\delta
\end{array} \right)$, where $\alpha,\beta,\gamma,\delta$ stand for  unknown positive numbers. Thus we see that $detJ=\alpha\delta+\beta\gamma>0$, and $trJ=-\alpha -\delta<0$, thus (see fig.\ref{fig8.2}) we have a stable equilibrium (stable node or stable spiral). On the basis of these  data we cannot say whether this equilibrium is a node or a spiral as we cannot compute whether  the discriminant of the Jacobian is positive or negative, as this depends on the exact values  of the partial derivatives. 

From fig.\ref{graph3_fig}b we find:
$J=\left(\begin{array}{lr}  -\alpha  & 0\\
                          0 &  -\delta
\end{array} \right)$. Thus  $detJ=\alpha\delta>0$, $trJ=-\alpha -\delta<0$ and $D=tr^2-4Det=(\alpha +\delta)^2-4\alpha\delta=\alpha^2+2\alpha\delta+\delta^2-4\alpha\delta =\alpha^2-2\alpha\delta+\delta^2=(\alpha -\delta)^2>0$, thus  we have a stable node.

From fig.\ref{graph3_fig}c we find:
$J=\left(\begin{array}{lr}  +\alpha  & 0\\
                          0 &  -\delta
\end{array} \right)$. Thus  $detJ=-\alpha\delta<0$, and we have a saddle.

From fig.\ref{graph4_fig}a we find:
$J=\left(\begin{array}{lr}  +\alpha  & +\beta\\
                          +\gamma &  -\delta
\end{array} \right)$, hence  $detJ=-\alpha\delta-\gamma\beta<0$, and 
we again have  a saddle.

From fig.\ref{graph4_fig}b we find:
$J=\left(\begin{array}{lr}  +\alpha  & +\beta\\
                          -\gamma &  -\delta
\end{array} \right)$.  $detJ=-\alpha\delta+\gamma\beta$, $trJ=\alpha-\delta$. Because we do not know the values of the coefficients $\alpha,\beta,\gamma,\delta$, just their signs, we do not know if the $detJ$ and $trJ$ are  positive or negative, thus we cannot determine the equilibrium type in this case using the graphical Jacobian.

Finally from \ref{graph4_fig}c we find:
$J=\left(\begin{array}{lr}  +\alpha  & +\beta\\
                          +\gamma &  +\delta
\end{array} \right)$. Thus  $detJ=\alpha\delta-\gamma\beta$,  $trJ=\alpha+\delta>0$. Again we cannot determine the sign of the $detJ$, but positive $trJ$ implies that the equilibrium will be unstable. We do not know its type: non-stable node, non-stable  spiral and a  saddle  are all possible.

We see that the method of graphical Jacobian is a useful tool for
finding the  equilibrium type from the vector field, but sometimes it is
not sufficient and we need to know not only the sign but also the values of the
 Jacobian coefficients.

\section{Exercises}

\ben

\item Sketch the phase portraits of the following systems using null-clines.
Try to find the type of equilibrium using the graphical Jacobian (section
\ref{gjac_sec}). If you are unable to find the equilibrium type using
the graphical Jacobian, find it using the method of sec.\ref{dettr_sec}.
Sketch a qualitative phase portrait (without computation of eigen
values and eigen vectors).
\ben

\item  $\left\{
\begin{array}{l}
{dx \over dt} = 3x+y\\ 
{dy \over dt}=-x-y
\end{array}
\right.$

\item  $\left\{
\begin{array}{l}
{dx \over dt} = x+2y \\ 
{dy \over dt}=-2x-2y
\end{array}
\right.$

\item   $\left\{
\begin{array}{l}
{dx \over dt} = y\\ 
{dy \over dt}=3x
\end{array}
\right.$

\item    $\left\{
\begin{array}{l}
{dx \over dt} = x-4y\\ 
{dy \over dt}=x+y
\end{array}
\right.$
\een

\item The figure below  shows the null-clines of three   systems of two non-linear differential equations  
\begin{figure}[H]
\centerline{
\psfig{type=pdf,ext=.pdf,read=.pdf,figure=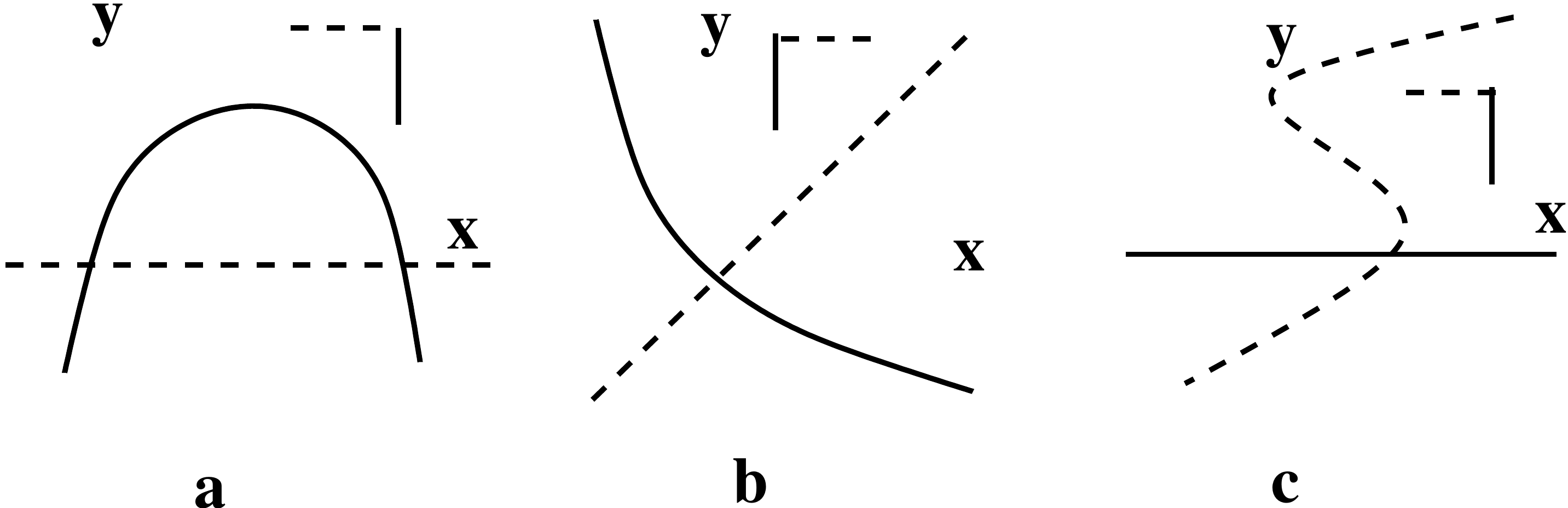,width=15.cm}
}
\end{figure}
 
\begin{enumerate}
 
\item Complete the figures by showing the direction of the vector field
in all regions on the plane and on the
null-clines.

\item Mark the equilibria  and find their types using the graphical Jacobian (if possible).

\end{enumerate}
\item Algae model using graphical Jacobian (see also problem \ref{ProbAlgae} from section \ref{exCh5}):
$$\left\{
\begin{array}{l}
{dx \over dt} =2x(1-y) \;\;x \geq 0;\\ 
{dy \over dt}=2-y-x^2 \;\; y \geq 0.
\end{array}
\right.$$ 
\ben
\item  Sketch the vector field for  the Algae  system using null-clines.

\item  Find equilibria.

\item Find the type of each equilibrium using the graphical Jacobian and draw qualitative local phase portraits.

\een

\item Lotka Volterra model using graphical Jacobian (see also problem \ref{ProbLV} from section \ref{exCh5}):
$$
\left\{
\begin{array}{l}
dN/dt=aN-bNP  \quad N \geq 0\\
  dP/dt=cNP-dP  \quad P \geq 0
\end{array}
\right.\;\;\;\;\;\;\;a,b,c,d,e >0
$$
\ben
\item  Sketch the vector field  of the  system using null-clines.

\item  Find equilibria.

\item Find the type of each equilibrium using the graphical Jacobian  and draw qualitative local phase portraits.

\item Does  the vector field  change  if we change the values of parameters $a,b,c,d$?
\een

\subsection*{Additional exercises}
\item Find the vector fields of these  systems using null-clines. Find equilibria. Determine  the type of each equilibrium using the graphical Jacobian  and draw qualitative local phase portraits.

\ben
\item $\left\{
\begin{array}{l}
{dx \over dt} =x-3y\\ 
{dy \over dt}=x+y
\end{array}
\right.$

\item  $\left\{
\begin{array}{l}
{dx \over dt} =y\\ 
{dy \over dt}=-x
\end{array}
\right.$

\item   $\left\{
\begin{array}{l}
{dx \over dt} =9x + y^2\\ 
{dy \over dt}=x-y
\end{array}
\right.$

\item  $\left\{
\begin{array}{l}
{dx \over dt} =y^2-x^2\\ 
{dy \over dt}=y-1
\end{array}
\right.$

\een

\een

\chapter{Plan of  qualitative analysis and examples}
\section{Plan}
Let us now formulate a plan to qualitatively  study  a  system of two differential equations with two variables and consider several examples.

We study the system:

\beq
\left\{
\begin{array}{l}
 {dx \over dt}=f(x,y) \\ {dy \over dt}=g(x,y)
\end{array}
\right.
\eeq
Our main aim is to plot the phase portrait of this system and then predict its dynamics. Based on methods which we have developed we will do it in two steps:
\ben 
\item [I]{\bf Null-cline analysis}

\item [II]{\bf Jacobian analysis}

\een
Where the Jacobian analysis can be either performed using the determinant-trace method from section \ref{dettr_sec}, or using  the graphical Jacobian from sec.\ref{gjac_sec}. Note, that the determinant-trace method always give us a definitive answer, while the graphical Jacobian method sometimes fails. In more details:

{\bf Null-cline analysis}

We assume that on the $Oxy$-plane the $x$-axis is the horizontal axis  and the $y$-axis is the vertical axis.

\ben
\item Draw the ${dx \over dt}=0$ null-clines from the equation $f(x,y)=0$ using    dashed lines and  the ${dy \over dt}=0$ null-clines from  the equation $g(x,y)=0$ using   solid lines.

\item  Choose a point in one of the regions on the $x,y$ plane and find
 the vector field for the $x$-component.  Denote the $x$ component as a dashed  '$\rightarrow$' if $f(x,y) >0$ and as a dashed '$\leftarrow$' if $f(x,y) < 0$.

\item  Find
 the vector field for the $y$-component  at the same point. Denote the $y$ component as a 
solid  $\uparrow$ if $g(x,y) >0$ and  as a solid $\downarrow$ if $g(x,y) < 0$

\item Find the vector field in the adjacent regions using the following rule:
\ben
\item change the direction of the dashed component of the vector
field 
 if to get to
the  adjacent region you cross the dashed null-cline.

\item change the direction of the solid component of the vector
field 
 if to get to
the  adjacent region you cross the solid null-cline.

\item show the direction of the vector field on the null-clines 
close to the equilibrium

\een

\een

{\bf Jacobian analysis using the determinant-trace  method}
\ben
\item
Find equilibria  from  equations:
\beq
\left\{
\begin{array}{l}
 f(x,y)=0 \\ g(x,y)=0
\end{array}
\right.
\eeq
\item For each equilibrium $(x^*,y^*)$, find the Jacobian at that equilibrium point
\beq
 J={\left(\begin{array}{lr}  {\partial f /  \partial x} & { \partial f /  \partial y }\\
                         { \partial g  /  \partial x} &  {\partial g /  \partial y } 
\end{array} \right)}_{(x^*,y^*)}
\eeq
Note: Do not forget to substitute $x=x^*,y=y^*$ into the Jacobian.
\item Determine the type of each equilibrium $(x^*,y^*)$. For this compute
\beq
detJ= {\left( {\partial f \over   \partial x}  * {\partial g \over  \partial y}
-  {\partial f \over   \partial y}  * {\partial g \over  \partial x } \right)}_{(x^*,y^*)}
\eeq
\beq
trJ=  {\left( {\partial f \over   \partial x} +  {\partial g \over  \partial y}  \right)}_{(x^*,y^*)}
\eeq

\beq
D=(trJ)^2-4detJ
\eeq

To find the type of equilibrium use fig.\ref{fig8.2}, or the following list:
\ben
\item  If  $detJ<  0$; the point is a saddle point

\item  If  $detJ>  0 ,trJ>0,D \geq 0$;  the point is a non-stable node

\item  If  $detJ>  0 ,trJ<0,D \geq 0$;  the point is a stable node

\item  If  $detJ>  0 ,trJ>0,D<0$;  the point is a non-stable spiral

\item  If  $detJ>  0 ,trJ<0,D<0$;  the point is a stable spiral

\item  If  $detJ>  0 ,trJ=0$;  the point is a center
\een

\item Draw local phase portraits using both knowledge on the equilibrium type  and the vector fields obtained using null-cline analysis

\item Connect local phase portraits to get the global picture and show attractors and their basins of attraction.
\een

{\bf Jacobian analysis using the graphical Jacobian}

\ben
\item For each equilibrium point (point of intersection of {\it different} null-clines) choose two points, one of which is located to the right and the other upward from the equilibrium. Find the components of  the Jacobian using the following rule:
the sign of the $x$ component of the vector field to the right of
the equilibrium point gives the sign of ${\partial f
\over \partial x}$ and the sign of the  $y$ component of the vector field to the right of  the equilibrium point gives  the sign of  ${\partial g
\over \partial x}$.  The
sign of the $x$ and $y$ vector field components upward of  the
equilibrium point give the sign of  ${\partial f \over \partial y}$
and ${\partial g \over \partial y}$.
\item Note that if one of the components is zero, then the corresponding derivative of the Jacobian will be zero. The latter can  happen if at a given equilibrium point one or both  null-clines are exactly 
horizontal or exactly vertical.

\item Try to compute the sign of the determinant and of the trace of the Jacobian and try to identify the type of equilibrium from fig.\ref{fig8.2}.

\item If the type of equilibrium cannot be identified, use  the analytical  determinant-trace  method formulated above.

\item Draw local and then global phase portraits of the system and make predictions about its dynamics.
\een

\section{Examples}

{\bf {Example 1}} Study system (\ref{PPour}), using null-clines and the determinant-trace method. 

\beq
\label{PPour2}
\left\{
\begin{array}{l}
{dx \over dt}=3x(1-x)-1.5xy\\ 
 {dy \over dt}=0.5xy-0.25y
\end{array}
\right.
\eeq

{\bf {Solution}} We have made studies of different aspects of this model throughout the reader, so let us collect the information.

I.Null-cline analysis.

The null-clines of this system are given in fig.\ref{2d7_fig}. We repeat it here  in fig\ref{plan1_fig}a:

\begin{figure}[hhh]
\centerline{
\psfig{type=pdf,ext=.pdf,read=.pdf,figure=fig2d7,width=4.5cm}
\psfig{type=pdf,ext=.pdf,read=.pdf,figure=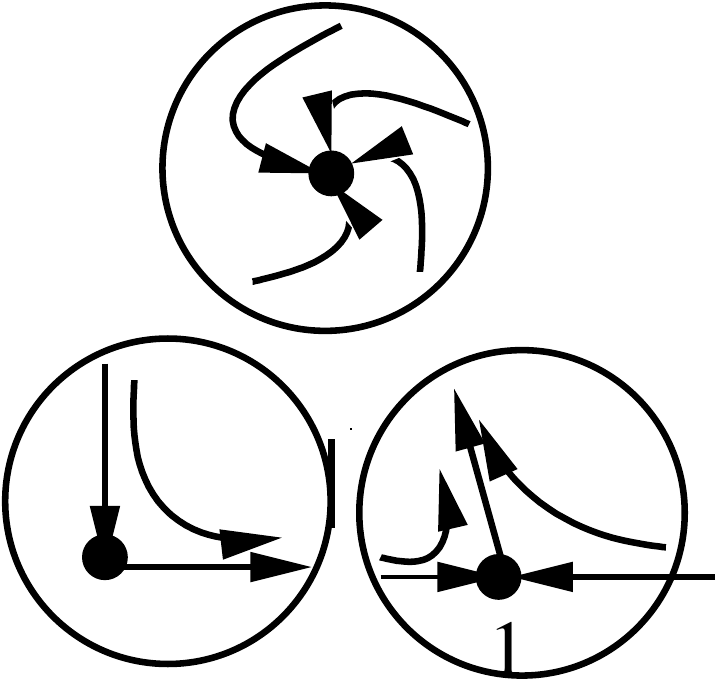,width=3.5cm}
\psfig{type=pdf,ext=.pdf,read=.pdf,figure=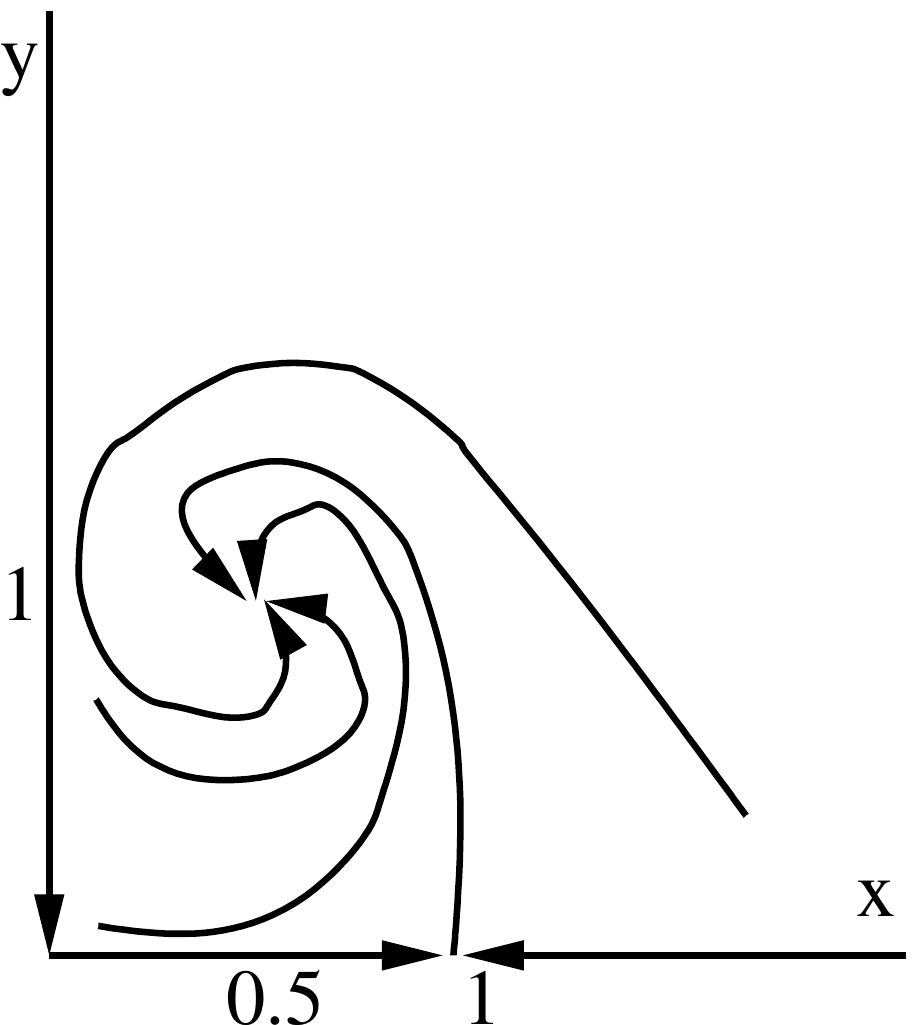,width=4.5cm}
}
\caption{\label{plan1_fig}}
\end{figure}

II.Jacobian analysis:
\ben
\item Equilibria. In example eq.(\ref{ex_eq}) 
we found that this system has three equilibria points: $(0,0),(1,0)$
and $(0.5,1)$. 

\item Jacobian. 
We compute the Jacobian as: $ {\partial f / \partial x}=3-6x-1.5y; \; {\partial f / \partial y}=-1.5x; \;
{\partial g / \partial x}=0.5y;\; {\partial g / \partial y}=0.5x-0.25$, thus:
$J=\left(\begin{array}{lr}   3-6x-1.5y & -1.5x\\
                         0.5y  &  0.5x-0.25
\end{array} \right)$. 

Let us find equilibria types from the Jacobian.

At  point $(0,0)$ the Jacobian is:
$J_1=\left(\begin{array}{lr}  3 & 0\\
                         0 &  -0.25
\end{array} \right)$, $det J_1=3*(-0.25)<0$, thus this is a saddle point.

At the point $(1,0)$ the Jacobian is:
$J_2=\left(\begin{array}{lr}  -3 & -1.5\\
                         0 &  0.25
\end{array} \right)$, $det J_2=(-3)*0.25<0$, thus this is a saddle point.

At the point $(0.5,1)$ the Jacobian is:
$J_3=\left(\begin{array}{lr}  -1.5 & -0.75\\
                         0.5 &  0
\end{array} \right)$, $det J_3=(-1.5)*0-0.5*(-0.75)<0.375$, thus we need to find $tr J_3=-1.5$ and $D=(-1.5)^2-4*0.375=0.75>0$ thus  this is a stable node.

\item Local phase portraits are presented in fig.\ref{plan1_fig}b. From the null-cline
analysis we find  the approximate locations of the manifolds of the  saddle points and the direction of the trajectories around the stable node.

\item Global picture. There are no general rules to draw the global picture.
However in the case of fig.\ref{plan1_fig}b we can expect that the  non-stable
manifold of the  saddle point at $(0,0)$ will end  at the other saddle point $(1,0)$, the non-stable manifold of the saddle point $(1,0)$ as well as most of  other trajectories should go to the stable node at $(0.5,1)$ as this is the only attractor here. Thus we get the phase portrait presented in fig.\ref{plan1_fig}c.
We see that  we will have  a stable global attractor with the basin of attraction the whole region $x>0,y>0$ (except two axes $x=0;  \; y=0$). 
\een
We see that the phase portrait which we got as a result of our study is qualitatively the same as the phase portrait of this system obtained using a computer (see fig.\ref{2d2_fig}b).

Let us try to apply for this problem the method of the graphical Jacobian

{\bf {Example 2}} Study the same system (\ref{PPour2}), using the graphical Jacobian method.

{\bf {Solution}}
From fig.\ref{plan1_fig}a we find  in the following components of the Jacobian:

\ben
\item Point $(0,0)$,
$J=\left(\begin{array}{lr}  +\alpha  & 0\\
                          0 &  -\delta
\end{array} \right)$. Thus  $detJ=-\alpha\delta<0$, thus we have a saddle. 

\item Point $(1,0)$,
$J=\left(\begin{array}{lr}  -\alpha  & -\beta\\
                          0 &  +\delta
\end{array} \right)$. Thus  $detJ=-\alpha\delta<0$, thus we have a saddle. 

\item Point $(0.5,1)$,
$J=\left(\begin{array}{lr}  -\alpha  & -\beta\\
                          +\gamma &  0
\end{array} \right)$. Thus  $detJ=\gamma \beta>0$, $trJ=-\alpha<0$, thus we have a stable equilibrium (stable node or stable spiral).

\item Drawing of the local and global phase portraits is exactly the same as with the previous method 
\een

Thus we see that in this case we were able to solve the  problem completely using the graphical Jacobian method. The only difference with the determinant-trace method is that we were not able to say that the stable equilibrium is a node. However, this has almost no consequences for the phase portrait.

\section{Exercises}

\begin{figure}[hhh]
\centerline{
\psfig{type=pdf,ext=.pdf,read=.pdf,figure=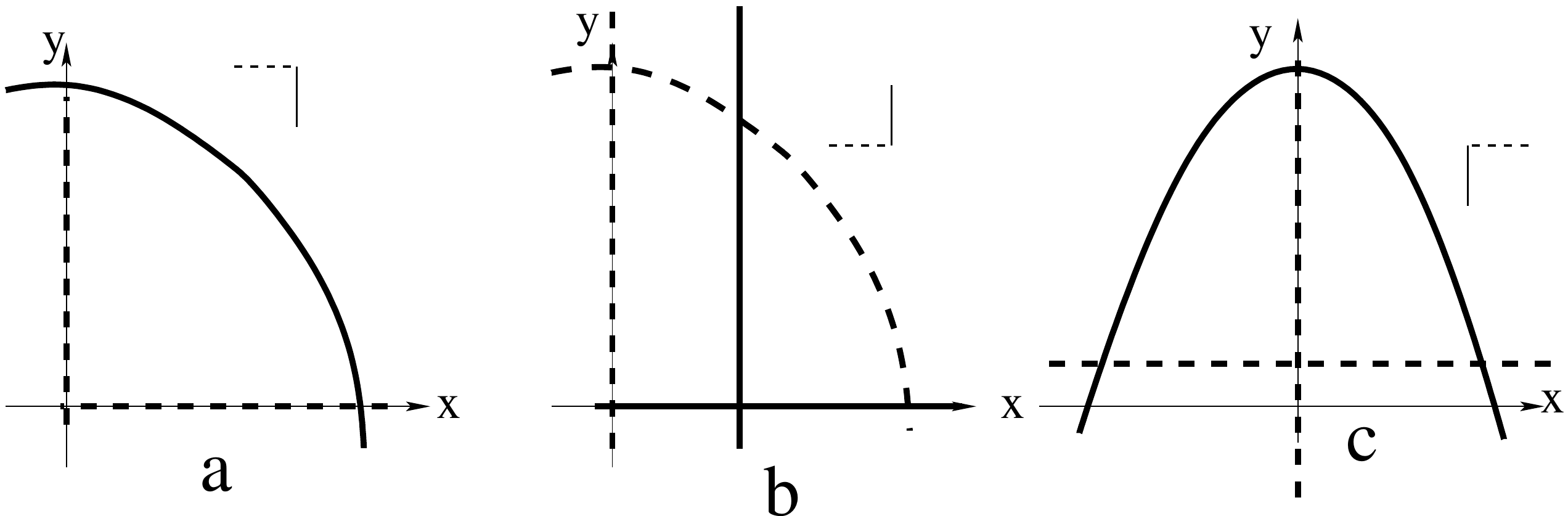,width=12cm}
}
\caption{\label{fig9.6}}
\end{figure}

\ben 
\item Complete the vector field approximations for the null-clines
shown in fig\ref{fig9.6}a,b,c.  Mark equilibria and determine their type
using the graphical Jacobian method. Draw local and then global phase
portraits. Try to  describe the dynamics of the system, by saying what happens with the variables $x$ and $y$ in the course of time.

\item  Draw the phase portrait of the following systems of differential equations. Explain their dynamics.
\ben
\item    $\left\{
\begin{array}{l}
{dx \over dt} = x(1-x-y)\\ 
{dy \over dt}=y(1-2x)
\end{array}
\right.  x \geq 0; \; y \geq 0 $

\item    $\left\{
\begin{array}{l}
{dN \over dt} = 2N-NP-N^2\\ 
{dP \over dt}=3P-2NP-P^2
\end{array}
\right.  N \geq 0; \; P \geq 0 $
\een

\item  In theoretical biology the following model has been used to study mRNA ($M$) protein ($P$) interaction:
\beq
\left\{
\begin{array}{l}
dP /dt=bM-d_PP\;\;\;\;\;\;\;a,b,d >0\\
dM/dt= a {K^2 b\over K^2+P^2} -d_M M   \qquad P>0; M>0
\end{array}
\right.
\eeq
Study this system using the graphical Jacobian approach, i.e. draw
null-clines, mark equilibria as points of intersection of
the null-clines, determine stability of these equilibria and sketch a
qualitative phase portrait.

\item  The following equations describe the dynamics
of  predator ($P$) and prey ($N$) populations:
 
$$
\left\{
\begin{array}{l}
dN/dt=rN (1-{ \displaystyle N \over \displaystyle K}) -bNP     \quad       \quad  N \geq 0  \quad  P \geq 0   \\
\\
dP /dt =bNP-2bP   \quad   \quad   \quad     \quad        b \geq 0 \; \; K \geq 0 \; \; r \geq 0 \; \;
\end{array}
\right.
$$
here $K,b,r$  are  parameters
 
This system of differential equations always has an equilibrium corresponding to an 
extinct  population of the predator and non-zero population of prey ($P=0;N \ne
0$).
 
\begin{enumerate}
 
\item
Find this ($P=0;N \ne 0$) equilibrium.

\item Find the Jacobian at this equilibrium point.
 
\item
For which values of the parameters the predator population can be driven to extinction
(i.e. to that equilibrium ).
 
\end{enumerate}
(N.B. Do not use the 'graphical Jacobian approach for this problem!).

\subsection*{Additional exercises}

\item  Draw the phase portrait of the following systems of differential equations. Explain their dynamics.
\ben
\item    $\left\{
\begin{array}{l}
{dx \over dt} = 2y\\ 
{dy \over dt}=x-x^2-0.5y
\end{array}
\right.$

\item    $\left\{
\begin{array}{l}
{dx \over dt} = x + y^2\\ 
{dy \over dt}=x+y
\end{array}
\right.$

\item 
 $\left\{
\begin{array}{l}
{dx \over dt} = x(25-x^2-y^2)\\ 
{dy \over dt}=y(x-3)
\end{array}
\right.  x \geq 0; \; y \geq 0 $

\item  $\left\{
\begin{array}{l}
{dx \over dt} = xy\\ 
{dy \over dt}=4-y-x^2 
\end{array}
\right.   \quad  - \infty < x  < \infty;  \quad - \infty < y <   \infty $

\een

\item  Lotka Volterra model with competition in the prey population:
$$
\left\{
\begin{array}{l}
dN/dt=aN-eN^2-bNP\\ 
  dP/dt=cNP-dP
\end{array}
\right.
$$
for $a=3,b=1.5,c=0.5,d=0.25,e=3$.

\item   A swimming pool is infested with algae whose population is
$N(t)$. The owner attempts to control the infestation with an algicidal
chemical, poured into the pool at a constant rate. In the absence of
algae, the chemical decays naturally; when algae are present it is metabolized by them and kills them. The equations of the rates of change of $N(t)$ and the concentration of the chemical in the pool, $C(t)$, are
$$
\left\{
\begin{array}{l}
dN/dt=aN-bNC\\ 
  dC/dt=Q-\alpha C-\beta NC
\end{array}
\right.
$$
where $a,b,Q,\alpha,\beta$ are positive constants. Discuss the meaning of each term in these equations.
\ben
\item Put $a=1,b=1,\alpha=1,\beta=1$ and show, that the system has two non-negative equilibria if $Q>1$, find them,
draw the phase portrait of this system and explain the dynamics.

\item What happens if $Q<1$?

\item* Find the conditions for control of the infestation for arbitrary positive values of  $a,b,Q,\alpha,\beta$.
\een

\een

\chapter{ Limit cycle}

\section{Stable and non-stable limit cycles}

\begin{figure}[h]
\centerline{
\psfig{type=pdf,ext=.pdf,read=.pdf,figure=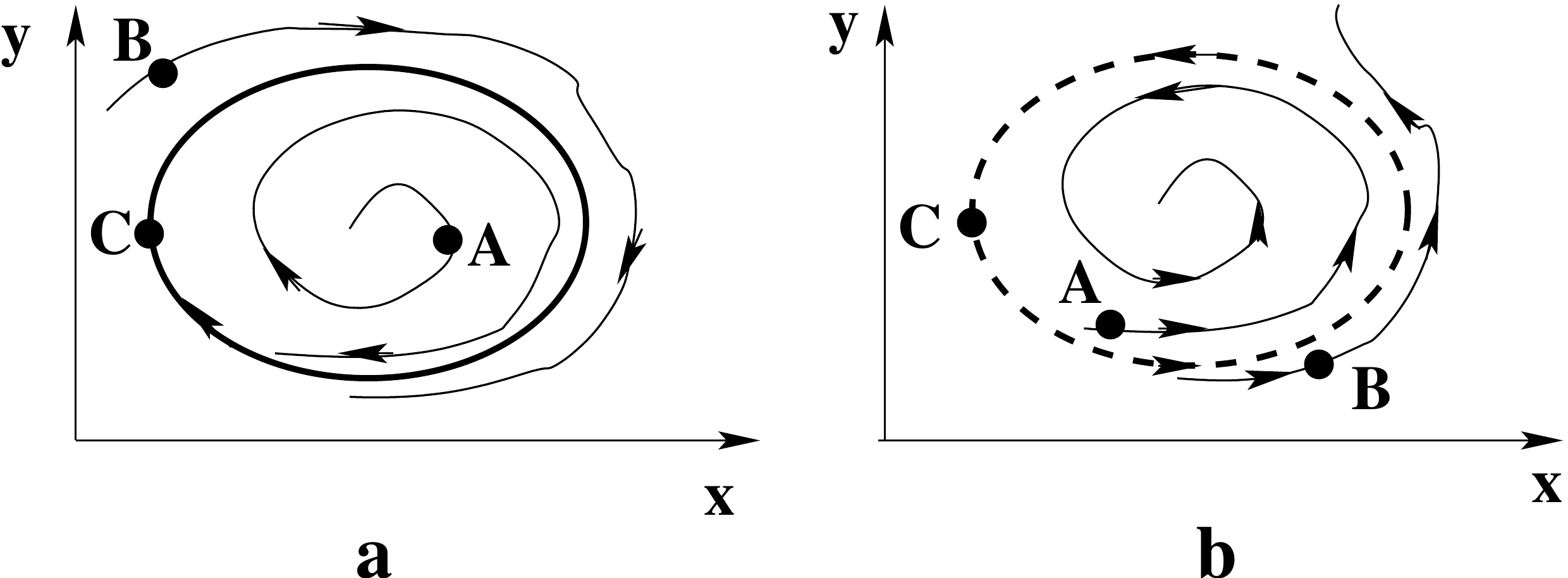,width=12cm}
}
\caption{ Phase portrait of  a system of two differential equations with a limit cycle. \label{flc1} }
\centerline{
\psfig{type=pdf,ext=.pdf,read=.pdf,figure=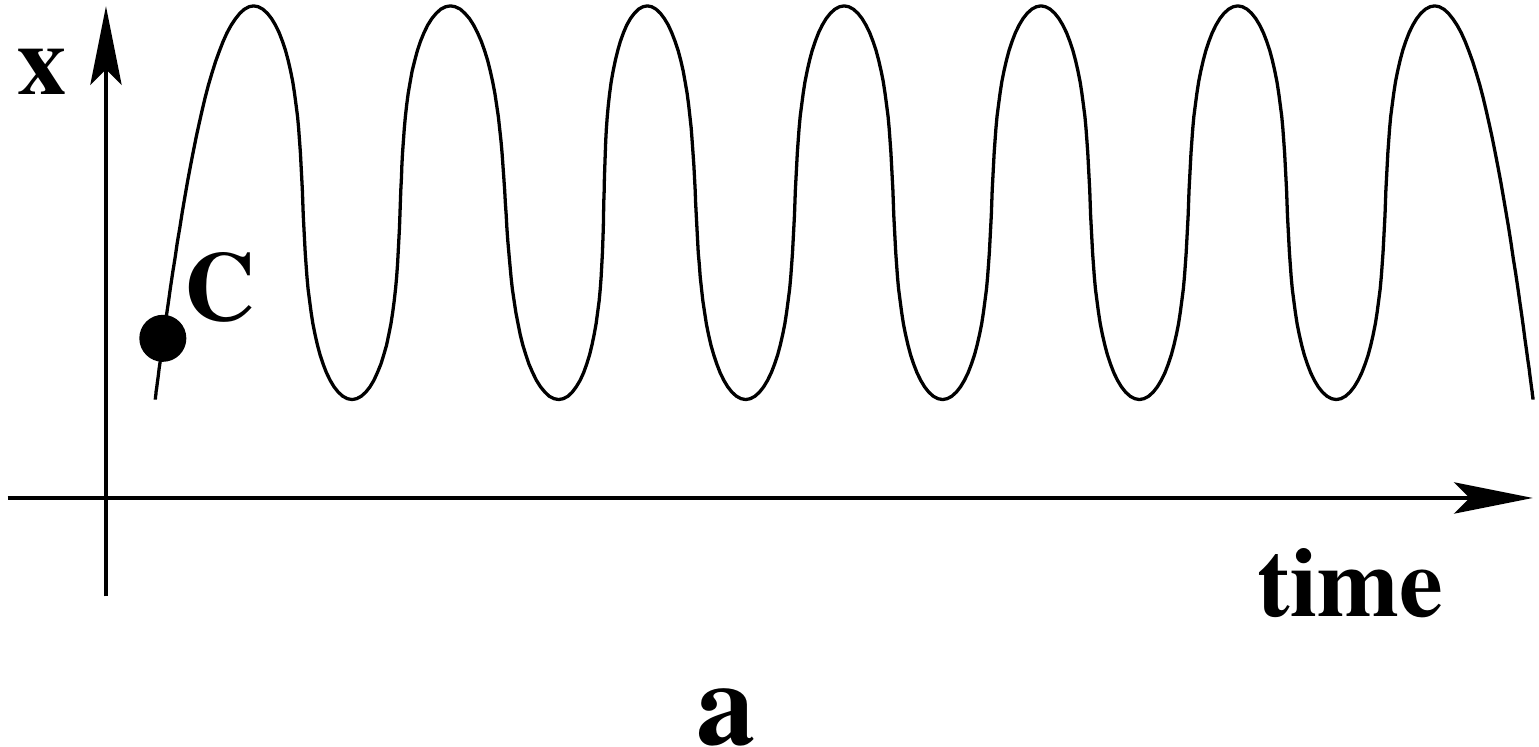,width=5cm}
\psfig{type=pdf,ext=.pdf,read=.pdf,figure=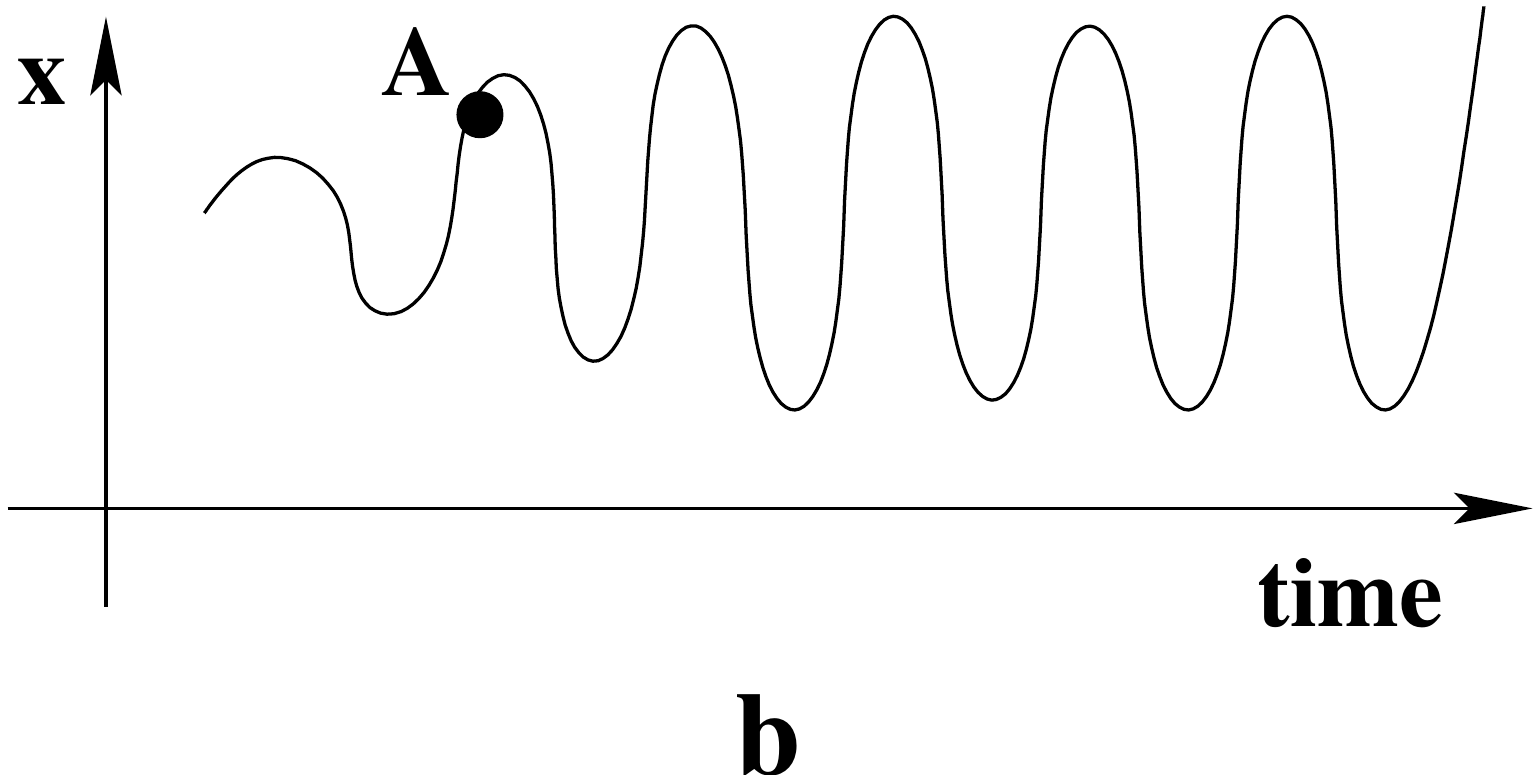,width=5cm}
\psfig{type=pdf,ext=.pdf,read=.pdf,figure=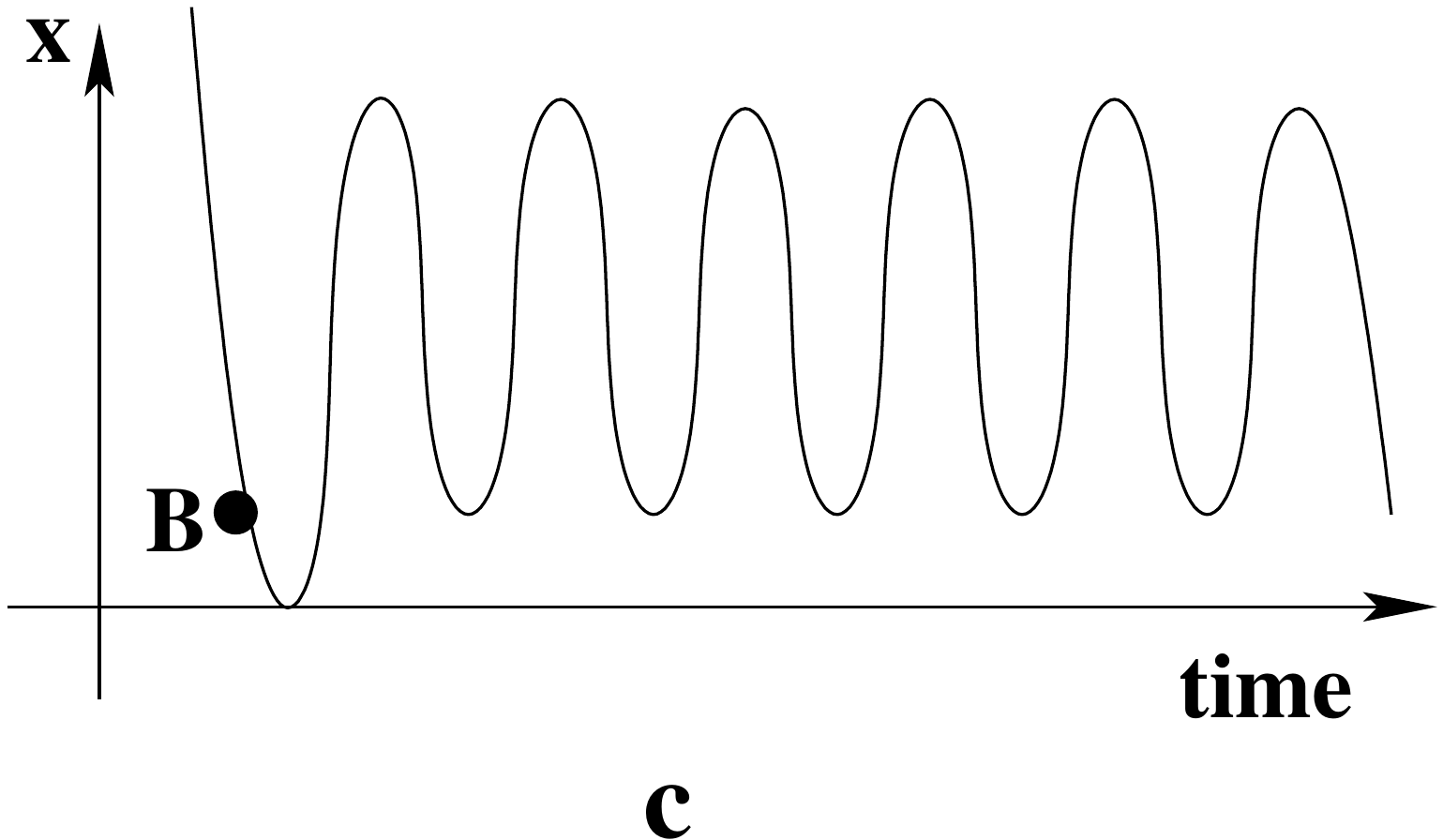,width=5cm}
}
\caption{ Dynamics of a system with a stable limit cycle from fig.\ref{flc1}a. \label{flc2} }
\centerline{
\psfig{type=pdf,ext=.pdf,read=.pdf,figure=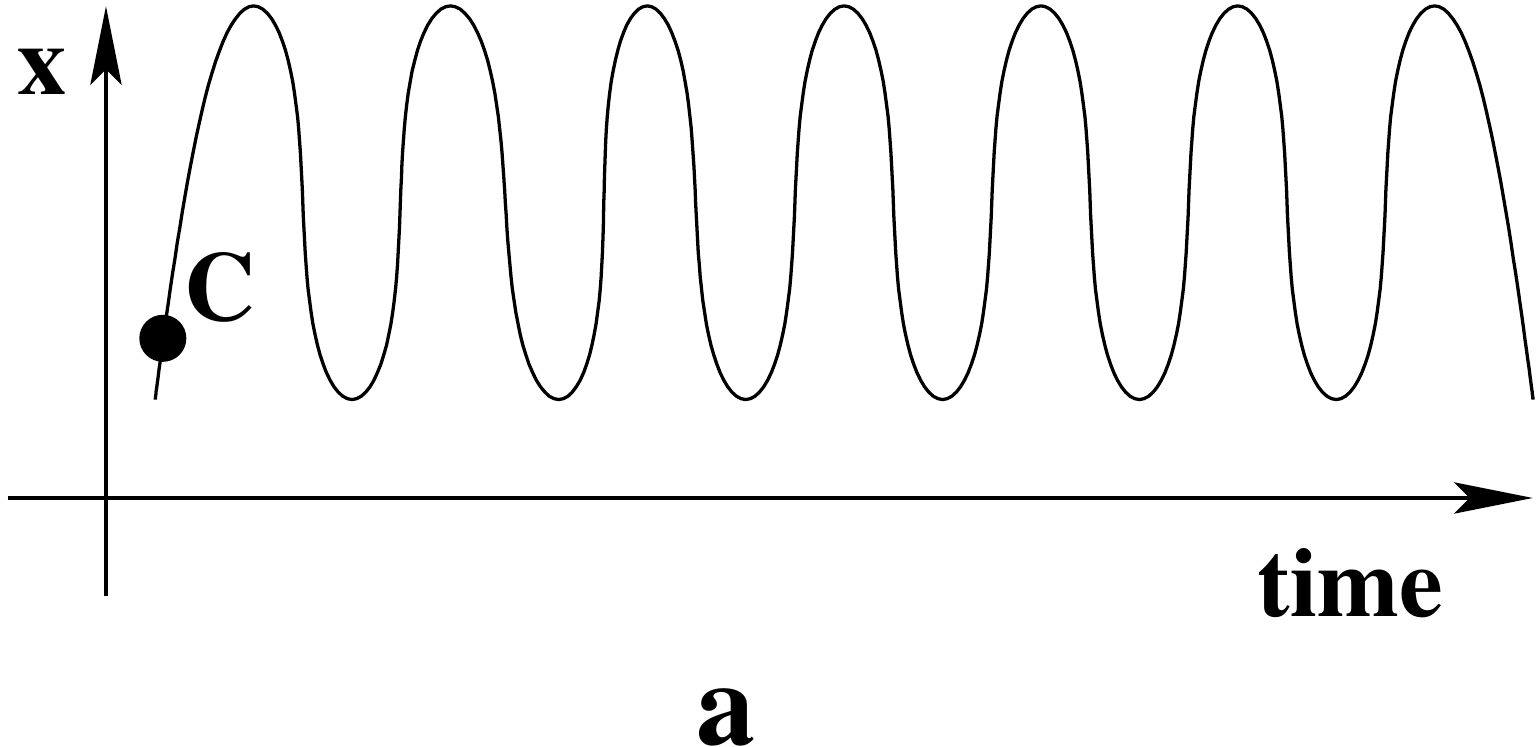,width=5cm}
\psfig{type=pdf,ext=.pdf,read=.pdf,figure=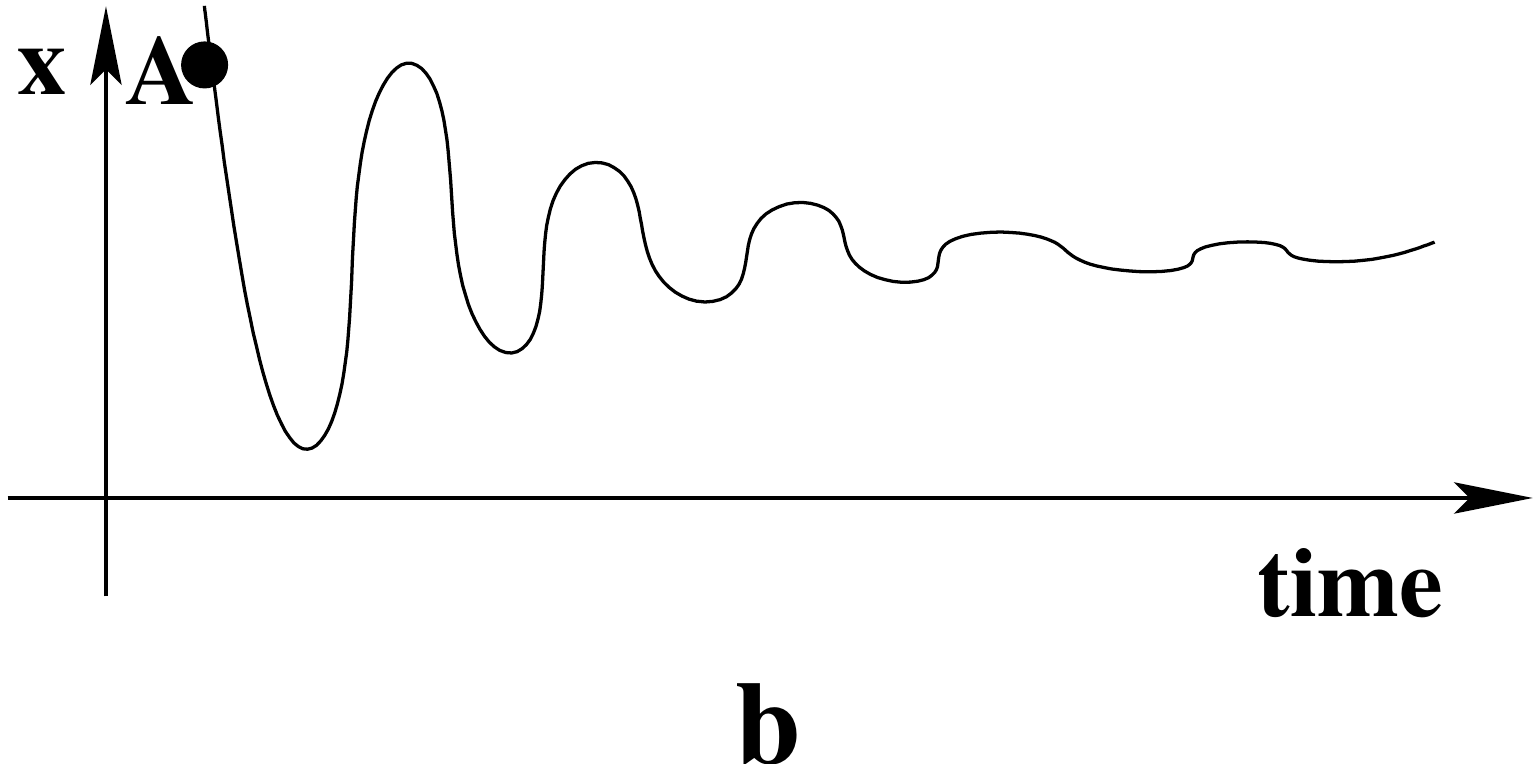,width=5cm}
\psfig{type=pdf,ext=.pdf,read=.pdf,figure=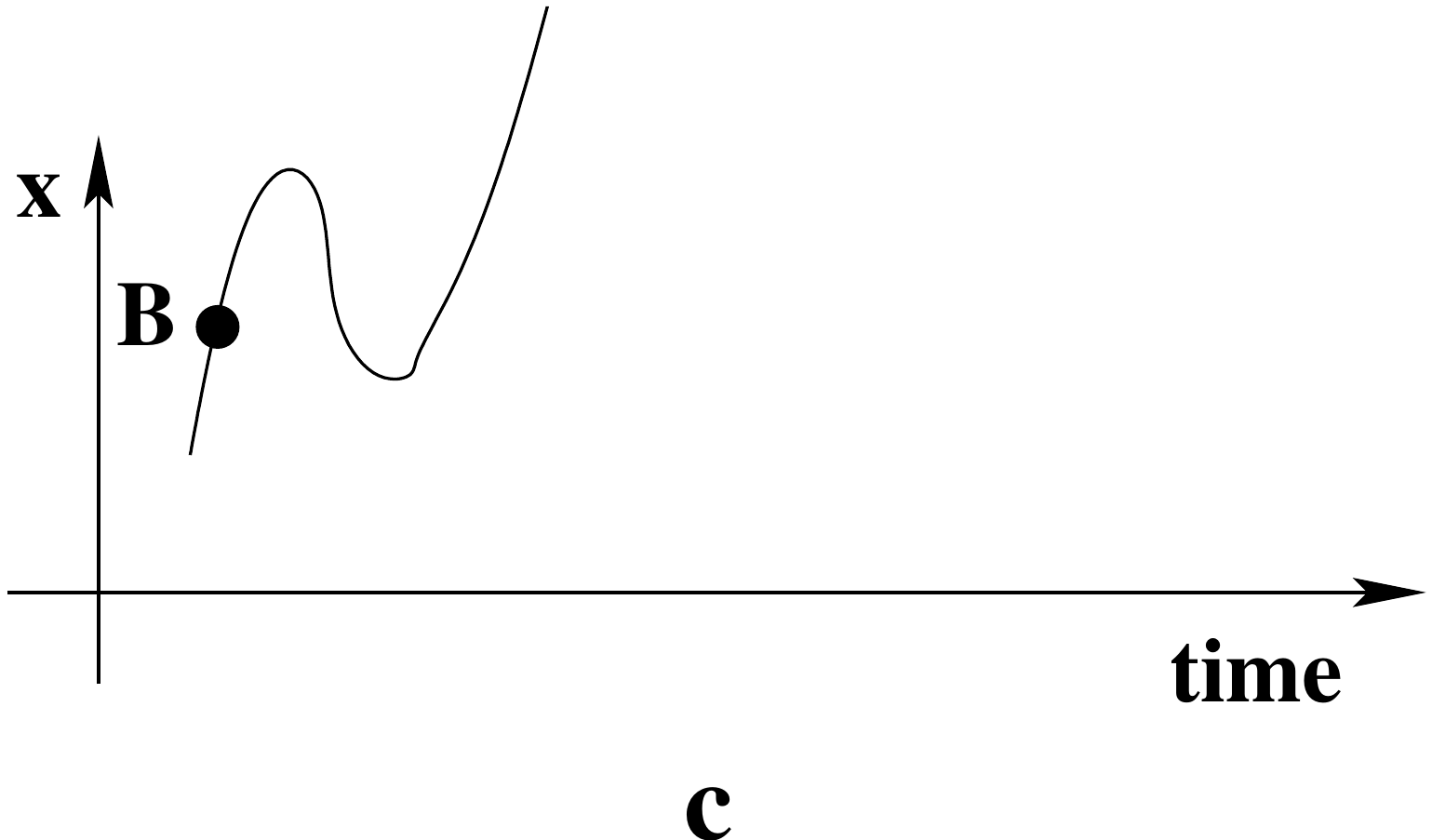,width=5cm}
}
\caption{ Dynamics of a system with an non-stable  limit cycle from fig.\ref{flc1}b. \label{flc3} }
\end{figure}

In previous chapters we found several possible types of equilibria: saddle, node, spiral and  center.
Some of these  equilibria can be attractors  which  determine the ultimate dynamics of a system.
 However, it turns out  that there exists  another  important attractor for  systems of two 
differential equations. It is called a {\bf limit cycle}. The geometrical image of a limit cycle on a phase portrait 
is a closed curve. 
The usual phase portrait of a system with a limit cycle is shown in fig.\ref{flc1}a. Here the limit cycle is shown as a bold
ellipse.  The figure  also shows two more trajectories: one inside the limit cycle and one outside the limit cycle.
Inside the limit cycle you can see a spiral. This is not unusual.
There is a theorem which states   that for systems of two equations there is always
at least one equilibrium point inside the limit cycle. In most of the cases such an equilibrium 
point is a spiral.
Outside the limit cycle we can see a trajectory which approaches the 
limit cycle and winds  onto  it.

The limit cycle which is shown in  fig.\ref{flc1}a is called a {\bf stable limit cycle}.
This is because if we start on a trajectory close to this limit cycle, this trajectory 
will approach  this  limit cycle in the course of time.

There also exists  another type of limit cycle called  a  {\bf non-stable limit cycle} (fig.\ref{flc1}b).
If we start on a trajectory close to  a non-stable limit cycle, it  will diverge
from  this  limit cycle. In order to distinguish the stable and non-stable limit cycles we will draw
the  non-stable limit cycle using a  dashed line.

The main questions  regarding the limit cycle are:
what will be the dynamics of systems with limit cycles and  how do limit cycles occur in systems
of two differential equations?  We will also consider an example of a biological system with a limit cycle.

\section{Dynamics of a system with a limit cycle.}
What will be the  dynamics of systems with a limit cycle?
In section \ref{sec_center} we studied an  equilibrium point called a ``center''. The phase portrait of that point
(fig.\ref{fig7.1}) was a set of embedded ellipses. The corresponding dynamics are oscillations of the 
$x$ and $y$ variables. Therefore,  the dynamics which corresponds to the trajectory which
starts on  the limit cycle will also be  oscillations. 
Fig.\ref{flc2}a, fig.\ref{flc3}a shows  an example of dynamics of the variable $x$ for the  trajectory 
which starts on the limit cycle (at the point {\bf C} in fig.\ref{flc1}).

The dynamics of trajectories originating  around  the limit cycles will depend on its  type.
If the limit cycle is stable, then a trajectory which starts inside the limit cycle will approach it in the course of time and
we obtain oscillations with  initially increasing amplitude  (fig.\ref{flc2}b). 
If the trajectory starts outside the limit cycle, then we obtain oscillations with an initial decrease of  amplitude 
 until
the trajectory reaches the limit cycle (fig.\ref{flc2}c).

For a non-stable limit cycle the dynamics are   different. Only the   trajectory which starts on a limit cycle will 
have oscillatory dynamics (fig.\ref{flc3}a).  Other trajectories will have different behavior. 
The trajectory which originates  inside the limit cycle will approach the stable equilibrium inside the limit cycle 
and we obtain oscillations with gradually decreasing amplitude  (fig.\ref{flc3}b).  The trajectory which starts outside the limit 
cycle blows up to infinity or to another equilibrium  (fig.\ref{flc3}c).
The oscillations  in fig.\ref{flc3}a with a non-stable limit cycle are highly improbable in real systems. This is 
because  for such oscillations the  trajectory must  start exactly on the  limit cycle and even small disturbances 
will switch our system either to the behavior of fig.\ref{flc3}b or  fig.\ref{flc3}c. Therefore,  for  real systems 
the  non-stable limit cycle  just determines the basin of attraction of the  stable equilibrium 
which is located inside this limit cycle.

\section{How do  limit cycles occur?}
\label{hopf}
\begin{figure}[h]
\centerline{
\psfig{type=pdf,ext=.pdf,read=.pdf,figure=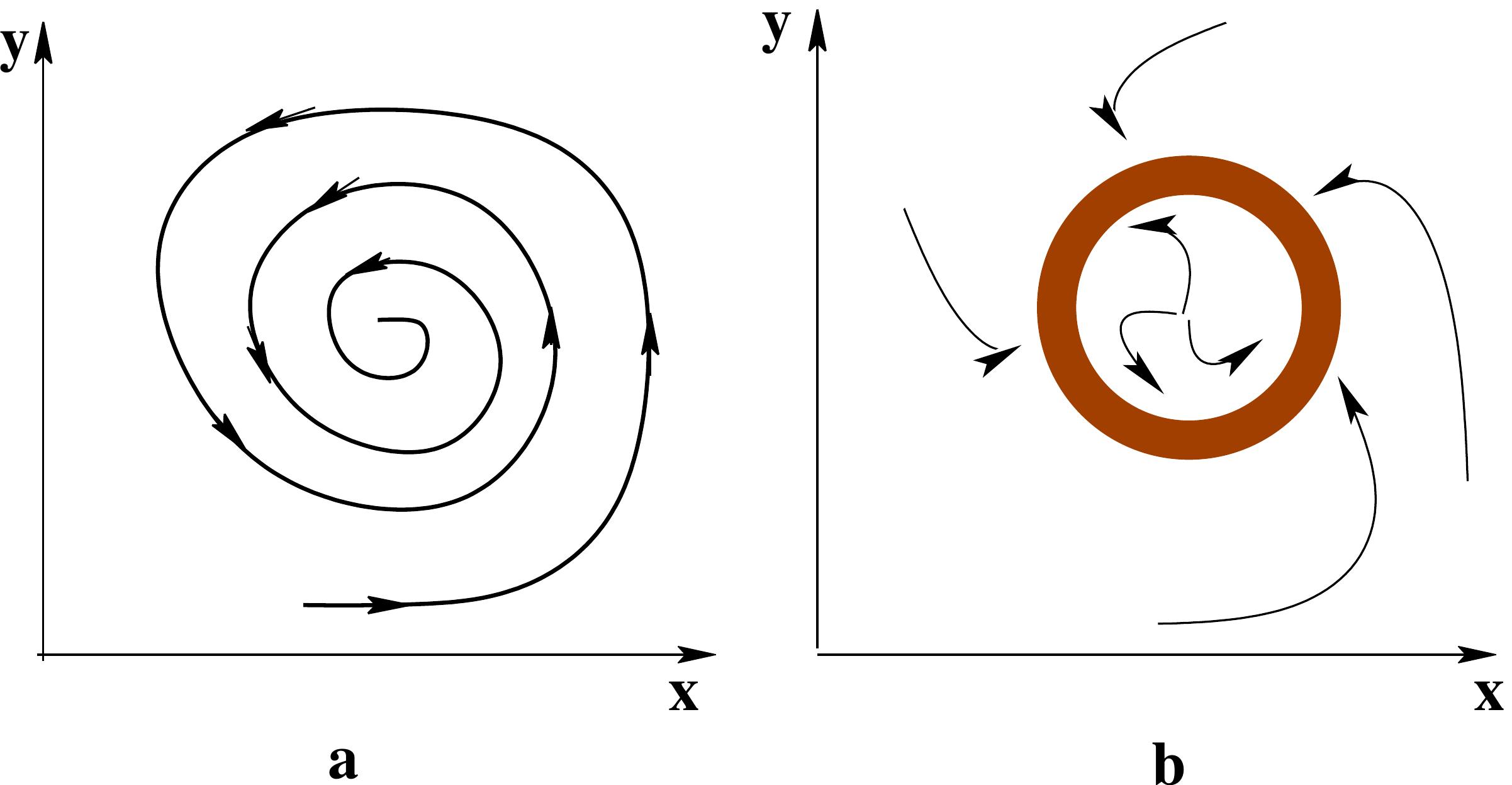,width=10cm}
}
\caption{ Appearance   of a limit cycle in a system with a parameter. \label{flc4} }
\end{figure}

In many cases  the limit cycle occurs  as a result of the changing of a parameter in a   system  of differential equations.
The most usual process of formation of a limit cycle is the following. Assume we have a system of two differential equations with a parameter $c$:
\beq
\label{lc1}
\left\{
\begin{array}{l}
{dx \over dt}=f(x,y,c) \\
{dy \over dt}=g(x,y,c)
\end{array}
\right.
\eeq
Assume that  at some parameter value $c=c_{1}$ system (\ref{lc1}) has a global  attractor  which is a stable spiral
(fig.\ref{flc4}a).  This means that  all the trajectories which originate close to or even far from this 
equilibrium approach it in the course of time.   Because our equilibrium point is a stable  spiral the Jacobian of the  system at this point
will have two complex eigen values $\lambda_{1,2}=\alpha \pm i\beta; \alpha < 0$. However,  because our system will now depend on   the parameter $c$  the eigen values will also depend on this parameter:

\beq
\lambda_{1,2}(c)=\alpha(c) \pm i\beta(c); 
\eeq
Because we have assumed that  at $c=c_{1}$ system (\ref{lc1}) has an equilibrium which is  a stable spiral, then
$\alpha(c_{1})<0$. When we change the parameter $c$   the value of $\alpha(c)$ will  change and  at some $c_2$ 
it can  
become a  positive number $\alpha(c_{2})>0$.
 This means  that  the stable spiral in  fig.\ref{flc4}a will become an  unstable spiral.  However, as we found in chapter 4 (fig.4.3), 
the  eigenvalues only give us
the dynamics  close to the equilibrium point  of our system.  Therefore, 
 it can happen that although close to the equilibrium point  our spiral becomes unstable,
 the global behavior far from the equilibrium remains  the same,
i.e. far from the equilibrium we still have a converging flow (fig.\ref{flc4}b). 
Hence,   in our phase portrait we have two types of flow: the  diverging flow from the unstable spiral around the equilibrium  and the  converging flow from the periphery. 
Simple geometrical consideration shows that these  two flows must be separated from each other.
The line of separation will be the  limit cycle in our system.  Therefore,  we will obtain  a phase portrait as in fig.\ref{flc1}a.

Note, that such a mechanism  of limit cycle formation  frequently occurs in systems of  two equations. It is called
 the Hopf bifurcation. From our analysis it follows that the Hopf bifurcation occurs when $\alpha(c)$ changes its sign,
i.e. at the parameter value $c^*$ where:
\beq
\label{ehopf1}
  \alpha(c^*)=0
\eeq

\section{Example of a system with a limit cycle}

\begin{figure}[h]
\centerline{
\psfig{type=pdf,ext=.pdf,read=.pdf,figure=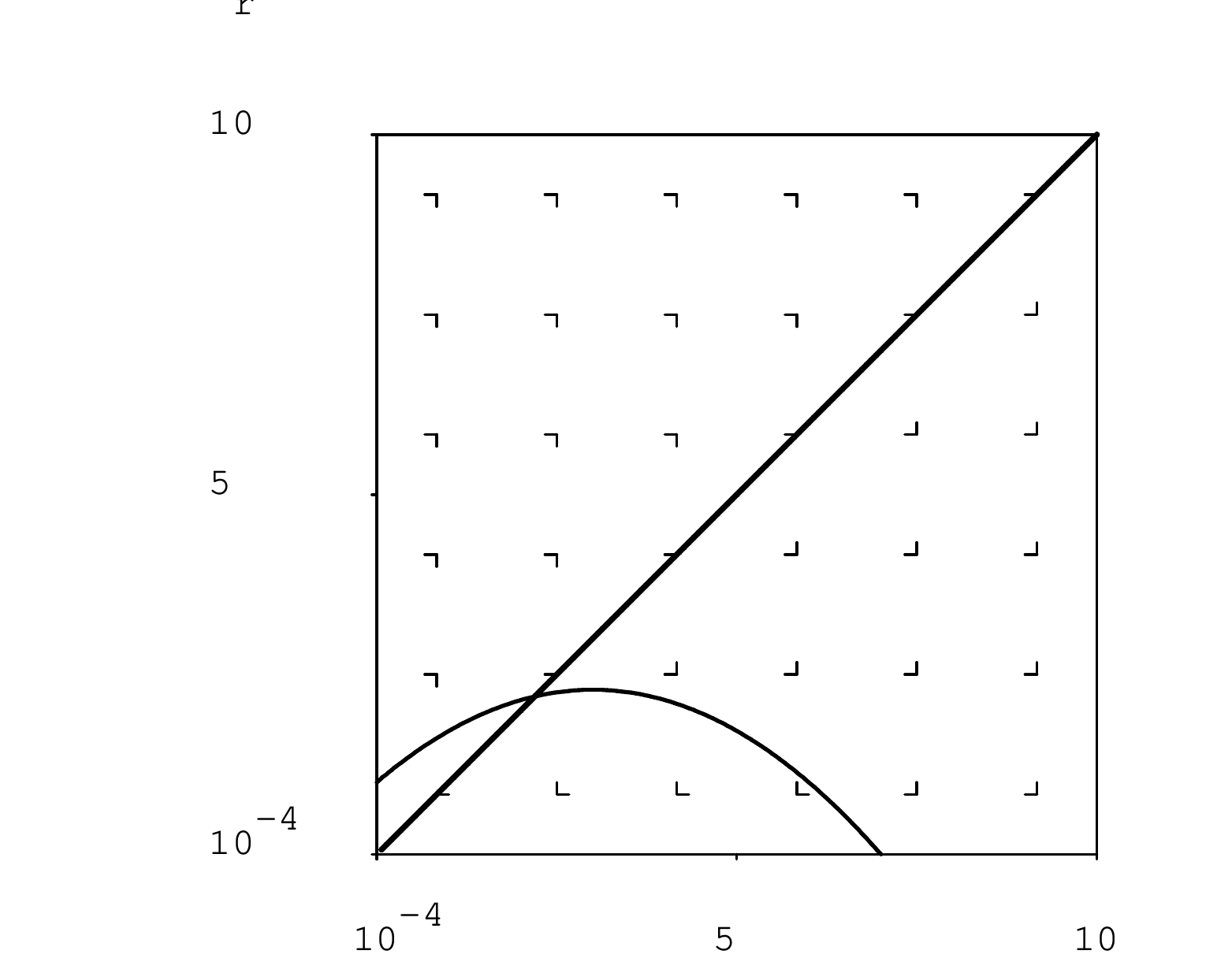,width=5cm}
\psfig{type=pdf,ext=.pdf,read=.pdf,figure=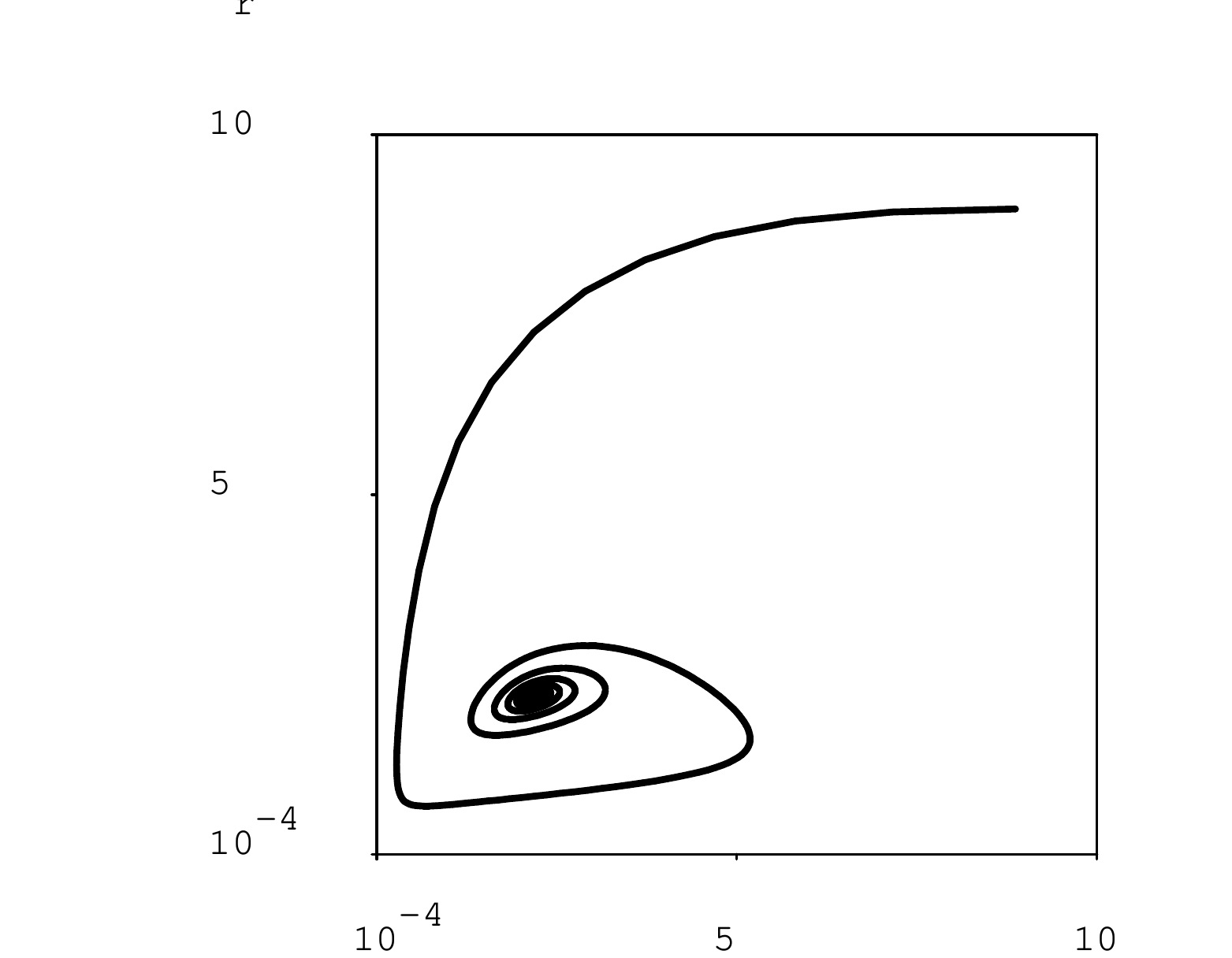,width=5cm}
\psfig{type=pdf,ext=.pdf,read=.pdf,figure=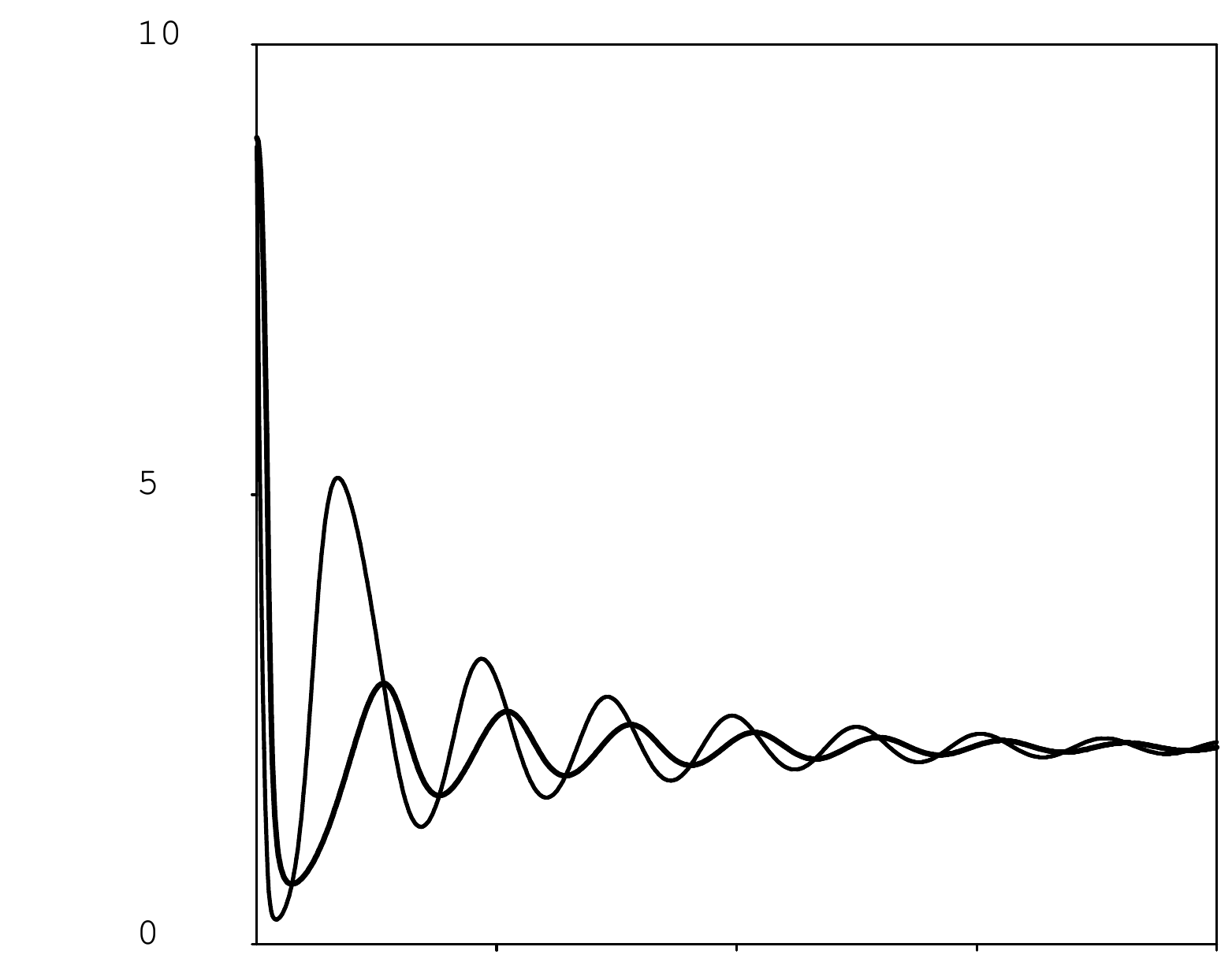,width=5cm}
}
\caption{ Dynamics of  system (\ref{lc2}) 
for  $K=0.7$. a-null-clines; b-an orbit; c-time-plot for the both variables for the orbit from fig.b
 \label{flc5} }
\end{figure}

Let us consider the Holling-Tanner  model for predator-prey interactions.
\beq
\label{lc2}
\left\{
\begin{array}{l}
dP /dt=rP(1-{P \over K} )-{aRP \over d+P} \\
dR/dt= bR(1-{R \over P})  \qquad P>0; R>0
\end{array}
\right.
\eeq
here $P$  denotes the prey and $R$ the predator population, the term $rP(1-P/K)$ describes the growth of the prey
in absence of predator, ${aRP \over d+P} $ accounts for the predator-prey interaction. 
At $a=1,b=0.2,r=1.,d=1.,K=0.7$ system (\ref{lc2}) has one equilibrium point (for $P>0,R>0$) which is a stable spiral.
The null-clines, phase portrait and dynamics of this system are  shown in fig.\ref{flc5}.

\begin{figure}[h]
\centerline{
\psfig{type=pdf,ext=.pdf,read=.pdf,figure=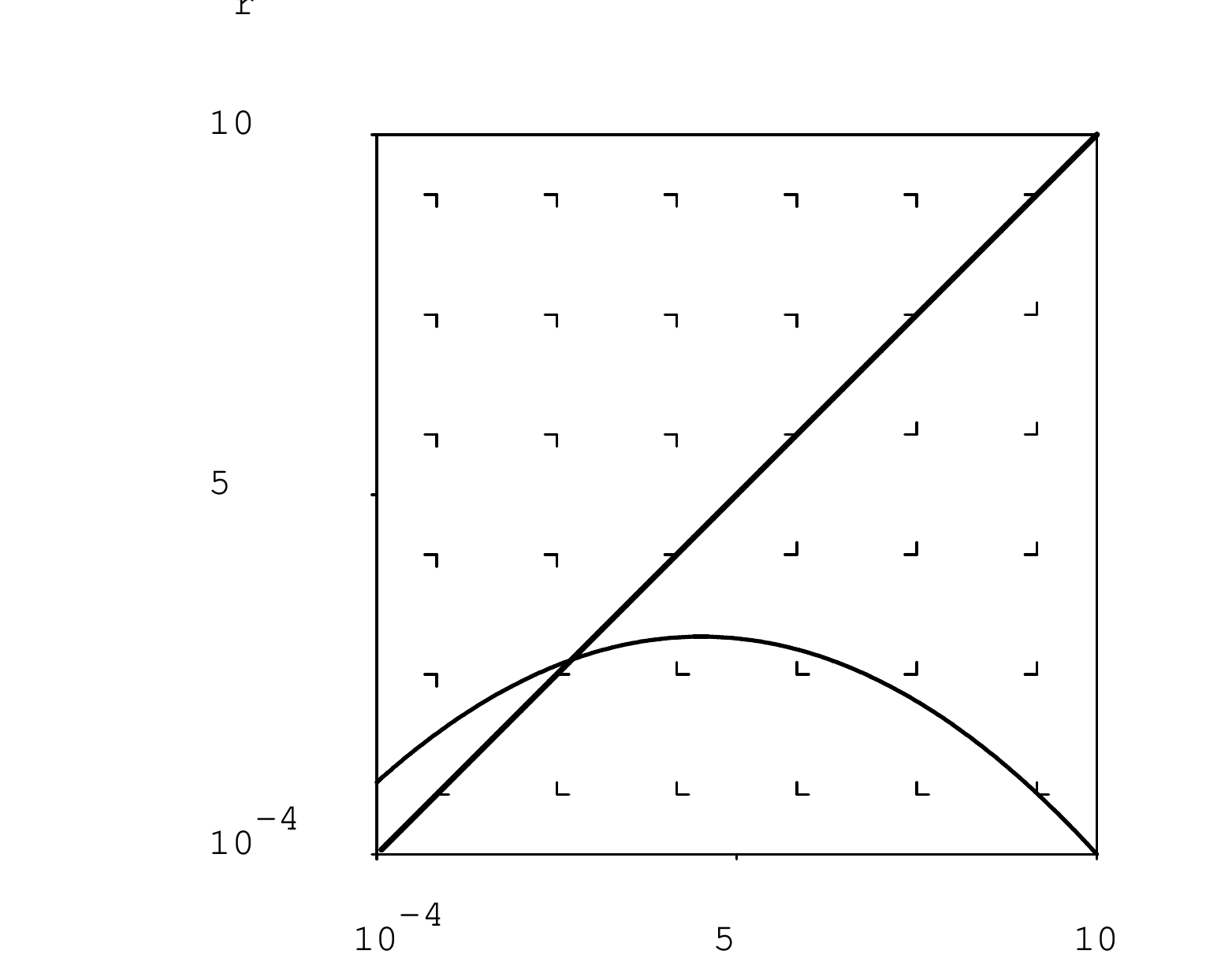,width=5cm}
\psfig{type=pdf,ext=.pdf,read=.pdf,figure=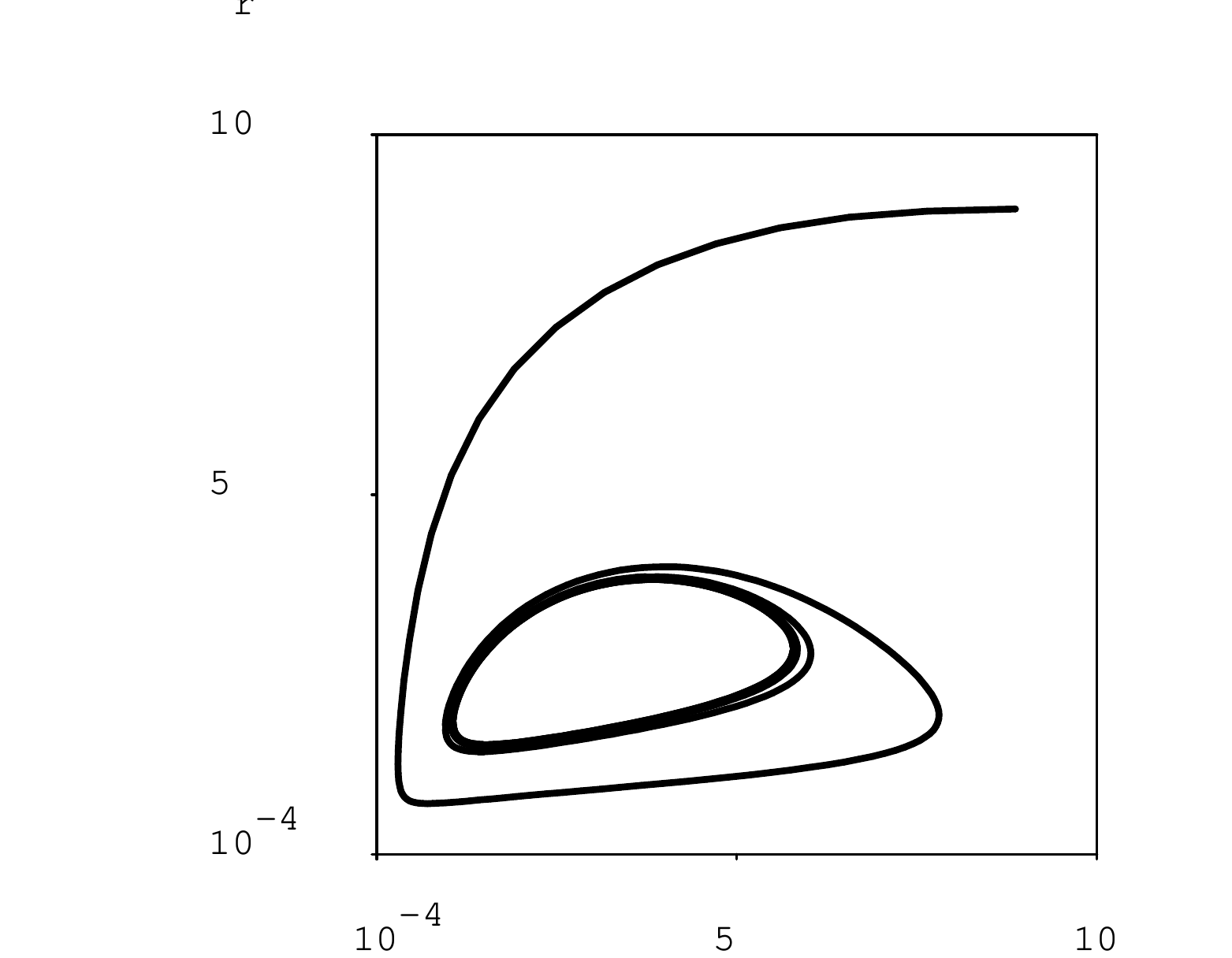,width=5cm}
\psfig{type=pdf,ext=.pdf,read=.pdf,figure=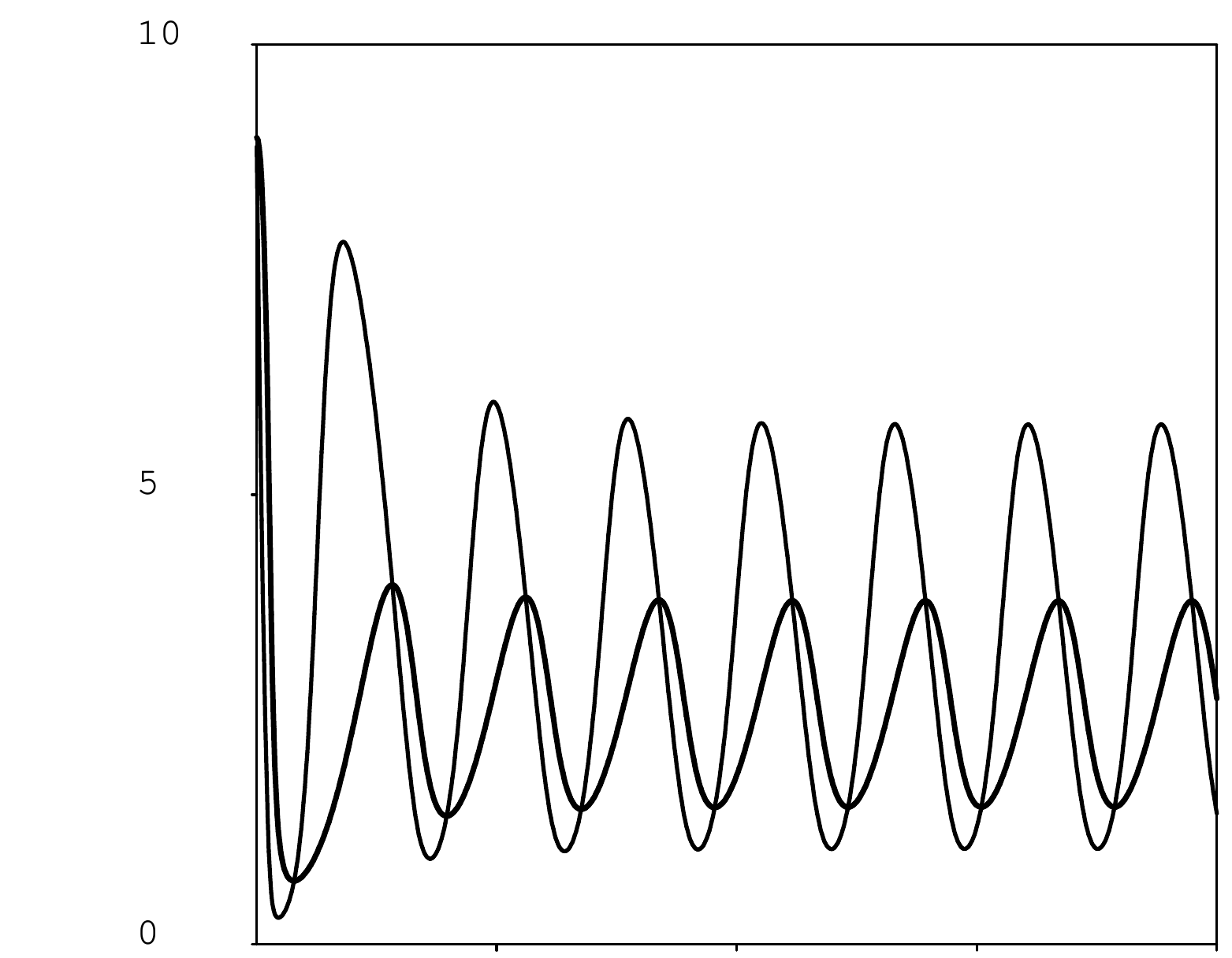,width=5cm}
}
\caption{ Dynamics of  system (\ref{lc2}) 
for  $K=1.0$. a-null-clines; b-an orbit; c-time-plot for  both variables for the orbit from fig.b
 \label{flc6} }
\end{figure}

If we increase the value of the parameter $K$ which accounts for the
carrying capacity ($K=1.0$) the type of equilibrium changes and
it becomes an unstable spiral. As we discussed in section \ref{hopf}
we expect the formation of a stable limit cycle. We can clearly see it
in fig.\ref{flc6}b.  The trajectory, which starts from the same
initial conditions as the trajectory from fig.\ref{flc5}b will now
approach some closed curve which is a limit cycle and the dynamics of
the system will be oscillatory (fig.\ref{flc6}c).

\begin{figure}[h]
\centerline{
\psfig{type=pdf,ext=.pdf,read=.pdf,figure=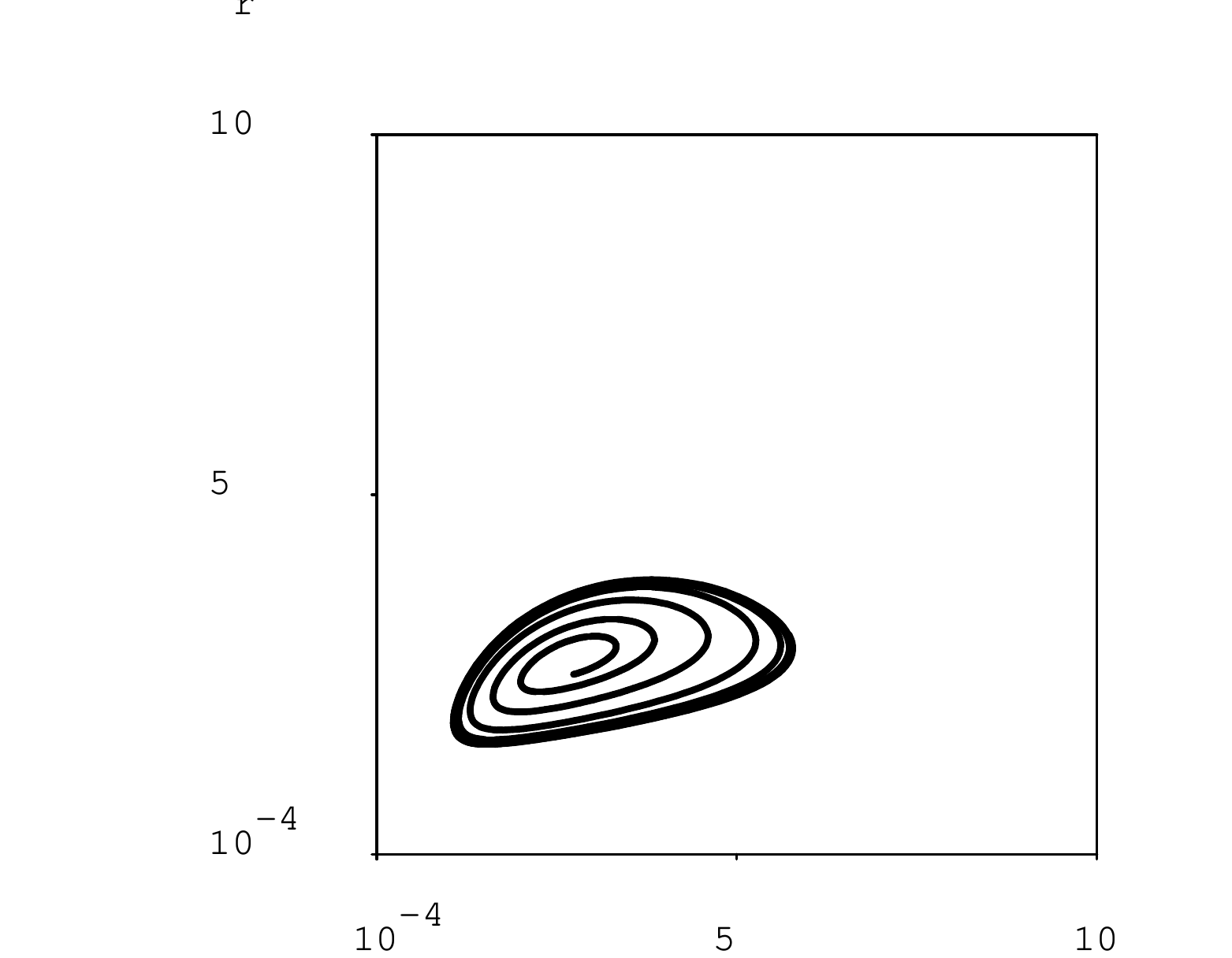,width=5cm}
\psfig{type=pdf,ext=.pdf,read=.pdf,figure=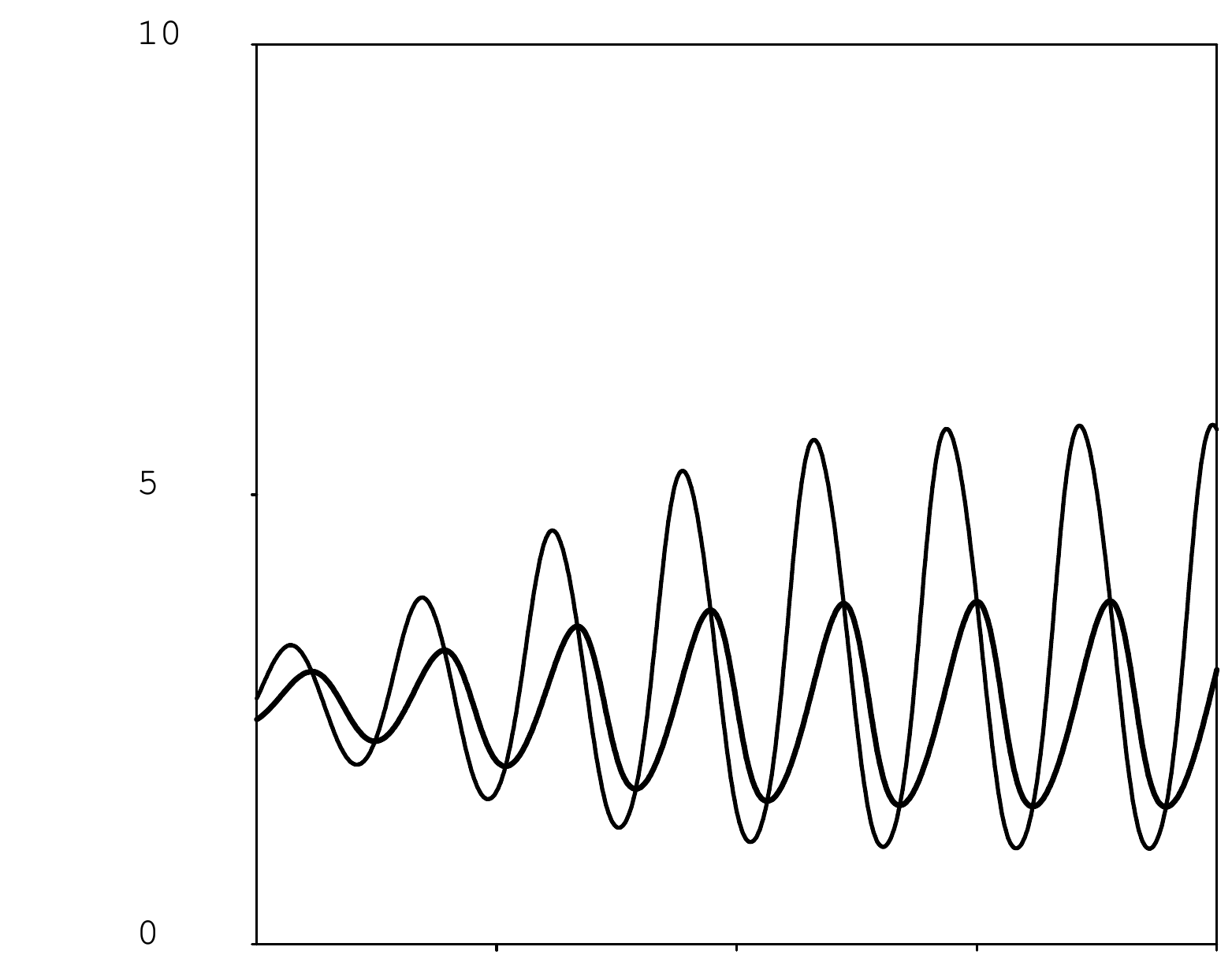,width=5cm}
\psfig{type=pdf,ext=.pdf,read=.pdf,figure=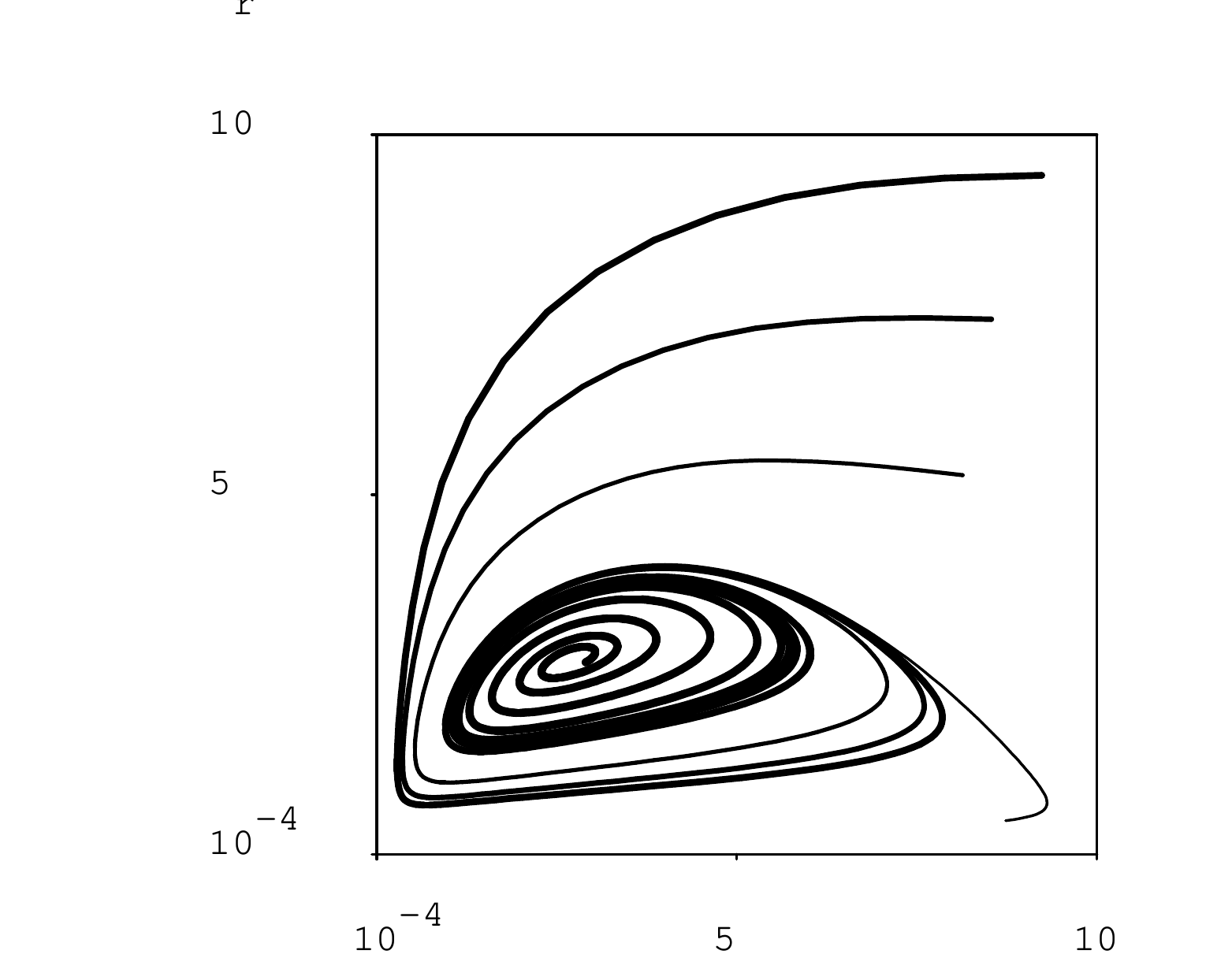,width=5cm}
}
\caption{ Dynamics of  system (\ref{lc2}) 
for  $K=1.0$. a-an orbit originating inside the limit cycle; b-time-plot for the both variables for the orbit from fig.a; c-phase portrait of system  (\ref{lc2}) 
 \label{flc7} }
\end{figure}

If we start a  trajectory inside the limit cycle  then,  as we predicted in  fig.\ref{flc1}a and fig.\ref{flc2}b,  the trajectory will approach the limit cycle and  the dynamics will
be oscillation with increasing amplitude fig.\ref{flc7}a,b.
The complete phase portrait of this system at $K=1$ is shown in fig.\ref{flc7}c.

You will learn more about the biological consequences  of this type of dynamics  in the course ``Theoretical Biology'' by R. de Boer. 

\begin{C}
A limit cycle is  a closed trajectory on the  phase portrait of a system of two differential equations. 

If  trajectories around the limit cycle converge onto it, then the limit cycle is called a stable limit cycle.
If  trajectories around the limit cycle diverge away from  it, then the limit cycle is called a non-stable limit cycle.

The dynamics of a system with a stable limit cycle is oscillatory.
The dynamics of a system with an non-stable limit cycle is either converging to the  equilibrium which is located within  the limit cycle
or diverging, possibly to infinity.

 Limit cycles can appear as a result of a  Hopf bifurcation, i.e. the process where the real part of the 
complex eigenvalues change their sign.
\end{C}

\section{Exercises}
\begin{figure}[hhh]
\centerline{
\psfig{type=pdf,ext=.pdf,read=.pdf,figure=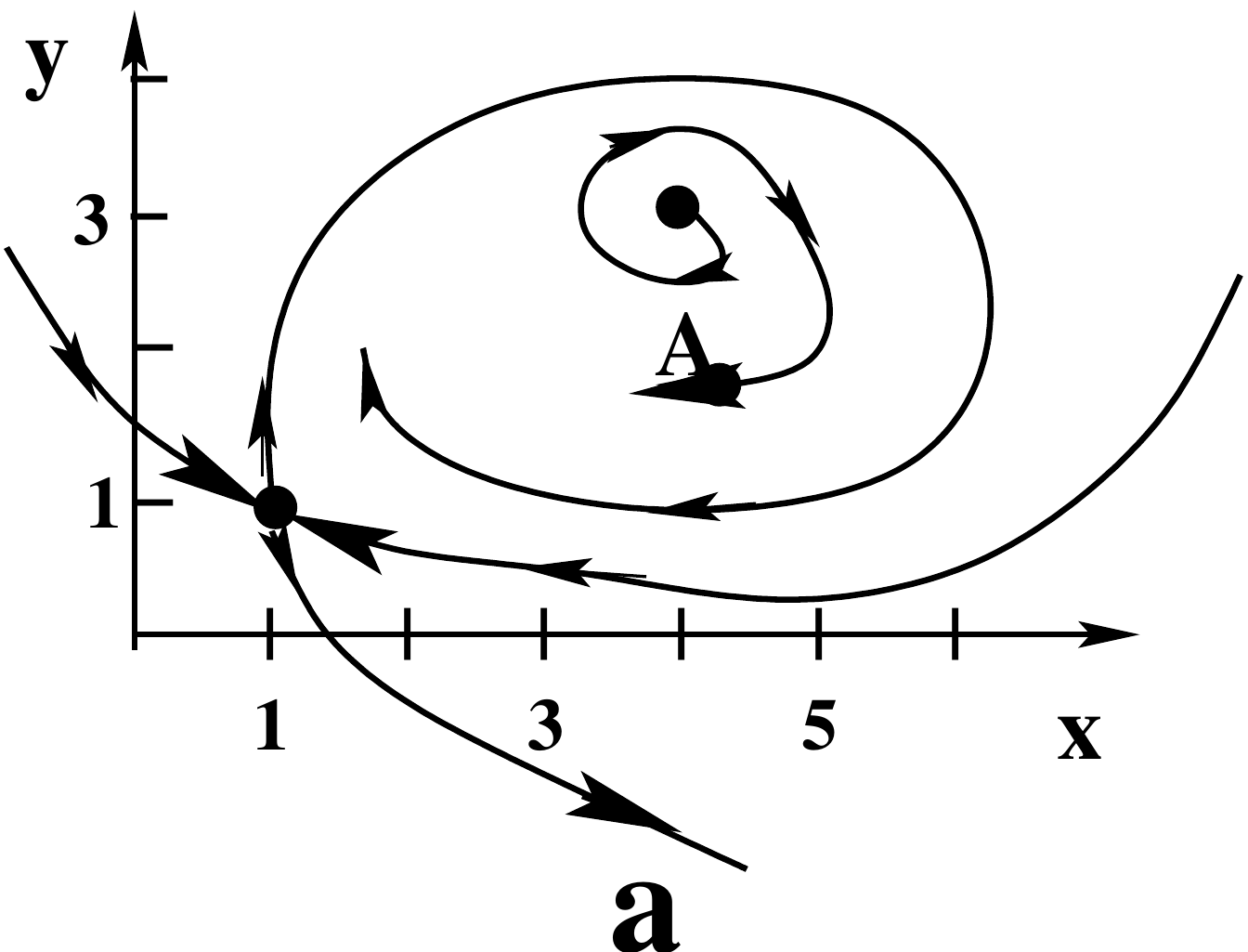,width=5cm}
\psfig{type=pdf,ext=.pdf,read=.pdf,figure=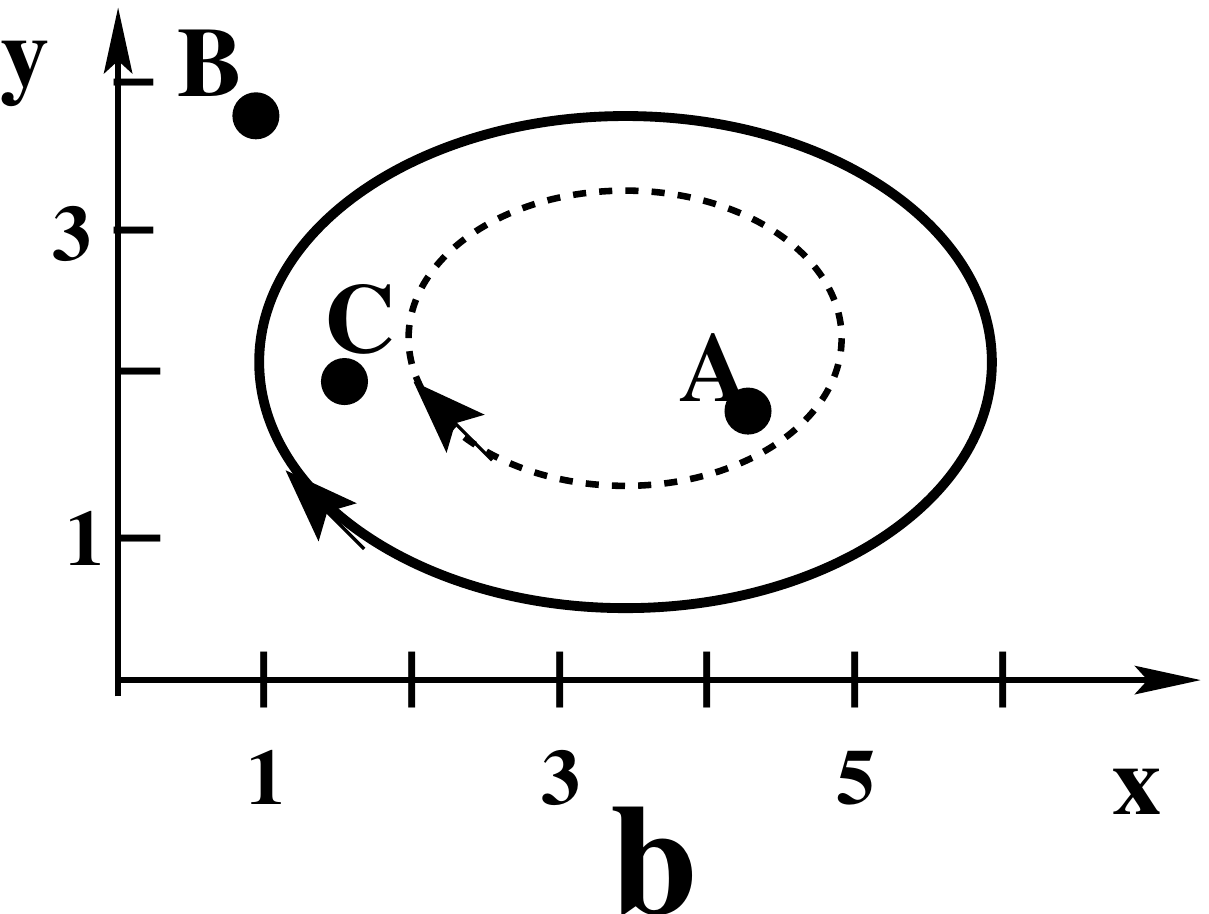,width=5cm}
\psfig{type=pdf,ext=.pdf,read=.pdf,figure=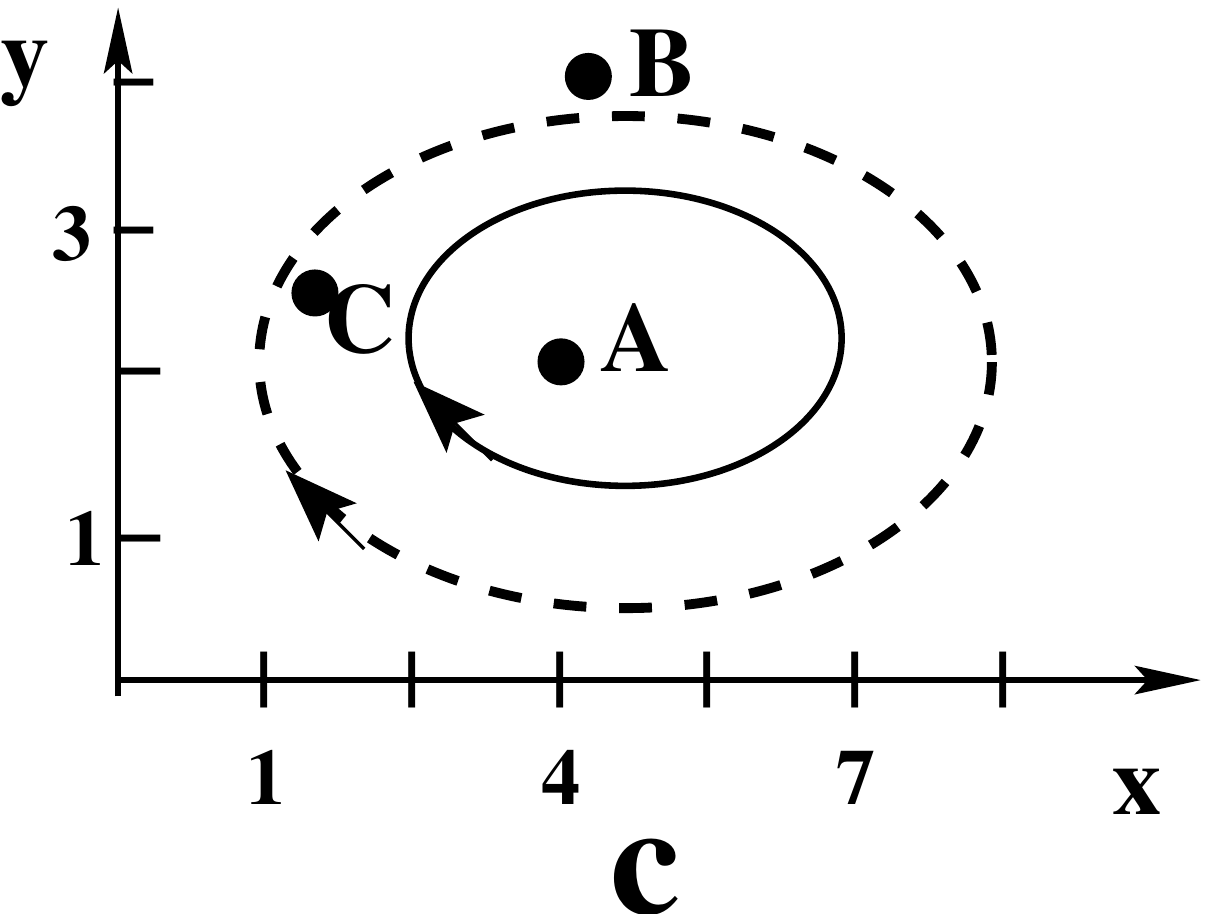,width=5cm}
}
\caption{\label{examsfig4}}
\end{figure}

\ben
\item  Assume that a system of two differential equations has two equilibria
which are a saddle point and a non-stable spiral.
The phase portrait of this system is partially shown in  fig.\ref{examsfig4}a 
What is missing here?

\begin{enumerate} 

\item Complete the phase portrait of this system 

\item Qualitatively sketch the dynamics of the variable $x$ ( dependence of the variable $x$ on time)
 for the initial conditions which are shown in fig.\ref{examsfig4} by
 the point {\bf A} in fig.\ref{examsfig4}a.

\end{enumerate}

\item Complete the phase portrait of this system shown in fig.\ref{examsfig4}b,c. Qualitatively sketch the dynamics of the variable $x$ ( dependence of the variable $x$ on time)
 for the initial conditions which are shown  by points {\bf
 A,B,C} in fig.\ref{examsfig4}b,c. (Note, that in fig.\ref{examsfig4}b,c a stable
 limit cycle is shown by the solid line and a non-stable limit cycle
 is shown by the dashed line)

\item As we discussed in (\ref{ehopf1}) the necessary condition for the appearance of a limit cycle
via a  Hopf bifurcation is  that the eigenvalues of the Jacobian matrix of a system at the equilibrium point are complex and the real part of the eigenvalues are  zero $\alpha(c^*)=0$. Prove that this condition can be rewritten in the following way using the $det$ and $tr$ of the Jacobian matrix:
\beq
\label{ehopf2}
trJ(c^*)=0; \qquad detJ(c^*)>0
\eeq

\item One of the classical models for oscillatory phenomena in biochemical systems  is a 
model called the Brusselator. In dimensionless form this model can be written as the following system of two differential equations:
\beq
\label{ehopf3}
\left\{
\begin{array}{l}
{dx \over dt}=a-(b+1)x+x^2y  \qquad x>0;y>0\\
{dy \over dt}=bx-x^2y \qquad a>0;b>0
\end{array}
\right.
\eeq
here $x$ and $y$ are concentrations of two biochemical species and $a$ and $b$ are parameters.
Study system (\ref{ehopf3}) for $a=1$. 
\ben
\item Find  a non-trivial equilibrium.
\item Determine the stability of this  equilibrium as a function of the parameter $b$.
\item  Find the value of $b$ when system  (\ref{ehopf3}) undergoes a  Hopf bifurcation.
  (Note: it can be helpful to use equations (\ref{ehopf2})).
\item For which values of $b$ do you  expect  oscillations in  system  (\ref{ehopf3})?
\een

\item Study the following predator-prey model
\beq
\label{lc_ex2}
\left\{
\begin{array}{l}
dP /dt=rP(1-{P \over K} )-{aRP \over h+P} \\
dR/dt= {caRP  \over h+P} -dR  \qquad P \geq 0; R \geq 0
\end{array}
\right.
\eeq
\ben
\item Draw null-clines of this system for $r=1,a=3,h=0.1,c=1.,d=2.5$ and two values of $K$, $K=0.8$ and $K=1.6$. (Hint: the maximum of the parabola $A(x-a)(x-b)$ is reached at the middle between its roots (i.e. at $x={a+b \over 2}$)).
\item Determine the stability of non-trivial equilibrium  in both cases using the graphical Jacobian.
\item For which value  of $K$ do you  expect  oscillatory behavior?

\item Extend your analysis for arbitrary positive values of the parameters ($r,a,h,c,d,K>0$) provided $ca>d$. Find for which values of $K$ the non-trivial equilibrium is stable. When do we expect oscillations? (Note, that critical values of $K$ should be a function of the other parameters of the system).
\een
\een

\chapter{Historical notes}
I developed this course in 1995-2010 for biology students from Utrecht University. The main idea was to develop  a course which will allow students with minimal mathematical  background to understands basic methods used in mathematical biology. Over the years the course was many times adjusted and modified. The current text contains the latest version of the course from 2010.
\chapter{Dictionary}

\begin{tabular}{ll}
absolute value  & absolute waarde \\
autonomous system & autonoom systeem \\
attractor  & attractor   \\
basin of attraction   & basin van attractie \\
bifurcation  & bifurcatie \\
carrying capacity  & draagkracht \\
center point  & centrumpunt \\
complex conjugate number  & complex toegevoegde \\
component of the vector  & component van de vector \\
determinant & determinant \\
derivative  & afgeleide \\
differential equation  & differentiaal vergelijking\\
direction field  & vectorveld \\
eigen value & eigenwaarde \\
eigen vector & eigenvector \\
equilibrium  & evenwicht \\
general solution  & algemene oplossing \\
harvesting & oogsten \\
imaginary part of complex number  & het imaginaire deel van een
complex getal \\
initial value problem  & beginwaarde probleem \\
Jacobian  & Jacobi-matrix \\
linear approximation   & lineaire benadering \\
modulus  & modulus \\
node  & knooppunt \\
non-stable manifold  & instabiele manifold  \\
\end{tabular}
\eject
\begin{tabular}{ll}
null-cline  & isocline \\
parameter  & parameter \\
partial derivative  & parti\"ele afgeleide \\
particular solution  & specifieke oplossing \\
phase space  & faseruimte \\
phase portrait  &  faseplaatje \\
real part of complex number  & het re\"eele deel van een complex
getal\\
saddle  & zadelpunt \\
spiral  & spiraalpunt \\
stability  & stabiliteit \\
stable manifold  & stabiele manifold   \\
system of differential equations  & stelsel differentiaal vergelijkingen\\
trajectory & trajectorie \\
trace of the matrix & spoor van de matrix \\
variable  & variabele \\
vector  & vector \\
\end{tabular}

\textheight=25.7cm
\textwidth=18.5cm
\topmargin=-2.0cm
\oddsidemargin=-0.9cm
\evensidemargin=-0.9cm
\renewcommand{\baselinestretch}{0.5}\normalsize
\small 
\chapter{Hints:}

\subsection{Solution of the initial value problem for a linear system \label{solutionIVP}}

 Let us illustrate how we can find a  solution for the  initial value problem, on  example of system (\ref{2dlinour}).

{\bf Problem:} Find the solution for the following initial value problem:

\beq
\label{2dlinour2}
\left( \begin{array}{c} 
{dx \over dt} \\  {dy \over dt}  \end{array} \right) = 
\left(\begin{array}{lr}  1 & 4\\
                         1 &  1
\end{array} \right),  \quad  
\left( \begin{array}{c} x
 \\  y  \end{array} \right)
\left( \begin{array}{c} x(0)
 \\  y(0)  \end{array} \right)=\left( \begin{array}{c} 4
 \\  6  \end{array} \right)
\eeq

{\bf Solution:}
We already found a general  solution of (\ref{2dlinour2}) as  formula (\ref{2dlinour_sol}):
\beq
\label{2dlinour_sol2}
\left( \begin{array}{c} 
x \\  y \end{array} \right) = 
C_1 \left( \begin{array}{c} -4
 \\  2  \end{array} \right)e^{-1* t }
+ 
C_2 \left( \begin{array}{c} -4
 \\  -2 \end{array} \right)e^{3* t }
\eeq

Using it we can find a  particular solution, corresponding to any given
initial conditions.

We proceed as follows. At time $t=0$   (\ref{2dlinour_sol2}) gives:
$$\left( \begin{array}{c} 
x \\  y \end{array} \right) = 
C_1 \left( \begin{array}{c} -4
 \\  2  \end{array} \right)e^{-1* 0 }
+ 
C_2 \left( \begin{array}{c} -4
 \\  -2 \end{array} \right)e^{3* 0 } = 
C_1 \left( \begin{array}{c} -4
 \\  2  \end{array} \right)
+ 
C_2 \left( \begin{array}{c} -4
 \\  -2 \end{array} \right).
$$
This expression will satisfy the initial conditions if:
$
C_1 \left( \begin{array}{c} -4
 \\  2  \end{array} \right)
+ 
C_2 \left( \begin{array}{c} -4
 \\  -2 \end{array} \right)=\left( \begin{array}{c} 4
 \\  6  \end{array} \right)$. This gives the following system of equations for unknowns $C_1$ and $C_2$:
$\left\{
\begin{array}{l}
-4C_1-4C_2= 4 \\ 2C_1-2C_2=6
\end{array}
\right.$

We now solve this system using the method described in section \ref{SysEq}.  From the second equation we find: \\ $2C_1=6+2C_2$, or $C_1=3+C_2$. After substitution of this expression to the first equation we find:\\ $-4(3+C_2)-4C_2=4$, i.e. $-12-4C_2-4C_2=4$, or  $-8C_2=4+12=16$ and  $C_2=16/(-8)=-2$.\\ We now find $C_1$ from our substitution  as $C_1=3+C_2=3-2=1$. 

Thus  $C_1=1,C_2=-2$ give the solution of our system satisfying  given  initial conditions. Let us rewrite it as:
$$\left( \begin{array}{c} 
x \\  y \end{array} \right) = 
1* \left( \begin{array}{c} -4
 \\  2  \end{array} \right)e^{-1* t }
-2* \left( \begin{array}{c} -4
 \\  -2 \end{array} \right)e^{3* t }= 
\left( \begin{array}{c} -4e^{-1* t }-2*(-4)e^{3* t }
 \\  2e^{-1* t } -2*(-2)e^{3* t } \end{array} \right), 
$$  thus the particular solution  is given by $x(t)=-4e^{-1* t }+8e^{3* t }$, and $y(t)=2e^{-1* t } +4e^{3* t }$.

\subsection{Equilibria/derivatives \label{solutionsCh5.1}}
 {\bf Problem:}(A) Find equilibria of the  following systems  $\left\{
\begin{array}{l}
{dx \over dt} =f(x,y)\\ 
{dy \over dt}=g(x,y)
\end{array}
\right.$ (see definition in section \ref{equilibrium_sec}).\\ (B) Find the following partial derivatives at each equilibrium point ($ {\partial f \over  \partial x}, { \partial f  \over  \partial y }, { \partial g   \over  \partial x},  {\partial g  \over  \partial y } $).\\
(a) $\left\{
\begin{array}{l}
{dx \over dt} =4x-2xy\\ 
{dy \over dt}=2xy-4y
\end{array}
\right.$

{\bf Solution:} (A) We rewrite the first equation as $4x-2xy=2x(2-y)=0$. It has solutions $x=0$ or $y=2$.\\ Let us substitute them to the second equation:\\
For  $x=0$, we get $2*0*y-4y=-4y=0$, thus $y=0$. Therefore we found one equilibrium $x=0,y=0$.\\ Now substitute  $y=2$. We get  $2*x*2-4*2=4x-8=0$, thus $x=2$. We found the second equilibrium  $x=2,y=2$.

(B) $f(x,y)=4x-2xy$, thus $ {\partial f \over  \partial x}=4-2y;  { \partial f  \over  \partial y }=-2x$.\\
 $g(x,y)=2xy-4y$, thus $ { \partial g   \over  \partial x}=2y;  {\partial g  \over  \partial y }=2x-4$.\\
At the equilibiria points:

Equilibrium  $x=0,y=0$:  $ {\partial f \over  \partial x}=4-2*0=4;  { \partial f  \over  \partial y }=-2*0=0;  { \partial g   \over  \partial x}=2*0=0;  {\partial g  \over  \partial y }=2*0-4=-4$.

Equilibrium  $x=2,y=2$:  $ {\partial f \over  \partial x}=4-2*2=0;  { \partial f  \over  \partial y }=-2*2=-4;  { \partial g   \over  \partial x}=2*2=4;  {\partial g  \over  \partial y }=2*2-4=0$.

\textheight=25.7cm
\textwidth=18.5cm
\topmargin=-2.0cm
\oddsidemargin=-0.9cm
\evensidemargin=-0.9cm
\renewcommand{\baselinestretch}{0.5}\normalsize
\twocolumn
\small 
\chapter{Answers for selected exercises and Formulas lists}
\begin{center} \subsection*{ Exercises Chapter 1}\end{center}
\ben

\item 
\ben

\item $3axy$
\item $ {6 \over r} -{5r \over 30r+5}={ -r^2+36r+6 \over r(6r+1)}$
\een

\item 
\ben
\item $ {\lim_{x \rightarrow \infty }} {a x + q \over c^2+x^2} =0 $
\item $ {\lim_{N \rightarrow \infty }} {a N^2 + q \over {b\over N}+ c^2+ dN^2}={a \over d}  $
\een

\item 
\ben
\item $3r + 2 - 5(r+1) = 6r+4$; $r=-{7 \over 8}$
\item $x + {4 \over x} =4$;
$x_{1,2}=2 \pm \sqrt{0}=2$
\item $(b-{N \over k})N=0$;
$N=0$, or $N=bk$.
\item 
$N=0$, or $b-d(1+{N \over k})=0$, thus $N={k(b-d) \over d}$.
\item $N=0$, or $N=h({b \over d}-1)$
\een

\item 
\ben  
\item  From  1st eq. $x=2y-5$, substitution to 2nd eq. gives $y=4$, thus $x=3$.
\item  from 1st eq.$x=-{ b \over a}y$ to 2nd eq. gives  $x=-{ b \over a}y=-{ b \over a}{-ba \over da-bc}={b^2 \over da-bc}$

\item From 2nd eq. $y(4-x)=0$, this $y=0$ or $x=4$. Substitution to the 1st eq. gives or $x=0$ and $x=0.5$. Now substitute  $x=4$, we get  $y=7$, thus all solutions are given: $(0,0),(0.5,0),(4,7)$.

\item From 2nd eq. $y(9-3x-y)=0$, thus $y=0$ or $y=9-3x$. 
Substitution to the 1st eq. gives:\\ $(0,0),(4,0),(0,9),(2.5,1.5)$.

\item  From 2nd eq. $N=0$, or $R= \delta $. Substitution to the 1st eq. gives: : $(0,0), (k(1-d),0), (\delta,{1-d \over a} - {\delta \over ak})$.
\een

\item 
\ben
\item $f'(x)={1 \over x^3}=-{3 \over x^4}$
\item 
$f(x)=e^{-5x}$ and $f'(x)=-5e^{-5x}$
\item $((4x-x^2)*(2x+3))'=10x+12-6x^2$
\item $y'={x^2+a^2 \over (a^2-x^2)^2}$
\een

\item 
\ben
\item $y$-intercept $y=-{b \over c^2}$, zeros $n= \pm \sqrt{b \over a}$, horizontal asymptote: $y=a$
\item  $y$-intercept $y={hr \over a}$, zeros $R=-h$ and $R=K$.
\een

\item 
\ben
\item fig.a
\item fig.b
\item horizontal asymptote $y=7$, intercept $y(0)=4$. (fig.c). Dependence on $a$. If $a$ increases, the horizontal asymptote $y=7$ will be approached at a slower rate. (fig.c)  

\begin{figure}[H]
\centerline{
\psfig{type=pdf,ext=.pdf,read=.pdf,figure=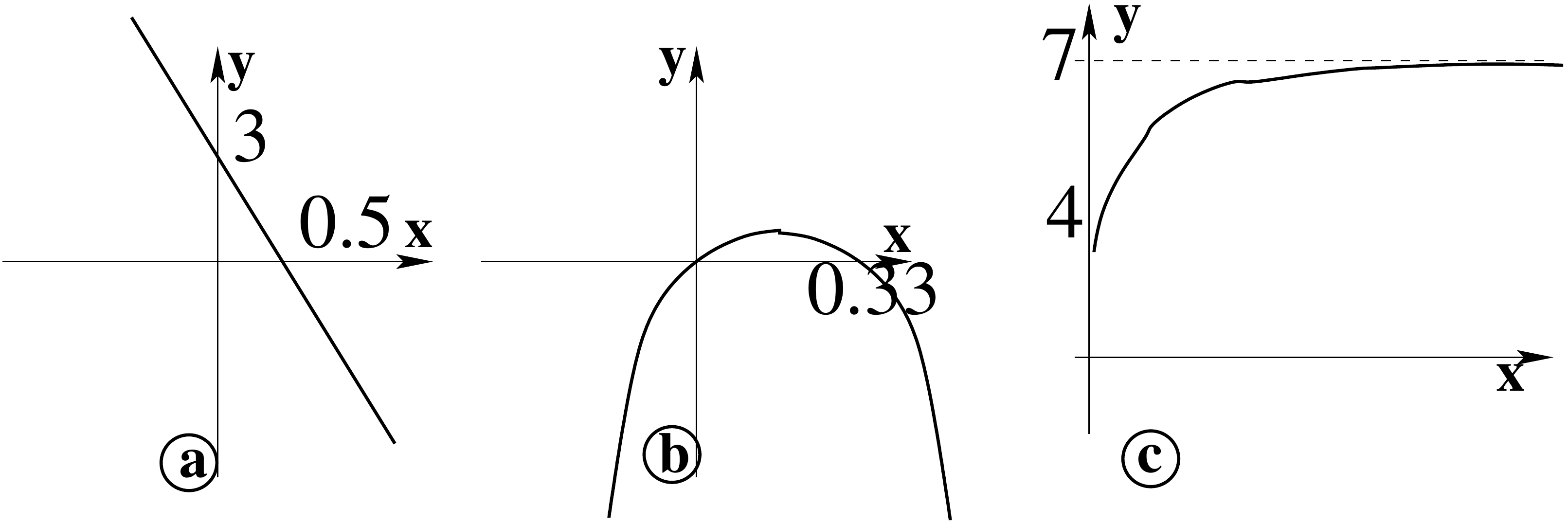,width=7cm}
}
\end{figure}

\een

\item 
\ben 
\item[(a)] (see below,fig.a)
\item[(b)] (see below,fig.b)
\begin{figure}[H]
\psfig{type=pdf,ext=.pdf,read=.pdf,figure=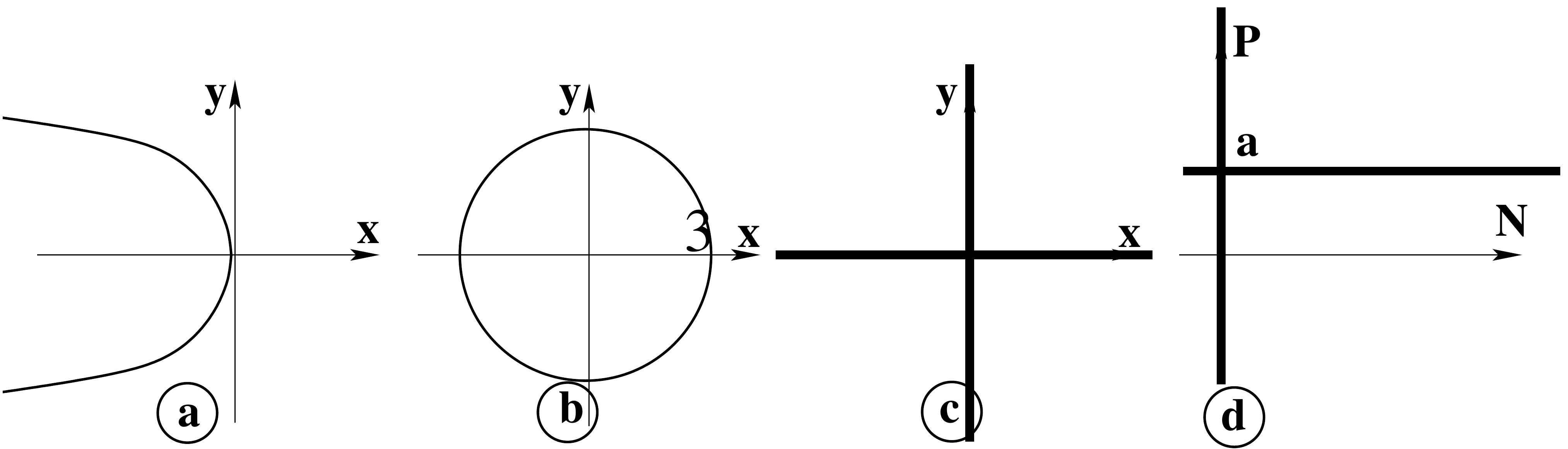,width=10cm}
\end{figure}

\item[(c)] (see above, fig.c)
\item[(d)]  (see above, fig.d)
\item[(e)] $N=0$ and  $P=-d(N+a)$ (see below, fig.a)
\item[(f)] $R=0$ and  ${d \over c}(b-R)(R+a)=N $ (see below, fig.b)
\begin{figure}[H]
\psfig{type=pdf,ext=.pdf,read=.pdf,figure=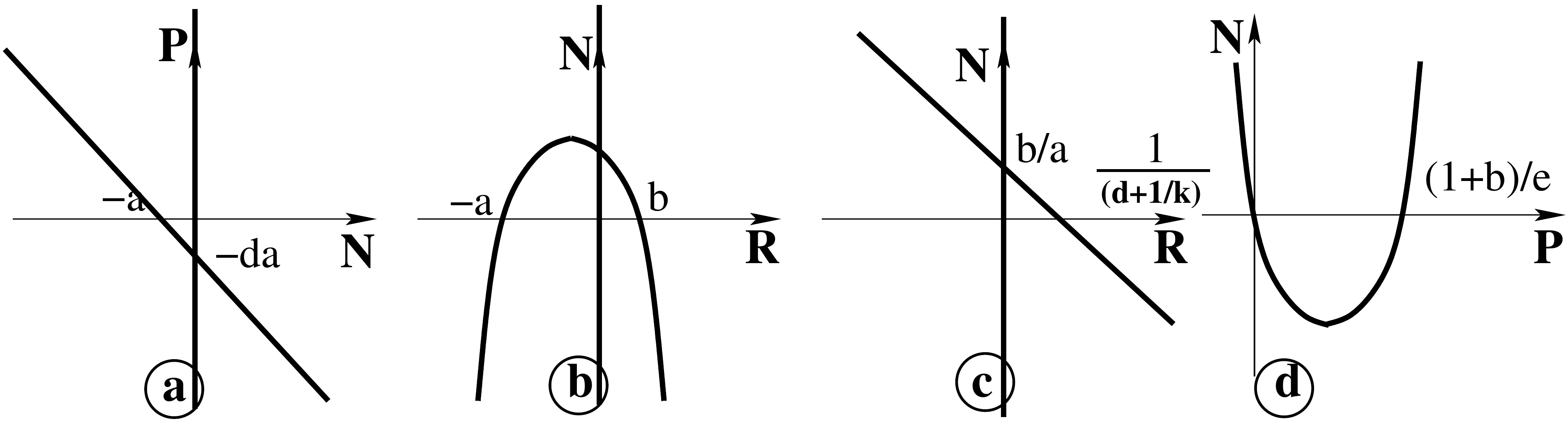,width=10cm}
\end{figure}

\item[(g)] $R=0$ and ${b \over a}(1-({1 \over k} +d)R) =N$ (see above,  fig.c).
\item[(h)] $N={P\over a}(eP-(1+b))$  (see above,  fig.d).
\een
\item
\ben
\item $((x-2y)*(y-2x)+2y^2)* {1 \over x}={5xy-2x^2 \over x}=5y-2x$
\item ${a -2 b \over 2p} : {4b - 2a \over \sqrt{p}}=-{ 1 \over 4 \sqrt{p}} $
\een

\item \no

\item
\ben
\item $2=e^{ln(2)}$, thus $(2^x)'=(e^{ln(2)x})'=ln(2)(e^{ln(2)x})=ln(2) 2^x$
\item  $\sqrt{1 \over x^3}=\sqrt{x^{-3}}=x^{-{3 \over 2}}$, thus $(x^{-{3 \over 2}})'=-{3 \over 2}x^{-{5 \over 2}}$
\item $f(x)=cos(x^2)$ and $f'(x)=-2xsin(x^2)$
\item $f(x)=cos^2(x)$ and $f'(x)=-2sin(x)cos(x)$
\item $(ax*e^{bx})'=ae^{bx}+abxe^{bx}$
\item  $y'={2x ( 2x^2-3x)-(4x-3)(x^2-5) \over ( 2x^2-3x)^2}=
{-3x^2+20x-15 \over ( 2x^2-3x)^2}
$
\item $y'={2ax(bx-c)-abx^2 \over (bx-c)^2}={abx^2 -2acx \over (bx-c)^2}$
\item $y'={(1+{x \over d})-{x \over d} \over  (1+{x \over d})^2}=
{1 \over  (1+{x \over d})^2}$
\item $y'={nx^{n-1}(x^n + a^n)-nx^{n-1}x^n \over (x^n + a^n)^2}=
{nx^{n-1}a^n \over (x^n + a^n)^2} $
\een

\item \no

\item
\ben
\item ${ df \over  dx}=3x^2$, ${df \over dt}=3x^2 {dx \over dt}$

\item  ${ df \over  dx}=-ae^{-ax}$, ${df \over dt}=-ae^{-ax}{dx \over dt}$

\item  ${df \over dt}={ df \over  dx}{dx \over dt}$

\een

\item 
\ben
\item $y$-intercept $y=-2$,zero $x=4$,horizontal asymptote $y=0$, as $y= {\lim_{x \rightarrow \infty }} {x -4\over x^2 -3x+2}={\lim_{x \rightarrow \infty }} {x \over x^2 }=0$, vertical asymptote is at $x^2 -3x+2=0$, i.e. $x=1$, and $x=2$. 

\item  $y= a:   {b \over x^3-c}= {a \over b}(x^3-c), $ thus $y$-intercept $y=-{ac \over b}$, one zero: $x=\sqrt[3]{c}$
\een
\item 
\ben
\item \no

\item $y=x^2+2x-3=(x-1)(x+3)$, as we have '+' at $x^2$ the parabola is opened upward.  (Graph fig.a)
\item horizontal asymptote $y=0$, vertical asymptote $x=-3$, no zeros (fig.b)

\begin{figure}[H]
\centerline{
\psfig{type=pdf,ext=.pdf,read=.pdf,figure=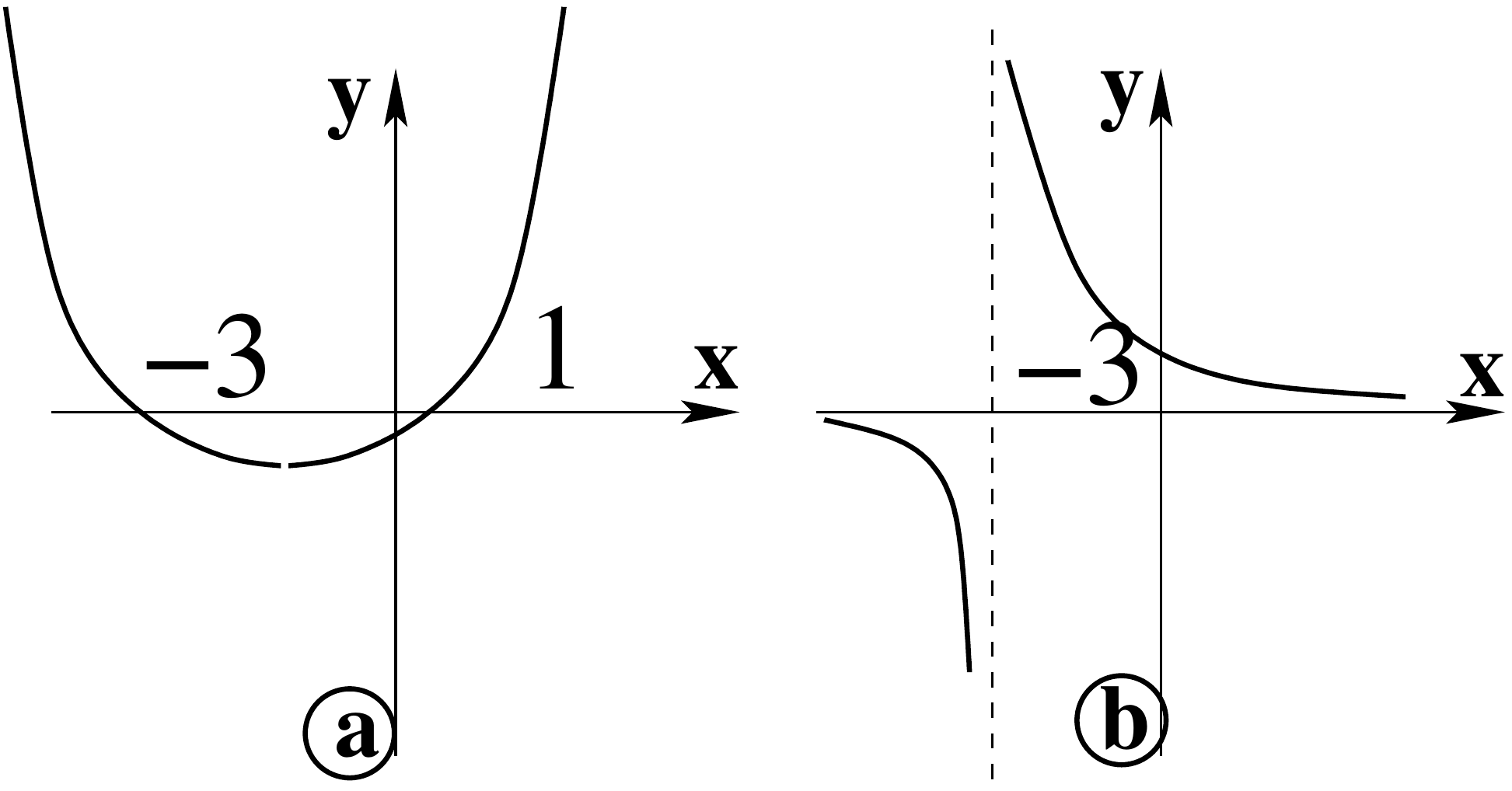,width=7cm}
}
\end{figure}

\item  (fig.a below). Parameter $b$ shifts the horizontal asymptote.

\item \no

\item  $h={kr \over 4}$.
\een
\begin{figure}[H]
\centerline{
\psfig{type=pdf,ext=.pdf,read=.pdf,figure=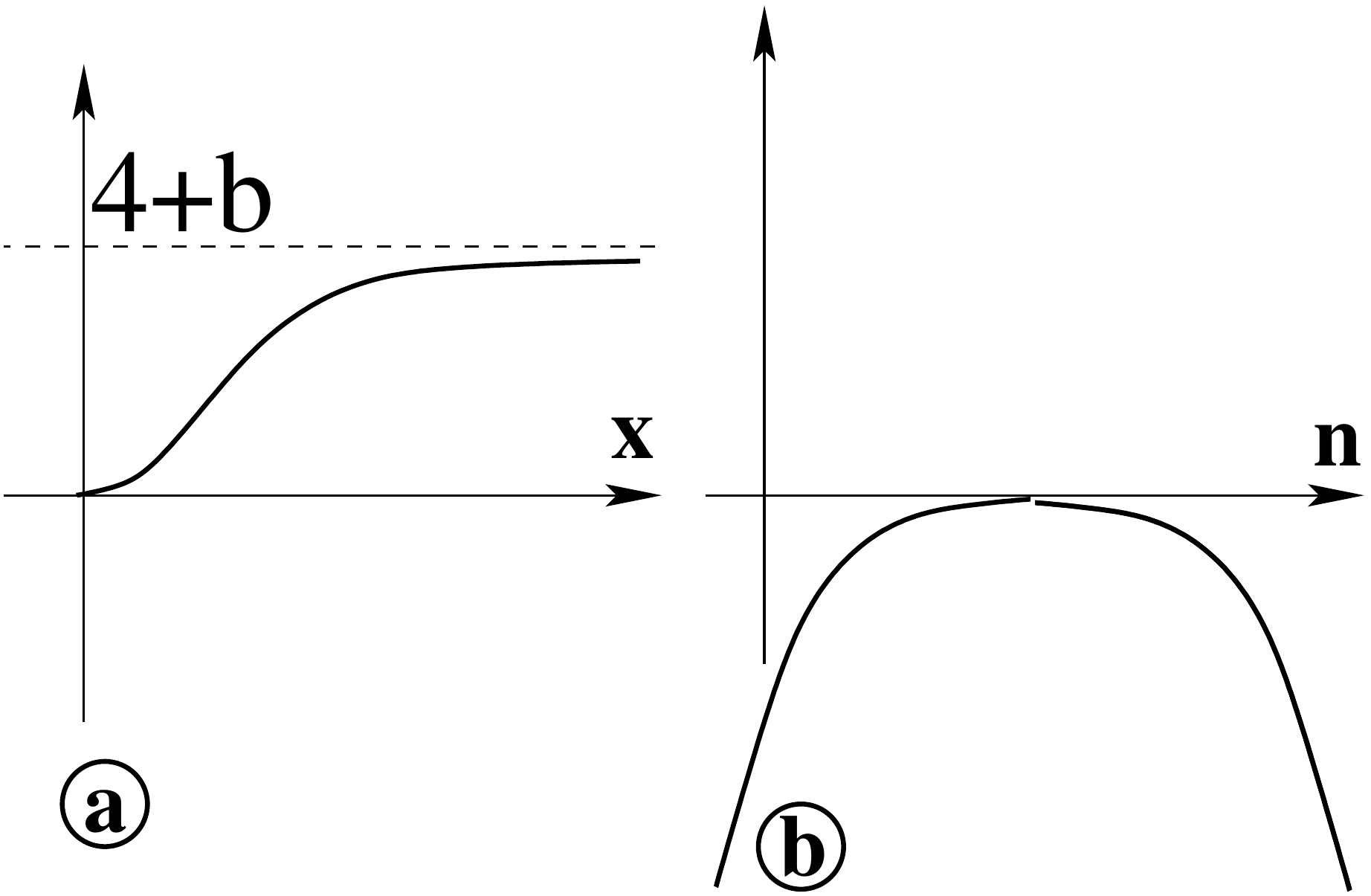,width=5cm}
}
\end{figure}

\item \no

\een

\begin{center} \subsection*{ Exercises Chapter 2} \end{center}

\ben
\item 
\ben
\item $x_{1,2}=-2 \pm i$.
\item $x_{1}=2, x_2=3$.
\een

\item 
\ben  
\item  $\left(\begin{array}{lr} 2 & -4\\ 1 & 1 \end{array} \right)
\left(\begin{array}{c}  x \\ y\end{array} \right)
=
\left(\begin{array}{c}  3\\ 1\end{array} \right)
$, $detA=2+4=6$.
\item $
 \left(\begin{array}{lr} a & b\\ c &  d \end{array} \right)
\left(\begin{array}{c}  x \\ y\end{array} \right)
=
\left(\begin{array}{c}  0\\ -b\end{array} \right)
$, $detA=ad-bc$.
\een

\item
\ben
\item $\lambda_1=-1$, $
\mathbf{v_1}=k\left( \begin{array}{c} -1 \\ -1 \end{array} \right)
$; $\lambda_1=-3$, $
\mathbf{v_2}=k\left( \begin{array}{c} -1 \\ 1 \end{array} \right)
$, where $k$ is an arbitrary number.
\item $\lambda_1=-1$, $
\mathbf{v_1}=k\left( \begin{array}{c} -4 \\ 2 \end{array} \right)
$; $\lambda_1=3$, $
\mathbf{v_2}=k\left( \begin{array}{c} -4 \\ -2 \end{array} \right)
$, where $k$ is an arbitrary number.
\item $\lambda_1=1+i$, $
\mathbf{v_1}=k\left( \begin{array}{c} -5 \\-2 -i \end{array} \right)
$; $\lambda_1=1-i$, $
\mathbf{v_2}=k\left( \begin{array}{c} -5 \\ -2+i \end{array} \right)
$, where $k$ is an arbitrary number.
\een

\item 
\ben 
\item ${\partial z \over \partial x}=2x$ at $(1,2)$ it is 2; ${\partial z \over \partial y}=2y$ at $(1,2)$ it is 4; 

\item  ${\partial z \over \partial x}=25-3x^2-y^2$, at $(3,4)$ it is -18;

\item  ${\partial z \over \partial N}=bR-d$, at $0,0$ it is $-d$; at $R={d \over b},N=1$ it is $0$;
${\partial z \over \partial R}=bN$,  at $0,0$ it is $0$; at $R={d \over b},N=1$ it is $b$.

\item ${\partial z \over \partial P}=-{a \over (1+P)^2}$, 
 ${\partial z \over \partial M}=-b$;

\item ${\partial z \over \partial N}=a-2eN-bP$, 
 ${\partial z \over \partial P}=-bN$;

\item ${\partial z \over \partial M}=L-{\nu A \over h+A}$, 
 ${\partial z \over \partial A}=-\delta  -{\nu Mh \over (h+A)^2}$;

\item ${\partial z \over \partial P_1}={2aP_1P_2  \over (h +P^2_1 + 2P_2)^2}$ and 
${\partial z \over \partial P_2}=-{a(h+P^2_1)  \over (h +P^2_1 + 2P_2)^2}$ 
\item  ${\partial z \over \partial N}=\frac{b^2NT(2+cN)}{(1+cN+bTN^2)^2} $
${\partial z \over \partial T}=\frac{b^2N^2(1+cN()}{(1+cN+bTN^2)^2} $
\een

\item
\ben
\item $\sqrt{3^2-90}=\sqrt{-81}= \pm 9 i$

\item $(-1+2i)+(4+7i)=3+9i$

\item $(4+5i)*(7+2i)=28+8i+35i+10i^2=18+43i$

\item ${1 \over i}={i \over i^2}=-i$

\een

\item
\ben
\item $x_{1,2}=\pm 11i$
\item $x_{1,2}=-1 \pm i \sqrt{2}$
\een

\item 
\begin{itemize}
\item $
A= \left(\begin{array}{lr} 1 & 2\\ 2 &  1 \end{array} \right)
$, $detA=-3$;
$
D_x= \left(\begin{array}{lr} 5 & 2\\ 4 &  1 \end{array} \right)
$, $detD_x=-3$ 
$
D_y= \left(\begin{array}{lr} 1 & 5\\ 2 &  4 \end{array} \right)
$, $detD_y=-6$, thus $x={-3 \over -3}=1; y={-3 \over -3}=2$.
\item By usual method: $x=5-2y$ thus after subs into 2nd eq. we get $10-4y+y=4$, or $y=2$ and hence $x=5-2*2=1$.
\end{itemize}
\item
\ben
\item $\lambda_1=2$, $
\mathbf{v_1}=k\left( \begin{array}{c} -6 \\ -3 \end{array} \right)
$; $\lambda_1=-5$, $
\mathbf{v_2}=k\left( \begin{array}{c} -6 \\ 4 \end{array} \right)
$, where $k$ is an arbitrary number.
\item $\lambda_1=3$, $
\mathbf{v_1}=k\left( \begin{array}{c} -1 \\ -1 \end{array} \right)
$; $\lambda_1=-5$, $
\mathbf{v_2}=k\left( \begin{array}{c} -1 \\ 7 \end{array} \right)
$, where $k$ is an arbitrary number.
\een

\item $f(x,y) \approx -2+2x+2y$

\een

\begin{center} \subsection*{Exercises Chapter 3} \end{center}

\ben
\item  at $t=4, n=30e^6 \approx 12102.86$, The double size at : $ t={ln(2) \over 1.5} \approx0.46$.  

\item $k={ln(2) \over 1200} \approx 0.58 \cdot 10^{-3} sec^{-1}$. 

\item 
\ben
\item  $ -15+ 8x -x^2=(3-x)(-5+x)$ Phase portrait in fig.a. (below), (zeros of parabola $x_1=3$ and $x_2=5$), attractor $x=5$, basin $x >3$. 

\item  Phase portrait in fig.a. (below),  attractor $x=4$, basin $x >1$. 
\item Phase portrait in fig.c. (below), attractor $x=-3$, basin $x <0$ and $x=2$, basin, $x>0$. 
\begin{figure}[H]
\psfig{type=pdf,ext=.pdf,read=.pdf,figure=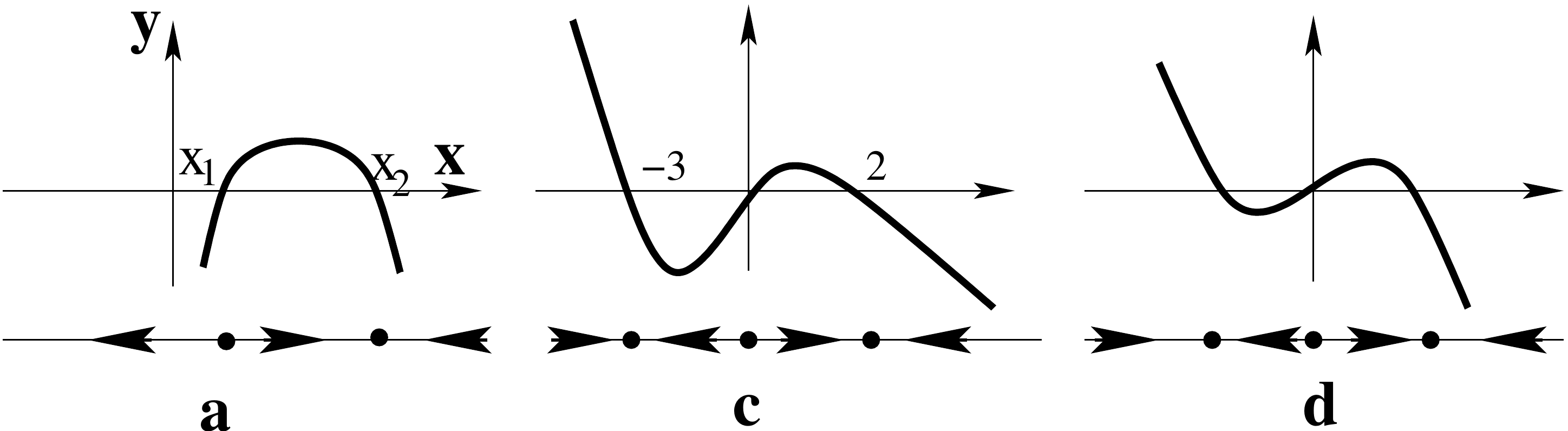,width=10cm}
\end{figure}

\item fig.d. (above), attractor $x=-2 \sqrt{2}$, basin $x <0$ and $x=2 \sqrt{2}$, basin, $x>0$. 

\item fig.a  (below) attractor $x=4$ basin $0 <x < 10$; fig.b  (below), attractor $x=0$, basin $x<1$ and attractor $x=6$, basin $x>1$.
\begin{figure}[H]
\psfig{type=pdf,ext=.pdf,read=.pdf,figure=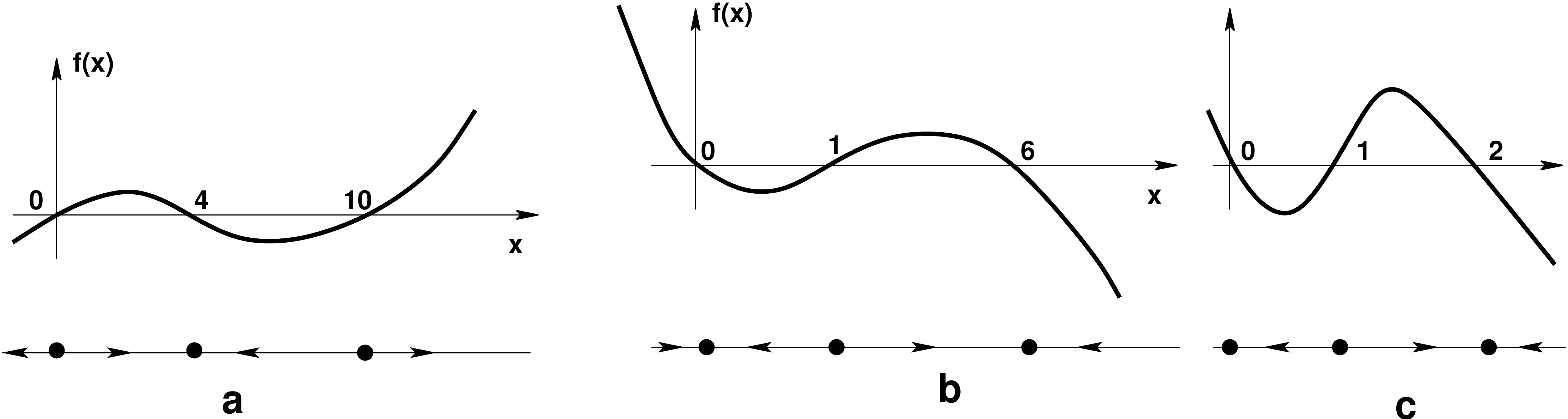,width=10cm}
\end{figure}

\item graph is shown in fig.c (above), attractor $x=0$, basin $x<1$ and attractor $x=2$, basin $x>1$.
\een

\item General consideration.
Typical graphs are shown in the figure. $s$ shifts the graph upward. If the total shift is less then the minimum of the graph, nothing will change. If the shift is more than the minimum (fig.b), the $x$ will go to the right equilibrium.

\begin{figure}[H]
\psfig{type=pdf,ext=.pdf,read=.pdf,figure=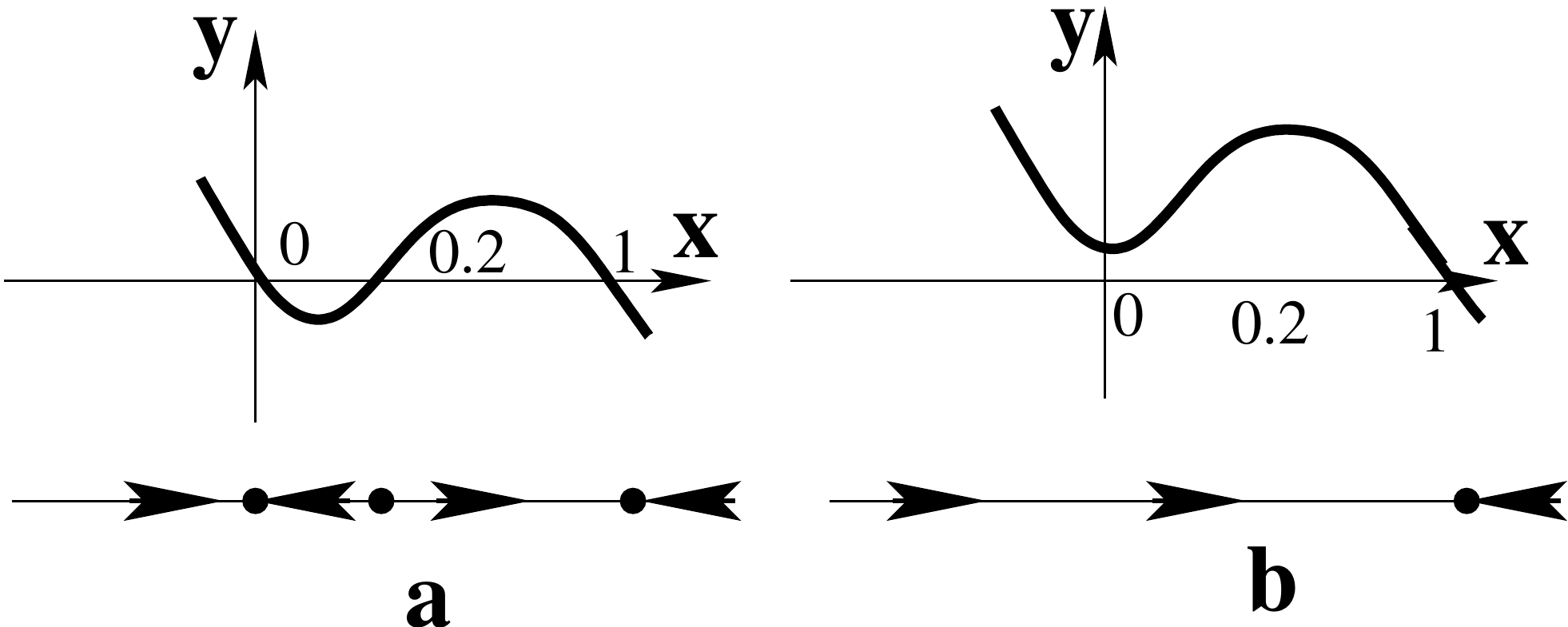,width=8cm}
\end{figure}

For questions (a) $f_{min} \approx -0.009$ (b) $x$ final is $x=0$, for question (c) $x$ final is $x=1$.
 (d) $s_{max}=0.009$

\item The maximal yield if $h$ equals maximum of $r*n*(1-n/k)$ (see the last section of this chapter), which gives $h_{max}=\frac{rk}{4}$

\item  (a) the general  solution $W=400/0.3+Ae^{-0.3t}$, 
$W=1333-1323e^{-0.3t}$;

 (b) $t=-ln(0.504)/0.3=2.29$; 

(c)   $t=-ln(0.9)/0.3=0.42$;

\item   the steady state value is  $m={\alpha \over\alpha+\beta}$ and the characteristic time $\tau={1 \over\alpha+\beta}$.

\item stable equilibrium at $n^*={k(r-h)\over r}$, yield is given by $yield=hn^*={kh(r-h)\over r}$, which has  maximal value $yield_{max}=\frac{kr}{4}$.

\item the  last strategy is better as population is more stable: small decrease in population size is OK for the last strategy, but for the other case small decrease  will result in the population extinction.

\een
\begin{center} \subsection*{Exercises chapter 4} \end{center}
\ben
\item
\ben
\item  $\left( \begin{array}{c} 
x \\  y \end{array} \right) = 
C_1 \left( \begin{array}{c} -1 \\  -1  \end{array} \right)e^{-t }
+ C_2 \left( \begin{array}{c} -1 \\  1 \end{array} \right)e^{-3* t }$

\item  $\left( \begin{array}{c} 
x \\  y \end{array} \right) = 
C_1 \left( \begin{array}{c} 1 \\  1  \end{array} \right)e^{2t }
+ C_2 \left( \begin{array}{c} 1 \\  -2 \end{array} \right)e^{5t }$

\een

\item 
$\left( \begin{array}{c} 
x \\  y \end{array} \right) = 
\frac{3}{2} \left( \begin{array}{c} 2 \\  -2  \end{array} \right)e^{3t }$

\item system is $\left\{
\begin{array}{l}
\frac{dC_1}{dt}=-0.01C_1+0.01C_2\\\frac{dC_2}{dt}= 0.04C_1-0.04C_2
\end{array}
\right.
$, the solution: 
 $\left( \begin{array}{c} 
C_1 \\  C_2 \end{array} \right) = 
-240 \left( \begin{array}{c} -0.01 \\  -0.01  \end{array} \right)
-60 \left( \begin{array}{c} -0.01 \\  0.04 \end{array} \right)e^{-0.05t }$, or $C_1=2.4+0.6e^{-0.05t }$, and $C_2=2.4-2.4e^{-0.05t }$.

\item
\ben
\item 
$det \left|\begin{array}{lr} 1-\lambda & 4\\ 2 & 3-\lambda
\end{array} \right|=(1-\lambda)( 3-\lambda)-8=\lambda^2-4\lambda-5=0$,
$\lambda_1=5$,$
\mathbf{v_1}=k\left( \begin{array}{c} -4 \\ -4 \end{array} \right)
$; $\lambda_2=-1$,$
\mathbf{v_2}=k\left( \begin{array}{c} -4 \\ 2 \end{array} \right)
$, thus this is a saddle point (the corresponding  phase portrait form the figure (a)  below). 
\item
characteristic  eq. is $\lambda^2-6\lambda+8=0$, $\lambda_1=2$,$
\mathbf{v_1}=k\left( \begin{array}{c} 1 \\ 3 \end{array} \right)
$; $\lambda_2=4$,$
\mathbf{v_2}=k\left( \begin{array}{c} 1 \\ 1 \end{array} \right)
$, i.e. a non-stable node (phase portrait see fig.b).

\item
characteristic  eq. is $\lambda^2-2\lambda+2=0$, $\lambda_1=1+i$, $\lambda_2=1-i$,(no eigen vectors required). A non-stable spiral  (qualitative phase portrait see fig.c).

\begin{figure}[H]
\centerline{
\psfig{type=pdf,ext=.pdf,read=.pdf,figure=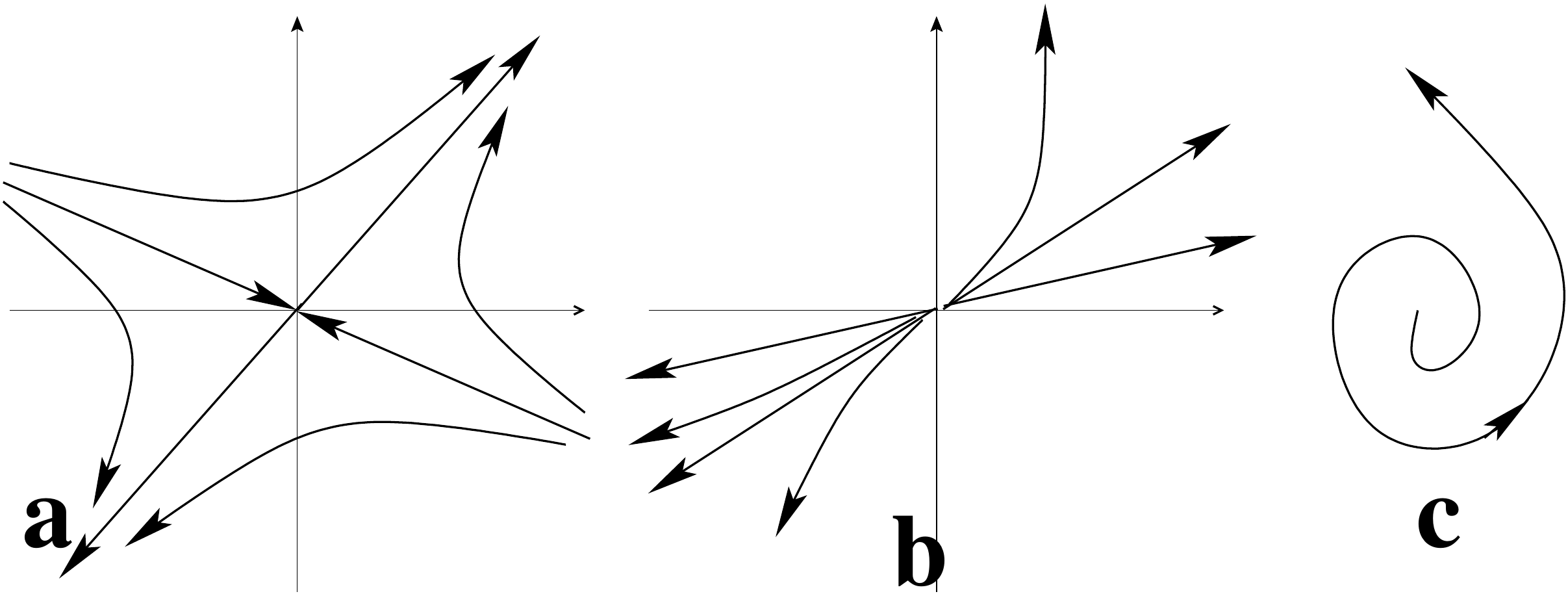,width=8cm}
}
\end{figure}

\item 
characteristic  eq. is $\lambda^2+4\lambda+3=0$, $\lambda_1=-3$,$
\mathbf{v_1}=k\left( \begin{array}{c} -1 \\ 1 \end{array} \right)
$, $\lambda_2=-1$,$
\mathbf{v21}=k\left( \begin{array}{c} -1 \\ -1 \end{array} \right)
$ i.e. a stable node  ( phase portrait see fig.a below).

\item 
characteristic  eq. is $\lambda^2+2\lambda+2=0$, $\lambda_1=-1+i$, $\lambda_2=-1-i$,(no eigen vectors required). A stable spiral  (qualitative phase portrait see fig.b).

\item 
characteristic  eq. is $\lambda^2+1=0$, $\lambda_1=+i$, $\lambda_2=-i$,(no eigen vectors required). A  center  (qualitative phase portrait see fig.c).

\item 
characteristic  eq. is $\lambda^2-1=0$, $\lambda_1=+1$,$
\mathbf{v_1}=k\left( \begin{array}{c} 1 \\ -3 \end{array} \right)
$; $\lambda_2=-1$,$
\mathbf{v_1}=k\left( \begin{array}{c} 1 \\ -1 \end{array} \right)
$, i.e. a  saddle  ( phase portrait see fig.d).

\begin{figure}[H]
\centerline{
\psfig{type=pdf,ext=.pdf,read=.pdf,figure=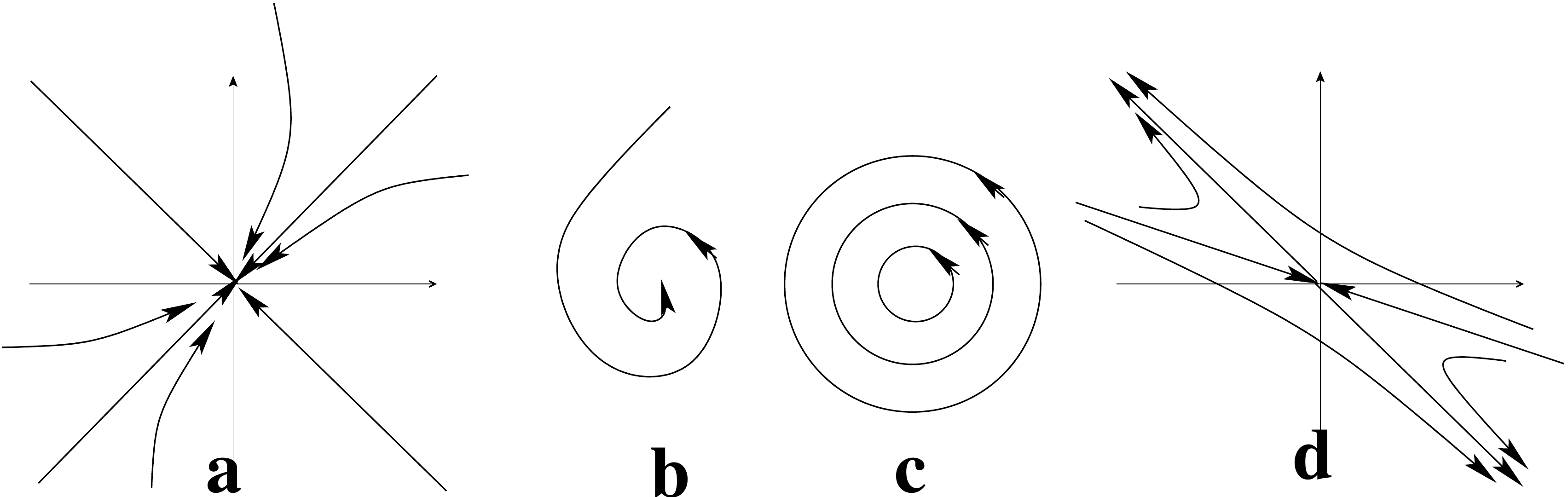,width=8cm}
}
\end{figure}

\een

\item Characteristic equation is given by\\ $det \left|\begin{array}{lr} -2-\lambda & -a\\ 3 & -1-\lambda
\end{array} \right|=(-2-\lambda)( -1-\lambda)+3a=\lambda^2+3\lambda+2+3a=0$,
$\lambda_{1,2}={ -3 \pm \sqrt{1-12a} \over 2}$.  If $1-12a<0$,
i.e. $a>{1 \over 12}$ we have complex roots. Because the real part for
complex roots is negative we have a stable spiral.  If the roots are
real ($a>{1 \over 12}$), then $\lambda_{2}={ -3 - \sqrt{1-12a} \over 2}$ is
always negative. The first root $\lambda_{1}={ -3 + \sqrt{1-12a} \over
2}$ will be positive is $-3 + \sqrt{1-12a} >0$, i.e. if the expression under the root
we have a value more than $9$ ( $\sqrt{9}=3$), this gives us $
1-12a>9$, or $a < {-8 \over 12}=-{2 \over 3}$. Thus if $a < -{2 \over
3}$, then $\lambda_{1}>0$ and we have a saddle point. In the interval $-{2
\over 3} < a < {1 \over 12}$ we have  $\lambda_{1}<0$, thus both  real roots negative, and we have  a stable node. Conclusion, if we increase $a$ from
$-\infty$, we will first have a  saddle point  (unstable eq.) till
$a=-{2 \over 3}$, then the saddle will become stable node (stable eq.),
until $a={1 \over 12}$ . If $a$ becomes  more than ${1 \over 12}$ we will
have a stable spiral (stable eq.). Qualitative phase portrait can be
plotted as in the previous exercise.

\item  In linear system oscillation occur if the equilibrium type is a
center. Characteristic
equation is given by $\lambda^2+(a+b)\lambda +ab+3-2a=0$. 
The center occurs if the roots  complex with a zero real part ($b=-a$).
These conditions give
$a^2+2a-3<0$, thus  oscillations occur, if $b=-a$ and $-3<a<1$.

\item (a) system is $\left\{
\begin{array}{l}
\frac{dx}{dt}=-(a+c)x+by \\\frac{dy}{dt}= ax-(b+e)y
\end{array}
\right.
$,\\ (b) system is $\left\{
\begin{array}{l}
\frac{dx}{dt}=-5x+2y \\\frac{dy}{dt}= 0.5x-5y
\end{array}
\right.
$, eigen values $\lambda_1=-4$, $\lambda_2=-6$, stable node, stable equilibrium. (c) formula for  eigen values in general case are:\\ $\lambda_{1,2}=\frac{-(a+b+c+e) \pm \sqrt{(a+b+c+e)^2-4*[(a+c)(b+e)-ab]}}{2}$. We see that  $(a+c)(b+e)-ab=ae+bc+ce>0$, thus expression under $\sqrt{}$ is less than $(a+b+c+e)^2$, thus we can have either stable node or stable spiral, both types of equilibria are stable.

\een

\begin{center} \subsection*{Exercises chapter 5} \end{center}
\ben
\item

(a) Equilibria:  $(0,0),(4,0)$.
$ {\partial f \over \partial x}, { \partial f \over \partial y }, {
\partial g \over \partial x}, {\partial g \over \partial y } $ are:\\
at  $(0,0)$:$0,-4,4,-0.5$;\\ at $(4,0)$: $0,-4,-4,-0.5$. 

\item[(b)]  Equilibria: $(0,0),(-9,-9)$.
Derivatives are:\\ at  $(0,0)$:$9,0,1,-1$;\\ at $(-9,-9)$: $9,-18,1,-1$. 
\item[(c)]   Equilibria:  $(0,0),(0.5,2).$
Derivatives are:\\ at $(0,0)$: $2,0,0,-1$; at $(0.5,2)$:
 $0,-0.5,4,1$. 
\item[(d)]  Equilibria: $(0,0), (0,0.5), (1,0), (0.25,0.25).$ 
Derivatives are: at $(0,0)$:$1,0,0,1$;\\ at $(0,0.5)$:
$-0.5,0,-1,-1$; at $(1,0)$:$-1,-3,0,-1$;\\ at $(0.25,0.25)$: $-0.25,-0.75,-0.5,-0.5$.  

\item  From 2nd $P=0$ or $N=\frac{d}{c}$, which after substitution to 1st equation  gives 3 equilibria $(0,0)$, $(\frac{a}{e},0)$,$(\frac{d}{c},\frac{a}{b}-\frac{ed}{bc})$. All non-negative if: $ac \geq ed$

\item From 2nd $M=\frac{d}{c}P$. Substitution to 1st gives one non-negative equilibrium: $P_1=\frac{-1+\sqrt{1+\frac{4ac}{bd}}}{2}$ and thus $M_1=\frac{d}{c}P_1$

\item From 2nd  $I=0$, or
$S={\alpha \over \beta }$. Substitution to 1st gives 
equilibria: ($I=0,S={B \over \mu}$) and ($I={B \over \alpha} -{\mu
\over \beta }, S={\alpha \over \beta }$.
The first equilibrium is always positive, the second one is positive
if $B \beta > \alpha \mu$.

\item 
\ben 
\item non-stable spiral, non-stable.
\item  saddle, non-stable.
 
\item  saddle, non-stable.
 
\item  center, neutrally stable.
 
\een

\item (a) 2 equilibria $(0,2)$, $(1,1)$ which are  stable node and  saddle; (d) see fig.a below. 
\begin{figure}[H]
\centerline{
\psfig{type=pdf,ext=.pdf,read=.pdf,figure=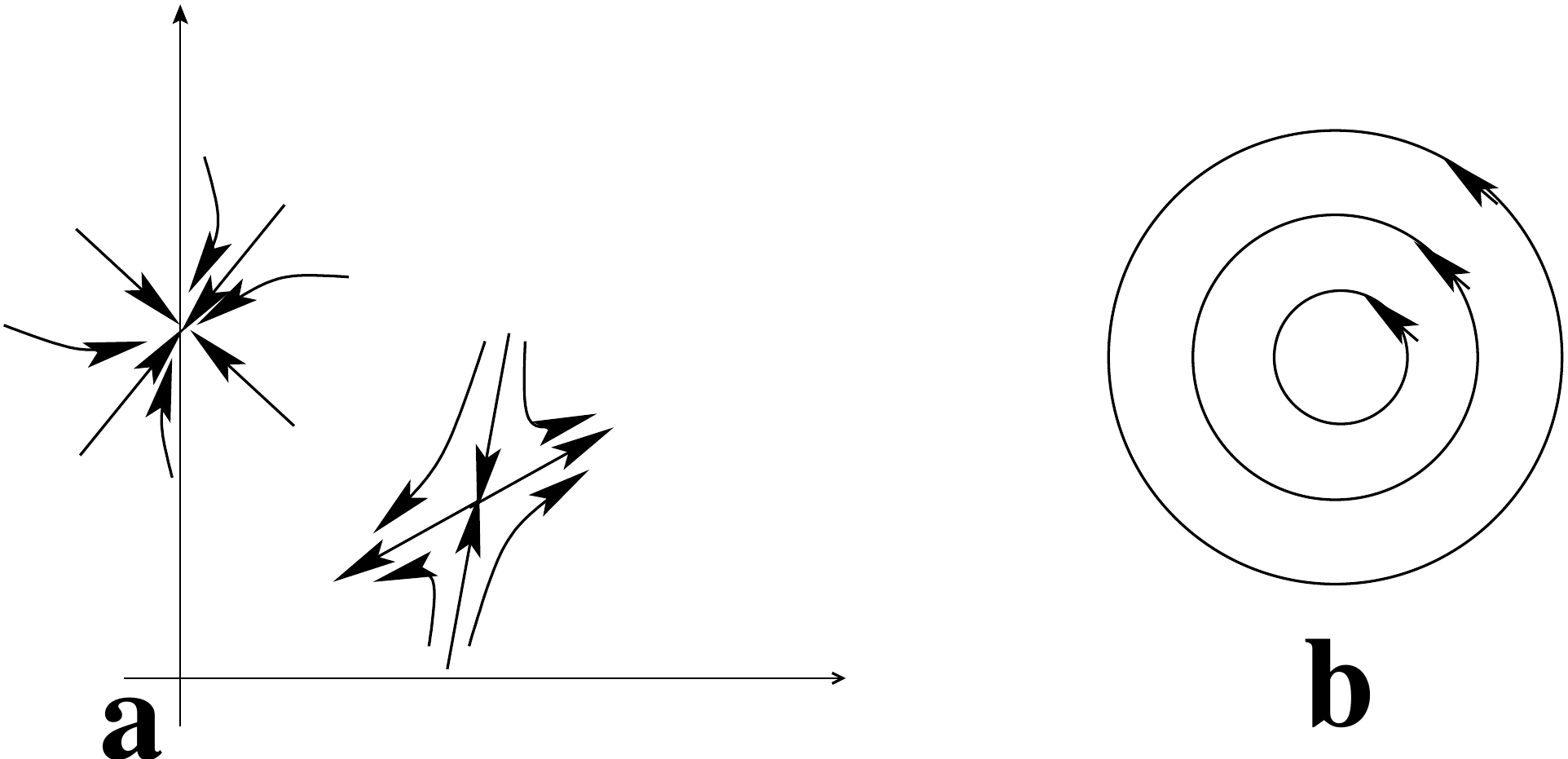,width=4cm}
}
\end{figure}

\item only one  equilibrium $(0,0)$. At this equilibrium, $detJ>0$,$trJ<0$, thus equilibrium is stable.

\item (a) $(\frac{d}{c},\frac{a}{b})$, linearization $\mathbf{\frac{dv}{dt}}=
J\mathbf{v}$, where $J=\left(\begin{array}{lr} 0 & -\frac{bd}{c}\\ \frac{ac}{b} & 0 \end{array} \right)$; (c) center point; (d) see fig.b above.

\item \no
\item \no
\item \no
\item \no

\item  Equilibria of system. From 2nd $e=0$, to 1st eq., $0-g=0$, i.e. equilibrium is $(0,0)$. Jacobian at the equilibrium is: ${\partial F \over \partial e}=
{\partial (-e^3+(1+a)e^2-ae-g) \over \partial e}= -3e^2+2(1+a)e-a$, at
$(0,0)$ it is $-a$, ${\partial F \over \partial g}=-1$, ${\partial G
\over \partial e}=\varepsilon$ ${\partial G \over \partial g}=0 $, thus $J=\left(\begin{array}{lr} -a
& -1\\ \varepsilon & 0 \end{array} \right)$ and we get
$detJ=\varepsilon>0, trJ=-a, D=a^2-4\varepsilon$, thus equilibrium is
always stable, and it is a node if $D>0$, i.e. $a^2>4\varepsilon$, and
a stable spiral if $a^2<4\varepsilon$.
\een

\begin{center} \subsection*{Exercises chapter 6} \end{center}

\ben
\item  (a) Null-clines are given below in fig.a, thus the graphical Jacobian on
basis of 'black' points is: $J=\left(\begin{array}{lr} \alpha & \beta
\\ -\gamma & -\delta \end{array} \right)$. This gives for
$detJ=-\alpha \delta + \beta\gamma$, and we do not know what is the
sign. Thus graphical Jacobian does not work here. The real Jacobian
here is: $J=\left(\begin{array}{lr} 3 & 1\\ -1 & -1 \end{array}
\right)$, that gives $detJ=-2<0$, thus we have a saddle point. Phase
portrait is in fig.b. 
\begin{figure}[H]
\centerline{
\psfig{type=pdf,ext=.pdf,read=.pdf,figure=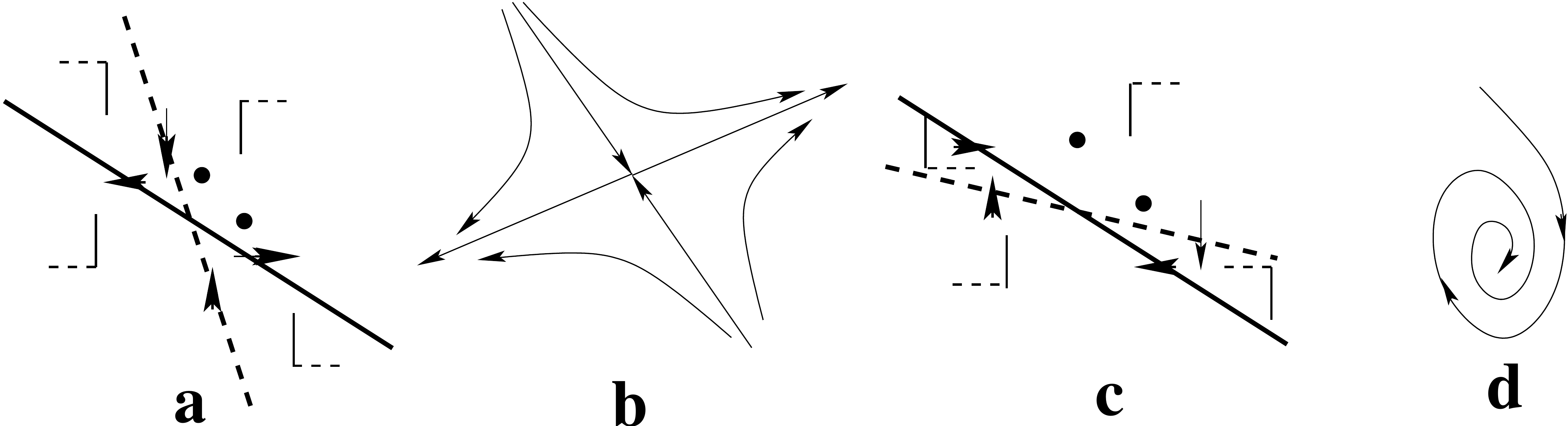,width=9cm}
}
\end{figure}

\item[(b)]  Null-clines are given above in fig.c, Graphical Jacobian does not work here. The real Jacobian gives  stable spiral. Phase
portrait is in fig.d. 

\item[(c)]  Null-clines are given below in fig.a, the graphical Jacobian gives  saddle point. Phase portrait  in fig.b.

\item[(d)]   
Null-clines are given below in fig.c, the graphical Jacobian gives  a non-stable equilibrium (node, spiral).  The real Jacobian gives a  non-stable spiral. Phase
portrait fig.d.
\begin{figure}[H]
\centerline{
\psfig{type=pdf,ext=.pdf,read=.pdf,figure=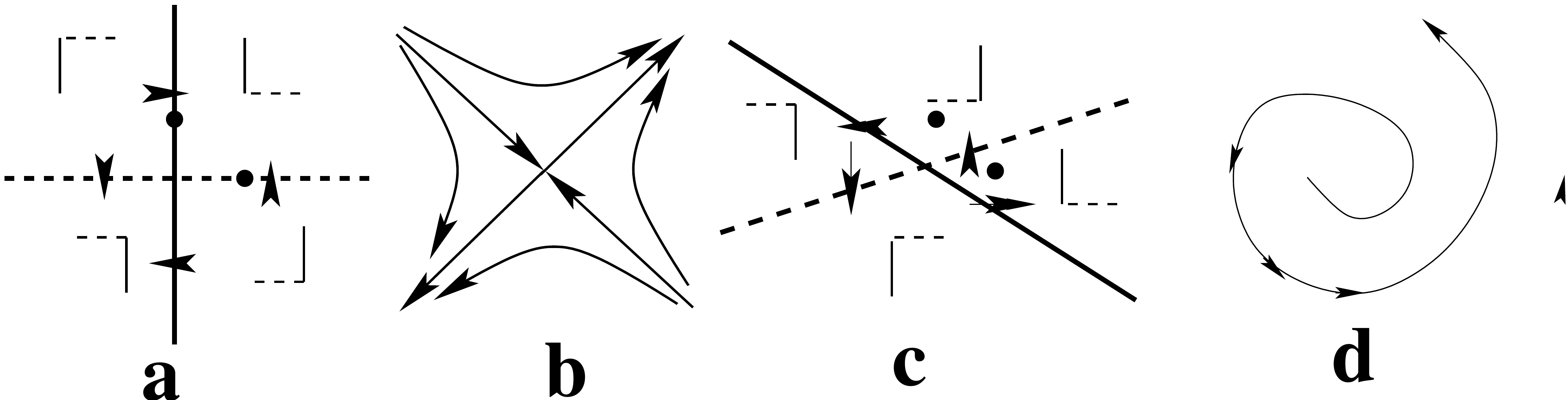,width=9cm}
}
\end{figure}

\item For fig.a below we have two equilibria (marked). For the left one the
graphical Jacobian gives a  stable equilibrium (stable node or stable spiral). The right equilibrium gives  a saddle.

\begin{figure}[H]
\centerline{
\psfig{type=pdf,ext=.pdf,read=.pdf,figure=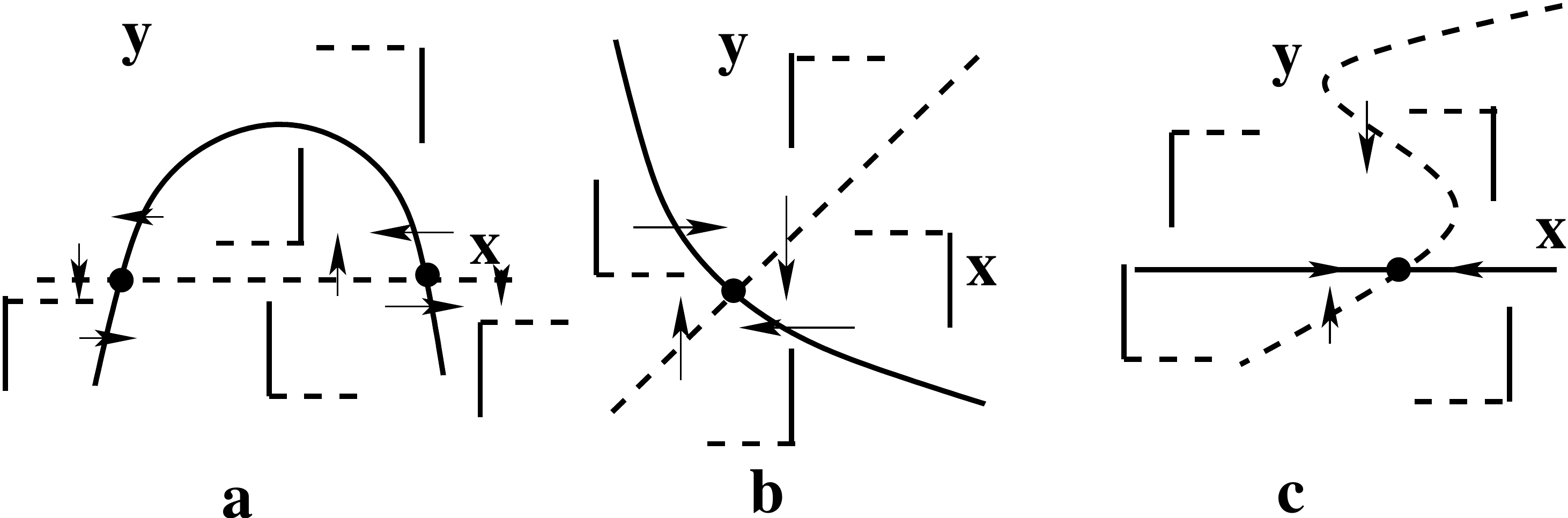,width=8.cm}
}
\end{figure}

For fig.b we have one equilibrium (marked). The
graphical Jacobian gives  a stable equilibrium (stable node or stable spiral).

For fig.c we have one equilibrium (marked). The
graphical Jacobian gives  a stable node.

\item (a) See fig.a below. (b,c) equilibria: $(0,2)$,graphical Jacobian gives a  stable node; $(1,1)$: saddle. For phase portrait see solution problem 6 chapter 5. 
\begin{figure}[H]
\centerline{
\psfig{type=pdf,ext=.pdf,read=.pdf,figure=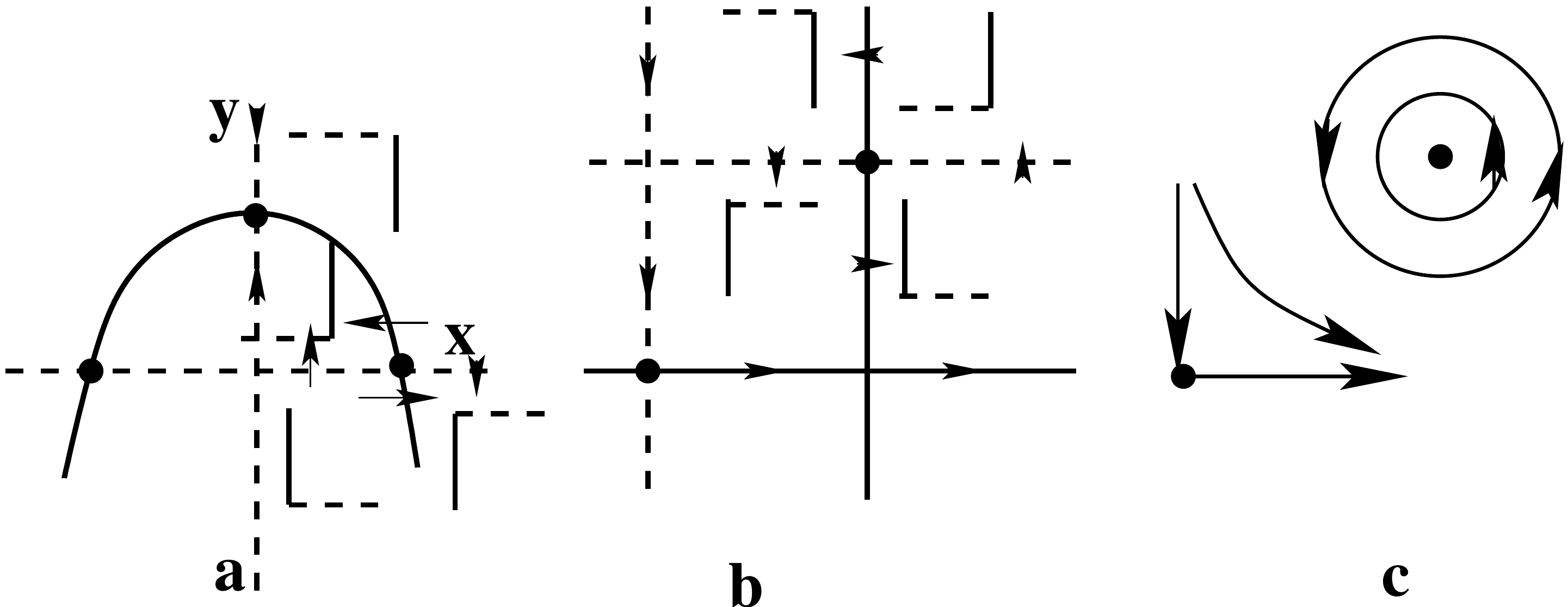,width=8.cm}
}
\end{figure}
 
\item (a) See fig.b above. (b,c)  equilibria: $(0,0)$,graphical Jacobian: saddle;
 2nd equilibrium $(\frac{d}{c},\frac{a}{b})$, graphical Jacobian: center. Phase portrait fig.c above. (d) no, for positive parameters.

\item (a), (b) \no
\item[(c)]vector field see fig.a; Equilibria $(0,0)$,graphical Jacobian: $J=\left(\begin{array}{lr} \alpha &  \beta  \\
\gamma & -\delta \end{array} \right)$, $detJ=-\alpha \delta -beta \gamma<0$,saddle;
 2nd equilibrium $(-9,-9)$,graphical Jacobian: $J=\left(\begin{array}{lr} \alpha &  -\beta  \\
\gamma & -\delta \end{array} \right)$, $detJ=-\alpha \delta +\beta \gamma$, sign unknown, $trJ=\alpha - \delta$, equilibrium type unknown; phase portrait fig.b.
\begin{figure}[H]
\centerline{
\psfig{type=pdf,ext=.pdf,read=.pdf,figure=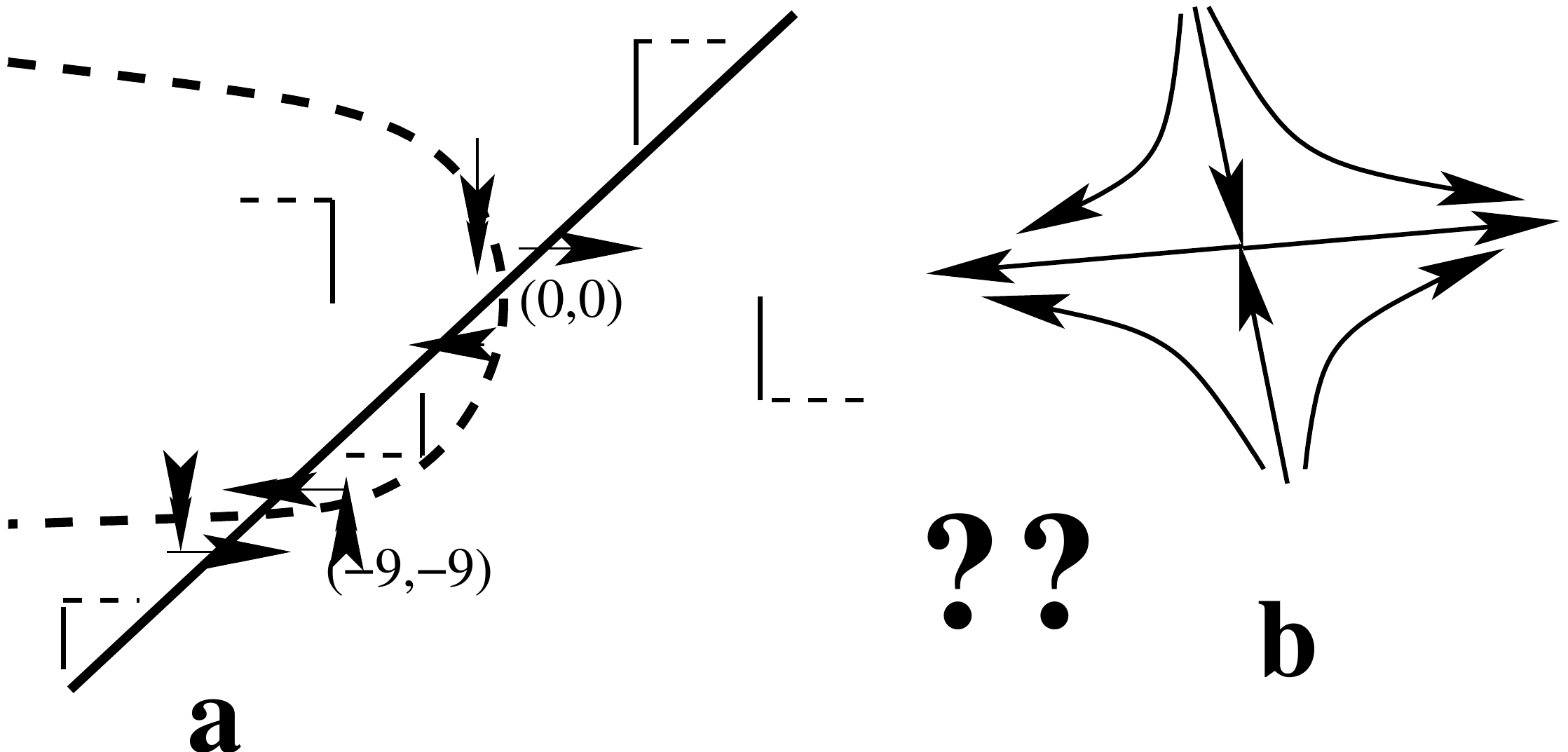,width=5cm}
}
\end{figure}

\item[(d)] \no

\een

\begin{center} \subsection*{Exercises chapter 7} \end{center}
\begin{figure}[H]
\centerline{
\psfig{type=pdf,ext=.pdf,read=.pdf,figure=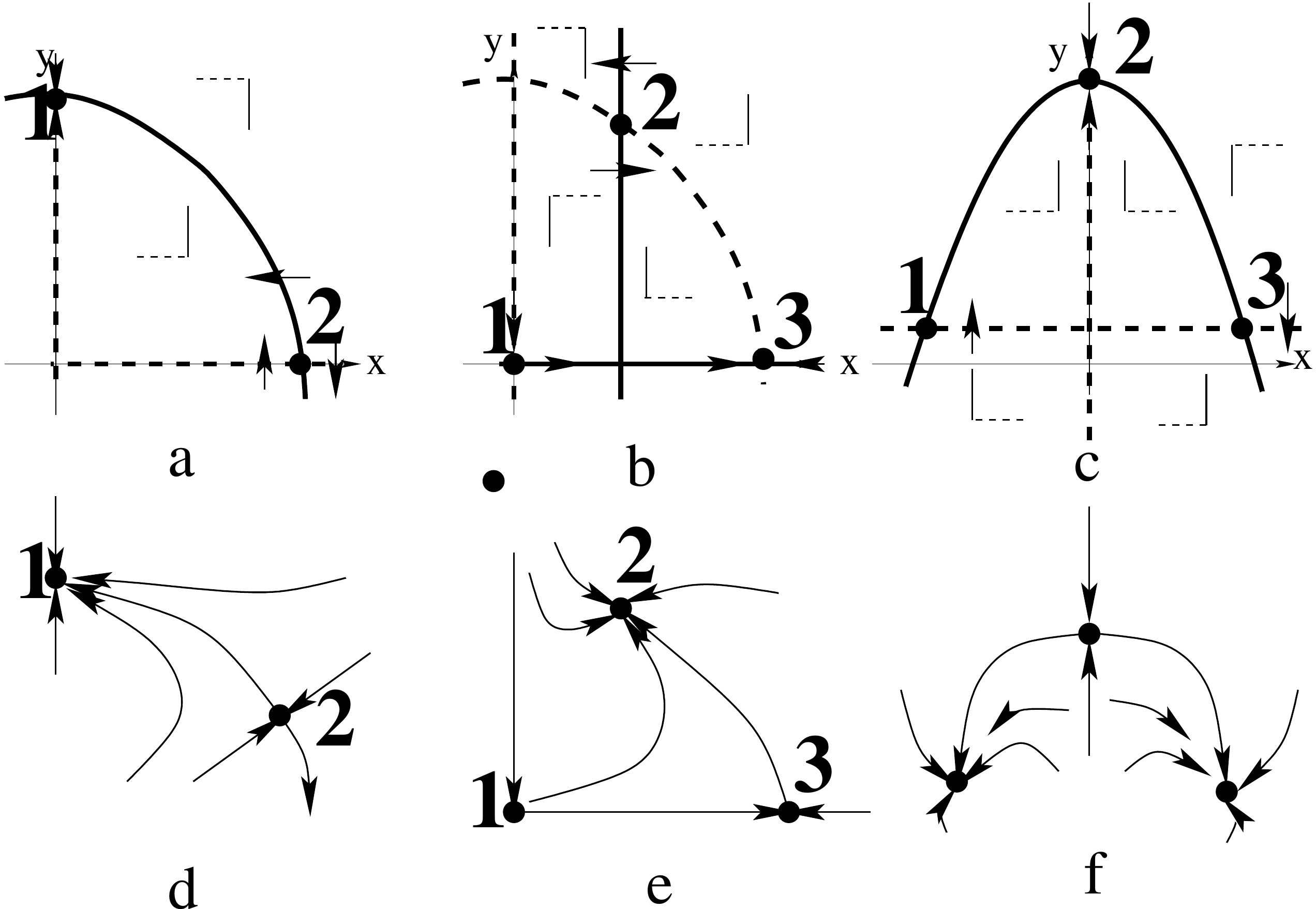,width=9cm}
}
\end{figure}
\ben
\item (a)  graphical Jacobian: at {\bf 1} stable node. At {\bf 2} saddle.
Null-clines fig.a above, phase portrait fig.d. 

\item[(b)] At {\bf 1} saddle. At {\bf 2}  not known, stable node/spiral.
At {\bf 3}  saddle. Null-clines fig.b above, phase portrait fig.e. 

\item[(c)] At {\bf 1}  not known, stable node/spiral.
At {\bf 2}  saddle. At {\bf 3}  not known, stable node/spiral.  Null-clines fig.c above, phase portrait fig.f. 

\item (a) Null-clines $x$: $x=0,y=1-x$, $y$: $y=0,x={1 \over 2}$, fig.a below.\\
Equilibria $(0,0),(1,0),({1 \over 2},{1 \over 2})$. Graphical Jacobian:
$(0,0)$
 unstable node. $(1,0)$:  stable node. $({1 \over 2},{1 \over 2})$ saddle. Phase portrait fig.b. One attractor $(1,0)$. Basin of attraction shaded.

\begin{figure}[H]
\centerline{
\psfig{type=pdf,ext=.pdf,read=.pdf,figure=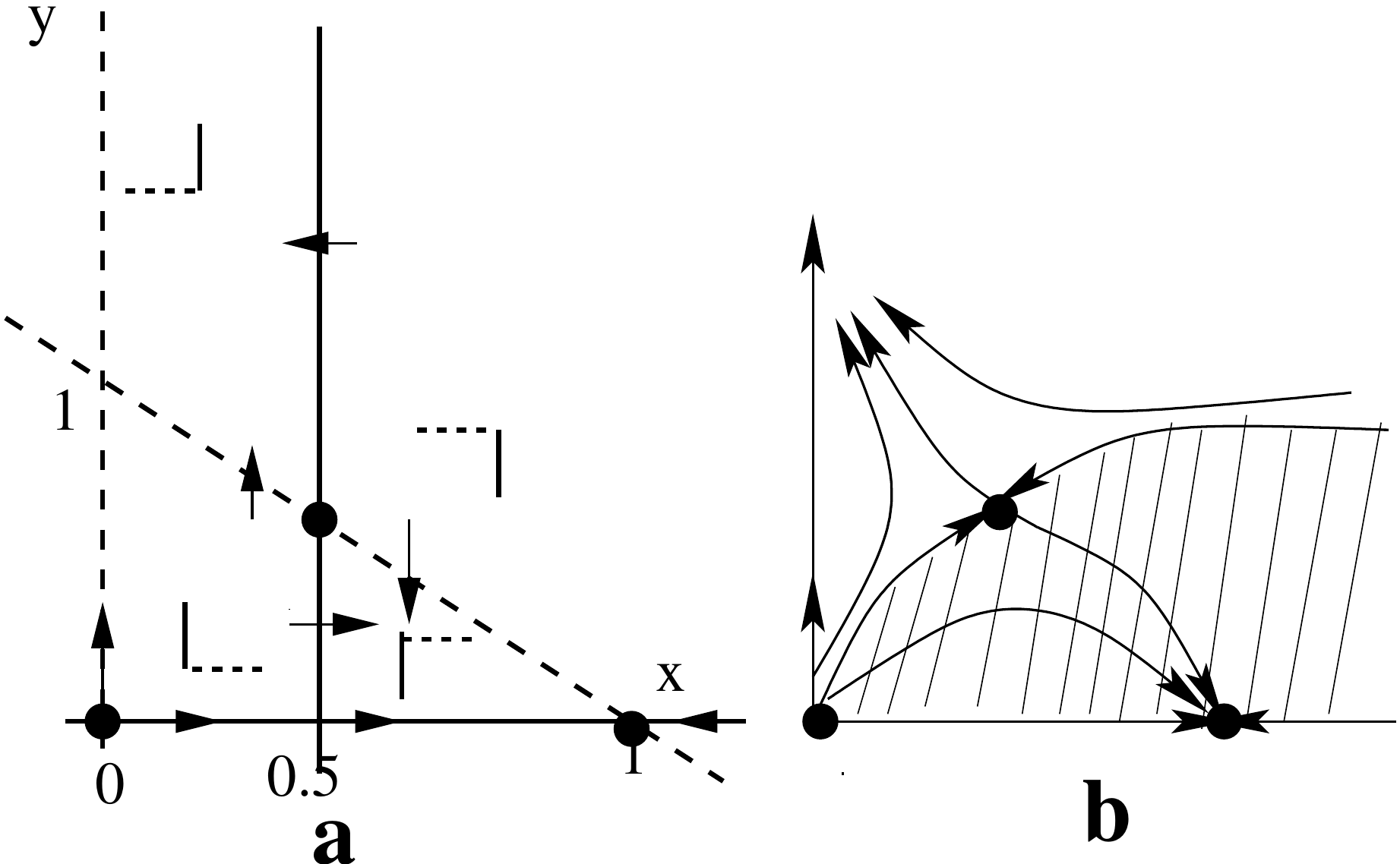,width=4cm}
\psfig{type=pdf,ext=.pdf,read=.pdf,figure=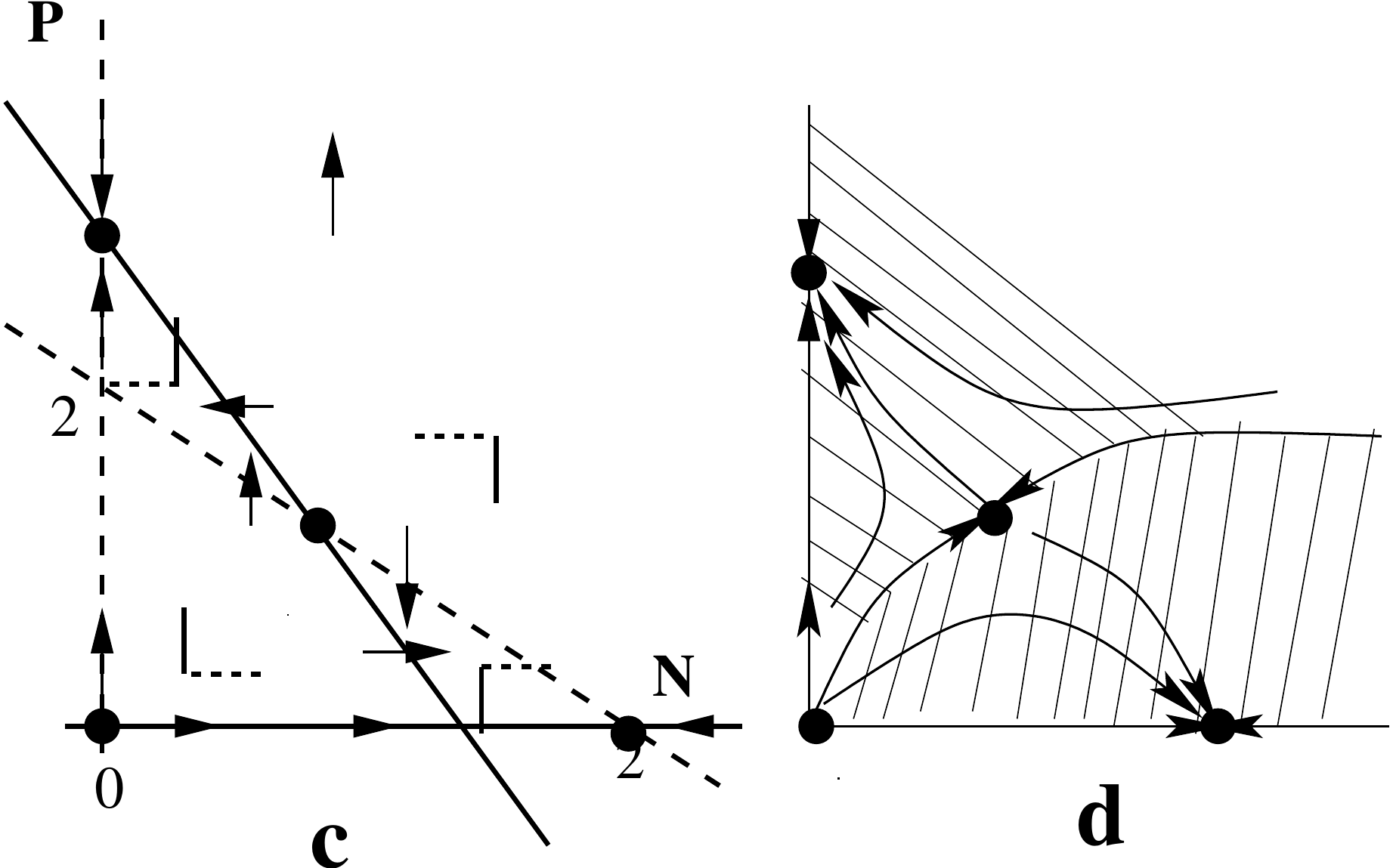,width=4cm}
}
\end{figure}

\item[(b)] Null-clines $N$: $N=0,P=2-N$, $P$: $P=0,P=3-2N$, fig.c above \\
Equilibria $(0,0),(2,0),(0,3),(1,1)$. Graphical Jacobian:
$(0,0)$ , unstable node; $(2,0)$ , stable node; $(0,3)$, stable node:
 $ (1,1)$ cannot determine equilibrium type using graphical Jacobian. Need 'real' Jacobian, which gives  saddle. Phase portrait fig.d. Two attractors $(2,0)$ and $(0,3)$, basins of attraction shaded.

\item Null-clines $P$: $M=\frac{d_P P}{b}$, $M$: $M=\frac{abK^2}{d_M(K^2+P^2)}$, fig.a below \\
One equilibrium. Graphical Jacobian gives stable equilibrium node/spiral. Phase portrait fig.b.  One attractor {\bf 1}. Basin of attraction shaded.

\begin{figure}[H]
\centerline{
\psfig{type=pdf,ext=.pdf,read=.pdf,figure=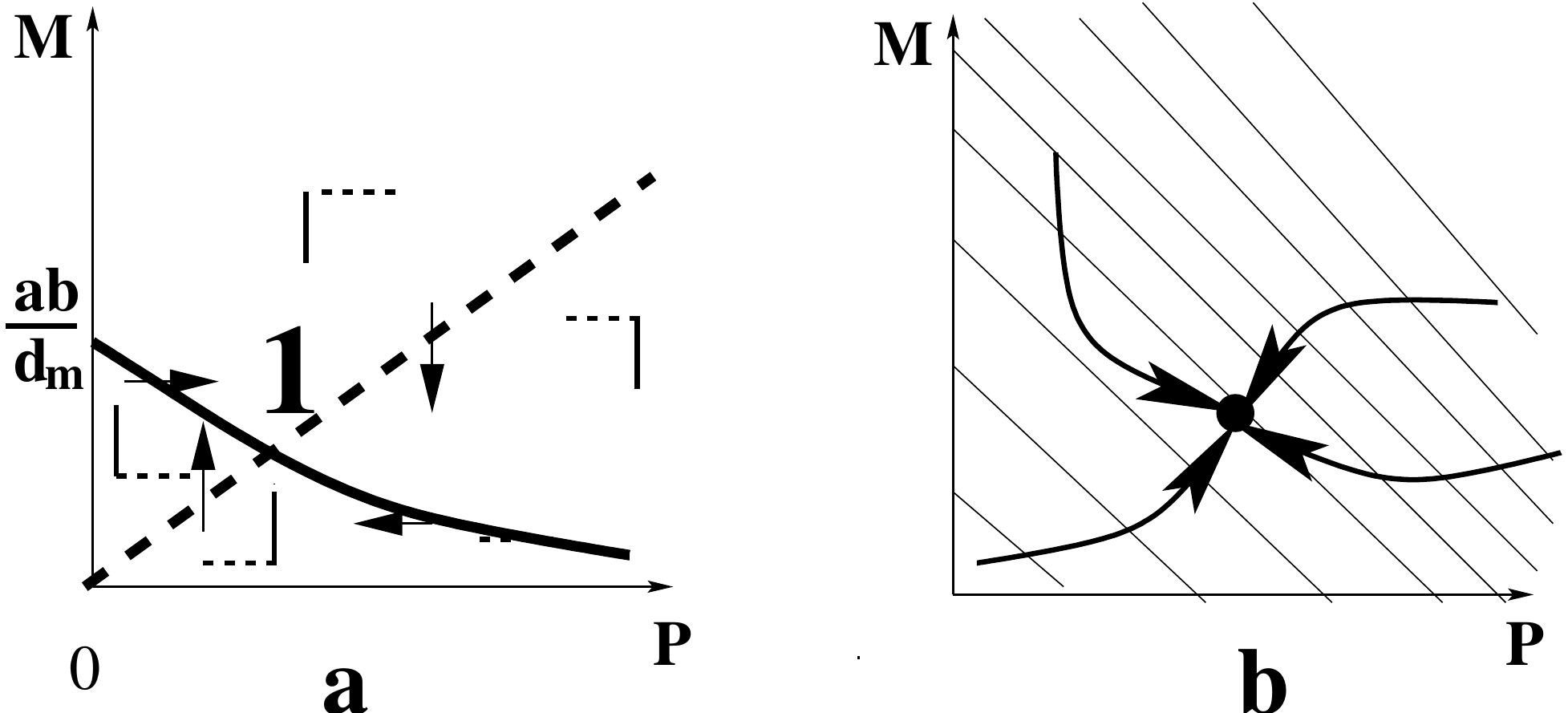,width=4cm}
\psfig{type=pdf,ext=.pdf,read=.pdf,figure=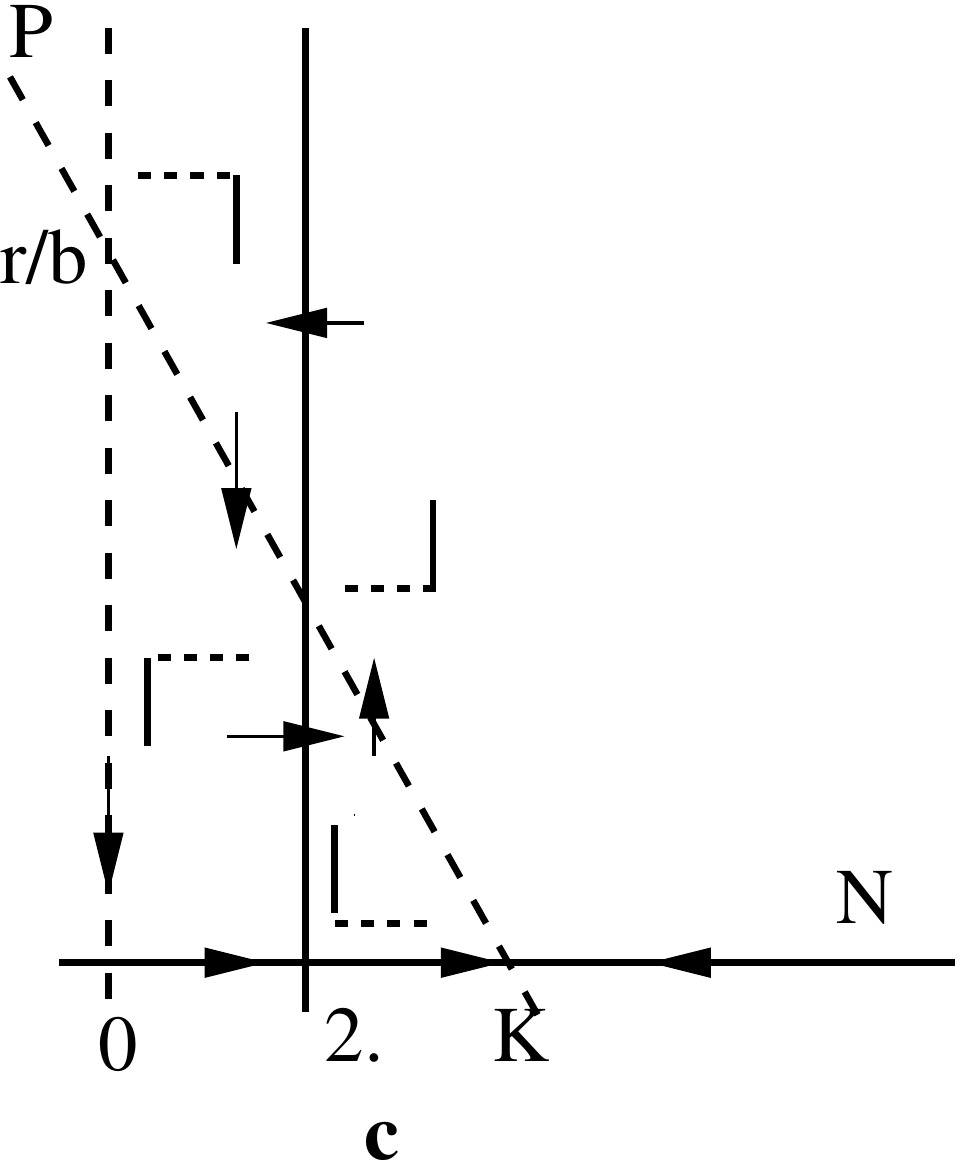,width=2.cm}
\psfig{type=pdf,ext=.pdf,read=.pdf,figure=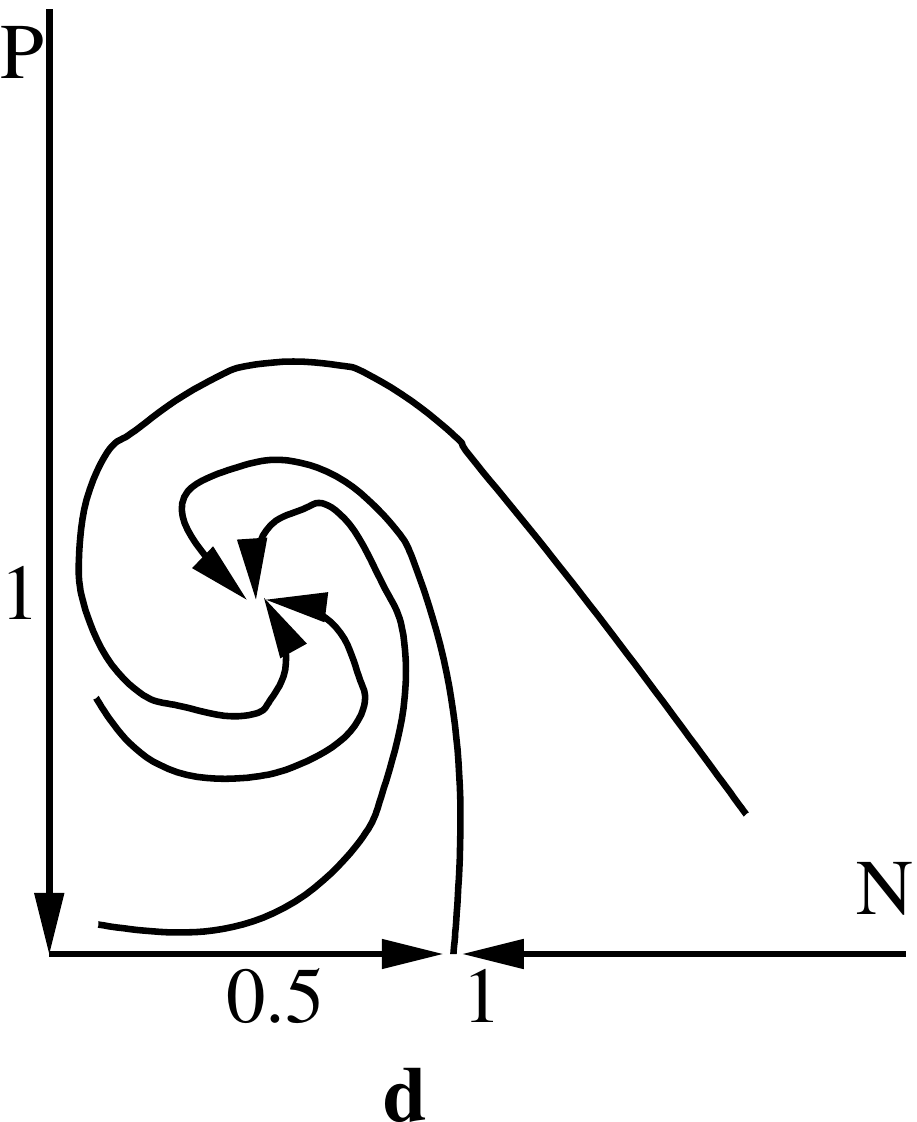,width=2.cm}
}
\end{figure}

\item (a)Equilibrium $(K,0)$; (b) $J=\left(\begin{array}{lr} r-\frac{2rN}{K} -bP &  -bN \\ bP & bN-2b \end{array} \right)=
\left(\begin{array}{lr} -r &  -bK \\ 0 & bK-2b \end{array} \right)$;
(c) Equilibrium is stable if $0<K<2$. 

\item (a)-(d) \no
\item Fig c,d above
\een
\begin{figure}[H]
\end{figure}

\begin{center} \subsection*{Exercises Chapter 8} \end{center}
\ben

\item see   figure below
\begin{figure}[hhh]
\psfig{type=pdf,ext=.pdf,read=.pdf,figure=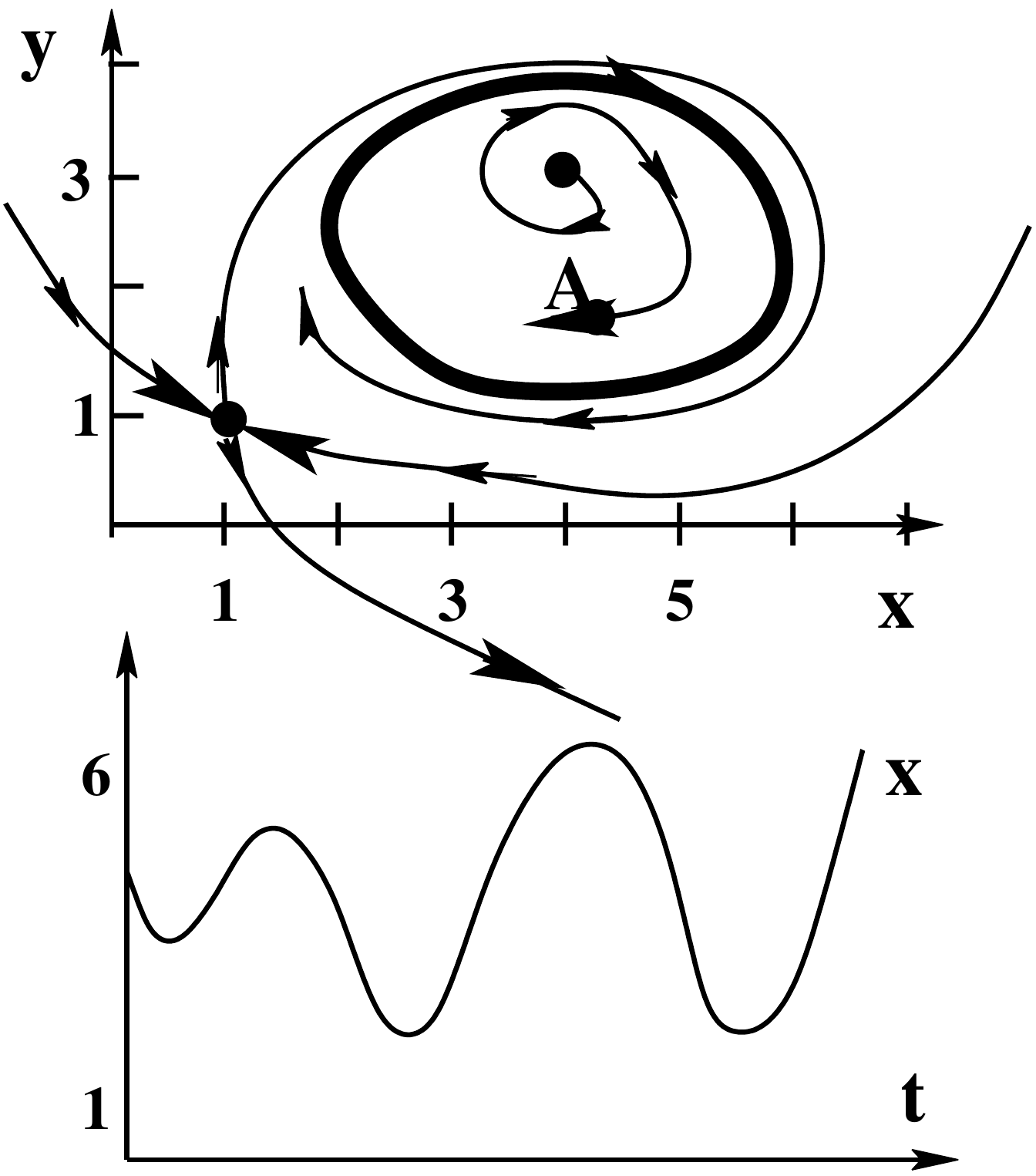,width=5cm}
\end{figure}
\item see figure below

\begin{figure}[hhh]
\psfig{type=pdf,ext=.pdf,read=.pdf,figure=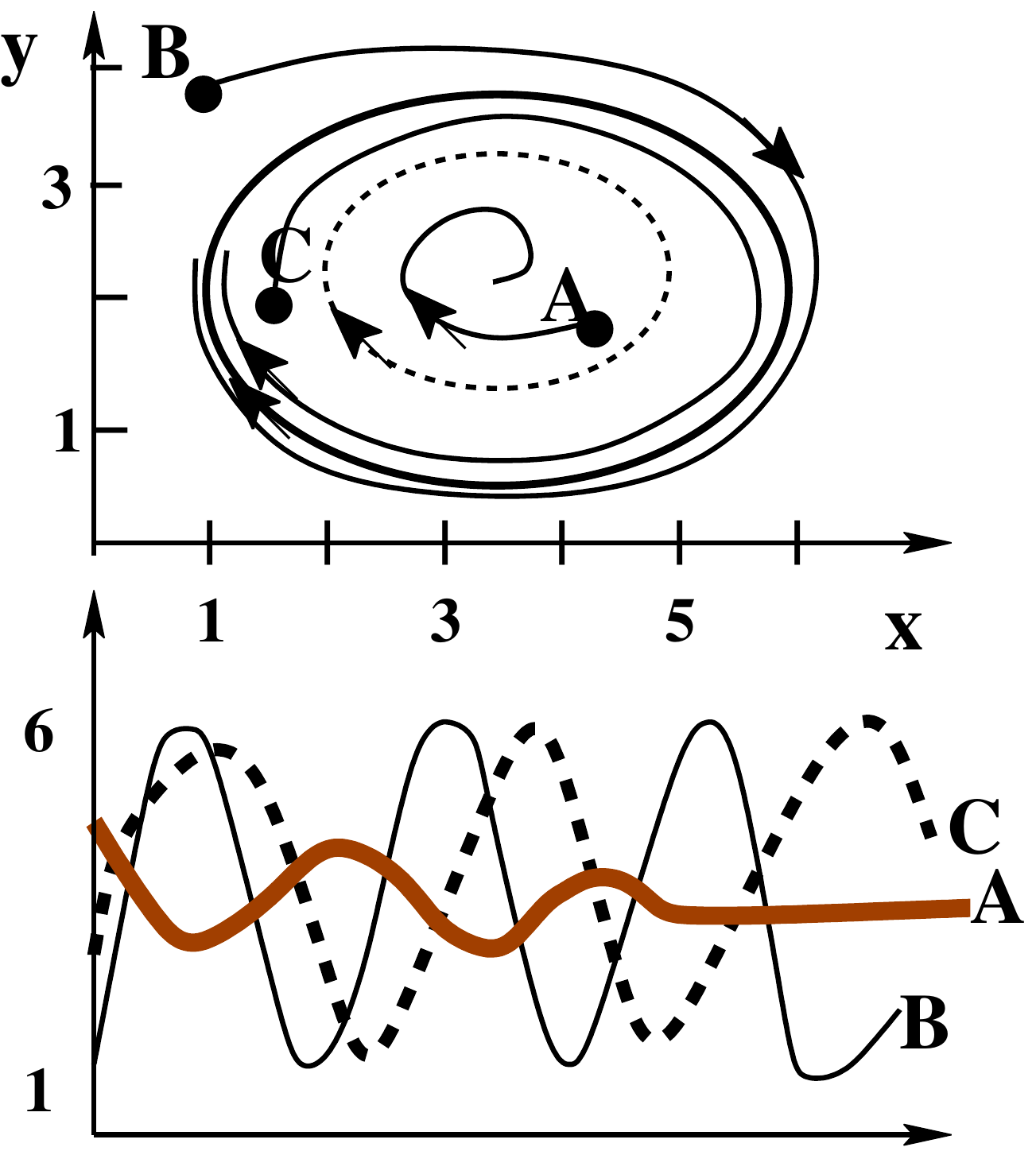,width=4cm}
\psfig{type=pdf,ext=.pdf,read=.pdf,figure=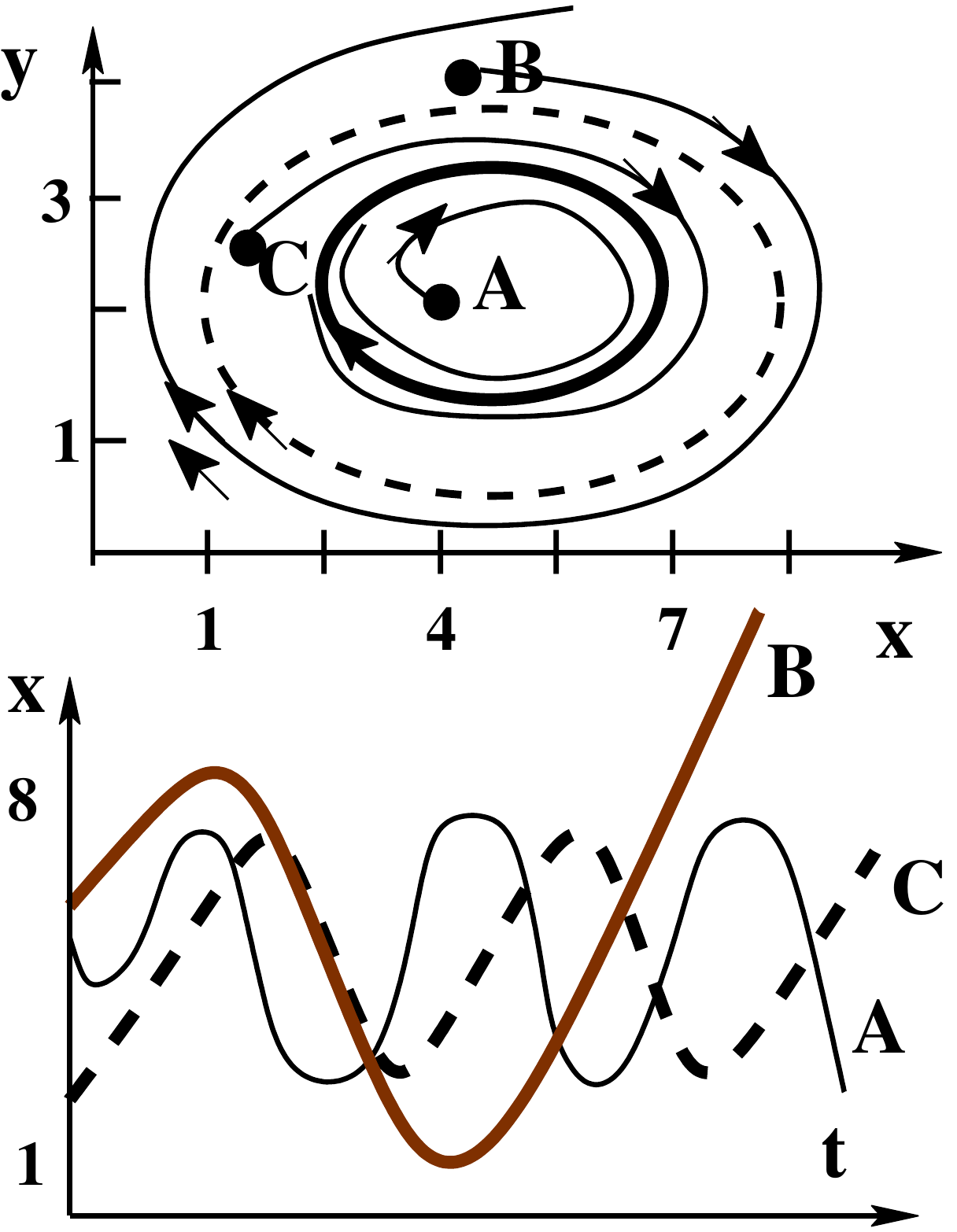,width=4cm}
\end{figure}
\een
\newpage
\onecolumn
\begin{figure}[hhh]
\vspace*{-3cm}\centerline{\psfig{type=pdf,ext=.pdf,read=.pdf,figure=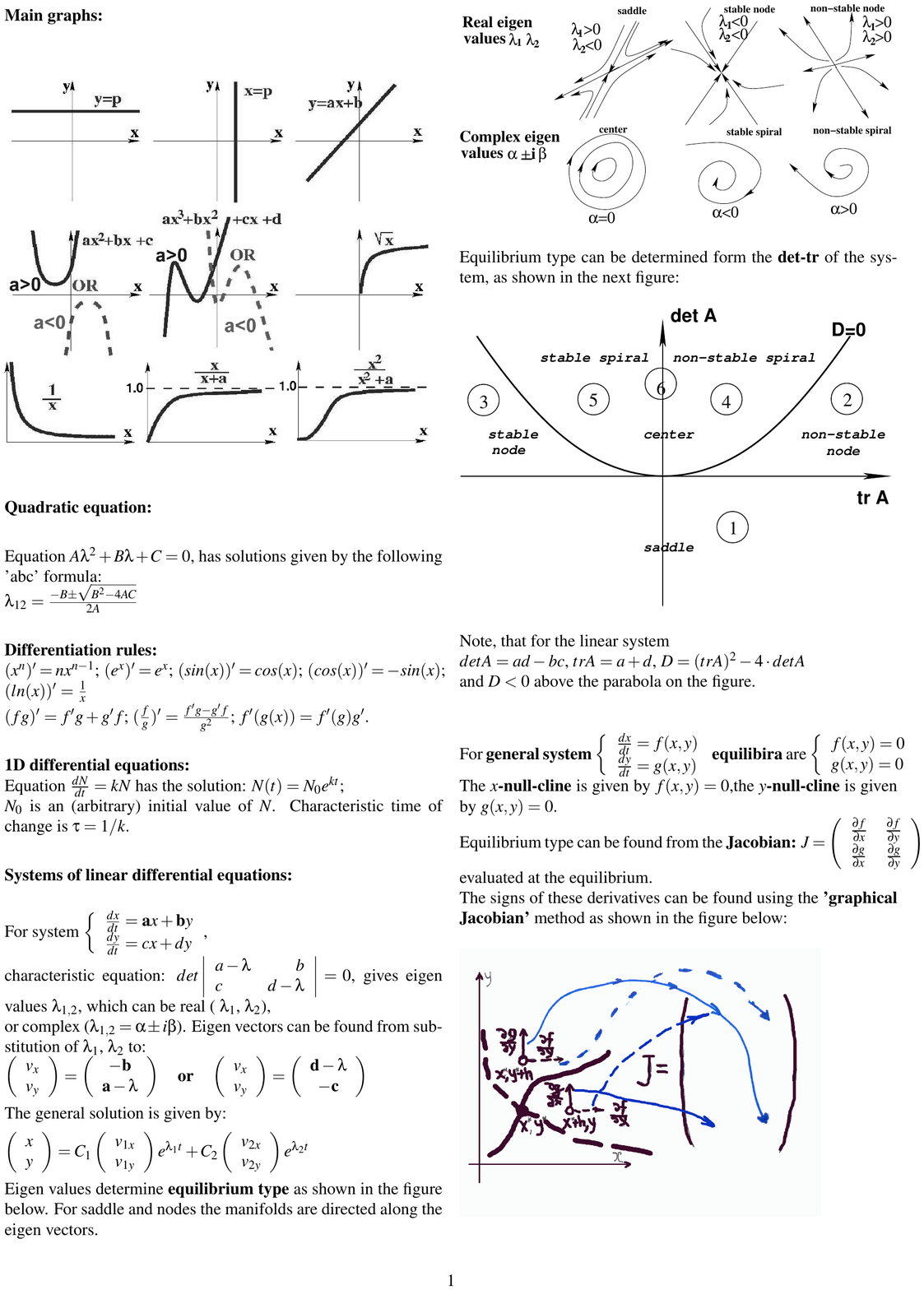,width=21cm}}
\end{figure}
\end{document}